\documentclass[10pt]{article}

\usepackage{template}
\newcommand{\m}{^{(l)}}
\renewcommand{\t}{^{\top}}
\newcommand{\inv}{^{-1}}
\newcommand{\uhat}{\mathbf{\hat U}}

\newcommand{\xtilde}{\mathbf{X}}
\newcommand{\xhat}{\mathbf{\hat X}}
\newcommand{\U}{\mathbf{U}}
\newcommand{\M}{^{(L)}}
\newcommand{\eps}{\varepsilon}
\newcommand{\one}{^{(1)}}
\newcommand{\A}{\mathbf{A}}

\newcommand{\E}{\mathbb{E}}
\newcommand{\wstar}{\mathbf{W}_*}
\newcommand{\ipq}{\mathbf{I}_{p,q}}
\newcommand{\yhat}{\mathbf{\hat Y}}
\newcommand{\ytilde}{\mathbf{ Y}}
\newcommand{\I}{\mathbf{I}}
\newcommand{\mi}{^{(l,-i)}}
\newcommand{\ycal}{\mathcal{Y}}
\newcommand{\yhatcal}{\mathcal{\hat Y}}
\newcommand{\snr}{\mathrm{SNR}}
\newcommand{\ave}{\mathrm{ave}}
\renewcommand{\l}{^{(l)}}

\newtheorem{proposition}{Proposition}[section]
\newtheorem{theorem}{Theorem}[section]
\newtheorem{corollary}{Corollary}[section]
\newtheorem{lemma}{Lemma}[section]

\theoremstyle{definition} 
\newtheorem{definition}{Definition}
\newtheorem{example}{Example}
\newtheorem{remark}{Remark}
\newtheorem{assumption}{Assumption}
\usepackage{natbib}

\usepackage[capitalize]{cleveref}
\crefname{assumption}{Assumption}{assumptions}
 \usepackage{thm-restate}

 \newcommand\numberthis{\addtocounter{equation}{1}\tag{\theequation}}
\renewcommand{\hat}{\widehat}
\newcommand{\specificthanks}[1]{\@fnsymbol{#1}}

\title{Joint Spectral Clustering in  Multilayer Degree-Corrected Stochastic Blockmodels\footnotetext{Corresponding author: Jes\'us Arroyo (Email: \href{mailto:jarroyo@tamu.edu}{jarroyo@tamu.edu}).}}

\author{Joshua Agterberg\thanks{Department of Statistics, University of Illinois Urbana-Champaign} 
\hspace{2em}  Zachary Lubberts\thanks{Department of Statistics, University of Virginia} 
\hspace{2em} Jes\'us Arroyo\thanks{Department of Statistics, Texas A\&M University}}

\date{\today}

\begin{document}

\maketitle

\begin{abstract} 
Modern network datasets are often composed of multiple layers, either as different views, time-varying  observations, or independent sample units, resulting in collections of networks over the same set of vertices but with potentially different connectivity patterns on each network. These data require models and methods that are flexible enough to capture local and global differences across the networks,
while at the same time being parsimonious and tractable to yield computationally efficient and theoretically sound solutions that are capable of aggregating information across the networks. 
This paper considers the multilayer degree-corrected stochastic blockmodel, where a collection of networks share the same community structure, but degree-corrections and block connection probability matrices are permitted to be different. We establish the identifiability of this model and propose a spectral clustering algorithm for community detection in this setting. Our theoretical results demonstrate that the misclustering error rate of the algorithm improves exponentially with multiple network realizations, even in the presence of significant layer heterogeneity with respect to degree corrections, signal strength, and spectral properties of the block connection probability matrices.
Simulation studies show that this approach improves on existing multilayer community detection methods in this challenging regime.
Furthermore, in a case study of US airport data through January 2016 -- September 2021, we find that this methodology identifies meaningful community structure and trends in airport popularity influenced by pandemic impacts on travel.

\end{abstract}

\tableofcontents

\section{Introduction}

Community detection, or the problem of clustering the vertices of a network into distinct groups (communities) in a coherent manner that somehow reflects the structure of the network, has become a fundamental tool for the analysis of network data, with many applications in fields such as neuroscience \citep{sporns2016modular}, biology \citep{luo2007modular}, social sciences \citep{conover2011political}, among others. 
In order to understand community detection in networks from a statistical perspective, a number of models have been proposed that characterize edge connectivity probabilities according to some notion of ground-truth communities. 

A workhorse community-based statistical model for networks is the \emph{stochastic blockmodel}, which posits that vertices belong to latent communities and that edges are drawn independently, with edge probability determined  by the community memberships of each vertex \citep{holland_stochastic_1983}.  A number of works have studied community detection from the lens of the stochastic blockmodel, including deriving information-theoretical limits \citep{zhang_minimax_2016} and phase transition phenomena \citep{abbe2017community}.  Of the various algorithms proposed for community detection in stochastic blockmodels, \emph{spectral clustering procedures} \citep{von_luxburg_tutorial_2007,rohe_spectral_2011,lei_consistency_2015}, which are collections of clustering techniques that use matrix factorizations such as eigendecompositions and singular value decompositions, have been shown to exhibit good performance both in practice and theoretically, including achieving perfect clustering down to the information-theoretical threshold \citep{lyzinski2014perfect,lei_unified_2019,abbe_entrywise_2020,su_strong_2020}.

One potential drawback of the stochastic blockmodel is that vertices are assumed to be ``equivalent'' within communities; i.e., edge probabilities are determined \emph{solely} by community memberships.  To relax this assumption, in the \emph{degree-corrected stochastic blockmodel} \citep{karrer2011stochastic} each vertex has associated to it a \emph{degree correction parameter} intended to shrink edge probabilities according to its magnitude. On the one hand, the degree-corrected stochastic blockmodel allows for vertex heterogeneity within communities, but on the other hand the model is more general than the stochastic blockmodel, often requiring more sophisticated procedures to recover communities. A number of variants of spectral clustering algorithms for community detection in this model have been considered \citep{lyzinski2014perfect,lei_consistency_2015,jin_fast_2015, gao_community_2018},  intended to ameliorate the ``nuisance'' degree correction parameters. {A strength of our approach is avoiding strong assumptions on the likelihood of the data, allowing for inference in a variety of settings without sacrificing generality.}

Many modern datasets deal with observations that consist of multiple networks on the same vertex set \citep{kivela2014multilayer,bazzi2020framework}, denoted as \emph{layers}, such as multiedges or multiview data, networks with time-varying structure, or multiple network observations. 
Community detection in these data presents additional challenges, as it is important to take advantage of a shared structure in the collection of graphs while respecting individual levels of idiosyncrasy. For these types of network data, which we refer to as \emph{multilayer networks}, perhaps the simplest community-based statistical model is the \emph{multilayer stochastic blockmodel} \citep{holland_stochastic_1983}. This model posits that communities are shared across networks but that edge probabilities change between networks.  

A key aspect of the multilayer stochastic blockmodel is that it allows for \emph{network heterogeneity} via the possibly changing edge probabilities.  However, as in the single network setting, vertices in the multilayer stochastic blockmodel are essentially equivalent; i.e., given their community memberships and the block probability matrices, their edge probabilities are entirely determined.  In the  \emph{multilayer degree-corrected stochastic blockmodel} that we consider in this work, individual vertices have network-specific degree correction parameters, so that there is \emph{global} network heterogeneity (via the connection probabilities), and \emph{local}  vertex heterogeneity (via the degree correction parameters). {In applications, degree heterogeneity is the rule, rather than the exception, as we can see in our real data analysis (\cref{sec:realdata}),} {where we observe changes in degree over time.
The inclusion of heterogeneous, network-specific degree corrections allows us to obtain accurate community estimates while simultaneously monitoring vertex-specific changes.}

{In this paper, we study the multilayer degree-corrected SBM, propose a spectral clustering algorithm to recover the communities, and study its misclustering error rate. To the best of our knowledge, this is the first paper to provide such an analysis for this model. More specifically, } our main contributions are as follows:
\begin{itemize}
    \item {We study the multilayer degree-corrected SBM, a flexible community model that allows for varying degree heterogeneity across layers.} We establish necessary and sufficient conditions for community identifiability of the multilayer degree-corrected stochastic blockmodel and propose a spectral clustering algorithm to estimate community memberships under this model.  Our necessary and sufficient conditions for identifiability also hold for the single network setting.
    \item We propose a spectral clustering procedure and obtain an expected misclustering error that improves exponentially with the number of networks, and we demonstrate perfect clustering under sufficient signal strength.  Our technical results rely only on signal strength conditions of each network and hold under severe degree heterogeneity within and between networks. {We also provide a lower bound on estimation, justifying our main technical assumptions on the network sparsity, and } {we also extend our result to settings with different community memberships across networks.} 
    \item In simulated data, we demonstrate that our method is competitive in multiple scenarios. Meanwhile, when there is  severe heterogeneity across the network layers, state-of-the-art community detection methods can fail in recovering the correct community structure of the model.
    \item We illustrate the flexibility of the model and methodology in a time series of United States flight network data from January 2016 to September 2021, identifying trends in airport popularity and the influence of COVID-19 on travel both at the local (vertex) and global (community) level.  
\end{itemize}
Our proposed algorithm consists of two stages: first, we compute individual (network-level) spectral embeddings, and then we compute a joint embedding by aggregating the output of the first stage.  
To prove our main technical results, we develop novel first-order entrywise expansions for each stage  of our algorithm that explicitly depend on all of the parameters of the model, including degree-corrections.  

The rest of this paper is structured as follows.  In the rest of this section  we consider closely related work and set notation.  We present our model, identifiability, and algorithm in \cref{s:model}.  The main results are presented in \cref{sec:mainresults}, and our simulations and real data analysis are presented in \cref{sec:sims} and \cref{sec:realdata} respectively.  We finish in \cref{S:discussion} with a discussion.   
The full proofs of all of our results are in the supplementary material.

\subsection{Related Work} \label{sec:relatedwork}

Community detection in the single network setting has received widespread attention in recent years \citep{abbe2017community,fortunato202220}. A number of works have studied community detection in the stochastic blockmodel, including consistency \citep{rohe_spectral_2011,zhao2012consistency,lei_consistency_2015}, phase transition phenomena \citep{abbe_entrywise_2020} and minimax rates \citep{gao_community_2018}. 
Beyond the stochastic blockmodel, a number of inference techniques have been considered for generalizations, such as the mixed-membership blockmodel \citep{airoldi_mixed_2008,mao_estimating_2021}, the random dot product graph \citep{athreya_statistical_2018} and generalised random dot product graph \citep{rubin-delanchy_statistical_2022}. This work is closely related to the literature on degree-corrected stochastic blockmodels \citep{karrer2011stochastic}.  The work \citet{jin_fast_2015} considered community detection in degree-corrected stochastic blockmodels using SCORE, or spectral clustering on ratios of eigenvectors, and several refinements, generalizations, and applications of this procedure have  been considered, including \citet{jin_estimating_2017,jin_improvements_2022,ke_optimal_2022} and \citet{fan_simple_2022}.  Our main results are perhaps most related to \citet{jin_improvements_2022}, who obtain an exponential error rate for spectral clustering with the SCORE procedure for a single network.  {We emphasize that the focus of our results is on the clustering error rate for \emph{multiple} networks, which presents its own different challenges.}

Turning to community detection in multilayer networks, several procedures have been considered for the multilayer stochastic blockmodel, including
spectral methods \citet{han_consistent_2015,bhattacharyya_spectral_2018,bhattacharyya2020consistent,huang_spectral_2020,lei_bias-adjusted_2022}, matrix factorization approaches 
\citep{paul_spectral_2020,lei_consistent_2020}, the 
expectation-maximization algorithm \citep{de2017community}, and 
efficient MCMC approaches
\citep{peixoto2015inferring,bazzi2020framework}. Extensions have also been considered, such as \citet{chen_global_2021}, which allows some members of each community to switch between networks.  Furthermore, \citet{jing_community_2021,pensky_clustering_2021}, and \citet{noroozi_sparse_2022} all consider generalizations of the multilayer stochastic blockmodel where there are a few different possible community configurations. 
Although spectral methods are competitive in terms of computation and accuracy, existing methods are limited in handling heterogeneous degree correction parameters. Both \citet{bhattacharyya2020consistent} and \citet{bhattacharyya_spectral_2018} consider  degree-corrections for each network, but they require that the degree-corrections remain the same across networks, making the analysis feasible.  Our work is perhaps most closely connected to the works \citet{arroyo_inference_2021} and \citet{zheng_limit_2022}, which consider the estimation of a common invariant subspace, but the model we consider in this paper is substantially different, and we provide finer theoretical results to analyze misclustering rates. {In particular, as the model we consider here permits degree heterogeneity, our algorithm requires an additional nonlinear normalization step, and the interplay between this step and the heterogeneous noise presents further challenges in the analysis.    Moreover, in \citet{arroyo_inference_2021} the authors only consider error rates in Frobenius norms, whereas we require a stronger characterization in the $\ell_{2,\infty}$ norm.   More detailed discussion of our proof techniques are provided in \cref{sec:proofoverview}.
}

From a technical point of view, our analysis is also closely related to the literature on entrywise eigenvector analysis of random matrices \citep{abbe_entrywise_2020,chen_spectral_2021}. {However, these works only focus on a single network and our analysis of the multilayer embedding is entirely novel. } 
Several authors have previously considered the entrywise analysis of the eigenvectors of a single degree-corrected stochastic blockmodel, such as \citet{lyzinski2014perfect,jin_estimating_2017,jin_improvements_2022,su_strong_2020} and \citet{ke_optimal_2022}. Here we provide an entrywise analysis of the \emph{scaled} eigenvectors of degree-corrected stochastic blockmodels, {which we empirically observe to perform better under unbalanced community sizes.} 

\subsection{Notation} \label{sec:notation}
We use bold or greek capital letters $\mathbf{M}$ or $\Lambda$ for matrices, and we let $\mathbf{M}_{i\cdot}$ and $\mathbf{M}_{\cdot j}$ denote the $i$'th row and $j$'th column respectively, where we view both as column vectors.  We let
$\|\mathbf{M}\|$ and $\|\mathbf{M}\|_{2,\infty}$ denote its spectral and $\ell_{2,\infty}$ norm, where the latter is defined as $\max_i \|\mathbf{M}_{i\cdot}\|$, where $\|\mathbf{M}_{i\cdot}\|$ is the usual (vector) Euclidean norm. For a vector $x$ we let $\|x\|_1, \|x\|_{\infty}$ denote its vector $\ell_1$ and $\ell_{\infty}$ norms respectively.  We let $\mathbf{I}_r$ denote the $r\times r$ identity. For two orthonormal matrices $\mathbf{U}$ and $\mathbf{V}$, we let $\|\sin\bTheta(\U,\mathbf{V})\|$ denote their (spectral) $\sin\bTheta$ distance, defined as $\|\sin\bTheta(\U,\mathbf{V})\| = \|(\mathbf{I} - \U \U\t)\mathbf{V}\|$.  We write $\mathbb{O}(r)$ to denote the set of $r\times r$ orthogonal matrices.  We also denote $e_i$ as the standard basis vector, and  we view $e_i\t \mathbf{M}$ as a column vector. We let $\mathbb{I}\{\cdot\}$ denote the indicator function, and $\real_+$ denote the strictly positive real numbers.  For two functions $f(n)$ and $g(n)$, we write $f(n) \lesssim g(n)$ if there exists some constant $C > 0$ such that $f(n) \leq C g(n)$, and we write $f(n) \ll g(n)$ if $f(n)/g(n) \to 0$ as $n \to \infty$.  We denote by $f(n) \asymp g(n)$ the case where both $f(n) \lesssim g(n)$ and $g(n) \lesssim f(n)$.  We also write $f(n) = O(g(n))$ if $f(n) \lesssim g(n)$.  We denote $[n] = \{1, 2, \dots, n \}$. 

\section{The Multilayer Degree-Corrected SBM}

\label{s:model}

Suppose one observes a collection of $L$ adjacency matrices $\bA^{(1)}, \ldots, \bA^{(L)}$ of size $n\times n$, with the vertices of the corresponding graphs aligned across the collection. For simplicity of the presentation and the theory,  we assume that the adjacency matrices represent simple undirected graphs, hence these matrices are symmetric with binary entries, and we allow the networks to have self-edges (loops), but the main results are not materially different if loops are not permitted. Much of the theory and methodology we consider here is also applicable in the settings of weighted or directed networks, but we focus on the binary and undirected setting since our primary concern in the present work is to quantify the misclustering error rate as a function of the degree parameters.

The model considered in this paper assumes a shared community structure across all the graphs, but allows for idiosyncrasy in the edge probabilities across the collection of graphs by letting the global and local individual parameters of each graph to be different. In particular, we consider a multilayer version of the degree-corrected stochastic blockmodel \citep{karrer2011stochastic}, in which both the block connectivity matrices and the vertex degree parameters can be different for each network. Some versions of this model have appeared in \cite{peixoto2015inferring,bazzi2020framework,bhattacharyya2020consistent,paul2021null}, but to be precise, we will use the following definition. 

\begin{definition}[Multilayer Degree-Corrected Stochastic Blockmodel]
A collection of $L$ graphs $\{\mathbf{A}\l\}_{l=1}^{L}$ on $n$ vertices are drawn from the \emph{multilayer degree-corrected stochastic blockmodel} (multilayer DCSBM) if:
\begin{itemize}
    \item each vertex $i$ belongs to one of $K$ communities.  Let $z: [n] \to [K]$ be the community membership function  satisfying $z(i) = r$ if vertex $i$ belongs to community $r$;
    \item $\theta\m_1, \ldots, \theta\m_n\in\real_+$ are \emph{degree correction parameters} associated to nodes in network $l$;
    \item $\bB^{(1)}, \ldots, \bB^{(L)}\in\real_+^{K\times K}$ are symmetric \emph{block connectivity matrices};
    \item the edges of the networks are mutually independent, and their expected values (probabilities) are described by
  $\mathbb{E}[\bA\m_{ij}] = \theta_i\m\theta_j\m\bB\m_{z(i), z(j)}, \quad\quad l\in[L],\ i,j\in[n], i\geq j.$
\end{itemize}
\end{definition}

The degree correction parameters denote a local connectivity component and the block connection probability matrices characterize a global connectivity component,  both of which can vary from graph to graph, while the community memberships remain constant. Since the edges are binary, the expected value also denotes the probability of the corresponding edge, but this definition can be used in other distributions (e.g. Poisson \citep{karrer2011stochastic}).

It is convenient to represent the multilayer DCSBM  using matrix notation. Denote the collection of matrices that encode the edge expectations by $\bP^{(1)}, \ldots, \bP^{(L)}\in[0,1]^{n\times n}$, such that $\e[\bA^{(l)}_{ij}] = \bP^{(l)}_{ij}$ for each $l\in[L]$ and $i,j\in[n], i\geq j$. Then we can write
\begin{align*}
    \bP^{(l)} = \bTheta^{(l)} \bZ \bB^{(l)} \bZ^\top\bTheta^{(l)}, \quad\quad l\in[L], \numberthis\label{eq:model}
\end{align*}
where $\bTheta^{(l)}\in\real^{n\times n}$ is a diagonal matrix with $\bTheta_{ii}^{(l)} = \btheta^{(l)}_i>0$, $\bZ\in\{0,1\}^{n\times K}$ is a binary matrix indicating community memberships ($\bZ_{ir}=1$ if $z(i)=r$, and $\bZ_{ir}=0$ otherwise), and $\bB^{(l)}\in \real_+^{K\times K}$ is a symmetric matrix. We assume that $\mathrm{rank}(\bB\l) = K_l$, and we allow $K_l$ to be less than $K$.


The multilayer DCSBM model is flexible enough to represent heterogeneous structures both at the vertex and the community levels, while retaining a joint community structure across the graphs. Due to these local and global idiosyncrasies, distinguishing between local and global graph structure at the single and multilayer level becomes important, as it is possible to formulate parameterizations of the model that give equivalent characterizations.
For instance, one may group high degree vertices in their own community according to degree correction parameters alone. To ensure identifiability and maintain a parsimonious model, we assume that the number of communities $K$ is the smallest possible that can represent the communities uniquely (up to label permutations). Our first result establishes the identifiability of the communities in the model.

\begin{theorem}[Community membership identifiability]
\label{prop:identifiability}
Suppose that $\{\mathbf{P}^{(l)}\}_{l=1}^{L} \in \mathbb{R}^{n\times n}$ are matrices such that
$\mathbf{P}^{(l)} = \bTheta^{(l)} \mathbf{Z} \mathbf{B}\m \mathbf{Z}\t \bTheta^{(l)},$ $l\in[L]$
where $\mathbf{Z} \in \{0,1\}^{n\times K}$ is a binary block membership matrix with at least one vertex in each community ( $\sum_{r=1}^{K} \mathbf{Z}_{ir}=1, i\in[n]$, and $\sum_{i=1}^n \bZ_{ir} \geq 1$, $r\in[K]$), $\{\mathbf{B}\m\}_{l=1}^L$ are symmetric matrices with entries in $\real_+$, and $\{\bTheta\m\}_{l=1}^L$ are diagonal matrices with positive entries on the diagonal. Let $\bB\l = \bV\l\bD\l(\bV\l)^\top$ be the eigendecomposition of $\bB\l$, with $\bV\l\in\real^{K\times K_l}$ a matrix with orthonormal columns and $\bD\l\in\real^{K_l}$ a diagonal matrix, and $\text{rank}(\bB\l) = K_l$. Write ${\bQ}\l$ as the matrix with normalized rows of $\bV\l$, i.e., $\bQ\l_{r\cdot} = \frac{1}{\|\bV\l_{r\cdot}\|}\bV\l_{r\cdot}$, and let $\bQ=[\bQ^{(1)}, \cdots, \bQ^{(L)}]$. The membership matrix $\mathbf{Z}$ is identifiable (up to label permutations) if and only if $\mathbf{ Q}$ has no repeated rows.
\end{theorem}
 The identifiability condition requires that the matrices $\{\bB^{(l)}\}$ have exactly $K$ jointly distinguishable rows, which determine the  community memberships.
  The condition $\bQ$ having no repeated rows implies that there are precisely $K$ unique \emph{directions} associated to the rows of $[\bV\one, \cdots, \bV^{(L)}]$. Therefore, if ${\bQ}\l$ is defined by normalizing the rows of $\bV\l$ in any other way, then as long as ${\bQ}$ has $K$ distinct rows, the communities will be identifiable. \cref{prop:identifiability}  also holds for $L = 1$, thereby establishing both  necessary and sufficient conditions for identifiability in the single network model.

The matrix $\bB\l$ is often assumed to be full rank \citep{qin_regularized_2013,jin_optimal_2022}, in which case there are exactly $K$ identifiable communities.  \cref{prop:identifiability} requires a milder condition to allow flexibility in modeling multiple networks, as
the number of identifiable communities in each layer may be smaller than $K$. The identifiable communities in the joint model are given by the different directions taken by the combined rows of $\bB\l$ across all the layers. 
Since this condition is also necessary for identifiability, this value of $K$ gives the most parsimonious representation in terms of the number of communities.
 
Identifiability of the degree correction and block connectivity parameters requires additional constraints, as it is otherwise possible to change their values up to a multiplicative constant. Multiple characterizations have been used previously for the single-network setting, and these immediately extend to the multilayer setting. For instance, if for all $i$, we have $\bB\m_{ii} = 1$, (e.g. \cite{jin_optimal_2022}) then the other model parameters are identifiable as well. We adopt this identifiability constraint to facilitate the presentation of the theoretical results in \cref{sec:mainresults}, as it allows us to isolate the effect of the degree correction parameters. Nevertheless, to ease interpretation, in \cref{sec:realdata} we adopt a different constraint, namely, that  the sum of degree corrections within each community is equal to 1. Both parameterizations are equivalent. 


\subsection{Degree-Corrected Multiple Adjacency Spectral Embedding}
\label{s:method}
In order to obtain a statistically principled, computationally efficient, and practical algorithm for community detection, we will consider a spectral clustering procedure.  General spectral clustering approaches for one network typically proceed as follows: first, using a few leading eigenvectors of the adjacency matrix (or related quantities, such as the graph Laplacian), obtain individual vertex representations by considering the rows of the matrices; we will refer to this first step as obtaining an \emph{embedding}. Then, the communities are estimated by clustering the rows of this matrix  using a clustering algorithm.  

For multilayer networks with shared community structure, the general procedure is similar, only now the requirement is to use all of the networks to obtain individual vertex representations in a low-dimensional space.  For the multilayer stochastic blockmodel,  a typical approach is to simply consider a few leading eigenvectors of the average adjacency matrix $\bar \bA = \frac{1}{L} \sum_{l} \bA\l$ \citep{tang2009clustering,han_consistent_2015}.  However, as discussed in e.g. \citet{paul_spectral_2020,lei_bias-adjusted_2022}, this procedure is only guaranteed to work when there is certain level of homogeneity in the block connectivity matrices, and it can fail if the $\mathbf{B}\l$ matrices are different.  \citet{lei_bias-adjusted_2022} proposed to rectify this by considering a bias-corrected version of the sum of the squared adjacency matrices.  Alternatively, one can look at an embedding obtained by aggregating  the projections onto the principal subspaces of each graph
\citep{paul_spectral_2020,arroyo_inference_2021}. 
In these situations, the population probability matrices $\{\bP\l\}$ share a common singular subspace, and running the relevant algorithm on those reveals the community memberships.  Unfortunately, this is not the case in the model considered herein, but with some modification, a certain matrix can be shown to have a left singular subspace that reveals the community memberships.

Our proposal to find an embedding is based on several observations concerning the joint spectral geometry of the matrices $\{\bP\l\}$, some of which have been considered before in the single-network literature.
\begin{itemize}
    \item \textbf{Observation 1}:  \textit{The rows of the $K_l$ (scaled) eigenvectors of $\bP\l$ are supported on at most $K$ different rays in $\real^{K_l}$, with each ray corresponding to a distinct community, and magnitude of each row determined by the magnitude of its corresponding degree-correction parameter. }
\end{itemize}
{By virtue of the clustering structure in the DCSBM, spectral embeddings of  $\bP\l$, such as scaled or unscaled eigenvectors \citep{lyzinski2014perfect,jin_fast_2015}, preserve these clusters.}
Suppose that each $\bP\l$ has eigendecomposition $\U\l \Lambda\l (\U\l)\t$, where $\U\l$ is an $n\times K_l$ orthonormal matrix and $\Lambda\l$ is the matrix of eigenvalues of $\bP\l$. Define
 $  \xtilde\l := \U\l |\Lambda\l|^{1/2}, 
$
where $|\cdot|$ is the entrywise absolute value.  It can be shown (see the proof of \cref{prop:algorithmcorrectness} below) that $\xtilde\l = \bTheta\l \mathbf{Z} \bM\l$, where $\bM\l \in \mathbb{R}^{K \times K_l}$ has $K'_l$ unique rows, with $K_l\leq K'_l\leq K$. {Note that in general $\bM\l$ may not have exactly $K_l$ unique rows unless $\bB\l$ is also assumed to be rank $K_l$.}
Explicitly, Observation 1 implies that each row $i$ of $\xtilde\l$ satisfies
\begin{align*}
    \xtilde\l_{i\cdot} &= \theta\l_i \bM\l_{z(i)\cdot}. \numberthis \label{xtilde}
\end{align*}
\begin{figure}[t]
    \centering
    \includegraphics[width=\textwidth]{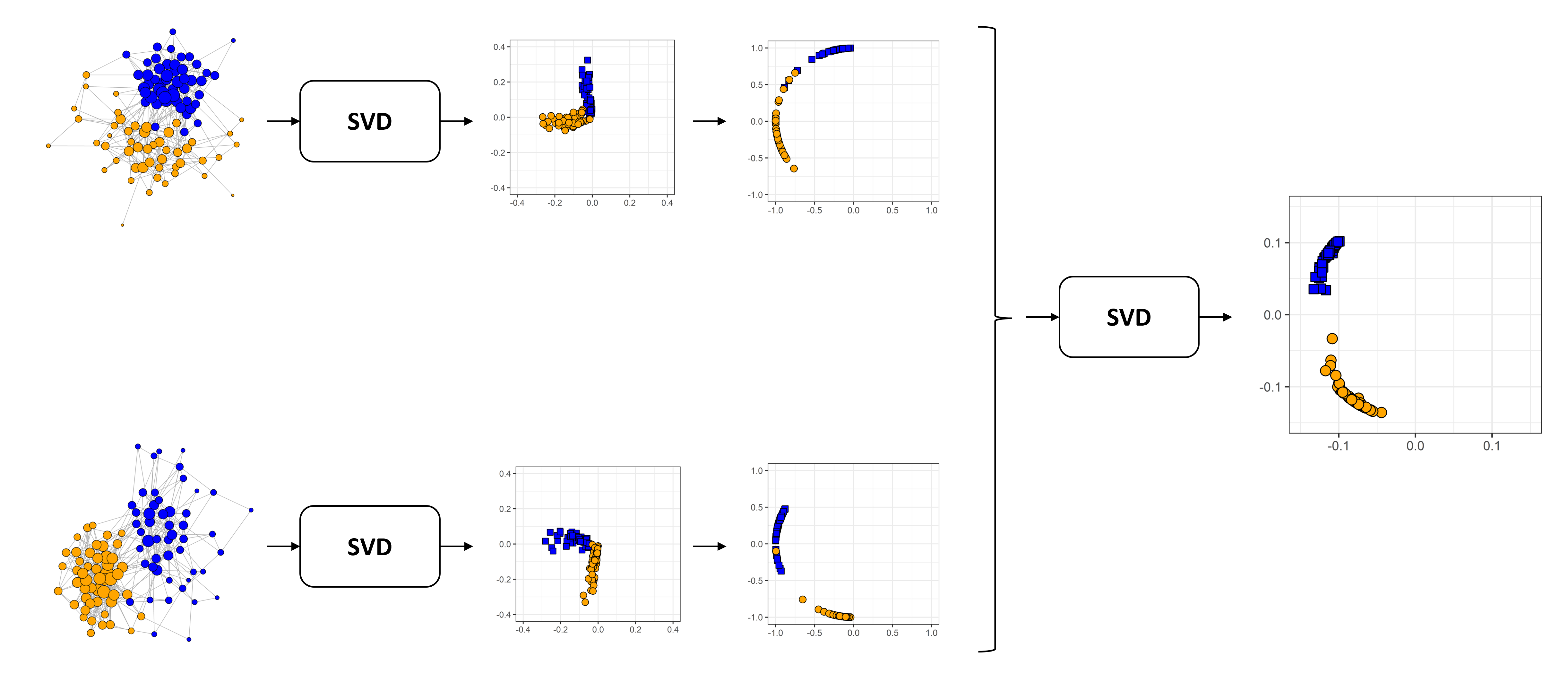}
    \caption{{\small Pictorial representation of \cref{alg:dcmase}.}}
    \label{fig:dc-mase}
\end{figure}
\vspace{-30pt}
\begin{itemize}
    \item \textbf{Observation 2}: \textit{ Projecting each row of $\xtilde\l$ to the sphere results in a matrix of at most $K$ unique rows, with each row corresponding to community membership.}
    \end{itemize}
    \noindent
{Different  normalization approaches have been proposed to handle degree heterogeneity  \citep{lei_consistency_2015,qin_regularized_2013,jin_fast_2015,zhang2020detecting}. In this work, we use the spherical normalization.} Define $\ytilde\l$ via 
$\ytilde\l_{i\cdot} = \frac{\xtilde\l_{i\cdot}}{\|\xtilde_{i\cdot}\l\|}.$
By \eqref{xtilde}, it holds that $\ytilde\l_{i\cdot} = \frac{\bM_{z(i)\cdot}\l}{\|\bM_{z(i)\cdot}\l\|}$.  In particular, there are only $K_l'\leq K$ unique rows of $\ytilde\l$, with each row corresponding to community membership.
\begin{itemize}
    \item \textbf{Observation 3}:  \textit{The left singular subspace of $\ycal= [ \ytilde\one, \ytilde^{(2)}, \cdots, \ytilde^{(L)}] \in \mathbb{R}^{n \times \sum_l K_l}$ reveals the community memberships}.
\end{itemize}
Suppose that $\ycal$ has singular value decomposition given by $\ycal = \U \Sigma \mathbf{V}\t$. It can be shown (see \cref{prop:algorithmcorrectness}) that under the condition of \cref{prop:identifiability}, $\text{rank}(\ycal)= \tilde{K} \leq K$ and $\U\in\real^{n\times \tilde{K}}$ satisfies $
    \U = \mathbf{Z}  \mathbf{M},
$
where $\mathbf{M} \in \mathbb{R}^{K\times \tilde{K}}$ is some matrix without repeated rows. 
Explicitly, this says that there are only $K$ unique rows of $\U$, with each row $i$ of $\U$ corresponding to community membership of vertex $i$.  Moreover, since $\U$ is obtained via the singular value decomposition of $\ycal$, it contains information from all the networks.  

The observations presented above lead to a joint spectral clustering algorithm applied to the sample adjacency matrices, summarized in \cref{alg:dcmase} and in \cref{fig:dc-mase}. Without the row-normalization step, one obtains the \emph{scaled} multiple adjacency spectral embedding (MASE) algorithm of \citet{arroyo_inference_2021}, who consider a model where each ``population'' network shares a common invariant subspace (which  includes the multilayer SBM as a special case). Due to the different degree correction parameters, the population matrices in the  multilayer DCSBM do not share a common invariant subspace, but our algorithm can  be viewed as a normalized version of the MASE algorithm, so we dub it our degree-corrected multiple adjacency spectral embedding (DC-MASE).
Introducing this normalization step is crucial in the presence of heterogeneous degree corrections and makes this methodology applicable to a  more flexible model. The following proposition formalizes the three arguments to construct the algorithm. 

 	\begin{algorithm} [t]
 		\caption{Degree-corrected multiple adjacency spectral embedding (DC-MASE)}
 		\begin{algorithmic} 
 			\Input Collection of adjacency matrices $\bA^{(1)},\ldots,\bA^{(L)}$; individual ranks $K_1, \ldots, K_L$, joint rank $\tilde{K}$, number of communities $K$.
 			\begin{enumerate}
 				\item For each graph $l\in[L]$, 
 				\begin{enumerate}
 				    \item Let $\widehat{\bX}\m\in\mathbb{R}^{n\times K_l}$ be defined $\xhat\l := \uhat\l |\widehat \Lambda\l|^{1/2}$, where $\uhat\l$ is the matrix containing the $K_l$  eigenvectors associated to the $K_l$ largest eigenvalues (in magnitude) of  $\bA\m$ and $\widehat \Lambda\l$ are the corresponding eigenvalues;
 				    \item let $\yhat\l\in\real^{n\times K_l}$ be the matrix containing the rows of $\xhat\l$ projected to the sphere, defined as 
 				      $ \yhat\l_{i\cdot} = \frac{\xhat\l_{i\cdot}}{\|\xhat\l_{i\cdot}\|}.$
 				\end{enumerate}
 		\item Form the matrix $\yhatcal = [\yhat\one, \cdots \yhat^{(L)}]$ by concatenating the row-scaled  matrices. 
 		\item Let $\uhat\in\mathbb{\real}^{n\times \tilde{K}}$ be the matrix  containing the $\tilde{K}$ leading left singular vectors of $\yhatcal$.
 		\item Assign  memberships as the clusters of the rows of $\uhat$ into $K$ groups via $K$-means.
 			\end{enumerate}
 			\Output Community memberships.
 		\end{algorithmic}
 		\label{alg:dcmase}
 	\end{algorithm}

\begin{proposition}\label{prop:algorithmcorrectness}
Under the conditions of \cref{prop:identifiability}, Algorithm \ref{alg:dcmase} applied to the collection of matrices $\mathbf{P}\one,$ $\dots,$ $\mathbf{P}\M$ recovers the community memberships exactly.
\end{proposition}

 \begin{remark}[Alternative Approaches to Embedding and Clustering]
  Variations of   Algorithm~\ref{alg:dcmase}  can be obtained by changing the initial embedding, row-normalization, or clustering procedures, for which we conjecture that similar results to \cref{prop:algorithmcorrectness} may hold, but we do not undertake a complete analysis of these different choices.
Other variations can be obtained by changing the embedding, for example, to unscaled eigenvectors or using the Laplacian matrix; the normalization procedure, for example, by using SCORE \citep{jin_fast_2015} or by changing the clustering procedure to $K$-medians \citep{lei_consistency_2015}. 
{In Appendix H  of the Supplementary Material, we discuss the implications of different choices of the embedding method. In practice, we have found that spherical normalization on the scaled eigenvectors is robust against severe degree heterogeneity and unbalanced communities, which is the reason we focus on this method for the theoretical analysis.}
 \end{remark}

\subsection{Estimating the Number of Communities}\label{sec:estimate-K}

Choosing the number of communities in the multilayer DCSBM via DC-MASE is an important yet challenging problem, as one is required to estimate the individual and joint embedding dimensions for each adjacency matrix, as well as the total number of communities in the joint model. Throughout this paper, we assume that these numbers are known or can be estimated appropriately, but we discuss here some approaches for choosing these parameters in practice.

The first step of Algorithm~\ref{alg:dcmase} requires the selection of $K_l$, which corresponds to the rank of the matrix $\bP\l= \e[\bA\m]$, and hence this corresponds to a rank estimation problem. A common practical approach is to look for an elbow in the scree plot of the eigenvalues of the adjacency matrix \citep{zhu2006automatic}. Similarly, to estimate $\tilde{K}$, one can look for elbows in the scree plot of the singular values obtained from the concatenated matrix $\widehat{\mathcal{Y}}$, as this matrix concentrates around a population matrix that has rank exactly equal to $\tilde{K}$. In simulations, we have  observed that overestimating these parameters typically does not have a significant effect on the performance of the clustering method.

The choice of $K$ is more important, as it controls the number of communities in the joint model. Several existing methods assume that the matrix $\bB\l$ has full rank, in which case the value of $K_l$ corresponds to the number of communities in the degree-corrected SBM for each network $l\in[L]$. A number of methods exist for estimating the communities in a single-layer DCSBM \citep{wang_likelihood-based_2017,ma2021determining,le2022estimating,li2020network}, 
including recent work by \cite{jin_optimal_2022}, who achieves the optimal phase transition under this assumption. Alternatively, one can use an appropriate criterion for choosing the number of clusters via $K$-means.

\section{Main Results} \label{sec:mainresults}
Having described our algorithm in detail, we are now prepared to discuss the associated community recovery guarantees.  In order to do so, we first must state some assumptions on the regularity of each network.  For simplicity of analysis and to facilitate interpretation, we assume that $\mathbf{B}\m_{rr}=1$ for all $r\in[K], l\in[L]$, and that
each $\mathbf{B}\m$ is rank $K$;  {
extensions are discussed at the end of this section.} 

\begin{assumption}[Regularity Conditions]
\label{ass:communityass}
Let $\mathcal{C}(r)$ denote the indices associated to community $r$; i.e., the set of $i$ such that $z(i) = r$.  It holds that $|\mathcal{C}(r)| \asymp |\mathcal{C}(s)|$ for $r\neq s$ and
$ K\|\theta\m_{\mathcal{C}(r)}\|^2 \asymp \|\theta\l\|^2$
for all $r \in [K]$.  
In addition, each matrix $\mathbf{B}\m$ is rank $K$ with unit diagonals; let $\lambda_{t}\m$ denote its ordered eigenvalues. Then $|\lambda_K\m| \geq \lambda_{\min}\m$ for some $\lambda_{\min}\m \in (0,1)$ and $\|\mathbf{B}\m\| = \lambda_1\l 
{\lesssim K}$. 
\end{assumption}
The first part of \cref{ass:communityass}  essentially requires that the communities and degree corrections within each community are balanced,  and it is commonly imposed in the analysis of the DCSBM \citep{jin_optimal_2022,su_strong_2020}, but it can be relaxed by keeping track of these constants. We also assume for simplicity that $\|\mathbf{B}\l\|\lesssim K$, which is not strictly required but facilitates analysis.  If $\|\mathbf{B}^{(l)} \| \gg K$, then $\mathbf{B}^{(l)}$ has a few very large entries, which makes clustering easier.
\newenvironment{psmallmatrix}
  {\left(\begin{smallmatrix}}
  {\end{smallmatrix}\right)}
  
We have also introduced the parameter $\lambda_{\min}\l$, which can be understood as a proxy for the community separation. For example, consider the matrix $\bB\l = \begin{psmallmatrix} 1 & 1 - \eta \\ 1- \eta & 1 \end{psmallmatrix}.$
Then it holds that $\lambda_{\min}\l = \eta$.  
We also assume for simplicity that $\lambda_{\min}\l \in (0,1);$ when this is not the case, the communities are well-separated, so the problem is qualitatively easier.

Next we introduce some assumptions on the individual network-level signal strengths and degree homogeneity.  Let $\theta_{\min}\m := \min_{i} \bTheta\m_{ii}$, and let $\theta_{\max}\m$ be defined similarly.  Define also the   average minimum eigenvalue parameter $\bar \lambda := \frac{1}{L}\sum_{l=1}^{L} \lambda_{\min}\m \in (0,1).$
The following is our main technical assumption on the individual network signal strengths.

\begin{assumption}[Network-Level Signal Strengths] \label{ass:networklevel}
There exist constants $C$ and $c$ (with $C$ depending on the community sizes) such that each network $l$ satisfies
\begin{align*}
&& 
 {C\bigg( \frac{\theta_{\max}\l}{\theta_{\min}\l} \bigg)  \frac{K^{8} \theta_{\max}\l \|\theta\l\|_1 \log(n)}{(\lambda_{\min}\l)^2\|\theta\l\|^4}}
 &\leq \bar \lambda; && \textit{(Signal Strength)}\\
  && \frac{\theta_{\min}\l}{\theta_{\max}\l} &\geq \sqrt{\frac{\log(n)}{n}} && \textit{(Degree Heterogeneity)}  \\
  &&  \theta_{\min}\m \|\theta\m\|_1 &\geq c\log(n). \qquad &&\textit{(Logarithmic Degree Growth)}.
\end{align*}
\end{assumption}
\noindent
To build intuition we consider several examples.

\begin{example}[Degree-Correction Heterogeneity] \label{ex:2}
We consider a setting with $\lambda_{\min}\asymp 1$, $K\asymp 1$ and we suppose that $\theta_i\l=a,$ for $1\leq i \leq \gamma n$ and $\theta_i\l=b>a$ for $\gamma n+1\leq i\leq n.$
It is easy to show that \cref{ass:networklevel} 
holds if $ 
    \frac{b^2(\gamma a + (1-\gamma)b)}{a(\gamma a^2 + (1-\gamma)b^2)^2}\lesssim  \frac{n}{\log(n)}
$
 and $a/b \gtrsim \sqrt{\frac{\log(n)}{n}}$.  For example, if $\gamma n = 1$ (an outlier model) and $b\gg a$, the first condition reduces to $ab \gtrsim \log(n)/n$. If $b = 1$,
$a=\sqrt{\log(n)/n}$ satisfies the degree heterogeneity assumption. On the other hand, when $\gamma n = n-1$ and $b=1$, $a\geq (\log(n)/n)^{1/4}$ is required.
\end{example}
\begin{example}[Close Communities with Homogeneous Degree Corrections] \label{ex:3}
We consider a setting with all $\theta_i\l\asymp \sqrt{\rho_n}$, $K\asymp 1$, and $\lambda_{\min}\l\asymp \lambda_{\min}$ for all $l$. Then we require $\lambda_{\min}^3\gtrsim \frac{\log(n)}{n\rho_n}.$ If only $o(L)$ networks have $\lambda_{\min}\l\asymp\lambda_{\min}$, and all others have $\lambda_{\min}\l\asymp 1$, then we have the weaker condition $\lambda_{\min}^2\gtrsim\frac{\log(n)}{n\rho_n}.$ Then so long as the majority of networks have strong signal, we can tolerate even weaker signal in the worst-behaved layers.
\end{example}

When $L=1$, our condition in \cref{ass:networklevel} is only slightly stronger than that of \citet{jin_improvements_2022} in terms of $\bar \lambda$ and slightly weaker in terms of $\frac{\theta_{\max}\l}{\theta_{\min}\l}$, though we include a more detailed comparison in Appendix H. To understand the intuition behind the signal-strength condition in terms of $\bar \lambda$ in Assumption \ref{ass:networklevel}, observe that when
$\bar \lambda$ is small, the average community separation is small, and hence the rays associated to each (unscaled) embedding $\xtilde\l$ (see \eqref{xtilde}) will be nearly colinear.  
Therefore, in order for the SVD step to succeed, we will require sufficient separation of the communities, which is why  \cref{ass:networklevel} concerns $\bar \lambda$. 

{
The assumption that $\mathbf{B}^{(l)}$ is rank $K$ ensures that the smallest nonzero eigenvalue of the matrix $\ycal \ycal\t$ is sufficiently large (Lemma A.2). In general, our main results will continue to hold as long as the communities are identifiable (\cref{prop:identifiability}) and
$\lambda_{\min}\big( \ycal\ycal\t) \gtrsim \frac{n}{K} L \bar \lambda$. 
For instance, suppose $\mathbf{B}^{(l)}$ has $K - K_l$ repeated rows,  which may occur if two communities  ``merge'' within one network but otherwise behave differently in other networks.
}

\subsection{Misclustering Error Rate and Perfect Clustering}
 
With these assumptions in hand, we are now prepared to state our main results.  For technical reasons we use $(1 + \eps)$ $K$-means. Let $\hat z$ denote the estimated clustering by applying $(1+ \eps)$ $K$-means to DC-MASE; i.e. $\hat z(i) =r $ if node $i$ is estimated to belong to community $r$.  Let $z$ denote the true clustering.  We define
\begin{align}
    \ell(\hat z, z) :&= \inf_{\text{Permutations } \mathcal{P}} \frac{1}{n} \sum_{i=1}^{n} \mathbb{I}\{ \hat z(i) \neq \mathcal{P}( z(i)) \}. \label{eq:misclusteringerrorrate}
\end{align}
In other words $\ell(\hat z,z)$ is the misclustering error up to label permutations.  The following theorem is our main technical result,  an upper bound on the  misclustering error.  

\begin{theorem} \label{thm:clusteringerror}
Suppose that \cref{ass:communityass} and  \cref{ass:networklevel} 
are satisfied, and suppose that $L \lesssim n^{5}$.  Define\begin{align}
    \mathrm{err}_{\ave}^{(i)} :&=  \frac{1}{L}\sum_{l} \frac{\|\theta\m\|_3^3}{\theta_i\m \| \theta\m\|^4 \lambda_{\min}\m} ; \label{eq:def-errave}\qquad 
    \mathrm{err}_{\max}^{(i)} := \max_{l} \frac{\theta_{\max}\m}{\theta_i\m \|\theta\m\|^2 (\lambda_{\min}\m)^{1/2}}. 
\end{align}
Then there exists a sufficiently small constant $c$ depending on the implicit constants in the assumptions such that the expected misclustering error is
\begin{align*}
    \mathbb{E} \ell( \hat z, z) &\leq \frac{2 K}{n} \sum_{i=1}^{n} \exp\Bigg( -c L \min\bigg\{ \frac{\bar \lambda^2}{K^4 \mathrm{err}_{\ave}^{(i)}}, \frac{\bar \lambda}{K^{2} \mathrm{err}_{\max}^{(i)}} \bigg\} \Bigg)  + O(n^{-10}).
\end{align*}
\end{theorem}
The assumption that $L \lesssim n^5$ is primarily for technical convenience; this is made so that we can take a union bound over all $L$ networks. If $L$ is larger but still polynomial in $n$, the result can still hold at the cost of increasing all of the implicit constants in the assumptions.  However, once $L$ is sufficiently large relative to $n$, the exponent can be made to be smaller than $e^{- c n}$ for some constant $c$, and hence Markov's inequality ensures that perfect community detection is possible. 
Therefore, while our theory only covers $L$ growing polynomially with $n$, for all practical purposes this assumption is irrelevant, as perfect clustering will be guaranteed once $L$ is larger than some polynomial of $n$.

\cref{thm:clusteringerror} makes precise the sense in which DC-MASE aggregates information across all of the networks.   In the bound there are two factors: one is the worst-case error for each network $\mathrm{err}_{\max}^{(i)}$, and one is the average-case error $\mathrm{err}_{\ave}^{(i)}$.  In order to further consider the rate of improvement relative to $L$, we also consider the following application in the regime that the signal strengths are comparable. 
\begin{corollary}[Network Homogeneity] \label{cor:homogeneous}
Instate the conditions of \cref{thm:clusteringerror}, and suppose that $\lambda_{\min}\l = \lambda_{\min}$ and $\theta_i\l = \theta_i$ for all $l$.  Then there exists a sufficiently small constant $c$ depending on the implicit constants in the assumptions such that
\begin{align*}
    \E \ell(\hat z, z) \leq \frac{2 K}{n} \sum_{i=1}^{n} \exp\bigg( -c L \theta_i \min\bigg\{ \frac{\|\theta\|^4 \lambda_{\min}^3}{K^4 \|\theta\|_3^3}, \frac{\|\theta\|^2 \lambda_{\min}^{3/2}}{K^{2}\theta_{\max}} \bigg\} \bigg) + O(n^{-10}).
\end{align*}
\end{corollary}


\cref{cor:homogeneous} further elucidates the sense in which DC-MASE aggregates information from multiple networks: the error rate includes a gain of $L$ but penalties of  $\lambda_{\min}$ (relative to which term is the minimizer in the rate).  In particular, if networks have extreme degree heterogeneity but well-separated communities, then the error rate for DC-MASE highly improves upon the corresponding rate for single networks.

In the homogeneous degree regime with $\theta_i\m \asymp \sqrt{\rho_n}$ for each $l$ as in \cref{ex:3}, this rate becomes $\exp\big( -c L n\rho_n \lambda_{\min}^3 \big) + O(n^{-10}),$ ignoring factors of $K$.
%
%
%
\citet{jin_improvements_2022} demonstrated that the SCORE clustering procedure with $L = 1$ yields the error rate of order $\exp( - c \lambda_{\min}^2 n\rho_n) + o(n^{-3})$.  
Therefore, we see that in this regime DC-MASE benefits whenever $\lambda_{\min} \gg \frac{1}{L}$,  \emph{even if each network is very sparse}. {However,  \cref{ass:networklevel} implies that we require that each network is sufficiently dense if the degrees are homogeneous, so in \cref{cor:homogeneous} we do not actually get to use the full strength of the exponent. 
 \cref{cor:homogeneous} still yields \emph{exponential} improvement with $L$ provided each network is sufficiently dense, and in \cref{sec:necessity} we demonstrate that an assumption similar to \cref{ass:networklevel} is inescapable in our setting. 
}

Our next result shows that under {an additional signal strength assumption 
DC-MASE yields perfect clustering with high probability. 
 
\begin{theorem}[Perfect Clustering] \label{cor:perfectclustering}
Suppose that the conditions of \cref{thm:clusteringerror} hold, and that 
\begin{align*}
\min_l  \frac{\theta_{\min}\l}{\theta_{\max}\l} \lambda_{\min}\l\|\theta\l\|^2 \geq C  \frac{K^{8}\log(n)}{L\bar \lambda^2}, \numberthis \label{snr} 
\end{align*}
 where $C$ is some sufficiently large constant. Then running $K$-means on the output of DC-MASE yields perfect recovery with probability at least $1 - O(n^{-9})$.  
\end{theorem}

\cref{cor:perfectclustering} demonstrates that if the layer-wise SNR is sufficiently strong relative to $\bar \lambda$, we achieve perfect clustering.  We note that \cref{ass:networklevel} already implies an assumption on the left hand side of  \eqref{snr} 
as well as imposing a lower bound on $\bar \lambda$.
If $(L \bar \lambda)\inv \lesssim \frac{\theta_{\max}\l\|\theta\l\|_1}{\|\theta\l\|^2\lambda_{\min}\l}$ for all $l$, then the condition in \cref{cor:perfectclustering} is already met.  Therefore, since the term $\frac{\theta_{\max}\l\|\theta\l\|_1}{\|\theta\l\|^2\lambda_{\min}\l}$ is always larger than one (by assumption), the condition in \cref{cor:perfectclustering} is only more stringent whenever $\bar \lambda \ll 1/L$, which can only happen in the moderate $L$ regime, since \cref{ass:networklevel} already imposes a lower bound on $\bar \lambda$.  At an intuitive level,  this condition further reflects the idea that the second SVD step may not perform as well when $\bar \lambda$ is small.  

 {In practice, the assumption of common community memberships across the layers may not hold exactly.  However, as the following result shows, as long as the fraction of nodes with different memberships is sufficiently small, the results continue to hold without significant modification. 
 \begin{theorem} \label{thm:newtheorem}
     Suppose that each network is given by $\bP^{(l)} = \bTheta^{(l)} \bZ^{(l)} \bB^{(l)} (\bZ^{(l)})\t \bTheta^{(l)}$, where the collection $\{\bZ^{(l)}\}$ satisfies $\max_{1\leq i\leq n} \frac{1}{L}\sum_{l=1}^{L} \mathbb{I}_{\{z_i^{(l)} \neq z_i\}} \leq \delta, $      where $z$ is the underlying ``ground truth'' communities.  
     If $\delta \ll \frac{\bar \lambda}{K}$ and $L \geq n/K^7$, then \cref{thm:clusteringerror} continues to hold.  
 \end{theorem}
 }

 {In words, \cref{thm:newtheorem} demonstrates that as long the fraction of networks that are different for a given vertex $i$ is at most $\delta$, then our main result remains unchanged. The maximal fraction $\delta$ is governed by the global signal strength $\bar \lambda$, and hence settings with more signal are permitted to have more ``errors'' (i.e., incorrect community assignments). Finally, while we assume that $L\geq n/K^7$ in the statement of the theorem, we believe this assumption to be a proof artifact.  }

\begin{remark}[Network outliers]
Consider the case that $L \delta_{\mathrm{out}}$ networks contribute no information at all or that \cref{ass:networklevel} is violated for these networks.  Let these networks be the \emph{outlier} networks, and let all other networks be \emph{inlier} networks.  Letting $\bar \lambda$ be defined only in terms of inlier networks, it is possible to show that the upper bound in \cref{thm:clusteringerror} continues to hold as long as $\delta_{\mathrm{out}} \ll \frac{\bar \lambda}{K}$, with the proviso that all of the quantities appearing on the right hand side of the misclustering error rate are replaced with inlier terms. We will not prove this to conserve space, as it is similar to the proof of \cref{thm:newtheorem}. 
\end{remark}

\subsection{Necessity of Individual-Network Signal Strength Condition} \label{sec:necessity}
\cref{ass:networklevel} imposes a minimal assumption on both the signal strength and degree homogeneity of each network.  When the networks have homogeneous degree corrections of order $\sqrt{\rho_n}$ and $\lambda_{\min}\l \asymp 1$, then this assumption is equivalent to the assumption that $n\rho_n \gtrsim \log(n)$.  In the setting of the multilayer stochastic blockmodel, it was shown in \citet{lei_bias-adjusted_2022} that a sufficient condition for  consistent community detection is that $\sqrt{L}n\rho_n \gg \sqrt{\log(n+L)}$ when $n\rho_n \lesssim 1$, which, to the best of our knowledge, is the weakest such condition from the literature with heterogeneous $\mathbf{B}^{(l)}$ matrices. Therefore, without heterogeneous degree corrections, \cref{ass:networklevel} may be stronger (by a factor of $\sqrt{L}$) than necessary.  In this section we study the necessity of this assumption in the presence of degree corrections.

Define the parameter space
\begin{align*}
    \mathcal{P}(\lambda_{\min},K,n,\theta,L) :&=   \bigg\{ \mathbf{P}^{(1)}, \cdots , \mathbf{P}^{(L)} \in [0,1]^{n\times n}: \mathbf{P}^{(l)} = \mathbf{\Theta}^{(l)} \mathbf{Z} \mathbf{B}^{(l)} \mathbf{Z}\t \mathbf{\Theta}^{(l)}; |\lambda_{\min}(\mathbf{B}^{(l)})| \geq \lambda_{\min};  \\
    &\qquad \max_{k \leq K} \frac{\|\theta^{(l)}_{\mathcal{C}_k}\|^2}{K} \leq C \min_{k\leq K }\frac{\|\theta^{(l)}_{\mathcal{C}_k}\|^2}{K}, c\frac{n}{K} \leq | \mathcal{C}_k | \leq C \frac{n}{K} \bigg\}.
\end{align*}
 Here $\theta = \{\theta_i\l\}$ and $L$ are allowed to depend on $n$, where for simplicity we focus on the regime that $K, \lambda_{\min}  \asymp 1$ to allow us to isolate the effect of degree correction parameters.   
The following result provides a lower bound on the estimation error.  
\begin{theorem} \label{thm:minimax}
Suppose that $K = O(1)$ and $L\lesssim n$.  Suppose further that there exists some constant $c_1$ such that $\theta$ satisfies
\begin{align*}
    c_1 \leq \frac{\theta_{\max}^{(l)} \|\theta^{(l)} \|_1 }{\|\theta^{(l)}\|^4 (\lambda_{\min}^{(l)})^2} \ll \sqrt{\frac{L}{\log(n+L)}}. \numberthis \label{leilintypething}
\end{align*}
{Assume further that  $\lambda_{\min}$ is fixed in $n$ and satisfies $\lambda_{\min} \geq c_2$ for some constant $c_2$.}
Then there exists some constant $c > 0$ such that $\inf_{\hat z} \sup_{\mathcal{P}(\lambda_{\min},K,n,\theta,L)} \mathbb{E} \ell(\hat z, z) \geq c. $
\end{theorem} 

{To ease intuition, we will compare the signal strength condition of Assumption~\ref{ass:networklevel} with the one of Theorem~\ref{thm:minimax} in the context of the two previous examples. Suppose the degree heterogeneous setting of Example~\ref{ex:2}, and $a/b\rightarrow 0$. Then the quantity in \cref{leilintypething} is given by 
$\frac{\theta_{\max}^{(l)} \|\theta^{(l)} \|_1 }{\|\theta^{(l)}\|^4 (\lambda_{\min}^{(l)})^2}= \frac{b(\gamma a +(1-\gamma)b)}{n(\gamma a^2+(1-\gamma)b^2)^2}.$ 
We will argue that \cref{leilintypething} and \cref{ass:networklevel} cannot simultaneously hold in the two extremes for $\gamma$.
\begin{itemize}
    \item \textbf{Case 1: $\gamma = 1/n$.} In this case the lower bound in \cref{leilintypething} gives $b\leq 1/\sqrt{n}$, in which case the signal strength assumption in Assumption~\ref{ass:networklevel} would require $a\geq \log(n)/\sqrt{n}>b$.  Since $a/b \to 0$, we see that \cref{ass:networklevel} cannot simultaneously hold. 
    \item \textbf{Case 2: $\gamma = 1 - 1/n$.}  When $\gamma=1-1/n$, the lower bound in \cref{leilintypething} gives $\frac{b((1-1/n)a+b/n)}{n((1-1/n)a^2+b^2/n)^2}=\frac{b((n-1)a+b)}{((n-1)a^2+b^2)^2}\geq c_1,$ so when $b=1$, $a\leq 1/n^{1/3}$ is needed in \cref{leilintypething}, as opposed to the condition $a\geq (\log(n)/n)^{1/4}$ given by the  Assumption~\ref{ass:networklevel} for this setting.
\end{itemize}
Now consider the degree homogeneous, weak-signal setting of \cref{ex:3}. The lower bound of \cref{leilintypething} is $\frac{\theta_{\max}^{(l)} \|\theta^{(l)} \|_1 }{\|\theta^{(l)}\|^4 (\lambda_{\min}^{(l)})^2}\asymp \frac{1}{n\rho_n \lambda_{\min}^2}\geq c_1,$ which implies that $\lambda_{\min}\lesssim 1/(n\rho_n)^{1/2}.$  Since $\lambda_{\min} \asymp 1$ by assumption, we see that we require that $n\rho_n \lesssim 1$, which is in contrast to \cref{ass:networklevel}, which essentially requires $n\rho_n \gtrsim \log(n)$. } 
Furthermore, in the homogeneous degree setting, the assumption \eqref{leilintypething} implies that $\sqrt{L}n\rho_n\gg \sqrt{\log(n+L)}$ but $n\rho_n \leq C$, which matches the sufficient condition from \citet{lei_bias-adjusted_2022}, so \cref{thm:minimax} can be understood as stating that the additional degree heterogeneity renders the problem significantly more difficult than its degree homogeneous counterpart. 
Therefore, 
\cref{thm:minimax} gives evidence that some minimal condition similar to \cref{ass:networklevel} is inescapable for consistent community detection in ML-DCSBMs.  
 Up to logarithmic terms and the factor of $\frac{\theta_{\max}\l}{\theta_{\min}\l}$, \cref{thm:minimax} shows that \cref{ass:networklevel} is necessary when $\lambda_{\min}\l \asymp 1$.

\subsection{Overview of the Proof of Theorem \ref{thm:clusteringerror}} \label{sec:proofoverview}
 This section gives a high-level overview and discusses the novelty of the proof of \cref{thm:clusteringerror}, though the full proof can be found in Appendix A. Our proof requires three key steps, each proved sequentially.

\textbf{Step 1: First Stage Asymptotic Expansion}. In Theorem A.1, we show that the initial estimates $\yhat\m$ satisfy
\begin{align*}
    \yhat\m \mathbf{W}_*\m - \mathbf{Y}\m &= \mathcal{L}(\mathbf{A}\m - \bP\m) + \mathcal{R}_{\mathrm{Stage \ I}}\m,
\end{align*}
where $\mathcal{L}(\cdot)$ is a linear function, $\mathcal{R}_{\mathrm{Stage \ I}}\m$ is a residual with small $\ell_{2,\infty}$ error, and $\wstar\m$ is a $K \times K$ orthogonal matrix. Unlike previous results of this type \citep{du_hypothesis_2021,fan_simple_2022}, our residual bounds depend explicitly on the degree corrections.  To prove these results we rely on the leave-one-out analysis technique established in \citet{abbe_entrywise_2020} and a Taylor expansion argument.

\textbf{Step 2: Second Stage $\sin\bTheta$ Perturbation Bounds}. 
We then prove Theorem A.2, applying Theorem A.1 to obtain concentration in $\sin\bTheta$ distance for the empirical singular vectors $\uhat$ to the true singular vectors $\U$ that reveal the community memberships.  In particular, by virtue of our first-order expansion, since $\mathcal{L}(\cdot)$ is linear in the noise, we are able to obtain stronger concentration for $\sin\bTheta$ distance than if one were to simply apply the na\"ive concentration using the triangle inequality, which would not yield improvement with $L$. This argument bears some resemblance to the concurrent work \citet{zheng_limit_2022}; however, in our analysis we also have to take into account the (nonlinear) transformation that normalizes the rows (i.e., projection to the sphere), and our second stage analysis requires several novel considerations for both the population and empirical versions of the algorithm.

\textbf{Step 3: Second Stage Asymptotic Expansion}.
The final step of our proof is our major technical contribution. In our final step we prove Theorem A.3, which shows that
\begin{align*}
    \uhat \mathbf{W}_* - \U &= \sum_{l} \mathcal{L}(\A\m - \bP\m ) (\ytilde\m)\t \U \Sigma^{-2} + \mathcal{R}_{\mathrm{Stage \ II}},
\end{align*}
where $\mathcal{R}_{\mathrm{Stage \ II}}$ is a smaller order term, $\mathcal{L}(\A\m - \bP\m)$ is the same linear operator as in the first step (Theorem A.1), and $\wstar$ is an orthogonal matrix.   To prove this result we use the asymptotic expansion established in the first step to obtain sharp concentration bounds for $\mathcal{R}_{\mathrm{Stage \ II}}$ in $\ell_{2,\infty}$ norm that takes into account the nonlinearity induced by the normalization procedure. The most similar work containing comparable technical results is the work \citet{zheng_limit_2022} analyzing the algorithm from \citet{arroyo_inference_2021} (which has no additional nonlinearity), and they make a number of simplifying assumptions that do not hold in our setting (such as that each network has comparable signal).  As an additional technical challenge, in our proofs we must also study the interplay between first and  second-order terms in a manner that is amenable to the different signal strengths within each network.

\section{Simulation Results}
\label{sec:sims}
We evaluate the performance of different methods for community detection in networks generated from the multilayer DCSBM. The experiments focus on the effect of the number of graphs $L$ for recovering the communities under different parameter setups. \footnote{An implementation of the code is available at \url{https://github.com/jesusdaniel/dcmase}}
 The performance measure reported in the experiments is the misclustering error rate as defined in \cref{eq:misclusteringerrorrate}, which is simply the proportion of nodes that are incorrectly clustered.

The benchmarks considered include spectral-based, optimization-based and likelihood-based clustering algorithms for multilayer networks. For spectral methods, the list comprises clustering on the embeddings defined as (i) the leading eigenvectors of  the aggregated sum of the adjacency matrices $\sum_l \bA\m$ \citep{han_consistent_2015,bhattacharyya2020consistent},  (ii) the leading eigenvectors of the bias-adjusted sum-of-squared (SoS) adjacency matrices of \cite{lei_bias-adjusted_2022}, and (iii) an estimate of the common invariant subspace of the adjacency matrices obtained via multiple adjacency spectral embedding (MASE) from
\cite{arroyo_inference_2021}. Existing methods and theoretical results for multilayer community detection with the aforementioned embedding procedures typically consider $K$-means clustering on the rows of these embeddings to obtain communities, but this clustering scheme is not expected to work well under high degree heterogeneity even for a single network. Thus, to isolate the performance of the embedding from the clustering method adopted, we employed spherical spectral clustering by normalizing the rows of the embeddings before performing $K$-means clustering \citep{lei_consistency_2015,bhattacharyya2020consistent}, as we observed better empirical performance compared to the unnormalized version.
We also consider the orthogonal linked matrix factorization (OLMF) of
\cite{paul_spectral_2020}, and an optimized
Monte Carlo Markov Chain approach \citep{peixoto2014efficient,peixoto2015inferring} implemented via the graph-tool package \citep{peixoto_graph-tool_2014}.

All the simulated graphs are generated using the multilayer DCSBM with $n=150$ vertices and $K=3$ equal sized communities, for which we assume that the membership matrix $\bZ$ is such that vertices in the same community have adjacent rows. We focus on studying the effect of number of graphs $L$ in the presence of different types of  parameter heterogeneity. For that goal, we consider scenarios in which the block connectivity matrices or the degree correction parameters are  the same or different across the collection of graphs. For the block connectivity matrices, we generate these parameters as follows:
\begin{itemize}
    \item \emph{Same connectivity matrices:} the matrices $\bB\m,l\in[L]$ are all set to be equal and defined as
    $\bB\m_{rr}=1, r\in[K]$, and $\bB\m_{rs}=0.4$, $r\neq s$.
    \item \emph{Different connectivity matrices:} each  $\bB\m, l\in[L]$ is generated independently with entries equal to $\bB\m_{rr}=p\m \sim\text{Unif}(0,1)$, for $r\in[K]$, and  $\bB\m_{rs}=q\m \sim\text{Unif}(0,1)$, $r\neq s$.
\end{itemize}
In terms of the degree correction parameters, we consider scenarios as follows:
\begin{itemize}
    \item \emph{Same degree corrections:} the diagonal entries of the matrices satisfy $\bTheta\m_{ii}=\theta_i$ and are generated from a shifted exponential distribution such that $\theta_1, \ldots, \theta_n\overset{\text{i.i.d.}}{\sim} \text{Exp}(1) + 0.2$.
    \item \emph{Different degree corrections:} the parameters are generated in a similar way, but now each network has its own parameters $\theta\m_1, \ldots, \theta\m_n\overset{\text{i.i.d.}}{\sim} \text{Exp}(1) + 0.2$..
    \item \emph{Alternating degrees:} the vertices within each community are split into two equal sized groups, and each group alternates between having low and high degrees on each network, that is, $\theta_i\m=0.8$ if either $l$ and $i$ are odd or $l$ and $i$ are even numbers, and $\theta_i\m=0.15$ otherwise.
\end{itemize}
The expected adjacency matrices are then defined as $\bP\m = \alpha\m\bTheta^{(l)}\bZ\bB^{(l)}\bZ\bTheta^{(l)}$ similar to Eq.~\eqref{eq:model}, and the constant $\alpha\m$ is introduced to keep the average expected degree equal to 10. For each parameter setup, the experiments are repeated 100 times, and the average results are reported.

The results are shown in Figure~\ref{fig:experiment-increasing-L}. As expected, the accuracy of the  methods generally improves with more graphs, and although there is no specific method that dominates in all the scenarios considered, we observe that DC-MASE is the only one that consistently improves its performance with $L$ until perfect clustering is achieved. When the degree correction parameters are the same (left column), most of the methods perform accurately, especially in the setting with the same connectivity matrices. In particular, spectral methods  perform  well due to the fact that the singular subspace is shared in the expected adjacency matrices, and the population version of the matrix in which the embedding is performed captures the community structure after further correcting  for degree heterogeneity via spherical normalization.  In the scenario with different but random degree corrections (middle column) several methods are still able to perform accurately even when the population matrix does not have the correct clustering structure, possibly due to an averaging effect of the degree-correction parameters generated independently at random for each graph. 
   Aggregation methods, such as the sum of the adjacency matrices, perform very well when the global structure of the graphs is the same, but are not able to identify the correct structure in the presence of severe parameter heterogeneity.  Notably, in the alternating degrees scenario (right column), DC-MASE is the only method that performs accurately, whereas other methods struggle to identify the model communities. 

\begin{figure}
    \centering
    \includegraphics[width=.7\textwidth,keepaspectratio]{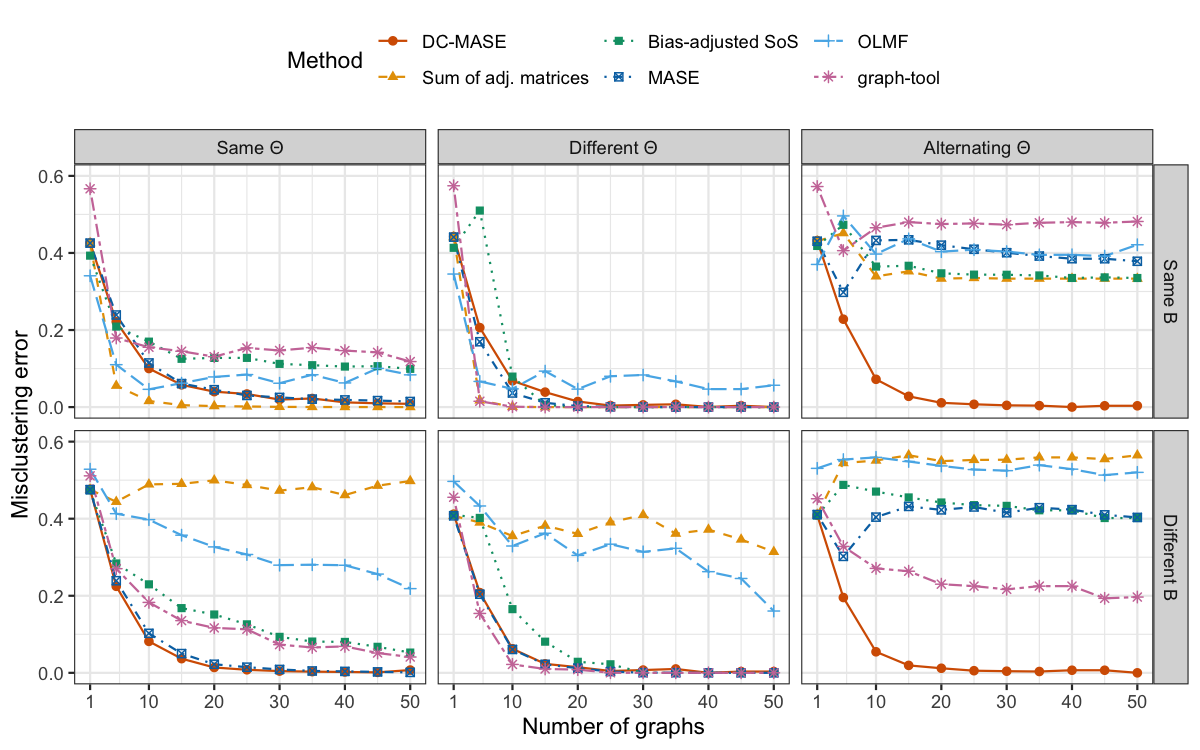}
    \caption{Community detection error of different methods (measured via {misclustering error}, averaged over 100 replications) as a function of the number of graphs. See \cref{sec:sims} for a discussion of the setups.}
    \label{fig:experiment-increasing-L}
\end{figure}

\section{Analysis of US Airport Network} \label{sec:realdata}

We evaluate the performance of the method in a time series of networks encoding the number of flights between airports in the United States within a given month for the period of January 2016 to September 2021. A multilayer degree-corrected SBM allows us to track the flight dynamics both at the airport and community levels to characterize the effect of the Covid 19 pandemic in flight connectivity. The data are publicly available and were downloaded from the US Bureau of Transportation Statistics \citep{usbts}.

The vertices of the networks correspond to some of the airports located within the 48 contiguous states in the US. For each network, the weighted edges contain the total number of flights  of class F (scheduled passenger/cargo service) between each pair of airports within a given month. We restricted the analysis to the vertices in the intersection of the largest connected components of all the networks, resulting in a total of $n=343$ airports. The period of the study contains 69 months (number of graphs). 

To identify communities of airports with similar connectivity patterns in the data, we apply DC-MASE to the collection of adjacency matrices. The number of communities was selected to be $K=4$ to facilitate interpretation and based on the scree plots of the individual network embeddings and
the concatenated matrix, as described in \cref{sec:estimate-K}. Figure~\ref{fig:airport} (left) shows the estimated community memberships of the airports. Three of the communities identified (communities 2, 3 and 4) appear to be related to the geographical area, (west, east and southwest, respectively), whereas community 1 contains most of the hub airports in the east side of the country, as well as other smaller airports that are mostly connected to these hubs. 

To characterize the dynamics in community and airport connectivity, we estimate the block connectivity matrices and degree correction parameters of the multilayer DCSBM. As the edges count the total number of flights between pairs of locations, the adjacency matrices are weighted, and thus, the parameters of the model describe the expected adjacency matrix $\mathbb{E}[\bA\m] = \bTheta\m\bZ\bB\m\bZ^\top\bTheta\m$. For ease of interpretation, we adopt a similar identifiability condition as in \cite{karrer2011stochastic} by constraining the sum of the degree correction parameters within each community to be equal to the size of the community, that is, if vertex $i$ is in community $r$ then
$\sum_{i\in\mathcal{C}(r)} \theta_i\m = |\mathcal{C}(r)|$ for $r\in[K], l\in[L].$
With this parameterization, we have the following relations. Let $d_i\m = \sum_{j=1}^n\e[\bA\m_{ij}]$ be the expected degree of node $i$ in network $l$. Then, for every $i\in[n]$, $r,s\in[K]$ and $l\in[L]$ we have
\begin{equation}
\theta\m_i =  \frac{d\m_i}{\frac{1}{|\mathcal{C}(r)|}\sum_{j\in\mathcal{C}(r)}d\m_j},\quad\quad\quad\bB\m_{rs} = \frac{1}{|\mathcal{C}(r)|\ |\mathcal{C}(s)|} \sum_{i\in\mathcal{C}(r), j\in\mathcal{C}(s)} \e[\bA\m_{ij}].
    \label{eq:plug-in-parameters}
\end{equation}
Under this parameterization, the degree correction parameters are on average equal to 1, and large values can be interpreted as higher individual connectivity of the corresponding vertex relative to other vertices in the community. Meanwhile, the block connectivity simply calculates the average number of edges within and between each pair of communities. When comparing the values of these parameters across time, this parameterization allows us to split global and local dynamics into the block connectivity matrices and degree corrections, respectively.
We obtain  plug-in estimates of the model parameters by using $\bA\m$ rather than $\e[\bA\m]$, and by using the estimated community memberships, which under certain edge distributions (e.g. Poisson) coincides with the maximum profile likelihood estimates given the fitted community memberships.

The multilayer DCSBM estimated parameters shown in Figure~\ref{fig:airport} (right panel) track the changes in airport connectivity at the community level, which are mostly related to regional dynamics. In contrast, Figure~\ref{fig:airport-theta} (left panel) also shows the individual airport popularity relative to airports within its community over time.  While the overall number of flights within and between communities decreased after the pandemic started, the impact on the airport traffic was not homogeneous, and this is captured by the changes in degree correction parameters. Figure~\ref{fig:airport-theta} (right) explores these changes in more detail for community 1, which includes some of the largest hubs, such as ATL, DFW and ORD. These became relatively more prominent with respect to other airports in their community at the start of the pandemic in the US. Meanwhile, the airports in the New York City area (EWR and LGA) were relatively more negatively affected, possibly due to the pandemic dynamics and related closures. This analysis illustrates the flexibility of the multilayer DCSBM model for tracking local and community-level dynamics with changes over time. 

{An additional analysis comparing the communities discovered by DC-MASE with the other spectral clustering algorithms is included in Appendix J of the Supplementary Material. In the absence of ground truth communities, the performance is measured via out-of-sample edge prediction accuracy. The results generally favor the communities discovered by DC-MASE, suggesting a better generalization error.}

\begin{figure}[ht!]
    \includegraphics[width=0.45\textwidth,keepaspectratio]{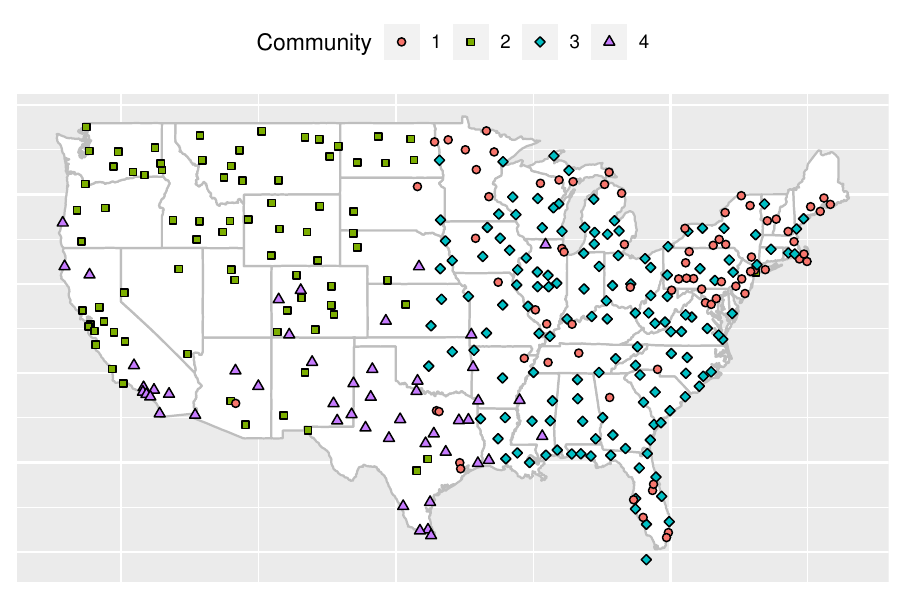}
    \hfill
    \includegraphics[width=0.45\textwidth,keepaspectratio]{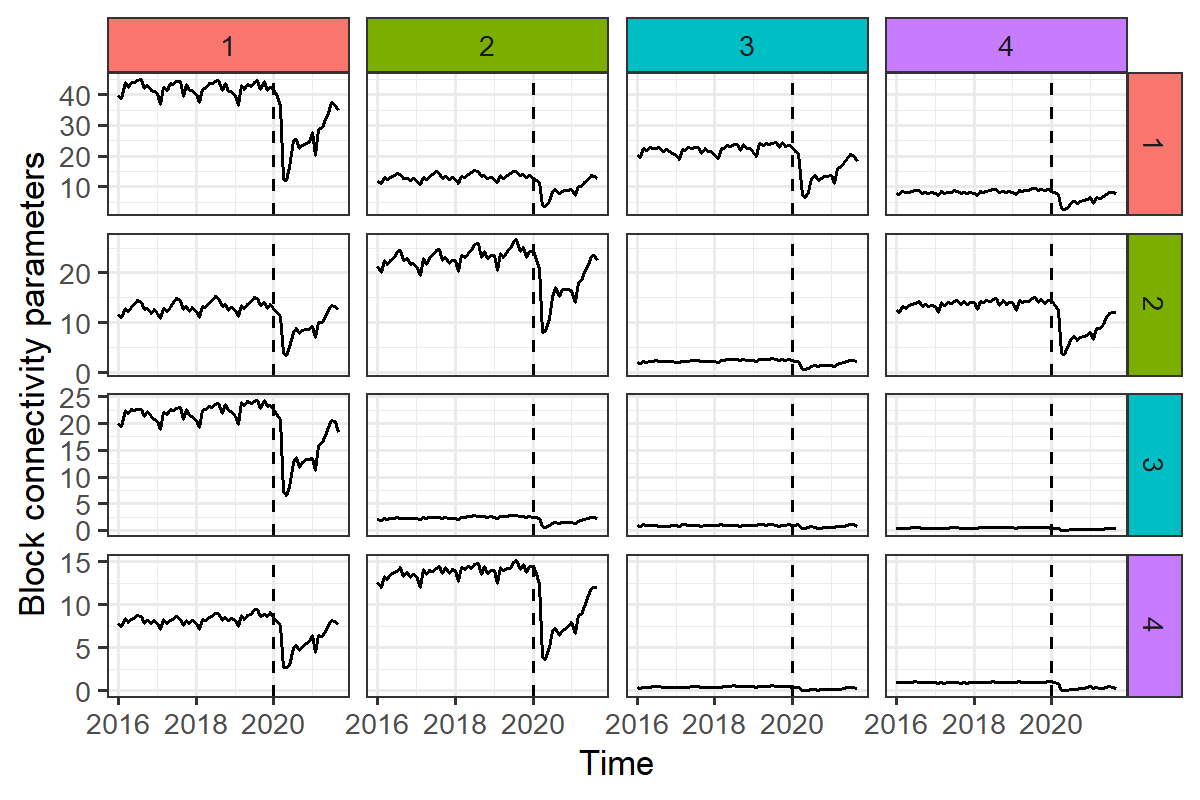}
    \caption{Map of US airports colored according to the communities discovered by DC-MASE (left) and time series  of the corresponding estimated block connectivity matrices in the model (right). Each cell in this plot represents an entry of these matrices over time. The vertical line indicates January 1st, 2020.}\label{fig:airport}
\end{figure}

\begin{figure}[ht!]
    \includegraphics[width=0.48\textwidth]{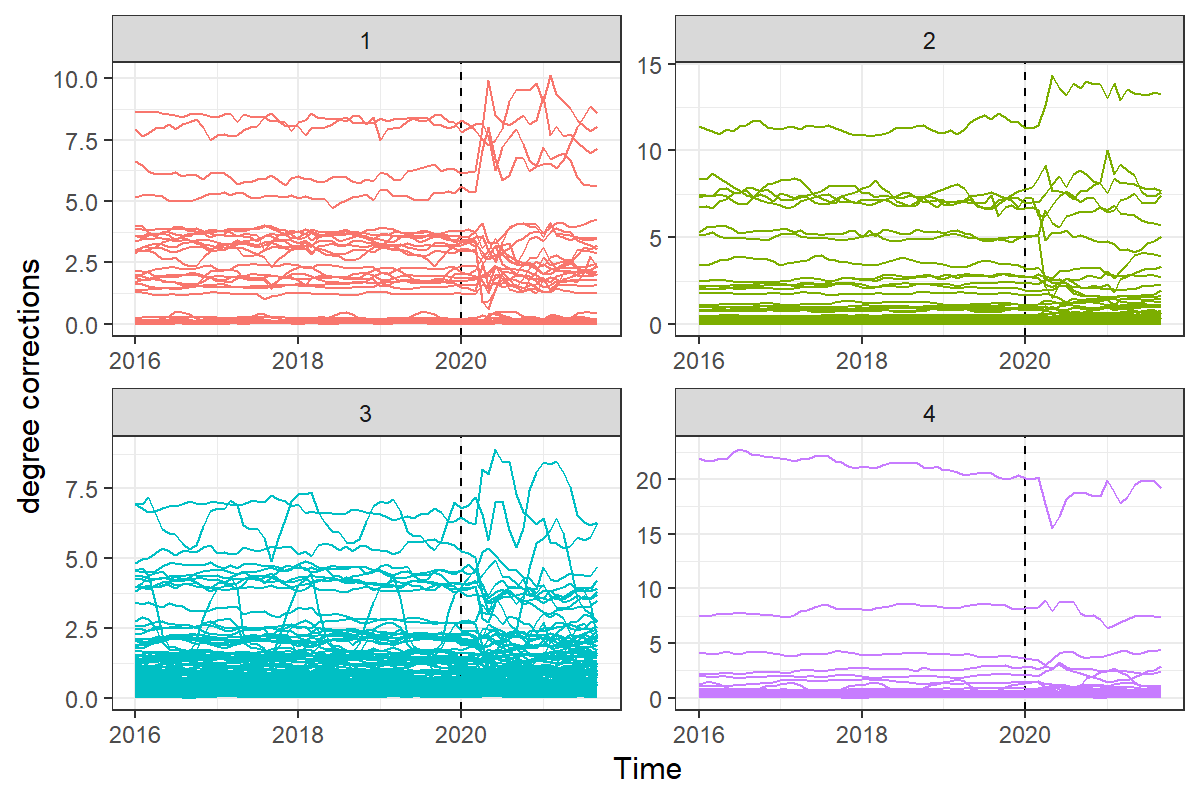}
    \includegraphics[width=0.48\textwidth]{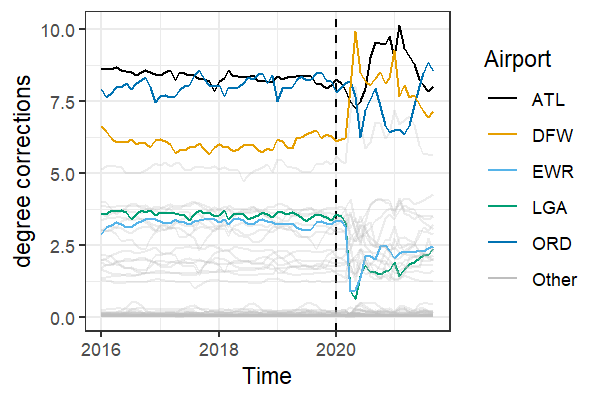}
    \caption{Degree correction parameter estimates in the US airport data divided by communities (left).  Each line corresponds to the parameter for some specific airport over time; the collection is divided according to the communities discovered by the algorithm.
The variability in the parameter estimates suggests the need for different degree correction parameters at each time point. In the right panel,  the results for community 1 are zoomed in, with some  major airports highlighted.   }\label{fig:airport-theta}
\end{figure}

\section{Discussion}
\label{S:discussion}
In this work we have considered the multilayer degree-corrected stochastic blockmodel, established its identifiability, and proposed a joint spectral clustering algorithm based on clustering the rows of a matrix that appropriately aggregates information about the communities in the model. The proposed method is simple and efficient, 
while the most expensive computations (required to estimate the leading eigenvalues and eigenvectors of each network) are able to be performed in parallel. This allows the methodology to scale to large datasets, both in terms of network size and in the number of graphs or layers.
Our main results demonstrate that the method can effectively leverage the information across the graphs to obtain an improvement in community estimation, particularly when the number of networks $L$ is large, even in the presence of significant vertex and layer heterogeneity. 
In our simulations, we observe that clustering with DC-MASE performs consistently well in various scenarios, and it is competitive with other state-of-the-art methods for multilayer community detection, particularly in situations with extreme degree heterogeneity.  In our flight data studies, we see that the multilayer DCSBM is a flexible but succinct model, allowing us to identify clusters, track degree corrections, and observe block connectivity over time.

Finally, while the multilayer DCSBM is a flexible model, our main results require an assumption on the amount of degree heterogeneity and signal strength within each network.  
The recent work \citet{ke_optimal_2022} demonstrates that the eigenvectors of the \emph{regularized Laplacian} can yield optimal mixed-membership estimation under extreme degree heterogeneity; it would be interesting to study the multilayer DCSBM in this regime.





\section*{Acknowledgements}
Joshua Agterberg acknowledges support from a fellowship from the Johns Hopkins Mathematical Institute of Data Science (MINDS) via its NSF TRIPODS award CCF-1934979, the Charles and Catherine Counselman Fellowship, and the Acheson J. Duncan Fund for the Advancement of Research in Statistics. Jes\'us Arroyo acknowledges support from the National Science Foundation under grant DMS-2413553.

\appendix

\appendix
\section{Proof Ingredients and Proof of Theorem~\ref{thm:clusteringerror}} 
\label{sec:proofingredients}
 This section elaborates on the informal results stated in \cref{sec:proofoverview}.  Recall that
%
%
 %
 we let $\yhat\l$ be defined in \cref{alg:dcmase}, and we let $\ytilde\l$ denote the corresponding matrix associated to the population matrix $\bP\l$.  We also recall $\yhatcal = [\yhat\one, \dots, \yhat^{(L)}]$, and we let $\ycal = [\ytilde\one, \dots, \ytilde^{(L)}]$.  Finally, we let $\U$ and $\Sigma$ denote the leading $K$ left singular vectors and singular values of $\ycal$, and we let $\uhat$ and $\hat \Sigma$ be defined similarly. For simplicity of notation, we assume that $\widehat{\bZ}$ and $\widehat{z}$ satisfy
 $$\|\widehat{\bZ}-\bZ\|_F = \min_{\bP}\|\widehat{\bZ} - \bZ\bP\|_F,$$
 $$ \sum_{i=1}^{n} \mathbb{I}\{ \hat z(i) \neq  z(i)\} = \min_{\mathcal{P}} \sum_{i=1}^{n} \mathbb{I}\{ \hat z(i) \neq \mathcal{P}( z(i)) \},$$
 where the minimum is taking among all permutations $\mathcal{P}$ and permutation matrices $\mathbf{P}$.

\subsection{First Stage Characterization}

In the first step of the proof, we derive the following asymptotic expansion result for the individual networks. Recall that $\xhat\m$ and $\xtilde\m$ denote the scaled eigenvectors of  $\A\m$ and $\E \A\l = \bP\m$, respectively, and we let $\uhat\l$ and $\U\l$ be the leading $K$ eigenvectors of $\A\m$ and $\bP\l$ respectively. 
We let $\ipq\m$ denote the diagonal matrix with elements $\pm 1$, where $1$ appears $p$ times and $-1$ appears $q$ times, with $p$ corresponding to the number of positive eigenvalues of $\bP\l$ and $q$ corresponding to the number of negative eigenvalues of $\bP\l$.  Equivalently, $p$ and $q$ count the number of positive and negative eigenvalues of $\mathbf{B}\m$.  We let $\Lambda\m$ denote the nonzero eigenvalues of $\bP\m$, and $\hat \Lambda\m$ denote the leading $p$ positive and $q$ negative eigenvalues of $\A \m$, arranged in decreasing order by magnitude after splitting according to positive and negative.   

The following result characterizes the rows of $\yhat\m$.

\begin{restatable}[Asymptotic Expansion: Stage I]{theorem}{firststep}
\label{thm:firststep} 
Suppose that \cref{ass:communityass} and \cref{ass:networklevel} hold. Fix a given $l \in [L]$.  Let $\wstar\l$ denote the orthogonal matrix satisfing
\begin{align*}
    \wstar\l :&= \argmin_{\mathbf{W} \in \mathbb{O}(K)} \| \uhat\l - \U\l \wstar\l \|_F.
\end{align*}
Then there is an event $\mathcal{E}_{\mathrm{Stage \ I}}\m$ with $\p(\mathcal{E}_{\mathrm{Stage \ I}}\m) \geq 1 - O(n^{-15})$ such that the following expansion holds:
\begin{align*}
     \yhat\m (\wstar\m)\t-  \ytilde\m &=  \mathcal{L}(\A\m - \bP\l ) +  \mathcal{R}_{\mathrm{Stage \ I}}\m,
\end{align*}
where the matrix $\mathcal{R}_{\mathrm{Stage \ I}}\m$ satisfies
\begin{align*}
    \|  \mathcal{R}_{\mathrm{Stage \ I}}\m \|_{2,\infty} &\lesssim  \frac{K^2 \theta_{\max}\m \|\theta\m\|_1}{\lambda_{\min}\m \|\theta\m\|^4} \bigg( \log(n) \frac{\theta_{\max}\m}{\theta_{\min}\m} + \frac{{\sqrt{K}}}{\lambda_{\min}\m} +\bigg(\frac{\theta_{\max}\m}{\theta_{\min}\m} \bigg)^{1/2}\frac{K^{5/2}\log(n)}{ (\lambda_{\min}\m)^{1/2}} \bigg),
\end{align*}
and the matrix $\mathcal{L}(\A\m - \bP\l )$ has rows given by
\begin{align*}
   \mathcal{L}(\A\m - \bP\l  )_{i\cdot} &= \frac{1}{\|\xtilde_{i\cdot}\m \|} \bigg( \mathbf{I} - \frac{ \xtilde_{i\cdot}\m( \xtilde_{i\cdot}\m)\t }{\|\xtilde_{i\cdot}\m\|^2} \bigg) \bigg( \big( \A\m -  \bP\l  \big) \U\m |\Lambda\m|^{-1/2} \ipq\m \bigg)_{i\cdot}.
\end{align*}
\end{restatable}
Explicitly, \cref{thm:firststep} provides an entrywise expansion for the rows of $\yhat\l$ about their corresponding population counterparts, up to the orthogonal transformation most closely aligning $\uhat\l$ and $\U\l$.  

We remark briefly how \cref{thm:firststep} is related to and generalizes several previous results for single network analysis. In \cite{du_hypothesis_2021}, the authors consider the rows of $\yhat$ to test if $\mathbf{Z}_{i\cdot} = \mathbf{Z}_{j\cdot}$ (under a mixed-membership model). To prove their main result, they establish a similar asymptotic expansion to \cref{thm:firststep}.  Our asymptotic linear term is the same as theirs, but our residual term exhibits a much finer characterization of the dependence on degree correction parameters, as they implicitly assume that $\theta_{\max} \asymp \theta_{\min},$ whereas we allow significant degree heterogeneity and extremely weak signals (\cite{du_hypothesis_2021} also implicitly assume that $\lambda_{\min}\l \asymp 1$).  Similarly, \cite{fan_simple_2022} consider the asymptotic normality of rows of the SCORE-normalized eigenvectors for testing equality of membership in degree-corrected stochastic blockmodels.  However, they also require that $\theta_{\max} \asymp \theta_{\min}$, which again eliminates the possibility of severe degree correction.  Moreover, our results also allow $K$ to grow and $\lambda_{\min}$ to shrink to zero sufficiently slowly, provided this is compensated for elsewhere in the signal strength, and previous results require much stronger conditions on these parameters.  Finally, a similar asymptotic expansion (with explicit degree corrections and dependencies) was used implicitly to prove the main result in \cite{jin_improvements_2022}, albeit for the SCORE normalization (as opposed to spherical normalization).   Therefore, our results complement theirs by providing an analysis of the spherical normalization often used in practice, and our result exhibits slightly different dependence on degree corrections and ccommunity separation.  We will also apply \cref{thm:firststep} in the proof of \cref{thm:singlenetwork} in \cref{sec:singlenetwork}, and we provide a detailed comparison of our assumptions to \citet{jin_improvements_2022} therein.

The following result will be used as an intermediate bound in the proof of \cref{thm:singlenetwork}, demonstrating a concentration inequality for $\|\yhat\l - \ytilde\l\wstar\l\|_{2,\infty}$.
\begin{corollary} \label{cor:step1twoinfty}
With probability at least $1 - O(n^{-15})$, it holds that 
\begin{align*}
     \| \yhat\m - \ytilde\m \wstar\m \|_{2,\infty} &\lesssim \bigg( \frac{\theta_{\max}\l}{\theta_{\min}\l} \bigg)^{1/2} \frac{K\sqrt{\log(n)}}{\|\theta\l\|(\lambda_{\min}\l)^{1/2} } \\
     &\quad + \frac{K^2 \theta_{\max}\m \|\theta\m\|_1}{\lambda_{\min}\m \|\theta\m\|^4} \bigg( \log(n) \frac{\theta_{\max}\m}{\theta_{\min}\m} + \frac{{\sqrt{K}}}{\lambda_{\min}\m} +\bigg(\frac{\theta_{\max}\m}{\theta_{\min}\m} \bigg)^{1/2}\frac{K^{5/2}\log(n)}{ (\lambda_{\min}\m)^{1/2}} \bigg).
\end{align*}
\end{corollary}
The proof follows from \cref{lem:lineartermtwoinfty} (see \cref{sec:firststageproofs}) and \cref{thm:firststep}.  

\subsection{Second Stage Characterization I: \texorpdfstring{$\sin\bTheta$}{sin(Theta)} Bound}

With the strong upper bounds for the first stage in \cref{thm:firststep}, we can apply this result to establish $\sin\bTheta$ perturbation for the output of DC-MASE. 

For convenience we will define the following signal-to-noise ratio parameter vector
\begin{align}
   \snr_l :&= \bigg( \frac{\theta_{\min}\l}{\theta_{\max}\l} \bigg)^{1/2} (\lambda_{\min}\l)^{1/2} \|\theta\l\|. \label{eq:def-snrl}
\end{align}
We will denote $\snr\inv$ as the entrywise inverse of the $\snr$ vector. 
When $\theta_{\max} \asymp \theta_{\min} \asymp \sqrt{\rho_n}$, it holds that $\snr_l \asymp  \sqrt{\lambda_{\min}\l n\rho_n}$.

\begin{restatable}[$\sin\bTheta$ Perturbation Bound]{theorem}{sintheta} \label{thm:step2sintheta} Suppose the conditions in \cref{thm:clusteringerror} hold. Define
\begin{align*}
    \alpha_{\max} &= \frac{K^2 \theta_{\max}\m \|\theta\m\|_1}{\lambda_{\min}\m \|\theta\m\|^4} \bigg( \log(n) \frac{\theta_{\max}\m}{\theta_{\min}\m} + \frac{{\sqrt{K}}}{\lambda_{\min}\m} +\bigg(\frac{\theta_{\max}\m}{\theta_{\min}\m} \bigg)^{1/2}\frac{K^{5/2}\log(n)}{ (\lambda_{\min}\m)^{1/2}} \bigg); 
\end{align*}
i.e., $\alpha_{\max}$ is the residual upper bound from \cref{thm:firststep}.  
Then with probability at least $1 - O(n^{-10})$, it holds that
\begin{align*}
  \| \sin\bTheta(\uhat,\U) \| &\lesssim    K^2 \sqrt{\log(n)} \frac{\big( \frac{1}{L} \| \snr\inv \|_2^2\big)^{1/2}}{\sqrt{L}\bar\lambda} + K^3 \log(n) \frac{\| \snr\inv \|_{\infty}^2}{\bar \lambda}  \\
    &\quad + K^2 \sqrt{ \log(n)} \frac{\alpha_{\max} \| \snr\inv \|_{\infty}}{\bar \lambda} + \frac{K \alpha_{\max}}{\bar \lambda}.
\end{align*}
In particular, under the conditions of \cref{thm:clusteringerror} it holds that
\begin{align*}
    \|\sin\bTheta(\uhat,\U) \| &\lesssim \frac{1}{K}.
\end{align*}
\end{restatable}
\noindent
We note that the first bound provided in \cref{thm:step2sintheta} may actually be much stronger than the upper bound of $\frac{1}{K}$, which is all that is needed for the proof of \cref{thm:clusteringerror}. First, by combining \cref{ass:networklevel} and the definition of $\snr_l$ in Equation~\ref{eq:def-snrl} it is straightforward to check that each term is smaller than one, since 
we require that
\begin{align*}
   C \frac{K^{8} \theta_{\max}\l \|\theta\l\|_1 \log(n)}{\|\theta\l\|^2 \snr_l^2} \leq  \bar \lambda,
\end{align*}
for some large constant $C$.  
Since $\frac{\theta_{\max}\l \|\theta\l\|_1}{\|\theta\l\|^2}$ is always larger than one, we see that \cref{ass:networklevel} is a stronger assumption than each term in \cref{thm:step2sintheta} being smaller than one.

For ease of interpretation, when all $l$ have $\lambda_{\min}\l \asymp 1$, and $\theta_{\max}\l \asymp \theta_{\min}\l \asymp \sqrt{\rho_n}$ and $K \asymp 1$, we have that
\begin{align*}
    \| \snr\inv \|_{\infty} &\lesssim \frac{1}{\sqrt{ n\rho_n}}; \\
    \alpha_{\max} &\lesssim \frac{\log(n)}{n\rho_n} ; \\
     \frac{1}{L} \| \snr\inv \|_2^2 &\lesssim \frac{1}{n\rho_n}.
\end{align*}
Therefore, the $\sin\bTheta$ bound simplifies to
\begin{align*}
    \| \sin\bTheta(\uhat, \U) \| &\lesssim 
    \frac{\sqrt{\log(n)}}{ \sqrt{Ln\rho_n}} + \frac{\log(n)}{n\rho_n}.
\end{align*}
This final bound shows that $\uhat$ concentrates in $\sin\bTheta$ distance about $\U$ as $n$ increases by a factor that improves with $\sqrt{L}$ when $L\lesssim n\rho_n/\log(n)$. For a single stochastic blockmodel without degree corrections, the $\sin\bTheta$ distance between $\uhat$ and $\U$ can be upper bounded as $\sqrt{\frac{\log(n)}{  n\rho_n}}$ \citep{lei_consistency_2015}.  Therefore, \cref{thm:step2sintheta}, which utilizes the information from all the networks and allows degree heterogeneity, already demonstrates improvement from multiple networks by a factor of $\max\{\frac{1}{\sqrt{L}},\frac{\sqrt{\log(n)}}{ \sqrt{n\rho_n}}\}$ relative to the single-network setting.  However, it is important to emphasize that a primary benefit of this second stage aggregation is to ameliorate degree heterogeneity, which is not reflected in the homogeneous degree setting.  

\subsection{Second Stage Characterization II: Asymptotic Expansion} \label{sec:secondstageasympexp}
 In essence, we require \cref{thm:step2sintheta} to demonstrate that the clusters are correctly identified (see the proof of \cref{thm:clusteringerror} in \cref{sec:mainresultproofs}), but it falls short of providing a fine-grained characterization for the rows of $\uhat$, which is what is needed for the exponential error rate.

The following result demonstrates a first-order asymptotic expansion for the singular vectors in the second stage of our algorithm. The proof is given in \cref{sec:secondstageasympexp}.


\begin{restatable}[Asymptotic Expansion: Stage II]{theorem}{asymptoticexpansion} \label{thm:step2asympexp}
Suppose the conditions of \cref{thm:clusteringerror} hold.  Define
\begin{align*}
    \wstar :&= \argmin_{\mathbf{W} \in \mathbb{O}(K)} \| \uhat - \U \mathbf{W} \|_F.
\end{align*}
There is an event $\mathcal{E}_{\mathrm{Stage \ II}}$ satisfying $\p\big( \mathcal{E}_{\mathrm{Stage \ II}} \big) \geq 1 - O(n^{-10})$ such that on this event, we have the asymptotic expansion 
\begin{align*}
    \uhat \wstar\t - \U &=  \sum_{l}\mathcal{L}(\A\m - \bP\m ) ( \ytilde\m)\t \U \Sigma^{-2} + \mathcal{R}_{\mathrm{Stage \ II}},
\end{align*}
where $\mathcal{L}(\cdot)$ is the operator from \cref{thm:firststep} and the residual satisfies
\begin{align*}
    \|  \mathcal{R}_{\mathrm{Stage \ II}} \|_{2,\infty} &\lesssim  \frac{K^3 \sqrt{\log(n)}}{nL \bar \lambda} \|\snr\inv\|_2 + \frac{K^4 \log(n)}{L^2 \sqrt{n} \bar \lambda^2} \|\snr\inv\|_2^2 \\
&\quad + \frac{K^{7/2} \log(n)}{\sqrt{n} \bar \lambda} \| \snr\inv\|_{\infty}^2 + \frac{\alpha_{\max}}{\sqrt{n} \bar \lambda}. 
\end{align*}
Here $\alpha_{\max}$ is as \cref{thm:step2sintheta}.  In particular, under the assumptions of \cref{thm:clusteringerror}, it holds that
\begin{align*}
      \|  \mathcal{R}_{\mathrm{Stage \ II}} \|_{2,\infty} &\leq \frac{1}{16 \sqrt{n_{\max}}}.
\end{align*}
\end{restatable}

\cref{thm:step2asympexp} establishes a first-order expansion for the rows of the difference matrix $\uhat \wstar\t - \U$, which is the main technical tool required to establish \cref{thm:clusteringerror}.  The proof of \cref{thm:step2asympexp} relies on both \cref{thm:firststep} and \cref{thm:step2sintheta}, but requires a number of additional considerations to bound the residual term $\mathcal{R}_{\mathrm{Stage \ II}}$ in  $\ell_{2,\infty}$ norm.

\subsection{Proof of Theorem~\ref{thm:clusteringerror} and Theorem~\ref{cor:perfectclustering}} 
\label{sec:mainresultproofs}

With all of these ingredients in place, we are nearly prepared to prove \cref{thm:clusteringerror}.  In the proof we will also require several results concerning the population parameters, which we state in the following two lemmas. 
The proofs can be found in \cref{sec:identifiabilityproofs}.  

\begin{restatable}[Population Properties: Stage I]{lemma}{populationspectralproperties} \label{lem:popprop} 
Suppose \cref{ass:communityass} holds, and let $\lambda_r\l(\bP\l)$ denote the eigenvalues of $\bP\l$ and let $\lambda_r(\bB\l)$ denote the eigenvalues of $\bB\l$.  Then for all $1 \leq r \leq K$,
\begin{align*}
     \theta_i\l 
     \lesssim \| \xtilde\m_{i\cdot} \| &\lesssim  \theta_i\m {\sqrt{K}} ;\\
    \| \U\m_{i\cdot} \| &\lesssim \sqrt{K} \frac{\theta\m_i}{\|\theta\m\|}; \\
    \lambda_r\l(\bP\l) &\asymp \frac{\|\theta\m\|^2}{K} \lambda_{r}(\mathbf{B}\m). 
\end{align*}
\end{restatable}

Next, the following result establishes the population properties of the the second stage; in particular demonstrating a lower bound on the smallest eigenvalue of the population matrix $\ycal\ycal\t$ in terms of $\bar \lambda$.

\begin{restatable}[Population Properties: Stage II]{lemma}{lemsteptwopopprop}\label{lem:step2popprop}
Suppose that $\ycal$ is rank $K$, and let $\ycal = \U \Sigma \mathbf{V}\t$ be its (rank $K$) singular value decomposition.  Then it holds that
\begin{align*}
    \U &= \mathbf{Z}  \mathbf{M},
\end{align*}
where $\mathbf{M} \in \mathbb{R}^{K\times K}$ is some invertible matrix satisfying
\begin{align*}
    \| \mathbf{M}_{r\cdot} - \mathbf{M}_{s\cdot} \| &= \sqrt{n_r\inv + n_s\inv}.
\end{align*}
In addition, when $n_{\min} \asymp n_{\max}$, it holds that
\begin{align*}
    \lambda_Y^2 :&= \lambda_{\min} \bigg( \sum_l \ytilde\m(\ytilde\m)\t \bigg) \gtrsim \frac{n}{K} L\bar \lambda. 
\end{align*}
\end{restatable}
Armed with these lemmas as well as Theorems \ref{thm:firststep}, \ref{thm:step2sintheta}, and \ref{thm:step2asympexp}, we are prepared to prove \cref{thm:clusteringerror}.

\begin{proof}[Proof of \cref{thm:clusteringerror}]
We follow the analysis technique developed in \cite{jin_improvements_2022} to derive an exponential rate for the output of $(1+\eps)$ $K$-means. 
First will use the the $\sin\bTheta$ bound (\cref{thm:step2sintheta})   together with Lemma 5.3 of \cite{lei_consistency_2015} to demonstrate a Hamming error of order strictly less than $\frac{n_{\min}}{4}$, so that each cluster has at a majority of its true members.  This allows us to associate each empirical cluster centroid to a true cluster centroid.  Next, we will study the empirical centroids of these clusters to show that they are strictly closer to their corresponding true cluster centroid than they are to each other.  Finally, we decompose the expected error into individual node-wise errors, where we apply the asymptotic expansion in \cref{thm:step2asympexp} to obtain the exponential error rate.  

In what follows, let $\mathcal{E}_{\sin\bTheta}$ denote the event
\begin{align*}
    \| \sin\bTheta(\uhat,\U) \| \leq \frac{\beta}{8 K \sqrt{C_{\eps}}},
\end{align*}
where $\beta \in (0,1)$ is such that $n_{\min} \geq \beta n_{\max}$ and $C_{\eps}$ is a constant to be defined in the subsequent analysis. We note that by \cref{thm:step2sintheta} the event $\mathcal{E}_{\sin\bTheta}$ holds with probability at least $1 - O(n^{-10})$. We also let $(\mathbf{\hat Z}, \mathbf{\hat M})$ denote the output of $(1+\eps)$ $K$-means on the rows of $\uhat$, where $\mathbf{\hat Z} \in \{0,1\}^{n\times K}$ and $\mathbf{\hat M} \in\mathbb{R}^{K\times K}$. \\ 
\\ 
\noindent
\textbf{Step 1: Initial Hamming Error} \\
First by \cref{lem:step2popprop} it holds that $\U = \mathbf{Z M}$ where $\mathbf{M}$ has $K$ unique rows satisfying 
\begin{align*}
    \frac{1}{\sqrt{n_{\max}} } \leq \| \mathbf{M}_{r\cdot} - \mathbf{M}_{s\cdot} \| \leq \frac{\sqrt{2}}{\sqrt{n_{\min}}}.
\end{align*}
Define the matrix
$\mathbf{\hat V} := \mathbf{\hat Z} \mathbf{\hat M}$.  Define 
$S_r := \{ i \in \mathcal{C}(r): \| \mathbf{W}_*\mathbf{\hat V}_{i\cdot} - \U_{i\cdot} \| \geq \delta_r/2\}$, where $\delta_r = \frac{1}{\sqrt{n_r}}$.  By Lemma 5.3 of \cite{lei_consistency_2015}, it holds that
\begin{align*}
    \frac{1}{n}  \sum_{i=1}^{n} \mathbb{I}\{ \hat z(i) \neq z(i) \}  &\leq \frac{1}{n} \sum_{r=1}^K |S_r | \leq \sum_{r=1}^K \frac{|S_r|}{n_r} = \sum_{r=1}^K |S_r| \delta_r^2 \\
    &\leq C_{\eps} \| \uhat \mathbf{W}_* - \U \|_F^2 \leq C_{\eps} \| \sin\bTheta(\uhat, \U) \|_F^2 \leq C_{\eps} K \| \sin\bTheta(\uhat,\U) \|^2.
\end{align*}
Therefore, on the event $\mathcal{E}_{\sin\bTheta}$, it holds that
\begin{align*}
   \sum_{i=1}^{n} \mathbb{I}\{ \hat z(i) \neq z(i) \}  &\leq C_{\eps} K n \frac{\beta^2}{64 K^2 C_{\eps}} 
    \leq \frac{\beta}{64} n_{\min},
\end{align*}
since $n \leq K n_{\max} \leq \frac{K}{\beta} n_{\min}$.  Therefore, since this error is strictly less than $ \beta n_{\min}/64\leq n_{\min}n_r/(64n_{\max})$, each cluster $r$ has at least $n_r -   \beta n_{\min} /64 \geq (1 - n_{\min} /(64n_{\max})) n_r\geq (63/64)n_r$ of its true members.  This implies that we can associate each empirical cluster to a true cluster -- let these empirical clusters be denoted $\mathcal{\hat C}(r)$.  Observe that we must have that $|\mathcal{\hat C}(r) | \geq (1 - \beta/64) n_{\min}$ and that $|\mathcal{\hat C}(r) \setminus \mathcal{C}(r) | \leq \beta n_{\min}/64$.
\\ \ \\
\noindent
\textbf{Step 2: Properties of Empirical Centroids}\\
Recall that the cluster centroid associated to $\mathcal{\hat C}(r)$ is equal to $\mathbf{\hat M}_{r\cdot}$.  Then by definition,
\begin{align*}
    \mathbf{\hat M}_{r\cdot} &= \frac{1}{|\mathcal{\hat C}(r)|} \sum_{i \in \mathcal{\hat C}(r)} \uhat_{i\cdot}.
\end{align*}
Recall that $\U$ consists of $K$ unique rows of $\mathbf{M}$. Without loss of generality assume that $\mathbf{M}_{r\cdot}$ is associated to $\mathcal{C}(r)$.  Then 
\begin{align*}
    \| \mathbf{W}_* \mathbf{\hat M}_{r\cdot} - \mathbf{M}_{r\cdot} \| &= \frac{1}{|\mathcal{\hat C}(r)|} \big\| \sum_{i \in \mathcal{\hat C}(r)} ( \wstar \uhat_{i\cdot} - \mathbf{M}_{r\cdot}) \| \\
        &\leq \frac{1}{|\mathcal{\hat C}(r)|} \big\| \sum_{i \in \mathcal{\hat C}(r)} (\wstar \uhat_{i\cdot} - \U_{i\cdot}) \big\| + \frac{1}{|\mathcal{\hat C}(r)|} \big\| \sum_{i \in \mathcal{\hat C}(r)}  (\U_{i\cdot} - \mathbf{M}_{r\cdot} )\big\|                    \\ &\leq \frac{1}{|\mathcal{\hat C}(r)|} \bigg\| \sum_{i \in  \mathcal{\hat C}(r)} (\mathbf{W}_* \uhat_{i\cdot} - \U_{i\cdot}) \bigg\| + \frac{1}{|\mathcal{\hat C}(r)|} \bigg\| \sum_{i \in \mathcal{\hat C}(r)\setminus\mathcal{C}(r)} (\U_{i\cdot} - \mathbf{ M}_{r\cdot}) \bigg\|.
\end{align*}
We observe that for $i \notin \mathcal{C}(r)$, it holds that
\begin{align*}
   \frac{1}{\sqrt{n_{\max}}}\leq  \| \U_{i\cdot} - \mathbf{ M}_{r\cdot} \| \leq \frac{\sqrt{2}}{\sqrt{n_{\min}}}
\end{align*}
by \cref{lem:step2popprop}.  Therefore,
\begin{align*}
    \| \mathbf{W}_* \mathbf{\hat M}_{r\cdot} - \mathbf{M}_{r\cdot} \| &\leq \frac{1}{|\mathcal{\hat C}(r)|} \bigg\| \sum_{i \in \mathcal{\hat C}(r)} (\mathbf{W}_* \uhat_{i\cdot} - \U_{i\cdot} )\bigg\| + \frac{1}{|\mathcal{\hat C}(r)|} \bigg\| \sum_{i \in \mathcal{\hat C}(r)\setminus\mathcal{C}(r)} (\U_{i\cdot} - \mathbf{ M}_{r\cdot})\bigg\| \\
    &\leq \frac{1}{|\mathcal{\hat C}(r)|^{1/2}} \| \uhat \mathbf{W}_*\t - \U \|_F + \frac{ | \mathcal{\hat C}(r) \setminus \mathcal{C}(r) |}{|\mathcal{\hat C}(r)|} \frac{\sqrt{2}}{\sqrt{n_{\min}}} \\
    &\leq \frac{\sqrt{2K}}{\sqrt{n_{\min}(1 - \beta/64)}}  \| \sin\bTheta(\uhat, \U) \| + \frac{ \beta n_{\min} }{64(1 - \beta/64) n_{\min} } \frac{\sqrt{2}}{\sqrt{n_{\min}}} \\
    &\leq \frac{\sqrt{2K}}{\sqrt{\beta n_{\max}(1 - \beta /64)}} \frac{\beta}{8 K C_{\eps}^{1/2}} + \frac{\beta}{64 (1 - \beta/64)} \frac{\sqrt{2}}{\sqrt{\beta n_{\max}}} \\
    &\leq \frac{1}{\sqrt{n_{\max}}} \bigg( \frac{\sqrt{2\beta }}{8 K^{1/2} C_{\eps}^{1/2}\sqrt{1 - \beta/64}} + \frac{\beta^{1/2}\sqrt{2}}{64 (1 - \beta/64)} \bigg) \\
    &\leq \frac{1}{8 \sqrt{n_{\max}}},
\end{align*}
since $n_{\min} \geq \beta n_{\max}$, $K \geq 1$ and $\beta < 1$, as well as the assumption $C_{\eps} \geq 4$.  Therefore, on the event $\mathcal{E}_{\sin\bTheta}$ it holds that
\begin{align*}
    \max_{1\leq r \leq K} \| \mathbf{W}_* \mathbf{\hat M}_{r\cdot} - \mathbf{M}_{r\cdot} \| \leq \frac{1}{8\sqrt{n_{\max}}}.
\end{align*}
\noindent
\textbf{Step 3: Applying The Asymptotic Expansion} \\
In this section we will use the previous bound on the cluster centroids  and \cref{thm:step2asympexp} to obtain the desired bound. Recall that by \cref{thm:step2sintheta},  $\mathbb{P}(\mathcal{E}_{\sin\bTheta}^c) = O(n^{-10})$.   It then holds that
\begin{align*}
    \mathbb{E}\ell(\hat z,z ) &= \frac{1}{n} \sum_{i=1}^n \p( \mathbf{Z}_{i\cdot} \neq \mathbf{\hat Z}_{i\cdot}  ) \\
    &\leq \frac{1}{n} \sum_{i=1}^n \p( \mathbf{Z}_{i\cdot} \neq \mathbf{\hat Z}_{i\cdot} , \mathcal{E}_{\sin\bTheta} )  + O(n^{-10}).
\end{align*}
Suppose that $\| (\uhat \mathbf{W}_*\t)_{i\cdot} - \U_{i\cdot} \| \leq \frac{1}{4\sqrt{n_{\max}}},$ and suppose the $i$'th node is in community $r$.  Then on the event $\mathcal{E}_{\sin\bTheta}$
\begin{align*}
    \| \mathbf{W}_* \uhat_{i\cdot} - \mathbf{W}_* \mathbf{\hat M}_{r\cdot} \| &\leq \| (\uhat \mathbf{W}_*\t)_{i\cdot} - \U_{i\cdot} \| + \| \U_{i\cdot} - \mathbf{\hat M}_{r\cdot} \| \\
    &\leq \frac{1}{4 \sqrt{n_{\max}}} + \max_{s} \| \mathbf{W}_* \mathbf{\hat M}_{s\cdot} - \mathbf{M}_{s\cdot} \| \\
    &\leq \frac{3}{8 \sqrt{n_{\max}}}.
\end{align*}
In addition, for any $s \neq r$, we have that
\begin{align*}
     \| \mathbf{W}_* \uhat_{i\cdot} - \mathbf{W}_* \mathbf{\hat M}_{s\cdot} \| &\geq \| \mathbf{M}_{r\cdot} - \mathbf{M}_{s\cdot} \| - \| \mathbf{W}_* \uhat_{i\cdot} - \U_{i\cdot} \| - \| \mathbf{W}_* \mathbf{\hat M}_{s\cdot} - \mathbf{M}_{s\cdot} \| \\
     &\geq \frac{1}{\sqrt{n_{\max}}} - \frac{1}{4 \sqrt{n_{\max}}} - \frac{1}{8 \sqrt{n_{\max}}} \\
     &\geq \frac{5}{8 \sqrt{n_{\max}}}.
\end{align*}
Therefore, node $i$ must belong to cluster $\mathcal{\hat C}(r)$, so that there is no error on node $i$.  Therefore,
\begin{align*}
    \p( \mathbf{Z}_{i\cdot} \neq \mathbf{\hat Z}_{i\cdot}, \mathcal{E}_{\sin\bTheta}) &\leq \p( \| (\uhat \mathbf{W}_*\t)_{i\cdot} - \U_{i\cdot}\| \geq \frac{1}{4 \sqrt{n_{\max}}} ) \\
    &\leq \p( \| (\uhat \mathbf{W}_*\t)_{i\cdot} - \U_{i\cdot}\| \geq \frac{1}{4 \sqrt{n_{\max}}}, \mathcal{E}_{\mathrm{Stage \ II}} ) + O(n^{-10}),
\end{align*}
where $\mathcal{E}_{\mathrm{Stage \ II}}$ is the event in \cref{thm:step2asympexp}.  On the event $\mathcal{E}_{\mathrm{Stage \ II}}$ it holds that 
\begin{align*}
    \uhat \mathbf{W}_*\t - \U &= \sum_{l} \mathcal{L}(\A\m - \bP\m ) (\ytilde\m)\t \U \Sigma^{-2} + \mathcal{R}_{\mathrm{Stage \ II}},
\end{align*}
with
\begin{align*}
    \| \mathcal{R}_{\mathrm{Stage \ II}}\|_{2,\infty} &\leq \frac{1}{16 \sqrt{n_{\max}}}. 
\end{align*}  
Therefore,
\begin{align*}
    \p( \| (\uhat \mathbf{W}_*\t)_{i\cdot} - \U_{i\cdot}\| \geq \frac{1}{4 \sqrt{n_{\max}}}, \mathcal{E}_{\mathrm{Stage \ II}} ) &\leq \p\bigg( \bigg\| e_i\t \sum_{l} \mathcal{L}(\A\m - \bP\m ) (\ytilde\m)\t \U \Sigma^{-2} \bigg\| \geq \frac{1}{8 \sqrt{n_{\max}}} \bigg)\\
    &\leq \p\bigg( \bigg\| e_i\t \sum_{l} \mathcal{L}(\A\m - \bP\m ) (\ytilde\m)\t \U \Sigma^{-2} \bigg\| \geq C \frac{\sqrt{K}}{\sqrt{n}} \bigg)\\
    &\leq K \max_k \p\bigg( \bigg| \sum_{l} e_i\t \mathcal{L}(\A\l - \bP\l) \big( \ytilde\m)\t \U \Sigma^{-2} e_k \bigg| \geq C \frac{1}{  \sqrt{n}}\bigg).
\end{align*}
We will apply the Bernstein inequality now.  We have that 
\begin{align*}
    \sum_{l} e_i\t \mathcal{L}(\A\m - \bP\m ) (\ytilde\m)\t \U \Sigma^{-2} e_k &= \sum_{l} \sum_{j} (\A\m - \bP\m)_{ij} \bigg(\U\m |\Lambda\m|^{-1/2} \ipq \mathbf{J}( \xtilde_{i\cdot})(\ytilde\m)\t \U \Sigma^{-2}\bigg)_{jk}.
\end{align*}
Using \cref{lem:popprop} and \cref{lem:step2popprop}, the variance $v$ of this quantity is upper bounded by
\begin{align*}
    v &\leq \sum_{l} \sum_{j} \theta_i\l \theta_j\l \|e_j\t \U\m |\Lambda\m|^{-1/2} \ipq \mathbf{J}( \xtilde_{i\cdot})(\ytilde\m)\t \U \Sigma^{-2} \|^2 \\
    &\leq \sum_{l} \sum_{j} \theta_i\l \theta_j\l \| e_j\t \U\m \|^{2} \| |\Lambda\m|^{-1/2} \|^2 \| \mathbf{J}(\xtilde_{i\cdot}) \|^2 \| (\ytilde\m)\t \U \Sigma^{-2} \|^2 \\
    &\leq C \sum_{l} \sum_{j} \theta_i\m \theta_j\m \frac{K (\theta_j\m)^2}{\|\theta\m\|^2} \frac{K}{\|\theta\m\|^2 \lambda_{\min}\m} \frac{1}{\|\xtilde_{i\cdot}\|^2} \frac{n K^2 }{n^2 L^2 \bar \lambda^2} \\
    &\leq C \frac{K^4}{n L^2 \bar \lambda^2} \sum_{l} \sum_{j} \theta_i\m \theta_j\m \frac{(\theta_j\m)^2}{\|\theta\m\|^4 \lambda_{\min}\m (\theta_i\m)^2} \\
    &\leq  C \frac{K^4}{n L^2 \bar \lambda^2} \sum_{l} \frac{\|\theta\m\|_3^3}{\theta_i\m\|\theta\m\|^4 \lambda\m_{\min}}.
\end{align*}
In addition, each term satisfies
\begin{align*}
    \max_{l,j} \| e_j\t \U\m |\Lambda\m|^{-1/2} \ipq \mathbf{J}( \xtilde_{i\cdot})(\ytilde\m)\t \U \Sigma^{-2} \| &\leq \max_{l,j} C\frac{\theta_j\m K}{ \theta_i\m \|\theta\m\|^2 (\lambda_{\min}\m)^{1/2}} \frac{K}{\sqrt{n} L \bar \lambda}\\
    &\leq C \frac{ K^2}{\sqrt{n} L \bar \lambda} \max_{l} \frac{\theta_{\max}\m}{\theta_i\m\|\theta\m\|^2 (\lambda_{\min}\m)^{1/2}}.
\end{align*}
By Bernstein's inequality,
\begin{align*}
  \p\bigg( &\bigg| \sum_{l} e_i\t \mathcal{L}(\A\m - \bP\m ) \big( \ytilde\m)\t \U \Sigma^{-2} e_k \bigg| \geq C \frac{1}{ \sqrt{n}} \bigg) \\
  &\leq 2 \exp\bigg( - \frac{ C^2\frac{1}{128 n}}{C_1 \frac{K^4}{n L^2 \bar \lambda^2} \sum_{l} \frac{\|\theta\m\|_3^3}{\theta_i\m\|\theta\m\|^4 \lambda\m_{\min}} + C_1 \frac{1}{\sqrt{n}} \frac{ K^2}{\sqrt{n} L \bar \lambda} \max_{l} \frac{\theta_{\max}\m}{\theta_i\m\|\theta\m\|^2 (\lambda_{\min}\m)^{1/2}}} \bigg) \\
  &\leq 2 \exp\bigg( - \frac{ C_1}{C_2 \frac{K^4}{ L^2 \bar \lambda^2} \sum_{l} \frac{\|\theta\m\|_3^3}{\theta_i\m\|\theta\m\|^4 \lambda\m_{\min}} + \frac{1}{24 } C_2 \frac{ K^{2}}{ L \bar \lambda} \max_{l} \frac{\theta_{\max}\m}{\theta_i\m\|\theta\m\|^2 (\lambda_{\min}\m)^{1/2}}} \bigg) \\
  &\leq 2 \exp\Bigg( - C_3 \min\bigg\{ \frac{\bar \lambda^2 L}{K^4} \bigg( \frac{1}{L}\sum_{l} \frac{\|\theta\m\|_3^3}{\theta_i\m \| \theta\m\|^4 \lambda_{\min}\m} \bigg)\inv, \frac{L \bar \lambda}{K^{2}} \min_{m}\frac{\theta_i\m \|\theta\m\|^2 (\lambda_{\min}\m)^{1/2}} {\theta_{\max}\m} \bigg\} \Bigg) \\
  &\leq 2 \exp\Bigg( -c L \min\bigg\{ \frac{\bar \lambda^2}{K^4 \mathrm{err}_{\ave}^{(i)}}, \frac{\bar \lambda}{K^{2} \mathrm{err}_{\max}^{(i)}} \bigg\} \Bigg),
\end{align*}
where  $\mathrm{err}_{\ave}^{(i)}$ and $ \mathrm{err}_{\max}^{(i)}$ are as defined in Eq~\eqref{eq:def-errave}. This completes the proof. 
\end{proof}

\subsubsection{Proof of Theorem~\ref{cor:perfectclustering}}
\begin{proof}[Proof of \cref{cor:perfectclustering}]
The proof proceeds from partway through the proof of \cref{thm:clusteringerror}.  We have already shown that on the event $\mathcal{E}_{\sin \bTheta}$  if $\| (\uhat \wstar\t )_{i\cdot} - \U_{i\cdot} \| \leq \frac{1}{4 \sqrt{n_{\max}}}$ then node $i$ must be classified correctly. 
By repeating the argument in step 3 of the proof of \cref{thm:clusteringerror}, it holds that \begin{align*}
    \p( \mathbf{Z}_{i\cdot} \neq  \mathbf{\hat{Z}}_{i\cdot} ) \leq 2 K \exp\bigg( - cL \min\bigg\{ \frac{\bar \lambda^2}{K^4 \mathrm{err}_{\ave}^{(i)}}, \frac{\bar \lambda}{K^{2} \mathrm{err}_{\max}^{(i)}} \bigg\} \bigg) + O(n^{-10}). 
\end{align*}
In order for the exponential to be strictly less than $O(n^{-10})$, we require that 
\begin{align*}
     \min\bigg\{ \frac{\bar \lambda^2}{K^4 \mathrm{err}_{\ave}^{(i)}}, \frac{\bar \lambda}{K^{2} \mathrm{err}_{\max}^{(i)}} \bigg\} \geq \frac{C \log(n)}{L},
\end{align*}
where $C$ is a sufficiently large constant.  Recalling the definitions of $\mathrm{err}_{\ave}^{(i)}$ and $\mathrm{err}_{\max}^{(i)}$, we see that we must have
\begin{align*}
    \frac{\bar \lambda}{K^{2}} &\geq \frac{C \log(n)}{L} \max_l \frac{\theta_{\max}\l}{\theta_{\min}\l}  \frac{1}{\|\theta\l\|^2 (\lambda_{\min}\l)^{1/2}}; \\
    \frac{\bar \lambda^2}{K^4} &\geq \frac{C\log(n)}{L} \bigg( \frac{1}{L} \sum_{l} \frac{\|\theta\l\|_3^3}{\theta_{\min}\l \|\theta\l\|^4 \lambda_{\min}\l} \bigg).
\end{align*}
Considering the first term and rearranging, we see that we require that
\begin{align*}
    \frac{\bar \lambda}{K^{2}} \min_l \bigg( \frac{\theta_{\min}\l}{\theta_{\max}\l} \bigg)  \|\theta\l\|^2 (\lambda_{\min}\l)^{1/2} \geq \frac{C \log(n)}{L}. 
\end{align*}
A sufficient condition is that 
\begin{align*}
     \min_l \snr_l^2 \geq \frac{C K^8 \log(n)}{L \bar \lambda} 
\end{align*}
As for the second term, by upper bounding $\|\theta\l\|_3^3 \leq \theta_{\max}\l \|\theta\l\|^2$, we see that it sufficient to have that
\begin{align*}
    \frac{\bar \lambda^2}{K^4} &\geq \frac{C\log(n)}{L} \bigg( \frac{1}{L} \sum_{l} \bigg( \frac{\theta_{\max}\l}{\theta_{\min}\l} \bigg)\frac{1}{ \|\theta\l\|^2 \lambda_{\min}\l} \bigg). 
    \numberthis 
    \label{tosnr2}
\end{align*}
Therefore, rearranging \eqref{tosnr2} yields the sufficient condition 
\begin{align*}
   \bigg( \frac{1}{L} \sum_{l} \frac{1}{\snr_l^2} \bigg)\inv \geq C \frac{K^8 \log(n)}{L \bar \lambda^2}.
\end{align*}
It is straightforward to check that the condition in \cref{cor:perfectclustering} is sufficient for the result to hold.
\end{proof}


\section{Proofs of Identifiability and Algorithm Recovery Results} \label{sec:identifiabilityproofs}
In this section we prove \cref{prop:identifiability} and \cref{prop:algorithmcorrectness}, as well as \cref{lem:popprop} and \cref{lem:step2popprop}.

\subsection{Proof of Theorem~\ref{prop:identifiability}}

\begin{proof}[Proof of \cref{prop:identifiability}] We first prove the ``if'' direction.  Suppose for contradiction that there is another block membership matrix $\widetilde{\bZ}\in\{0,1\}^{n\times K'}$ with at least one vertex assigned to each community, and  positive diagonal matrices 
$\{\tilde{\bTheta}\m\}_{l=1}^L$ and symmetric matrices $\{\tilde{\bB}\m\}_{l=1}^L$
such that
$${\bTheta}\m{\bZ}{\bB}\m{\bZ}^\top {\bTheta}\m = \tilde{\bTheta}\m\tilde{\bZ}\tilde{\bB}\m\tilde{\bZ}^\top \tilde{\bTheta}\quad\quad\quad\text{for each $l\in[L]$.}$$
Equivalently, since the matrices $\tilde\bTheta\l$ have positive diagonal, for all $l\in[L]$ it holds that
\begin{align}
    \tilde{\bZ}\tilde{\bB}\m\tilde{\bZ}^\top = & [\tilde{\bTheta}\m]^{-1}{\bTheta}\m {\bZ}{\bB}\m{\bZ}^\top [\tilde{\bTheta}\m]^{-1}{\bTheta}\m \nonumber\\
    := & \bGamma\m {\bZ} {\bB}\m (\bGamma\m\bZ)^\top\nonumber\\
    = &  \bGamma\m {\bZ} {\bV}\m\bD\m (\bGamma\m\bZ\bV\m)^\top. \label{eq:proof-identifiability}
\end{align}
For any vertex index $i\in[n]$,  denote by $z(i)$ and $\tilde{z}(i)$  the community memberships  according to $\bZ$ and $\tilde\bZ$. We will show that $K'\geq K$ and if $K'=K$ then $z(i) = z(j)$ if and only if $\tilde{z}(i) = \tilde{z}(j)$.

By the RHS of Eq.~\eqref{eq:proof-identifiability}, the column space of $\bGamma\l\bZ\bV\l$ should be contained within the column space of $\tilde\bZ$ (as these two matrices are full rank by construction), and hence, there is a matrix $\bM\l\in\real^{K'\times K_l}$ such that
\begin{equation}
    \bGamma\m\bZ\bV\l = \tilde\bZ \bM\m,\quad\quad\text{for all }l\in[L].
    \label{eq:proof-ident-tildeZ-Z}
\end{equation}
In particular, this implies that  for any $i\in[n]$,
\begin{equation}
    \bGamma\l_{ii}\bV\l_{z(i)\cdot} = \bM\l_{\tilde{z}(i) \cdot}\quad\quad\text{for all }l\in[L].
    \label{eq:propidentif-equalityrows}
\end{equation}
If $\tilde{z}(i) = \tilde{z}(j)$ then 
$$\bGamma\l_{ii}\bV_{z(i)\cdot}\l = \bGamma\l_{jj}\bV_{z(j)\cdot}\l\quad\quad\text{for all }l\in[L].$$
This equation implies that the normalized rows are the same, i.e., $\bQ_{z(i)\cdot}\l = \bQ_{z(j)\cdot}\l$ for all $l\in[L]$, and hence 
$\bQ_{z(i)\cdot}=\bQ_{z(j)\cdot}$, which is only possible if $z(i) = z(j)$ according to the condition in the proposition.

Now, take a set of vertices $\mathcal{T}\subset [n]$ such that each vertex is in a different community according to $\tilde{\bZ}$. Without loss of generality, suppose that $\tilde{\bZ}_{\mathcal{T} \cdot} = I$, and hence, Eq.~\eqref{eq:proof-ident-tildeZ-Z} implies
$$(\bGamma\l\bZ\bV\l)_{\mathcal{T}\cdot} = \bGamma\l_{\mathcal{T}\cdot}\bV\l_{z(\mathcal{T})\cdot} = \bM\l\quad\quad\text{for all } l\in[L].$$
If there are two indexes $i,j\in\mathcal{T}$ such that $z(i) = z(j)$, then the corresponding rows of $\bM\l$ are proportional, that is $\bM_{z(i)\cdot} = \bGamma_{ii}\l \bV\l_{z(i)\cdot}$ and $\bM_{z(j)\cdot} = \bGamma_{jj}\l \bV\l_{z(i)\cdot}$ for all $l\in[L]$. If $K'=K$, this implies that  $\bM\l$ can only have at most $K-1$ different rows that are not proportional, and these are the same for all $l\in[L]$. Hence, by Eq.~\eqref{eq:propidentif-equalityrows} the matrix $\bQ$ has at most $K-1$ different rows, which contradicts the assumption. Note that this is also the case if $K'<K$. If $K'>K$, then it is still possible to have $z(i) = z(j)$, but then $\bZ$ can fit the same model with fewer communities.

We now prove the ``only if'' direction.  Suppose for contradiction that $\bQ$ has repeated rows; we will construct $\tilde{\bZ}$ and $\bB\l$ that yield the same $\bP\l$ matrices.  Without loss of generality we may assume that rows one and two are repeated, since communities are identifiable up to permutation. Furthermore, without loss of generality we can have $\bQ\l = \bV\l$. Indeed, for $i \in \mathcal{C}(r)$, we can rescale $\theta_{i}\l$ via
$\theta_i\l \mapsto \theta_i\l \| \bV\l_{r\cdot} \|$, which still yields the same matrix $\bP\l$ since $$\bP\l_{ij} = \theta_i\l \theta_j\l (\bV\l \bD\l \bV\l)\t_{z(i)z(j)} =\big(  \theta_i\l \|\bV\l_{z(i)\cdot}\| \big) \big( \theta_j\l \|\bV\l_{z(j)\cdot} \| \big) \frac{ (\bV\l \bD\l \bV\l)\t_{z(i)z(j)}}{\|\bV\l_{z(i)\cdot}\| \|\bV\l_{z(j)\cdot} \|}.$$
Therefore, the first two rows of $\bV\l$ are repeated for all $l$.  However, this implies that
\begin{align*}
    \bB\l_{1r} = \big( \bV\l \bD\l (\bV\l)\t \big)_{1r} = \sum_{s=1}^{K_l} \bV\l_{1s} \bD\l_{ss} \bV\l_{rs} = \sum_{s=1}^{K_l} \bV\l_{2s} \bD\l_{ss} \bV\l_{rs} = \bB\l_{2r},
\end{align*}
which shows that the first row and column of $\bB\l$ is repeated.  Therefore, we can collapse the first two communities into one community, creating a new matrix $\tilde{\bB}\l$ with $K-1$ communities (with the first two communities merged).  Then we have that
\begin{align*}
    \bP\l_{ij} = \theta_i\l \theta_j\l \bB\l_{z(i)z(j)} = \theta_i\l \theta_j\l \bB\l_{\tilde{z}(i) \tilde{z}(j)},
\end{align*}
which shows that $\bZ$ is not identifiable unless $\bQ$ has no repeated rows.
\end{proof}


\subsection{Proof of Proposition~\ref{prop:algorithmcorrectness}}

\begin{proof}[Proof of \cref{prop:algorithmcorrectness}]
We will demonstrate that the left singular vectors obtained immediately before clustering contain exactly $K$ unique rows, for which the final result follows. We will analyze each stage separately.
\\ \ \\ \noindent
\textbf{First Stage (individual network embedding):}  First, suppose that $\bP\l=\U\l \Lambda\l (\U\l)\t$ with $\bU\in\real^{n\times K_l}$ is the eigendecomposition of $\bP\l$, and let $\bB\l=\bV\l \bD\l (\bV\l)\t$ be the eigendecomposition of $\bB\l$, with $\bV\in\real^{K\times K_l}$ a matrix with orthogonal columns and $\bD\l\in\real^{K_l\times K_l}$ a diagonal matrix with non-zero elements in the diagonal.  From this factorization it is evident that
\begin{align*}
    \bP\l &= \bTheta\l \bZ \bB\l \bZ\t \bTheta\l = \bTheta\l \bZ \bV\l\bD\l (\bV\l)\t\bZ\t \bTheta\l.
\end{align*}
Since $\bU\l$ and $\bTheta\l\bZ\bV\l$ are full rank matrices, they have the same column space, so
\begin{align*}
    \U\l &= \bTheta\l \bZ \bV\l\mathbf{H}\l,
\end{align*}
where $\mathbf{H}\l\in\real^{K_l\times K_l}$  is a full rank matrix. From this decomposition it is immediate that $\U\l$ consists of rows of $\bV\l\mathbf{H}\l$ with each row of $\U\l$ scaled by $\theta_i\l$.  Let $\xi_r\l$ denote the $r$'th row of $\bV\mathbf{H}|\Lambda\l|^{1/2}$.  Then  if $z(i) = r$,
\begin{align*}
    (\U\l |\Lambda\l|^{1/2})_{i\cdot} &= \theta_i\l \xi_r\l = \theta_i\l\bV_{r\cdot}\l \bH\l,
\end{align*}
and hence
\begin{align*}
    \ytilde\l_{i\cdot} &= \frac{\theta_i\l \xi_r\l}{\|\theta_i\l \xi_r\l\|} = \frac{1}{\|\xi_r\l\|} \xi_r\l = 
    \frac{1}{\|\bV_{r\cdot}\l\bH\l\|} \bV_{r\cdot}\l\bH\l,
\end{align*}
which does not depend on  $\theta_i\l$.
\\ 

\noindent \textbf{Second Stage (joint network embedding)}: We now consider the left singular vectors of the matrix $\ycal$ defined as
\begin{align*}
    \ycal &= [ \ytilde\one, \cdots ,\ytilde^{(L)}].
\end{align*}
Observe that the leading $K$ left singular vectors $\U$ of $\ycal$ are given by the leading $K$ eigenvectors of the matrix $\ycal \ycal\t$, which can equivalently be written as
\begin{align*}
    \ycal\ycal\t &= \sum_{l=1}^{L} \ytilde\l (\ytilde\l)\t.
\end{align*}
Consider $i$ and $j$ in community $r$ and $s$ respectively.  Then from the analysis in the previous step,
\begin{align*}
    \big(\ycal\ycal\t\big)_{ij} &= \sum_{l=1}^{L} \frac{\langle \xi_r\l, \xi_s\l \rangle}{\|\xi_r\l\| \|\xi_s\l\|}.
\end{align*}
Consequently, this shows that $\ycal\ycal\t$ is a matrix of the form
\begin{align*}
    \ycal\ycal\t &= \bZ \bigg( \sum_{l=1}^{L} (\Xi\l) (\Xi\l)\t  \bigg) \bZ\t = (\bZ \bXi)(\bZ\bXi)^\top,
\end{align*}
 where $\Xi\l$ is the matrix whose rows are $\xi_r\l/\|\xi_r\l\|$ and $\bXi = [\Xi^{(1)} \cdots \Xi^{(l)}]$. Next observe that
 $$\Xi\l = \bD\l \bQ\l \bH\l = \bQ\l \bM\l$$
 for some matrix $\bM\l$ that is full rank, where $\bQ\l$ is as in \cref{prop:identifiability}.  Since $\bQ$ has $K$ different rows (by assumption),  $\bXi$ has $K$ different rows, and hence $\bXi \bXi\t$ is a $K \times K$ block matrix.  Let $\U$ denote the leading $\tilde K$ eigenvectors  of $\ycal\ycal\t$, where $\tilde K$ is the rank of $\ycal$. Let $\mathbf{V \Gamma V}\t$ denote the eigendecomposiion of $(\bZ\t \bZ)^{1/2} \bXi \bXi\t (\bZ\t \bZ)^{1/2}$.  Then it is straightforward to see that $\U = \bZ (\bZ\t \bZ)^{-1/2} \bV$ since they both  have orthonormal columns.  It suffices to argue that $\bV$ does not have repeated rows.  Assuming this for the moment, by taking $\bM = (\bZ\t \bZ) \bV$, it holds that $\U = \bZ \bM$, with $\bM$  having no repeated rows, whence the result is proven.
 
 It remains to argue that $\bV$ does not have repeated rows.  Under the conditions of \cref{prop:identifiability} we have already shown that $\bXi$ does not have repeated rows.  Hence $\bXi \bXi\t$ is a block matrix with no repeated rows and columns, and hence $(\bZ\t \bZ)^{1/2} \bXi \bXi\t (\bZ\t \bZ)^{1/2}$ is also a block matrix with no repeated rows and columns.  Now assume for contradiction that $\bV$ has repeated rows.  This implies that $\bV = \tilde{\bZ} \tilde{\bV}$ for some matrices $\tilde{\bZ} \in \{0,1\}^{K \times \tilde K}$ and $\tilde{\bV} \in \mathbb{R}^{\tilde K\times \tilde K}$ a full rank matrix.  Suppose that $\bV$ has rows $r$ and $r'$ repeated, and without loss of generality suppose that row is the first row of $\tilde{\bV}$ (or else permute $\tilde{\bV}$), so that $\tilde{\bZ}_{r1} = \tilde{\bZ}_{r'1} =1$. Then from the equation $(\bZ\t \bZ)^{1/2} \bXi \bXi\t (\bZ\t \bZ)^{1/2} = \tilde{\bZ} \tilde{\bV} \Gamma \tilde{\bV}\t \tilde{\bZ}\t$, it holds that for all $1\leq s\leq K$,
 \begin{align*}
     \frac{1}{\sqrt{n_r n_s}} \langle \xi_{r\cdot}, \xi_{s\cdot} \rangle &= \langle \big( \tilde{\bZ} \tilde{\bV} \Gamma^{1/2} \big)_{r\cdot},\big( \tilde{\bZ} \tilde{\bV} \Gamma^{1/2} \big)_{s\cdot} \rangle \\
     &= \langle \big(\tilde{\bV} \Gamma^{1/2} \big)_{1\cdot},\big( \tilde{\bZ} \tilde{\bV} \Gamma^{1/2} \big)_{s\cdot} \rangle \\
    &= \langle \big( \tilde{\bZ} \tilde{\bV} \Gamma^{1/2} \big)_{r'\cdot},\big( \tilde{\bZ} \tilde{\bV} \Gamma^{1/2} \big)_{s\cdot} \rangle \\ 
     &= \frac{1}{\sqrt{n_{r'} n_s}} \langle \xi_{r'\cdot}, \xi_{s\cdot} \rangle.
 \end{align*}
 Consequently, since the above identity holds for all $s$, this shows that the $r$ and $r'$'th rows and columns of $(\bZ\t \bZ)^{1/2} \bXi \bXi\t (\bZ\t \bZ)^{1/2}$ are identical.  However, this is a contradiction, which completes the proof.
\end{proof}

\subsection{Proof of Lemma~\ref{lem:popprop}}

We will restate \cref{lem:popprop} for convenience.

\populationspectralproperties*

\begin{proof}[Proof of \cref{lem:popprop}]
Define the matrix
$$\mathbf{G}\m := K \|\theta\m\|^{-2} \diag( \|\theta_{\mathcal{C}(1)}\m\|, \dots, \|\theta_{\mathcal{C}(K)}\m\| ) \mathbf{B}\m\diag( \|\theta_{\mathcal{C}(1)}\m\|, \dots, \|\theta_{\mathcal{C}(K)}\m\| ).$$ 
Letting $\lambda_r(\cdot)$ denote the eigenvalues of a matrix, by Ostrowski's Theorem (Theorem 4.5.9 of \cite{horn2012matrix}) and \cref{ass:communityass}, the eigenvalues of $\bG\l$ satisfy $\lambda_r(\bG\l) \asymp \lambda_r(\bB\l)$.  Since the eigenvalues of $\bP\l = \bTheta\l \bZ \bB\l \bZ\t \bTheta\l$ are the same as the eigenvalues of the matrix $$\big( \bZ\t (\bTheta\l)^2 \bZ \big)^{1/2} \bB\l \big( \bZ\t (\bTheta\l)^2 \bZ \big)^{1/2},$$ we have 
\begin{align*}
    \lambda_r\bigg( \big( \bZ\t (\bTheta\l)^2 \bZ \big)^{1/2} \bB\l \big( \bZ\t (\bTheta\l)^2 \bZ \big)^{1/2} \bigg) \asymp \frac{\|\theta\l\|^2}{K} \lambda_r\big( \bG\l \big) \asymp \frac{\|\theta\l\|^2}{K} \lambda_r(\bB\l).
\end{align*}
To prove the other two assertions, we first observe that
\begin{align*}
    \bP\l &= \bTheta\l \bZ \big(\bZ\t (\bTheta\l)^2 \bZ  \big)^{-1/2} \bigg[ \big(\bZ\t (\bTheta\l)^2 \bZ  \big)^{1/2}  \bB\l \big(\bZ\t (\bTheta\l)^2 \bZ  \big)^{1/2} \bigg] \big(\bZ\t (\bTheta\l)^2 \bZ  \big)^{-1/2}\bZ\t \bTheta\l.
\end{align*}
Suppose that the matrix $\big(\bZ\t (\bTheta\l)^2 \bZ  \big)^{1/2}  \bB\l \big(\bZ\t (\bTheta\l)^2 \bZ  \big)^{1/2}$ has eigendecomposition $\tilde{\U} \Lambda^{(l)} \tilde{\U}\t$, which is permissible as both matrices share the same eigenvalues.
Then it holds that
\begin{align*}
    \bP\l = \bTheta\l \bZ \big(\bZ\t (\bTheta\l)^2 \bZ  \big)^{-1/2} \tilde{\U} \Lambda^{(l)} \tilde{\U}\t\big(\bZ\t (\bTheta\l)^2 \bZ  \big)^{-1/2} \bZ \bTheta\l.
\end{align*}
However, since the columns of the matrix $\bTheta\l \bZ \big(\bZ\t (\bTheta\l)^2 \bZ \big)^{-1/2} \tilde{\bU}$ are orthonormal, the decomposition above is a valid
eigenvector-eigenvalue decomposition for $\bP\l$.  
In particular, this shows that without loss of generality, we may take $\U\l$ to be
\begin{align*}
    \U\l = \bTheta\l \bZ  \big(\bZ\t (\bTheta\l)^2 \bZ \big)^{-1/2} \tilde{\bU}.
\end{align*}
We immediately obtain the bound
\begin{align*}
    \| \U_{i\cdot}\l \| &= \theta_i\l \frac{1}{\|\theta\l_{\mathcal{C}(z(i))}\|} \lesssim \theta_i\l \frac{\sqrt{K}}{\|\theta\l\|},
\end{align*}
where we have used the fact that $\|\theta_{\mathcal{C}(r)}\|^2 \asymp \frac{\|\theta\l\|^2}{K}$ for all $r$. Similarly, it holds that
\begin{align*}
    \| \xtilde_{i\cdot}\l \| \leq \| \U\l_{i\cdot} \| \| | \Lambda\l |^{1/2} \| \lesssim \theta_i\l {\sqrt{K}}.
\end{align*}
It remains to provide a lower bound on $\xtilde_{i\cdot}$. We have that
\begin{align*}
     \U\l |\Lambda\l|^{1/2} &=  \bTheta\l \bZ  \big(\bZ\t (\bTheta\l)^2 \bZ \big)^{-1/2} \tilde{\bU}  |\Lambda\l|^{1/2}.
 \end{align*}
Observe that $\big(\bZ\t (\bTheta\l)^2 \bZ \big)^{-1/2} \tilde{\bU} = \bB\l \big(\bZ\t (\bTheta\l)^2 \bZ \big)^{1/2} \tilde{\U} (\Lambda^{(l)})\inv$, which shows that
\begin{align*}
    \U\l |\Lambda\l|^{1/2} &= \bTheta\l \bZ \bB\l \big(\bZ\t (\bTheta\l)^2 \bZ \big)^{1/2} \tilde{\U} (\Lambda\l)\inv |\Lambda^{(l)}|^{1/2} \\
    &= \bTheta\l \bZ \bB\l \big(\bZ\t (\bTheta\l)^2 \bZ \big)^{1/2} \tilde{\U} |\Lambda\l|^{-1/2} \mathbf{W}, 
\end{align*}
where $\mathbf{W}$ is the diagonal matrix of signs of $\Lambda\l$.  
Consider a given row $i$ and suppose that $z(i) = r$.  Then by Ostrowki's Theorem again,
\begin{align*}
    \| \xtilde\l_{i\cdot} \| &\geq \theta_i\l \|\bB\l_{r\cdot} \|\sigma_{\min}\bigg(  \big(\bZ\t (\bTheta\l)^2 \bZ \big)^{1/2} \tilde{\U} |\Lambda\l|^{-1/2} \mathbf{W} \bigg) \\
    &\geq  \theta_i\l \|\bB\l_{r\cdot} \|\sigma_{\min} \big( \bZ\t (\bTheta\l)^2 \bZ \big) ^{1/2} \sigma_{\min} \big( |\Lambda\l|^{-1/2} \big) \\
    &\geq \theta_i\l \|\bB\l_{r\cdot} \| \min_{r} \|\theta_{\mathcal{C}(r)}\| \sigma_{\min} \big( |\Lambda\l|^{-1/2} \big) \\
    &\geq \theta_i\l \|\bB\l_{r\cdot} \|\min_{r} \|\theta\l_{\mathcal{C}(r)}\| \frac{\sqrt{K}}{\|\theta\l\|} \\
    &\gtrsim \theta_i\l \| \bB\l_{r\cdot} \| \\
    &\gtrsim \theta_i\l,
\end{align*}
where the final line follows from the assumption that $\bB\l$ has unit diagonals. 
This completes the proof.
\end{proof}

\subsection{Proof of Lemma~\ref{lem:step2popprop}}
We restate \cref{lem:step2popprop} for convenience.

\lemsteptwopopprop*

\begin{proof}[Proof of \cref{lem:step2popprop}]
The first part of the proof holds by Lemma 2.1 of \cite{lei_consistency_2015} applied to the matrix $\ycal \ycal\t$, which is a block matrix.  See also the proof of \cref{prop:algorithmcorrectness}.

For the second part we proceed as follows.  First recall by the proof of \cref{prop:algorithmcorrectness} that we can write the matrix $\ycal\ycal\t$ as the matrix
\begin{align*}
    \ycal\ycal\t &= \sum_{l=1}^{L} \mathbf{\Xi}\l (\mathbf{\Xi}\l)\t,
\end{align*}
where the matrix $\mathbf{\Xi}\l$ is defined as follows.  First, let $\mathbf{Q}\l$ be the matrix such that
\begin{align*}
    \U\l = \bTheta\l \bZ \mathbf{Q}\l.
\end{align*}
Then the rows of $\Xi\l$ are equal to the rows of $\mathbf{Q}\l|\Lambda\l|^{1/2}$ normalized by their magnitude. It was discussed in the proof of \cref{prop:algorithmcorrectness} that the entries of $\mathbf{Q}\l$ are of order $\frac{1}{\|\theta\l\|}$.  Observe that we can write $\ytilde\l = \mathbf{Z} (\tilde{\bD}\l)\inv \mathbf{Q}\l |\Lambda\l|^{1/2}$, where $\tilde{\bD}\l$ is the $K \times K$ diagonal matrix of row norms of $\mathbf{Q}\l |\Lambda\l|^{1/2}$.  Observe that 
\begin{align*}
    \| \tilde{\bD}\l \|^2 &= \max_i \| \big( \mathbf{Q}\l |\Lambda\l|^{1/2} \big)_{i\cdot} \|^2 \\
    &= \max_i \sum_{r=1}^{K} (\mathbf{Q}\l_{ir})^2 |\lambda_r|  \\
    &= \max_i \sum_{r=2}^{K} \frac{C}{\|\theta\l\|^2} \frac{\|\theta\l\|^2}{K} |\lambda_r(\bB\l)|  + \frac{C}{\|\theta\l\|^2} \|\theta\l\|^2 \\
    &\lesssim 1,
\end{align*}
 where we have applied \cref{lem:popprop} to observe that $\lambda_r \asymp \frac{\|\theta\l\|^2}{K} \lambda_r(\bB\l)$ for $2 \leq r \leq K$ and $\lambda_1 \asymp \|\theta\l\|^2$, since by \cref{ass:communityass} that the largest eigenvalue of $\bB\l$ is upper bounded by {$CK$}.  Therefore, we have that
 \begin{align*}
     \lambda_{\min} \bigg( \sum_{l} \ytilde\l (\ytilde\l)\t \bigg)  &= \lambda_{\min} \bigg( \sum_{l} \bZ (\tilde{\bD}\l)\inv \mathbf{Q}\l |\Lambda\l| (\mathbf{Q}\l)\t (\tilde{\bD}\l)\inv \bZ\t \bigg) \\
     &= \lambda_{\min} \bigg( \bZ\t \bZ \big( \sum_{l} (\tilde{\bD}\l)\inv (\mathbf{Q}\l |\Lambda\l| (\mathbf{Q}\l)\t (\tilde{\bD}\l)\inv \big) \bigg) \\
     &\geq \lambda_{\min}( \bZ\t \bZ) \lambda_{\min} \bigg( \sum_{l} (\tilde{\bD}\l)\inv \mathbf{Q}\l |\Lambda\l| (\mathbf{Q}\l)\t (\tilde{\bD}\l)\inv \bigg) \\
     &\gtrsim \frac{n}{K} L\bigg( \frac{1}{L} \sum_{l} \lambda_{\min} \bigg[  (\tilde{\bD}\l)\inv \mathbf{Q}\l |\Lambda\l| (\mathbf{Q}\l)\t (\tilde{\bD}\l)\inv \bigg] \bigg)
 \end{align*}
 where we have used the fact that $\lambda_{\min}(\bZ\t\bZ) = n_{\min} \asymp n/K$ and that the term inside the sum is rank $K$ and hence invertible. 
 Consequently, it suffices to show that
 \begin{align*}
     \lambda_{\min} \bigg[  (\tilde{\bD}\l)\inv \mathbf{Q}\l |\Lambda\l| (\mathbf{Q}\l)\t (\tilde{\bD}\l)\inv \bigg]  \gtrsim \lambda_{\min}\l.
 \end{align*}
 However, by the argument in \cref{lem:popprop}, it holds that $\mathbf{Q}\l (\mathbf{Q}\l)\t =  (\bZ\t(\bTheta\l)^2 \bZ)\inv$.
Set $\mathbf{G}\l := K\inv \|\theta\l\|^2 (\bZ\t(\bTheta\l)^2 \bZ).$ By \cref{ass:communityass}, $\|(\mathbf{G}\l)\inv\| \leq C$ and $\mathbf{Q}\l (\mathbf{Q}\l)\t = K \|\theta\l\|^{-2} (\mathbf{G}\l)\inv.$  
 Consequently,
 \begin{align*}
      \lambda_{\min} \bigg[  (\tilde{\bD}\l)\inv \mathbf{Q}\l |\Lambda\l| (\mathbf{Q}\l)\t (\tilde{\bD}\l)\inv \bigg]  &\gtrsim   \lambda_{\min}( (\tilde{\bD}\l)\inv)^2     \lambda_{\min} \big( \mathbf{Q}\l |\Lambda\l| (\mathbf{Q}\l)\t \big) \\
      &\gtrsim \lambda_{\min} \big((\mathbf{Q}\l)\t \mathbf{Q}\l |\Lambda\l|  \big) \\
      &\gtrsim \lambda_{\min} \big((\mathbf{Q}\l)\t \mathbf{Q}\l\big) \lambda_{\min}( |\Lambda\l|)\\
      &\gtrsim K \|\theta\l\|^{-2} \lambda_{\min} (\mathbf{G}\l) \frac{\|\theta\l\|^2}{K} \lambda_{\min}\l \\
      &\gtrsim \lambda_{\min}\l,
 \end{align*}
 where we have used \cref{ass:communityass} and \cref{lem:popprop} implicitly.  
This completes the proof.
\end{proof}


\section{Proof of First Stage Characterization (Theorem~\ref{thm:firststep})} \label{sec:firststageproofs}
This section contains the full proof of \cref{thm:firststep}.  First, we will restate \cref{thm:firststep} here for convenience.

\firststep*

As an immediate application of \cref{thm:firststep}, we can obtain a spectral norm concentration bound for the residual, which will be useful in subsequent steps.
\begin{lemma}\label{lem:step2resspectral}
The residual term $\mathcal{R}_{\mathrm{Stage \ I}}\m$ satisfies
\begin{align*}
    \| \mathcal{R}_{\mathrm{Stage \ I}}\m \| &\lesssim  \sqrt{n} \frac{K^2 \theta_{\max}\l \|\theta\l\|_1}{\lambda_{\min}\l \|\theta\l\|^4} \bigg( \log(n) \frac{\theta_{\max}\l}{\theta_{\min}\l} + \frac{{\sqrt{K}}}{\lambda_{\min}\l} + \bigg( \frac{\theta_{\max}\l}{\theta_{\min}\l} \bigg)^{1/2} \frac{K^{5/2} \log(n)}{(\lambda_{\min}\l)^{1/2}} \bigg).
\end{align*}
with probability at least $1 - O(n^{-15})$.
\end{lemma}
The proof of this result follows immediately by noting that $\|\cdot\|\leq \sqrt{n} \|\cdot\|_{2,\infty}$ and the bound in \cref{thm:firststep}.

We will also use an $\ell_{2,\infty}$ bound for the linear term appearing in \cref{thm:firststep} in the proof of \cref{thm:step2asympexp}.
\begin{lemma} \label{lem:lineartermtwoinfty}
The linear term in \cref{thm:clusteringerror} satisfies, with probability at least $1- O(n^{-15}),$
\begin{align*}
    \| \mathcal{L}(\A\m - \bP\m ) \|_{2,\infty} &\lesssim   \bigg( \frac{\theta_{\max}\l}{\theta_{\min}\l} \bigg)^{1/2} \frac{K\sqrt{\log(n)}}{(\lambda_{\min}\l)^{1/2} \|\theta\l\|}. 
\end{align*}
\end{lemma}

\begin{proof}[Proof of \cref{lem:lineartermtwoinfty}]
 Throughout this proof we suppress the dependence of $\bTheta\l,\Lambda\l$ and $\U\l$ on  the index $l$, and we denote $\lambda$ via $\frac{1}{\lambda} = \| (\Lambda\l)\inv \|.$   Define $\mathbf{E} := \A\m - \bP\m$, so that $\mathbf{E}$ is a mean-zero random matrix.

We will apply the Matrix Bernstein inequality to each row separately.  To wit, by Corollary 3.3 of \cite{chen_spectral_2021}, we have that with probability at least $1 - O(n^{-16})$ it holds that
\begin{align*}
     \|  e_i\t \mathcal{L}(\mathbf{E}) \| &\lesssim \sqrt{v \log(n)} +  w \log(n),
\end{align*}
where
\begin{align*}
    v &=\max\bigg\{ \bigg\| \sum_{j}\bigg( \U |\Lambda|^{-1/2} \ipq \mathbf{J}(\mathbf{X}_{i\cdot}) \bigg)_{j\cdot }\t  \mathbb{E} \mathbf{E}_{ij}^2 \bigg( \U |\Lambda|^{-1/2} \ipq \mathbf{J}(\mathbf{X}_{i\cdot}) \bigg)_{j\cdot } \bigg\|, \\
    &\qquad \qquad \bigg\| \sum_{j} \E \mathbf{E}_{ij}  \bigg( \U |\Lambda|^{-1/2} \ipq \mathbf{J}(\mathbf{X}_{i\cdot}) \bigg)_{j\cdot }  \bigg( \U |\Lambda|^{-1/2} \ipq \mathbf{J}(\mathbf{X}_{i\cdot}) \bigg)_{j\cdot }\t  \mathbf{E}_{ij} \bigg\| \bigg\}; \\
    w &= \max_{j} \bigg\| \mathbf{E}_{ij} \bigg( \U |\Lambda|^{-1/2} \ipq \mathbf{J}(\mathbf{X}_{i\cdot}) \bigg)_{j\cdot } \bigg\|.
\end{align*}
Since $\mathbb{E}(\mathbf{E}_{ij}^2)$ is a scalar, we have that by \cref{lem:popprop},
\begin{align*}
    v &\leq \sum_{j} \mathbb{E} \mathbf{E}_{ij}^2 \bigg\| \bigg(\U |\Lambda|^{-1/2} \ipq \mathbf{J}(\mathbf{X}_{i\cdot}) \bigg)_{j\cdot} \bigg\|^2 \\
    &\leq \sum_{j} \theta_i \theta_j \| e_j\t \U \|^2 \frac{1}{\lambda} \| \mathbf{J}(\mathbf{X}_{i\cdot}) \|^2 \\
    &\lesssim \sum_{j} \theta_i \theta_j \frac{\theta_j^2 K}{\|\theta\|^2} \frac{K}{\|\theta\|^2 \lambda_{\min}}  \frac{1}{\|\mathbf{X}_{i\cdot}\|^2} \\
    &\lesssim \sum_{j} \theta_i \theta_j \frac{\theta_j^2 K}{\|\theta\|^2} \frac{K}{\|\theta\|^2 \lambda_{\min}}\frac{1}{\theta_i^2 } \\ 
    &\lesssim \sum_{j} \frac{K^2}{\lambda_{\min} \|\theta\|^4} \bigg( \frac{\theta_j}{\theta_i} \bigg) \theta_j^2 \\
    &\lesssim \frac{K^2}{\lambda_{\min} \|\theta\|^4}\bigg( \frac{\theta_{\max}}{\theta_{\min}} \bigg) \|\theta\|^2 \\
    &\lesssim \frac{K^2}{\lambda_{\min} \|\theta\|^2}\bigg( \frac{\theta_{\max}}{\theta_{\min}} \bigg).
\end{align*}
Similarly,
\begin{align*}
    w &\leq \max_{j} \bigg\|\bigg( \U |\Lambda|^{-1/2} \ipq \mathbf{J}(\mathbf{X}_{i\cdot}) \bigg)_{j\cdot } \bigg\| \\
    &\lesssim \frac{\theta_{\max} \sqrt{K}}{\|\theta\|} \frac{\sqrt{K}}{\|\theta\| \lambda_{\min}^{1/2}}  \frac{1}{\theta_i } \\
    &\lesssim \bigg( \frac{\theta_{\max}}{\theta_{\min}}\bigg) \frac{K}{\|\theta\|^2 \lambda_{\min}^{1/2}}.
\end{align*}
Therefore, with probability at least $1 - O(n^{-16})$, we have that
\begin{align*}
     \|  e_i\t \mathcal{L}(\mathbf{E}) \| &\lesssim \bigg( \frac{\theta_{\max}}{\theta_{\min}} \bigg)^{1/2} \frac{K\sqrt{\log(n)}}{\lambda_{\min}^{1/2} \|\theta\|} + \bigg( \frac{\theta_{\max}}{\theta_{\min}} \bigg) \frac{K \log(n)}{\lambda_{\min}^{1/2} \|\theta\|^2} \\
     &\lesssim \bigg( \frac{\theta_{\max}}{\theta_{\min}} \bigg)^{1/2} \frac{K\sqrt{\log(n)}}{\lambda_{\min}^{1/2} \|\theta\|} \max\bigg\{1 , \bigg( \frac{\theta_{\max}}{\theta_{\min}} \bigg)^{1/2} \frac{\sqrt{\log(n)}}{\|\theta\|} \bigg\}.
\end{align*}
We now show that \cref{ass:networklevel} implies that 1 is the maximum above.  \cref{ass:networklevel} states that
\begin{align*}
  C  \frac{\theta_{\max}}{\theta_{\min}} \frac{K^{8} \theta_{\max}\|\theta\|_1 \log(n)}{\|\theta\|^4 (\lambda_{\min})^2} \leq \bar \lambda.
\end{align*}
Since $\bar \lambda \leq 1 $ by assumption and $\|\theta\|^2 \leq \theta_{\max}\|\theta\|_1$, it is straightforward to verify that \cref{ass:networklevel} implies that  
\begin{align*}
    \bigg( \frac{\theta_{\max}}{\theta_{\min}}\bigg)\frac{ \log(n)}{ \|\theta\l\|^2} \lesssim 1.
\end{align*}
Taking square roots reveals that
\begin{align*}
     \bigg( \frac{\theta_{\max}}{\theta_{\min}}\bigg)^{1/2}\frac{ \sqrt{\log(n)}}{ \|\theta\l\|} \lesssim 1,
\end{align*}
which shows that one is the dominant term in the maximum, as long as $C$ is larger than some universal constant.
Taking a union bound over all the rows completes the proof.
\end{proof}

\subsection{Preliminary Lemmas}
Throughout this section and its proof we suppress the dependence on $l$ in all terms.  We also let $\lambda$ denote the absolute value of the smallest nonzero eigenvalue of $\bP$. 
In what follows, we will assume that $\lambda \gtrsim \sqrt{\theta_{\max} \|\theta\|_1 \log(n)}$, which by \cref{lem:popprop} holds under \cref{ass:networklevel}. We will verify this explicitly at the beginning of the proof of \cref{thm:firststep}.  

The following result shows a form of spectral norm concentration.  

\begin{lemma}[Spectral Norm Concentration for One Graph]\label{lem:spectralnormconcentration}
When $\theta_{\max} \|\theta\|_1 \geq \log(n)$, it holds that
\begin{align*}
    \| \A - \bP  \| &\lesssim \sqrt{\theta_{\max} \| \theta\|_1}; \\
    \| \U\t (\A - \bP) \U \| &\lesssim \sqrt{K} + \sqrt{\log(n)},
\end{align*}
with probability at least $1 - O(n^{-20})$.
\end{lemma}

\begin{proof}
See Lemma C.1 of \cite{jin_estimating_2017}, or directly apply  Remark 3.13 from \cite{bandeira_sharp_2014}. 
The other part follows from a straightforward $\eps$-net argument.
\end{proof}

The following lemma demonstrates good concentration for several residual terms, showing that several terms ``approximately commute."  
\begin{lemma}[Approximate Commutation] \label{lem:approximatecommutation}
When $\lambda \gtrsim \sqrt{\theta_{\max} \|\theta\|_1 \log(n)}$ and $\min_{i} \theta_i \|\theta\|_1 \gtrsim\log(n)$, the following bounds hold with probability at least $1 - O(n^{-20}):$
\begin{align}
\| \wstar - \U\t \uhat \| &\lesssim \frac{\theta_{\max} \|\theta\|_1}{\lambda^2} \label{wstarapprox} \\
\| \uhat\t \U |\Lambda|^{1/2} - |\hat \Lambda|^{1/2} \uhat\t \U  \| &\lesssim \frac{K^2 
}{\lambda^{1/2}} \bigg( \sqrt{K \log(n) }+  \frac{\theta_{\max} \|\theta\|_1}{\lambda}  \bigg)  \label{approxcommute4} \\
\| \uhat\t \U |\Lambda|^{-1/2} \ipq - |\hat \Lambda|^{-1/2} \ipq \uhat \t \U \| &\lesssim \frac{K^2 
}{\lambda^{3/2}} \bigg( \sqrt{K \log(n) }+  \frac{\theta_{\max} \|\theta\|_1}{\lambda}  \bigg).\label{approxcommute3} 
\end{align}
\end{lemma}

\begin{proof}[Proof of \cref{lem:approximatecommutation}]
For \eqref{wstarapprox}, the argument follows since $\wstar$ is the product of the orthogonal matrices in the singular value decomposition of $\U\t \uhat$ and hence
\begin{align*}
    \| \wstar - \U\t \uhat \| &= \| \mathbf{I} - \cos\bTheta \| \\
    &\leq \| \sin\bTheta(\U, \uhat) \|^2 \\
    &\lesssim \frac{ \| \A - \bP\|^2}{\lambda^2} \\
    &\lesssim \frac{\theta_{\max} \|\theta\|_1}{\lambda^2}
\end{align*}
which holds with probability at least $1 - O(n^{-20})$ by \cref{lem:spectralnormconcentration}.

For all the following terms, we first show that $|\hat \Lambda| (\ipq \uhat\t \U - \uhat\t \U \ipq)$ is sufficiently small by modifying a similar argument to \citet{rubin-delanchy_statistical_2022}.  Observe that \begin{align*}
    \ipq \uhat\t \U - \uhat\t \U \ipq &= \begin{pmatrix} 0 & 2 \uhat_{+}\t \U_{-} \\ -2 \uhat_{-}\t \U_+ & 0 \end{pmatrix},
\end{align*}
where $\U_{+}$ denotes the eigenvectors of $\U$ corresponding to the positive eigenvectors (and similarly for $\U_{-}$, $\uhat_{+}$, and $\uhat_{-}$ respectively). Let $\mathbf{u}_j^+$ and $\mathbf{\hat u}_j^-$ denote the $j$'th columns of $\U_+$ and $\uhat_-$ respectively.  Then the $j,i$ entry of $\uhat_-\t \U_+$ is simply $(\mathbf{u}_{i}^{+})\t \mathbf{\hat u}_{j}^{-}$, and hence by the eigenvector-eigenvalue equation,
\begin{align*}
  (\mathbf{u}_{i}^{+})\t \mathbf{\hat u}_{j}^{-} &= \frac{(\mathbf{u}_i^+)\t (\A - \bP) \mathbf{\hat u}_{j,-}}{\hat \lambda_{j,-} - \lambda_{i,+}} \\
  &= \frac{(\mathbf{u}_i^+)\t (\A - \bP) \U_- \U_-\t \mathbf{\hat u}_{j,-}}{\hat \lambda_{j,-} - \lambda_{i,+}} + \frac{(\mathbf{u}_i^+)\t (\A - \bP)\big( \mathbf{I} -  \U_- \U_-\t\big) \mathbf{\hat u}_{j,-}}{\hat \lambda_{j,-} - \lambda_{i,+}},
  \end{align*}
  where $\lambda_{i,+}$ denotes the $i$'th largest in magnitude eigenvalue of $\bP$ (and similarly for $\hat \lambda_{j,-}$ for the negative eigenvalues of $\A$). It is straightforward to check that the $j,i$ entry of the matrix $|\hat \Lambda_+| (\ipq \uhat_+\t \U_{-} - \uhat_+\t \U_{-}\ipq )$ is given by
\begin{align*}
    \frac{|\hat \lambda_{j+}|}{\hat \lambda_{j,+} - \lambda_{i,-}} (\mathbf{u}_{i}^+)\t (\A - \bP) \mathbf{\hat u}_{j,-} &= \frac{1}{1 - \frac{\lambda_{i,-}}{\hat \lambda_{j,+}}} (\mathbf{u}_{i}^+)\t (\A - \bP) \U_+ \U_+\t \mathbf{\hat u}_{j,-} \\
    &\quad + \frac{1}{1 - \frac{\lambda_{i,-}}{\hat \lambda_{j,+}}} (\mathbf{u}_{i}^+)\t (\A - \bP)(\mathbf{I} - \U_+ \U_+\t) \mathbf{\hat u}_{j,-}
\end{align*}
Since $\lambda_{i,-}$ is negative and $\hat \lambda_{j,+}$ is positive with high probability, $1 - \frac{\lambda_{i,-}}{\hat \lambda_{j,+}}$ is strictly larger than one.  A similar argument holds for the entries with the ``$+$" changed to a ``$-$".

  Without loss of generality, consider the term corresponding to the negative eigenvalues.
  We can write the matrix as follows.  Denote $\mathbf{M}$ as the matrix whose $i,j$ entry is $\frac{|\hat \lambda_{j,-}|}{\hat \lambda_{j,-} - \lambda_{i,+}}$.  Then we have the equality
  \begin{align*}
      \U_+\t \uhat_- |\hat \Lambda_-| &= \mathbf{M} \circ \bigg( \U_+\t (\A - \bP) \U_- \U_-\t \uhat_- + \U_+\t (\A - \bP) \big(\mathbf{I} - \U_- \U_-\t \big) \uhat_-\bigg).
  \end{align*}
  Therefore,
  \begin{align}
      \| \U_+\t \uhat_- |\hat \Lambda_-| \| &\leq \| \mathbf{M} \| \bigg( \| \U_+\t (\A - \bP) \U_- \| + \| \U_+\t (\A - \bP) \big(\mathbf{I} - \U_- \U_-\t \big) \uhat_- \| \bigg) \nonumber \\
      &\leq \| \mathbf{M} \| \bigg( \| \U\t (\A - \bP) \U \| + \| \A - \bP \| \|  \big(\mathbf{I} - \U_- \U_-\t \big) \uhat_- \| \bigg), \label{jberg1000}
  \end{align}
  where we have used the fact that $\U_+\t (\A - \bP ) \U_-$ is a submatrix of $\U\t (\A - \bP) \U$.  We now note that 
  \begin{align*}
\|  \big(\mathbf{I} - \U_- \U_-\t \big) \uhat_- \|  &= \| \sin\bTheta(\U_-, \uhat_- ) \|.
  \end{align*}
  In addition, the eigenvalues corresponding to $\uhat_-$ are all negative, and the eigengap condition is satisfied since the eigenvalues corresponding to $(\I - \U_- \U_-) \bP$ are either all zero or positive.  Consequently, the eigengap satisfies
\begin{align*}
    \min_{\lambda_i > 0}  \lambda_i - \max_{p+1 \leq i \leq n} \hat \lambda_i \gtrsim \lambda
\end{align*}
by applying Weyl's inequality to the negative eigenvalues and the bottom $n - K$ eigenvalues separately.  We can therefore apply the Davis-Kahan Theorem to obtain
\begin{align}
    \| \sin\bTheta(\U_-, \uhat_- ) \| &\lesssim \frac{\| \A - \bP \|}{\lambda} \nonumber \\
    &\lesssim \frac{\sqrt{\theta_{\max} \|\theta\|_1}}{\lambda} \label{sinthetabound}
\end{align}
with probability at least $1 - O(n^{-20})$.  In addition, observe that the matrix $\mathbf{M}$ satisfies 
\begin{align}
    \|\mathbf{M} \| &\lesssim K. 
    \label{Mbound}
\end{align}
Finally, by \cref{lem:spectralnormconcentration}, we have that $\|\U\t (\A - \bP ) \U\| \lesssim \sqrt{K} + \sqrt{\log(n)}$ with high probability.  Plugging in this estimate, \eqref{Mbound}, and \eqref{sinthetabound} into \eqref{jberg1000} yields
\begin{align*}
      \| \U_+\t \uhat_- \| &\lesssim K \bigg( \sqrt{K} + \sqrt{\log(n)} +  \frac{\theta_{\max} \|\theta\|_1}{\lambda}  \bigg).
\end{align*}
Therefore, by applying a similar argument to $\U_-\t \uhat_+$, we obtain
\begin{align*}
    \| |\hat \Lambda| ( \ipq \uhat\t \U - \uhat\t \U \ipq) \| &\lesssim  K \bigg( \| \U\t (\A - \bP) \U \| + \| \A - \bP \| \| \sin\bTheta( \U_{-},\uhat_{-}) \| \\
    &\qquad + \| \A - \bP \| \| \sin\bTheta( \U_{+},\uhat_{+}) \| \bigg) \\
    &\lesssim K \bigg( \sqrt{K} + \sqrt{\log(n)} + \frac{\theta_{\max} \|\theta\|_1}{\lambda} \bigg), \numberthis \label{hatlambdaipqbound}
\end{align*}
which holds with probability at least $1 - O(n^{-20})$. 

We now bound \eqref{approxcommute4}.  First, note that we have
\begin{align*}
    \|  \uhat\t \U |\Lambda|^{1/2} - |\hat \Lambda|^{1/2}  \uhat \t \U \| &= \|  \uhat\t \U |\Lambda|^{1/2} \ipq - |\hat \Lambda|^{1/2}  \uhat \t \U \ipq \| \\
    &= \|  \uhat\t \U \ipq |\Lambda|^{1/2}  - |\hat \Lambda|^{1/2}  \uhat \t \U \ipq \|,
\end{align*}
where the first line follows since $\ipq$ is orthogonal and the second line follows since diagonal matrices commute. We observe that the $k,l$ entry of the matrix above can be written as
\begin{align*}
   \bigg(  \uhat\t \U \ipq|\Lambda|^{1/2}  - |\hat \Lambda|^{1/2} \uhat \t \U \ipq \bigg)_{kl} &= \langle \uhat_{\cdot k}, \U_{\cdot l} (\ipq)_{ll} \rangle \bigg( |\lambda_l|^{1/2} - |\hat \lambda_k|^{1/2}  \bigg) \\
   &= \langle \uhat_{\cdot k}, \U_{\cdot l}  \rangle (\ipq)_{ll} \frac{ |\lambda_l| - | \hat \lambda_k| }{|\lambda_l|^{1/2} + |\hat \lambda_k|^{1/2}  }.
\end{align*}
Define the matrix $\mathbf{H}$ via $\mathbf{H}_{kl} := \frac{1}{|\lambda_l|^{1/2} + |\hat \lambda_k|^{1/2}  }$.  Then the matrix $\uhat\t \U \ipq|\Lambda|^{1/2}  - |\hat \Lambda|^{1/2} \uhat \t \U \ipq$ can be written as
\begin{align*}
    \uhat\t \U \ipq|\Lambda|^{1/2}  - |\hat \Lambda|^{1/2} \uhat \t \U \ipq &= \mathbf{H} \circ \bigg(\uhat\t \U \ipq |\Lambda| -  |\hat \Lambda|  \uhat\t \U \ipq \bigg) \\
    &= \mathbf{H} \circ \bigg(\uhat\t \U \ipq |\Lambda| -  |\hat \Lambda| \ipq  \uhat\t \U  \bigg) \\
    &\quad + \mathbf{H} \circ \bigg(|\hat \Lambda| \bigg( \ipq  \uhat\t \U  - \uhat\t \U \ipq \bigg) \bigg) \\
    &= \mathbf{H} \circ \bigg(\uhat\t \U \Lambda -  \hat \Lambda  \uhat\t \U  \bigg) \\
    &\quad + \mathbf{H} \circ \bigg(|\hat \Lambda| \bigg( \ipq  \uhat\t \U  - \uhat\t \U \ipq \bigg) \bigg).
\end{align*}
where $\circ$ denotes the Hadamard product. It is straightforward to observe that $\|\mathbf{H} \| \lesssim \frac{K}{\lambda^{1/2}}$. Consequently, we have that
\begin{align*}
    \| \uhat\t \U \ipq | \Lambda |^{1/2} &- |\hat \Lambda|^{1/2} \uhat\t \U \ipq \|  \\
    &\leq \| \mathbf{H} \| \|\uhat\t \U \Lambda - \hat \Lambda \uhat\t \U \| + \|\mathbf{H} \| \| |\hat \Lambda|\big( \ipq \uhat\t \U - \uhat\t \U\ipq \big)\|  \\
    &\lesssim \frac{K}{\lambda^{1/2}} \| \uhat\t\U \Lambda - \hat \Lambda \uhat\t \U \| + \frac{K^2}{\lambda^{1/2}} \bigg( \sqrt{K} + \sqrt{\log(n)} + \frac{\theta_{\max} \|\theta\|_1}{\lambda} \bigg) \\
    &\lesssim \frac{K}{\lambda^{1/2}} \| \uhat\t \bP \U - \uhat\t\A \U \| + \frac{K^2}{\lambda^{1/2}}  \bigg( \sqrt{K} + \sqrt{\log(n)} + \frac{\theta_{\max} \|\theta\|_1}{\lambda} \bigg) \\
    &\lesssim \frac{K}{\lambda^{1/2}} \| \uhat\t (\bP - \A ) \U \| + \frac{K^2}{\lambda^{1/2}} \bigg( \sqrt{K} + \sqrt{\log(n)} + \frac{\theta_{\max} \|\theta\|_1}{\lambda} \bigg). \numberthis \label{1003}
\end{align*}
We note that
\begin{align*}
    \|\uhat\t (\bP - \A)  \U  \|  &\lesssim \| \uhat\t \U \U\t (\bP - \A) \U \| + \| \uhat\t (\I - \U \U\t) (\bP - \A) \U \| \\
    &\lesssim \sqrt{K} + \sqrt{\log(n)} + \| \sin\bTheta(\uhat,\U) \| \| \A - \bP \| \\
    &\lesssim \sqrt{K} + \sqrt{\log(n)} + \frac{\theta_{\max} \|\theta\|_1}{\lambda}.
\end{align*}
Plugging this into our bound \eqref{1003}, we obtain that
\begin{align*}
 \|  \uhat\t \U \ipq |\Lambda|^{1/2} - |\hat \Lambda|^{1/2} \uhat\t \U \ipq \| 
 &\lesssim \frac{K}{\lambda^{1/2}} \bigg( \sqrt{K} + \sqrt{\log(n)} + \frac{\theta_{\max} \|\theta\|_1}{\lambda} \bigg) \\
 &\quad + \frac{K^2}{\lambda^{1/2}}  \bigg( \sqrt{K} + \sqrt{\log(n)} + \frac{\theta_{\max} \|\theta\|_1}{\lambda} \bigg) \\
 &\lesssim \frac{K^2}{\lambda^{1/2}} \bigg( \sqrt{K\log(n)}  + \frac{\theta_{\max} \|\theta\|_1}{\lambda} \bigg).
\end{align*}
This proves \eqref{approxcommute4}. 

We now consider the term \eqref{approxcommute3}.  Since diagonal matrices commute,
\begin{align*}
     \| \uhat\t \U |\Lambda|^{-1/2} \ipq &- |\hat \Lambda|^{-1/2} \ipq \uhat \t \U \|\\ &= \| \uhat\t \U \ipq |\Lambda|^{-1/2} - |\hat \Lambda|^{-1/2} \ipq \uhat\t \U \| \\
    &= \|| \hat \Lambda|^{-1/2} \ipq \bigg( |\hat \Lambda|^{1/2} \ipq \uhat\t \U - \uhat\t \U \ipq |\Lambda|^{1/2} \bigg) \ipq |\Lambda|^{-1/2} \| \\
    &\lesssim \frac{1}{\lambda} \| |\hat \Lambda|^{1/2} \ipq \uhat\t \U - \uhat\t \U \ipq |\Lambda|^{1/2} \| \\
      &\lesssim \frac{1}{\lambda} \bigg(  \| |\hat \Lambda |^{1/2} (\ipq \uhat\t \U - \uhat\t \U \ipq ) \| + \| |\hat \Lambda|^{1/2} \uhat\t \U \ipq - \uhat\t \U \ipq |\Lambda|^{1/2} \| \bigg) \\
      &\lesssim \frac{1}{\lambda^{3/2}}  \| |\hat \Lambda| (\ipq \uhat\t \U - \uhat\t \U \ipq ) \| + \frac{1}{\lambda}  \| |\hat \Lambda|^{1/2} \uhat\t \U \ipq - \uhat\t \U \ipq |\Lambda|^{1/2} \| \\
      &\lesssim \frac{K}{\lambda^{3/2}} \bigg( \sqrt{K} + \sqrt{\log(n)} + \frac{\theta_{\max} \|\theta\|_1}{\lambda} \bigg) \\
      &\quad + \frac{K^2}{\lambda^{3/2}} \bigg( \sqrt{K\log(n)} + \frac{\theta_{\max} \|\theta\|_1}{\lambda} \bigg) \\
      &\lesssim \frac{K^2}{\lambda^{3/2}} \bigg( \sqrt{K\log(n)} + \frac{\theta_{\max} \|\theta\|_1}{\lambda} \bigg).
\end{align*}
where we have implicitly used the bound \eqref{approxcommute4} and \eqref{1003}. 
This bound holds cumulatively with probability at least $1 - O(n^{-20})$, which completes the proof. 
\end{proof}

The following lemma characterizes the row-wise concentration of terms that involve $\uhat$.  However, this proof requires the use of leave-one-out sequences, so we defer its proof to \cref{sec:rowwiseproof} after the proof of \cref{thm:firststep}.  

\begin{restatable}[Row-Wise Concentration I]{lemma}{rowwiseerror}\label{lem:rowwiseerror1}
When $\lambda \gtrsim \sqrt{\theta_{\max} \|\theta\|_1 \log(n)}$ and $\min_{i} \theta_i \|\theta\|_1 \gtrsim\log(n)$, it holds that
\begin{align*}
\| e_i\t (\mathbf{A} - \E \mathbf{A}) \uhat \| &\lesssim \sqrt{\theta_i \| \theta\|_1 \log(n)} \| \uhat\|_{2,\infty} 
\end{align*}
\end{restatable}

The following result demonstrates that $\uhat$ is sufficiently close to  $\U \U\t \uhat$ in $\|\cdot\|_{2,\infty}$.

\begin{restatable}[Closeness of $\uhat$ to $\U$]{lemma}{uhatucloseness}\label{lem:uhatucloseness}
When $\lambda \gtrsim \sqrt{\theta_{\max} \|\theta\|_1 \log(n)}$ and $\min_{i} \theta_i \|\theta\|_1 \gtrsim\log(n)$, the following bounds holds with probability at least $1 - O(n^{-19}):$
\begin{align*}
    \| \uhat - \U \U\t \uhat \|_{2,\infty} &\lesssim  \frac{ \sqrt{\theta_{\max} \|\theta\|_1 \log(n)}}{\lambda} \| \U \|_{2,\infty}\\
    \| \uhat - \U \wstar \|_{2,\infty} &\lesssim  \frac{ \sqrt{\theta_{\max} \|\theta\|_1 \log(n)}}{\lambda} \| \U \|_{2,\infty}\\
    \| \uhat \uhat\t \U - \U \|_{2,\infty} &\lesssim  \frac{ \sqrt{\theta_{\max} \|\theta\|_1 \log(n)}}{\lambda} \| \U \|_{2,\infty}.
\end{align*}
\end{restatable}
The bound above matches the bound in \citet{jin_estimating_2017}, Lemma D.2.  

\begin{proof}[Proof of \cref{lem:uhatucloseness}]
Observe that since $\U$ are the eigenvectors of $\bP$ and $\bP$ is rank $K$,
\begin{align*}
    e_i\t \big( \uhat - \U \U\t \uhat \big) &= e_i\t ( \mathbf{I} - \U \U\t )\uhat \\
    &= e_i\t ( \mathbf{I} - \U \U\t ) \A \uhat \hat \Lambda\inv \\
    &= e_i\t ( \mathbf{I} - \U \U\t ) ( \A - \bP) \uhat \hat \Lambda\inv \\
    &= e_i\t (\A - \bP) \uhat \hat \Lambda\inv - e_i\t \U \U\t (\A - \bP) \uhat \hat \Lambda\inv.
\end{align*}
Taking norms reveals that
\begin{align*}
    \| e_i\t \big( \uhat - \U \U\t \uhat \big) \| &\leq \| e_i\t \big( \A - \bP \big) \uhat \| \|\hat \Lambda\inv \| + \| e_i\t \U \| \| \A - \bP \| \|\hat \Lambda\inv \|.
\end{align*}
By \cref{lem:spectralnormconcentration}, we have that $\| \A - \bP \| \lesssim \sqrt{\theta_{\max} \|\theta\|_1}$.  In addition, Weyl's inequality implies that $\| \hat \Lambda\inv\| \lesssim \lambda\inv$.  Therefore, combining these bounds with \cref{lem:rowwiseerror1}, we see that with probability at least $1 - O(n^{-20})$ that
\begin{align*}
     \| e_i\t \big( \uhat - \U \U\t \uhat \big) \| &\lesssim \frac{ \sqrt{\theta_i \|\theta\|_1 \log(n)}}{\lambda} \|\uhat\|_{2,\infty} + \| e_i\t \U \| \frac{\sqrt{\theta_{\max} \|\theta\|_1}}{\lambda} \\
     &\lesssim \frac{ \sqrt{\theta_{\max} \|\theta\|_1 \log(n)}}{\lambda} \|\uhat\|_{2,\infty} + \|  \U \|_{2,\infty} \frac{\sqrt{\theta_{\max} \|\theta\|_1}}{\lambda}.
\end{align*}
This bound is independent of row $i$, so taking a union bound reveals that with probability at least $1 - O(n^{-19})$ that
\begin{align*}
    \| \uhat - \U \U\t \uhat \|_{2,\infty} &\lesssim  \frac{ \sqrt{\theta_{\max} \|\theta\|_1 \log(n)}}{\lambda} \|\uhat\|_{2,\infty} + \|  \U \|_{2,\infty} \frac{\sqrt{\theta_{\max} \|\theta\|_1}}{\lambda}.
\end{align*}
By \cref{lem:approximatecommutation}, it holds that
\begin{align*}
    \| \wstar - \U\t \uhat \| &\lesssim \frac{ \theta_{\max} \|\theta\|_1}{\lambda^2}.
\end{align*}
Therefore,
\begin{align*}
    \| \uhat - \U \wstar \|_{2,\infty} &\leq \| \uhat - \U \U\t \uhat \|_{2,\infty} + \| \U \|_{2,\infty} \| \wstar - \U\t \uhat \| \\
    &\lesssim \frac{ \sqrt{\theta_{\max} \|\theta\|_1 \log(n)}}{\lambda} \|\uhat\|_{2,\infty} + \|  \U \|_{2,\infty} \frac{\sqrt{\theta_{\max} \|\theta\|_1}}{\lambda} + \frac{ \theta_{\max} \|\theta\|_1}{\lambda^2} \| \U\|_{2,\infty} \\
    &\lesssim \frac{ \sqrt{\theta_{\max} \|\theta\|_1 \log(n)}}{\lambda} \|\uhat\|_{2,\infty}   + \|  \U \|_{2,\infty} \frac{\sqrt{\theta_{\max} \|\theta\|_1}}{\lambda}.
\end{align*}
As a byproduct, this also reveals that
\begin{align*}
    \|\uhat\|_{2,\infty} &\leq \| \uhat - \U \wstar \|_{2,\infty} + \| \U \|_{2,\infty} \\
    &\leq \frac{1}{2} \| \uhat \|_{2,\infty} + \frac{3}{2} \| \U \|_{2,\infty},
\end{align*}
as long as $\lambda \geq C \sqrt{\theta_{\max} \|\theta\|_1 \log(n)}$ for some sufficiently large constant $C$ (which we verify at the beginning of the proof of \cref{thm:firststep}, and which holds under \cref{ass:networklevel}).  By rearranging, it holds that $\|\uhat \|_{2,\infty} \lesssim \| \U \|_{2,\infty}$.  Plugging this in yields
\begin{align*}
    \| \uhat - \U \wstar \|_{2,\infty} &\lesssim \frac{ \sqrt{\theta_{\max} \|\theta\|_1 \log(n)}}{\lambda} \| \U \|_{2,\infty}; \\
    \| \uhat - \U \U\t \uhat \|_{2,\infty} &\lesssim \frac{ \sqrt{\theta_{\max} \|\theta\|_1 \log(n)}}{\lambda} \| \U \|_{2,\infty}.
\end{align*}
The final inequality holds since
\begin{align*}
    \| \uhat \uhat\t \U - \U \|_{2,\infty} &\leq \| \uhat \uhat\t \U - \uhat \wstar\t \|_{2,\infty} + \| \uhat - \U \wstar \|_{2,\infty} \\
    &\leq \| \uhat \|_{2,\infty} \| \uhat \t \U - \wstar\t \| + \| \uhat - \U \wstar \|_{2,\infty} \\
    &\lesssim \|\U \|_{2,\infty} \| \wstar - \U\t \uhat \| + \| \uhat - \U \wstar\|_{2,\infty}.
\end{align*}
The proof is completed by plugging in the previous bounds.
\end{proof}

The following result establishes finer control over the rows of the estimated eigenvectors.  We relegate the proof of his result to \cref{sec:rowwiseproof}, since it requires the use of leave-one-out sequences.

\begin{restatable}[Row-wise Concentration II]{lemma}{rowwiseerrortwo}\label{lem:rowwiseerror2} When $\lambda \gtrsim \sqrt{\theta_{\max} \|\theta\|_1 \log(n)}$ and $\min_{i} \theta_i \|\theta\|_1 \gtrsim\log(n)$, with probability at least $1 - O(n^{-19})$, it holds that
\begin{align*}
    \| e_i\t (\mathbf{A} - \E \mathbf{A}) (\uhat \uhat \t \U - \U) \| &\lesssim \sqrt{\theta_i \|\theta\|_1 \log(n)}\|\frac{\sqrt{\theta_{\max} \|\theta\|_1 \log(n)}}{\lambda} \| \U \|_{2,\infty}.
\end{align*}
\end{restatable}

\subsection{Proof of Theorem~\ref{thm:firststep}}
\begin{proof}[Proof of \cref{thm:firststep}]
Throughout the proof we suppress the dependence of these terms on the index $l$. Our proof proceeds in several steps: first, we express $\xhat \wstar - \xtilde$ as a linear term plus a residual term, where the residual term obeys a strong row-wise concentration bound.  Next, we demonstrate that the rows of $\yhat$ (i.e. the normalized rows of $\xhat$) concentrate about the corresponding rows of $\ytilde$.   Before embarking on the proof, we make note of several preliminary facts.  By \cref{lem:popprop}, we have that
\begin{align*}
    \lambda &\gtrsim \frac{\|\theta\|^2 \lambda_{\min}}{K}; \\
    \| e_i\t \U \| &\lesssim \frac{\sqrt{K}\theta_i}{\|\theta\|}; \\
   \theta_i \lesssim  \| e_i\t \xtilde \| &\leq \theta_i {\sqrt{K}}.
\end{align*}
We will use these bounds repeatedly without reference when simplifying our results.  

In addition, many of the previous lemmas require that $\lambda \gtrsim \sqrt{\theta_{\max} \|\theta\|_1 \log(n)}$.  We verify that this condition holds under \cref{ass:networklevel} now.
\cref{ass:networklevel} requires that 
\begin{align*}
   C \bigg( \frac{\theta_{\max}}{\theta_{\min}} \bigg) \frac{K^{8} \theta_{\max} \|\theta\|_1 \log(n)}{\|\theta\|^4 \lambda_{\min}^2} \leq  \bar \lambda. \numberthis \label{jjj}
\end{align*}
By \cref{lem:popprop} it holds that
\begin{align*}
    \lambda \gtrsim \frac{ \|\theta\|^2}{K} \lambda_{\min}.
\end{align*}
Consequently, it suffices to argue that \eqref{jjj} implies the condition
\begin{align*}
     \frac{ \|\theta\|^2}{K} \lambda_{\min} \gtrsim \sqrt{\theta_{\max} \|\theta\|_1 \log(n)},
\end{align*}
or equivalently,
\begin{align*}
    \frac{K \sqrt{\theta_{\max} \|\theta\|_1 \log(n)}}{\|\theta\|^2 \lambda_{\min}} \lesssim 1.
\end{align*}
Squaring both sides yields the condition
\begin{align*}
    \frac{K^2 \theta_{\max} \|\theta\|_1 \log(n)}{\|\theta\|^4 \lambda_{\min}^2} \lesssim 1.
\end{align*}
This is weaker than \eqref{jjj} as $\lambda_{\min}, \bar \lambda \in (0,1)$ by assumption and $K \geq 1$, as long as $C$ is larger than some universal constant.
\\ \ \\ \noindent
\textbf{Step 1: First-Order Approximation of $\xhat$:}\\ 
At the outset we recall that $\wstar$ is the Frobenius-optimal matrix aligning $\uhat$ and $\U$.  Moreover, by the concentration inequality in \cref{lem:spectralnormconcentration} and the assumption on the eigenvalue $\lambda$ above, we have that $\|\hat \Lambda\inv\| \lesssim \lambda\inv$ with probability at least $1 - O(n^{-20})$.  We now expand via:
\begin{align*}
\xhat \wstar\t - \xtilde &= (\mathbf{A} - \E \mathbf{A}) \mathbf{U} |\Lambda|^{-1/2} \ipq + \mathbf{R};\\
\mathbf{R}&=\mathbf{R}_1\wstar\t+\mathbf{R}_2\wstar\t+\mathbf{R}_3\wstar\t+\mathbf{R}_4+\mathbf{R}_5+\mathbf{R}_6;\\
\mathbf{R}_1 :&= - \U \U\t (\mathbf{A} - \E \mathbf{A}) \uhat |\hat \Lambda|^{-1/2} \ipq; \\
\mathbf{R}_2 :&= \U (\U\t \uhat |\hat \Lambda|^{1/2}  -  |\Lambda|^{1/2} \U\t \uhat ); \\
\mathbf{R}_3 :&= \U |\Lambda|^{1/2}(\U\t \uhat - \wstar); \\
\mathbf{R}_4 :&= (\mathbf{A} - \E \mathbf{A}) (\uhat \uhat \t \U - \U) |\Lambda|^{-1/2} \ipq; \\
\mathbf{R}_5 :&= -(\mathbf{A} - \E \mathbf{A})\uhat (\uhat \t \U |\Lambda|^{-1/2} \ipq - |\hat \Lambda|^{-1/2} \ipq \uhat \t \U); \\
\mathbf{R}_6 :&= (\mathbf{A} - \E \mathbf{A}) \uhat |\hat \Lambda|^{-1/2} \ipq (\wstar\t - \uhat \t \U).
\end{align*}
We now bound each residual in turn.   We will also use \cref{lem:approximatecommutation}, \cref{lem:rowwiseerror1}, \cref{lem:uhatucloseness}, \cref{lem:rowwiseerror2} repeatedly without reference; the cumulative probability will be at least $1 - O(n^{-18})$. 
\\ \ \\ \noindent
\textbf{The term $\mathbf{R}_1$:}\\
First, we note that
\begin{align*}
    \| e_i\t \mathbf{R}_1\| &\leq \| e_i\t \U\| \| \U\t ( \A - \bP ) \uhat |\hat \Lambda|^{-1/2} \| \\
    &\lesssim \frac{\| e_i\t \U \|}{\lambda^{1/2}} \bigg( \| \U\t ( \A - \bP) \U \| + \| \A - \bP \| \| \U_{\perp}\t \uhat \| \bigg) \\
    &\lesssim \frac{\| e_i\t \U \|}{\lambda^{1/2}}  \bigg( \| \U\t ( \A - \bP) \U \| + \frac{\| \A - \bP \|^2}{\lambda} \bigg).
\end{align*}
By \cref{lem:spectralnormconcentration}, we have that $\|\U\t (\A - \bP) \U \| \lesssim \sqrt{K} + \sqrt{\log(n)}$ with probability at least $1 - O(n^{-20})$.  Consequently,
\begin{align*}
     \| e_i\t \mathbf{R}_1\| &\lesssim \frac{\| e_i\t \U \|}{\lambda^{1/2}}  \bigg( \sqrt{K} + \sqrt{\log(n)} + \frac{\theta_{\max} \|\theta\|_1}{\lambda} \bigg).
\end{align*}
By \cref{lem:popprop}, we have that $\| e_i\t \U \| \lesssim \frac{\sqrt{K} \theta_i}{\|\theta\|}$ and that $\lambda \gtrsim \frac{\|\theta\|^2}{K} \lambda_{\min}$.  Putting it together, we arrive at the bound
\begin{align}
    \| e_i\t \mathbf{R}_1 \| &\lesssim \frac{K \theta_i}{\|\theta\|^2 \lambda_{\min}^{1/2}}  \bigg( \sqrt{K} + \sqrt{\log(n)} + \frac{K \theta_{\max} \|\theta\|_1}{\|\theta\|^2 \lambda_{\min}} \bigg) \nonumber \\
    &\lesssim \frac{K \theta_i}{\|\theta\|^2 \lambda_{\min}^{1/2}}  \bigg( \sqrt{K\log(n)} + \frac{K \theta_{\max} \|\theta\|_1}{\|\theta\|^2 \lambda_{\min}} \bigg) . \label{eq:r1bound}
\end{align}
\textbf{The term $\mathbf{R}_2$:} \\
We have
\begin{align}
     \| e_i\t \mathbf{R}_2 \| &\lesssim \| e_i\t \U \| \| \U\t\uhat |\hat \Lambda|^{1/2} - |\Lambda|^{1/2} \U\t \uhat \| \nonumber \\
     &\lesssim \frac{\sqrt{K} \theta_i}{\|\theta\|} \|  \U\t\uhat |\hat \Lambda|^{1/2} - |\Lambda|^{1/2} \U\t \uhat \| \nonumber \\
     &\lesssim \frac{\sqrt{K} \theta_i}{\|\theta\|} \frac{K^2}{\lambda^{1/2}} \bigg( \sqrt{K\log(n)} + \frac{\theta_{\max} \|\theta\|_1}{\lambda} \bigg) \nonumber \\
     &\lesssim \frac{\sqrt{K}\theta_i}{\|\theta\|} \frac{K^{5/2}}{\lambda_{\min}^{1/2} \|\theta\|} \bigg( \sqrt{K\log(n)} + \frac{K \theta_{\max} \|\theta\|_1}{ \|\theta\|^2\lambda_{\min}} \bigg) \nonumber \\
     &\asymp \frac{K^3 \theta_i}{\|\theta\|^2 \lambda_{\min}^{1/2}} \bigg( \sqrt{K\log(n)} + \frac{K\theta_{\max} \|\theta\|_1}{\|\theta\|^2 \lambda_{\min}} \bigg). \label{eq:r2bound}
\end{align}

\noindent
\textbf{The term $\mathbf{R}_3$:}\\
Following similarly as the previous step, we have that
\begin{align}
    \| e_i\t \mathbf{R}_3\| &\lesssim \| e_i\t \xtilde \| \| \U\t \uhat - \mathbf{W}_* \| \nonumber \\
    &\lesssim \theta_i  \frac{\theta_{\max} \|\theta\|_1}{\lambda^2} \nonumber \\
    &\lesssim \theta_i  {\sqrt{K}}  \frac{K^2\theta_{\max} \|\theta\|_1}{\|\theta\|^4 \lambda_{\min}^2}. 
   \label{eq:r3bound}
\end{align}
\textbf{The term $\mathbf{R}_4$:} \\
By \cref{lem:rowwiseerror2}, we have
\begin{align}
    \| e_i\t \mathbf{R}_4 \| &\lesssim  \| e_i\t (\A - \bP) (\uhat \uhat\t \U - \U) \|  \| |\Lambda|^{-1/2} \| \nonumber \\
    &\lesssim \frac{\sqrt{\theta_i \|\theta\|_1 \log(n)}}{\lambda^{1/2}} \frac{\sqrt{\theta_{\max} \|\theta\|_1 \log(n)}}{\lambda} \| \U \|_{2,\infty} \nonumber \\
    &\lesssim \frac{\sqrt{ K\theta_i \|\theta\|_1 \log(n)}}{\lambda_{\min}^{1/2} \|\theta\|} \frac{K\sqrt{\theta_{\max} \|\theta\|_1 \log(n)}}{\|\theta\|^2 \lambda_{\min}} \frac{\sqrt{K} \theta_{\max}}{\|\theta\|} \nonumber \\
    &\asymp \frac{ \theta_i^{1/2} K^{2} \theta_{\max}^{3/2} \|\theta\|_1 \log(n)}{\lambda_{\min}^{3/2} \|\theta\|^4}. \label{eq:r4bound}
\end{align}
\\ \ \\ \noindent
\textbf{The term $\mathbf{R}_5$:}\\
By \cref{lem:rowwiseerror1} and \cref{lem:approximatecommutation}, we have that
\begin{align*}
    \| e_i\t \mathbf{R}_5 \| &\lesssim \| e_i\t (\A - \bP ) \uhat \| \| \uhat\t \U |\Lambda|^{-1/2} \ipq - | \hat \Lambda|^{-1/2} \ipq \uhat\t \U \| \\
    &\lesssim \| e_i\t (\A - \bP ) \uhat \| \frac{K^2 
    }{\lambda^{3/2}} \bigg( \sqrt{K \log(n) }+  \frac{\theta_{\max} \|\theta\|_1}{\lambda}  \bigg) \\
    &\lesssim \sqrt{\theta_i \|\theta\|_1 \log(n)} \|\uhat\|_{2,\infty} \frac{K^2 
    }{\lambda^{3/2}} \bigg( \sqrt{K \log(n) }+  \frac{\theta_{\max} \|\theta\|_1}{\lambda}  \bigg).
    \end{align*}
By \cref{lem:uhatucloseness}, we have that $\|\uhat\|_{2,\infty} \lesssim \|\U\|_{2,\infty}$ as long as $\lambda \gtrsim \sqrt{\theta_{\max} \|\theta\|_1 \log(n)}$, which is true by \cref{ass:networklevel}.  Therefore,
\begin{align*}
    \| e_i\t \mathbf{R}_5 \| &\lesssim \sqrt{\theta_i \|\theta\|_1 \log(n)} \|\U\|_{2,\infty}\frac{K^2
    }{\lambda^{3/2}} \bigg( \sqrt{K \log(n) }+  \frac{\theta_{\max} \|\theta\|_1}{\lambda}  \bigg) \nonumber \\ 
    &\lesssim \sqrt{\theta_i \|\theta\|_1 \log(n)} \frac{\sqrt{K} \theta_{\max}}{\|\theta\|} \frac{K^2
    }{\lambda^{3/2}} \bigg( \sqrt{K \log(n) }+  \frac{\theta_{\max} \|\theta\|_1}{\lambda}  \bigg) \nonumber \\
        &\asymp \sqrt{\theta_i \|\theta\|_1 \log(n)} \frac{\sqrt{K} \theta_{\max}}{\|\theta\|}  \frac{K^{7/2}}{\lambda_{\min}^{3/2} \|\theta\|^3 } \bigg( \sqrt{K\log(n)} + \frac{K \theta_{\max} \|\theta\|_1}{\lambda_{\min} \|\theta\|^2}\bigg) \nonumber \\
    &\asymp \frac{\theta_i^{1/2}  \sqrt{\|\theta\|_1 \log(n)} \theta_{\max} K^{4}}{\lambda_{\min}^{3/2} \|\theta\|^4} \bigg( \sqrt{K\log(n)} + \frac{K \theta_{\max} \|\theta\|_1}{\lambda_{\min} \|\theta\|^2} \bigg). \numberthis
    \label{eq:r5bound}
\end{align*}  
    
\noindent
\textbf{The term $\mathbf{R}_6$:}\\
Similarly to the previous term, we obtain
\begin{align}
    \|e_i\t\mathbf{R}_6 \| &\lesssim \frac{\| e_i\t (\A - \bP) \uhat \|}{\lambda^{1/2}} \| \uhat\t \U - \wstar \| \nonumber \\
    &\lesssim \frac{\sqrt{\theta_i \|\theta\|_1 \log(n)} \| \U \|_{2,\infty} }{\lambda^{1/2}} \frac{\theta_{\max} \|\theta\|_1}{\lambda^2} \nonumber \\
    &\lesssim \frac{\sqrt{\theta_i \|\theta\|_1 \log(n)} \sqrt{K} \theta_{\max} }{\lambda^{1/2} \|\theta\|} \frac{\theta_{\max} \|\theta\|_1}{\lambda^2} \nonumber \\
    &\asymp \frac{\theta_i^{1/2} \| \theta\|_1^{3/2} \theta_{\max}^{2} \sqrt{\log(n)}K^3}{\lambda_{\min}^{5/2} \|\theta\|^6} \label{eq:r6bound}
\end{align}
\textbf{Putting it together:}\\
By \eqref{eq:r1bound}, \eqref{eq:r2bound}, \eqref{eq:r3bound}, \eqref{eq:r4bound}, \eqref{eq:r5bound}, and \eqref{eq:r6bound}, we obtain that
\begin{align*}
    \| e_i\t \mathbf{R}_1 \| &\lesssim \theta_i \frac{K^{3/2}}{\|\theta\|^2 \lambda_{\min}^{1/2}}  \bigg( \sqrt{K\log(n)} + \frac{K \theta_{\max} \|\theta\|_1}{\|\theta\|^2 \lambda_{\min}} \bigg); \\
    \| e_i\t \mathbf{R}_2 \| &\lesssim \frac{K^3 \theta_i}{\|\theta\|^2 \lambda_{\min}^{1/2}} \bigg( \sqrt{K\log(n)} + \frac{K\theta_{\max} \|\theta\|_1}{\|\theta\|^2 \lambda_{\min}} \bigg); \\
    \| e_i\t \mathbf{R}_3\| &\lesssim \theta_i {\sqrt{K}} \frac{K^2\theta_{\max} \|\theta\|_1}{\|\theta\|^4 \lambda_{\min}^2}; \\
    \| e_i\t \mathbf{R}_4 \| &\lesssim \theta_i^{1/2} \frac{ K^{2} \theta_{\max}^{3/2} \|\theta\|_1 \log(n)}{\lambda_{\min}^{3/2} \|\theta\|^4}; \\ 
    \| e_i\t \mathbf{R}_5 \| &\lesssim \frac{\theta_i^{1/2}  \sqrt{\|\theta\|_1 \log(n)} \theta_{\max} K^{4}}{\lambda_{\min}^{3/2} \|\theta\|^4} \bigg( \sqrt{K\log(n)} + \frac{K \theta_{\max} \|\theta\|_1}{\lambda_{\min} \|\theta\|^2} \bigg); \\
    \| e_i\t \mathbf{R}_6 \| &\lesssim  \theta_i^{1/2}\frac{K^3 \| \theta\|_1^{3/2} \theta_{\max}^{2} \sqrt{\log(n)}}{\lambda_{\min}^{5/2} \|\theta\|^6}.
\end{align*} 
We now group these terms for simplicity.  First, observe that the bound for $\| e_i\t \mathbf{R}_1\|$ is no more than the bound for $\|e_i\t \mathbf{R}_2\|$ since $\lambda_{\min} < 1$ and $K \geq 2$.  Therefore,
\begin{align*}
    \| e_i\t \mathbf{R}_1 \| + \| e_i\t \mathbf{R}_2 \| + \| e_i\t \mathbf{R}_3 \| &\lesssim \theta_i \bigg( \frac{K^{7/2} \sqrt{\log(n)}}{\|\theta\|^2 \lambda_{\min}^{1/2}} + \frac{K^4 \theta_{\max} \|\theta\|_1}{\|\theta\|^4 \lambda_{\min}^{3/2}} + {\sqrt{K}}\frac{K^2 \theta_{\max} \|\theta\|_1}{\|\theta\|^4 \lambda_{\min}^2} \bigg) 
\end{align*}
  We now simplify the remaining terms; i.e., the terms $\mathbf{R}_4$ through $\mathbf{R}_6$.  We observe that
\begin{align*}
    \| e_i\t \mathbf{R}_4\|& + \| e_i\t \mathbf{R}_5\| +\| e_i\t \mathbf{R}_6\| \\ &\lesssim 
  \theta_i^{1/2} \frac{K^2 \theta_{\max}^{3/2} \|\theta\|_1 \log(n)}{\lambda_{\min}^{3/2} \|\theta\|^4} + \theta_i^{1/2} \frac{K^{4} \sqrt{\|\theta\|_1 \log(n)} \theta_{\max} }{\lambda_{\min}^{3/2} \|\theta\|^4} \bigg( \sqrt{K\log(n)} + \frac{K \theta_{\max} \|\theta\|_1}{\lambda_{\min} \|\theta\|^2} \bigg) \\
    &\qquad + \theta_i^{1/2}\frac{K^3 \| \theta\|_1^{3/2} \theta_{\max}^{2} \sqrt{\log(n)}}{\lambda_{\min}^{5/2} \|\theta\|^6} \\
    &\lesssim \theta_i^{1/2} \frac{K^2 \theta_{\max}^{3/2} \|\theta\|_1 \log(n)}{\lambda_{\min}^{3/2} \|\theta\|^4} + \theta_i^{1/2} \frac{K^{9/2} \sqrt{\|\theta\|_1} \theta_{\max} \log(n)}{\lambda_{\min}^{3/2} \|\theta\|^4} \\
    &\quad + \theta_i^{1/2} \frac{K^5 \|\theta\|_1^{3/2} \theta_{\max}^2 \sqrt{\log(n)}}{\lambda_{\min}^{5/2} \|\theta\|^6} + \theta_i^{1/2} \frac{K^3 \|\theta\|_1^{3/2} \theta_{\max}^2 \sqrt{\log(n)}}{\lambda_{\min}^{5/2} \|\theta\|^6} \\
    &\lesssim \theta_i^{1/2} \frac{K^2 \theta_{\max}^{3/2} \|\theta\|_1 \log(n)}{\lambda_{\min}^{3/2} \|\theta\|^4} + \theta_i^{1/2} \frac{K^{9/2} \sqrt{\|\theta\|_1} \theta_{\max} \log(n)}{\lambda_{\min}^{3/2} \|\theta\|^4}  + \theta_i^{1/2} \frac{K^5 \|\theta\|_1^{3/2} \theta_{\max}^2 \sqrt{\log(n)}}{\lambda_{\min}^{5/2} \|\theta\|^6} \\
    &\lesssim (\theta_i \theta_{\max})^{1/2}\bigg( \frac{K^2 \theta_{\max} \|\theta\|_1 \log(n)}{\lambda_{\min}^{3/2} \|\theta\|^4} +  \frac{K^{9/2} \sqrt{\|\theta\|_1 \theta_{\max}} \log(n)}{\lambda_{\min}^{3/2} \|\theta\|^4}  + \frac{K^5 \|\theta\|_1^{3/2} \theta_{\max}^{3/2} \sqrt{\log(n)}}{\lambda_{\min}^{5/2} \|\theta\|^6} \bigg)\\ 
    &\lesssim \theta_i \bigg( \frac{\theta_{\max}}{\theta_{\min}} \bigg)^{1/2} \bigg(  \frac{K^2 \theta_{\max} \|\theta\|_1 \log(n)}{\lambda_{\min}^{3/2} \|\theta\|^4} +  \frac{K^{9/2} \sqrt{\|\theta\|_1 \theta_{\max}} \log(n)}{\lambda_{\min}^{3/2} \|\theta\|^4}  + \frac{K^5 \|\theta\|_1^{3/2} \theta_{\max}^{3/2} \sqrt{\log(n)}}{\lambda_{\min}^{5/2} \|\theta\|^6} \bigg)\\
    &\lesssim \theta_i \bigg( \frac{\theta_{\max}}{\theta_{\min}} \bigg)^{1/2} \bigg( \frac{K^{9/2} \theta_{\max} \|\theta\|_1 \log(n)}{\lambda_{\min}^{3/2} \|\theta\|^4} \bigg),
\end{align*}
where we have used the fact that $\lambda_{\min} \|\theta\|^2 \gtrsim  K\sqrt{\theta_{\max} \|\theta\|_1 \log(n)}$ and $\theta_{\max} \|\theta\|_1 \gtrsim \log(n)$, the first of which we verified at the beginning of this proof and the second by \cref{ass:networklevel}.  Putting these together, we arrive at
\begin{align*}
       \| e_i\t \mathbf{R} \| &\lesssim \theta_i \bigg( \frac{K^{7/2} \sqrt{\log(n)}}{\|\theta\|^2 \lambda_{\min}^{1/2}} + \frac{K^4 \theta_{\max} \|\theta\|_1}{\|\theta\|^4 \lambda_{\min}^{3/2}} +{\sqrt{K}}\frac{K^2 \theta_{\max} \|\theta\|_1}{\|\theta\|^4 \lambda_{\min}^2} \bigg) \\
       &\quad + \theta_i \bigg( \frac{\theta_{\max}}{\theta_{\min}} \bigg)^{1/2} \bigg( \frac{K^{9/2} \theta_{\max} \|\theta\|_1 \log(n)}{\lambda_{\min}^{3/2} \|\theta\|^4} \bigg) \\
       &\lesssim \theta_i\frac{K^{7/2} \sqrt{\log(n)}}{\|\theta\|^2 \lambda_{\min}^{1/2}} + \theta_i {\sqrt{K}}\frac{K^2 \theta_{\max} \|\theta\|_1}{\|\theta\|^4 \lambda_{\min}^2} + \theta_i \bigg( \frac{\theta_{\max}}{\theta_{\min}} \bigg)^{1/2} \bigg( \frac{K^{9/2} \theta_{\max} \|\theta\|_1 \log(n)}{\lambda_{\min}^{3/2} \|\theta\|^4} \bigg)
\end{align*}

\noindent
Consequently, we see that with probability at least $1 - O(n^{-18})$, each row $i$ of $\xhat$ satisfies
\begin{align*}
    e_i\t(\xhat \wstar\t - \xtilde) &= e_i\t (\A - \bP) \U |\Lambda|^{-1/2} \ipq + e_i\t \mathbf{R},
\end{align*}
where $\mathbf{R}$ satisfies 
\begin{align}
    \|e_i\t \mathbf{R} \| &\lesssim \theta_i\frac{K^{7/2} \sqrt{\log(n)}}{\|\theta\|^2 \lambda_{\min}^{1/2}} + \theta_i{\sqrt{K}} \frac{K^2 \theta_{\max} \|\theta\|_1}{\|\theta\|^4 \lambda_{\min}^2} + \theta_i \bigg( \frac{\theta_{\max}}{\theta_{\min}} \bigg)^{1/2} \bigg( \frac{K^{9/2} \theta_{\max} \|\theta\|_1 \log(n)}{\lambda_{\min}^{3/2} \|\theta\|^4} \bigg) \nonumber \\
    &\lesssim  \theta_i {\sqrt{K}}\frac{K^2 \theta_{\max} \|\theta\|_1}{\|\theta\|^4 \lambda_{\min}^2} + \theta_i \bigg( \frac{\theta_{\max}}{\theta_{\min}} \bigg)^{1/2} \bigg( \frac{K^{9/2} \theta_{\max} \|\theta\|_1 \log(n)}{\lambda_{\min}^{3/2} \|\theta\|^4} \bigg)\label{rbound} 
\end{align}
In what follows, denote
\begin{align*}
    \alpha_{\mathbf{R}} :&= {\sqrt{K}}\frac{K^2 \theta_{\max} \|\theta\|_1}{\|\theta\|^4 \lambda_{\min}^{2}} + \bigg( \frac{\theta_{\max}}{\theta_{\min}} \bigg)^{1/2} \bigg( \frac{K^{9/2} \theta_{\max} \|\theta\|_1 \log(n)}{\|\theta\|^4 \lambda_{\min}^{3/2}} \bigg),
     \numberthis\label{rbound1} 
\end{align*}
so that $\| e_i\t \mathbf{R} \| \lesssim \theta_i  \alpha_{\mathbf{R}}.$ 
\\ \ \\ \noindent
\textbf{Step 2: First Order Approximation of $\yhat$:}\\
Now, we note that
\begin{align*}
    e_i\t (\A - \bP) \U |\Lambda|^{-1/2} \ipq &= \sum_{j=1}^{n} (\A_{ij} - \bP_{ij}) (\U |\Lambda|^{-1/2} \ipq)_{j\cdot}
\end{align*}
is a sum of $n$ independent random matrices.  Bernstein's inequality shows that this is less than or equal to
\begin{align*}
    \frac{\sqrt{\theta_i \|\theta\|_1 \log(n)}}{\lambda^{1/2}} \| \U \|_{2,\infty} &\lesssim \frac{K\sqrt{\theta_i \|\theta\|_1 \log(n)}}{\lambda_{\min}^{1/2} \|\theta\|} \frac{\theta_{\max}}{\|\theta\|} \\
    &\asymp \frac{K \sqrt{\theta_i \|\theta\|_1 \log(n)} \theta_{\max}}{\|\theta\|^2 \lambda_{\min}^{1/2}} \\
    &\asymp \theta_i  \bigg( \frac{\theta_{\max}}{\theta_i} \bigg)^{1/2} \bigg[ \frac{K \sqrt{\theta_{\max} \|\theta\|_1 \log(n)}}{\|\theta\|^2 \lambda_{\min}^{1/2}} \bigg]. 
\end{align*}
Consequently, we obtain that
\begin{align*}
    \| e_i\t \xhat\wstar\t - e_i\t \xtilde \| &\lesssim \theta_i  \bigg( \frac{\theta_{\max}}{\theta_i} \bigg)^{1/2} \bigg[ \frac{K \sqrt{\theta_{\max} \|\theta\|_1 \log(n)}}{\|\theta\|^2 \lambda_{\min}^{1/2}} \bigg] + \theta_i\lambda_{\min}^{1/2} \alpha_{\mathbf{R}}  \\
    &= \theta_i \
    \Bigg\{  \bigg( \frac{\theta_{\max}}{\theta_i} \bigg)^{1/2} \bigg[ \frac{K \sqrt{\theta_{\max} \|\theta\|_1 \log(n)}}{\lambda_{\min}^{1/2}\|\theta\|^2} \bigg] +  \alpha_{\mathbf{R}} \Bigg\}\\
    &\leq \frac{1}{64}\| \xtilde_i \|,
\end{align*}
since $\|\xtilde_i\|\gtrsim \theta_i$, as long as $\alpha_{\mathbf{R}} \lesssim 1$ and that 
\begin{align*}
    \bigg( \frac{\theta_{\max}}{\theta_{\min}} \bigg)^{1/2} \frac{K \sqrt{\theta_{\max} \|\theta\|_1 \log(n)}}{\lambda_{\min}^{1/2} \|\theta\|^2} \lesssim 1, \numberthis \label{jjjj}
\end{align*}
both of which are guaranteed \cref{ass:networklevel}, which we will verify now. First, a direct comparison of $\alpha_{\mathbf{R}}$ with \cref{ass:networklevel} shows that $\alpha_{\mathbf{R}} \leq \frac{\bar \lambda}{C\sqrt{K}}$, which is strictly less than one.  In addition, by squaring \eqref{jjjj}, we see that we require that
\begin{align*}
    \frac{\theta_{\max}}{\theta_{\min}} \frac{K^2 \theta_{\max} \|\theta\|_1 \log(n)}{\lambda_{\min} \|\theta\|^4} \lesssim 1,
\end{align*}
but this is of smaller order than the first term in $\alpha_{\mathbf{R}}$.
Consequently, we are free to apply Taylor's Theorem to the function $x \mapsto x/\|x\|$ in a neighborhood of at most constant radius of $\xtilde_{i\cdot}$ not containing zero to obtain
\begin{align*}
     \big(\yhat \wstar\t\big)_{i\cdot} - \ytilde_{i\cdot} &= \frac{ \big(\xhat \wstar\t\big)_{i\cdot} }{\|\xhat_{i\cdot}\|} - \frac{ \xtilde_{i\cdot}}{\| \xtilde_{i\cdot}\|} \\
    &= \mathbf{J}( \xtilde_{i\cdot} ) \big(   (\xhat \wstar\t)_{i\cdot} - \xtilde_{i\cdot} \big) + \big( \mathbf{\tilde R}_Y\big)_{i\cdot},
\end{align*}
where
\begin{align*}
    \| e_i\t \mathbf{\tilde R}_Y \| &\lesssim r^2 \max_{|\alpha| = 2} \sup_{\|c - \xtilde_i\| \leq r} \| \mathbf{D}^{\alpha}(c) \|,
\end{align*}
where $\mathbf{D}^{\alpha}$ denotes the partial derivatives of the function $x \mapsto \frac{x}{\|x\|}$, and $r$ satisfies
\begin{align}
    r \leq C \theta_i 
    \Bigg\{  \bigg( \frac{\theta_{\max}}{\theta_i} \bigg)^{1/2} \bigg[ \frac{K \sqrt{\theta_{\max} \|\theta\|_1 \log(n)}}{\lambda_{\min}^{1/2}\|\theta\|^2} \bigg] +  \alpha_{\mathbf{R}} \Bigg\}, \label{eq:residual-taylor-bound}
\end{align}
for some constant $C > 0$.  We also have used the notation
\begin{align*}
    \mathbf{J}(  \xtilde_{i\cdot} ) &= \frac{1}{\|\xtilde_{i\cdot}\|}\bigg( \mathbf{I} - \frac{ \xtilde_{i\cdot} \xtilde\t_{i\cdot}}{\|\xtilde_{i\cdot}\|^2} \bigg),
\end{align*}
which is the Jacobian of the mapping $x\mapsto \frac{x}{\|x\|}$.  Expanding further, we have that
\begin{align*}
    \big( \yhat \wstar\t\big)_{i\cdot} - \ytilde_{i\cdot} &= \frac{1}{\|\xtilde_{i\cdot} \|} \bigg( \mathbf{I} - \frac{ \xtilde_{i\cdot} \xtilde\t_{i\cdot}}{\|\xtilde_{i\cdot}\|^2} \bigg) \big( (\xhat \wstar)_{i\cdot} - \xtilde_{i\cdot} \big) + \big(\mathbf{\tilde R}_Y\big)_{i\cdot} \\
    &= \frac{1}{\|\xtilde_{i\cdot} \|} \bigg( \mathbf{I} - \frac{ \xtilde_{i\cdot} \xtilde\t_{i\cdot}}{\|\xtilde_{i\cdot}\|^2} \bigg) \bigg( \big( \A - \bP \big) \U |\Lambda|^{-1/2} \ipq \bigg)_{i\cdot} + \frac{1}{\|\xtilde_{i\cdot} \|} \bigg( \mathbf{I} - \frac{ \xtilde_{i\cdot} \xtilde\t_{i\cdot}}{\|\xtilde_{i\cdot}\|^2} \bigg)\big( \mathbf{R} \big)_{i\cdot} + \big( \mathbf{\tilde R}_Y\big)_{i\cdot}.
\end{align*}
This justifies the linear part of the expansion, where we define
\begin{align*}
     \big( \mathcal{R}_{\mathrm{Stage \ I}} \big)_{i\cdot}:&= \frac{1}{\|\xtilde_{i\cdot} \|} \bigg( \mathbf{I} - \frac{ \xtilde_{i\cdot} \xtilde\t_{i\cdot}}{\|\xtilde_{i\cdot}\|^2} \bigg)\big( \mathbf{R} \big)_{i\cdot} + \big(\mathbf{\tilde R}_Y\big)_{i\cdot}.
\end{align*}
Therefore, it remains to bound this residual.  Recall that we already have the bound
\begin{align*}
    \| e_i\t \mathbf{R} \| &\lesssim\theta_i  \alpha_{\mathbf{R}} 
\end{align*}
with probability at least $1 - O(n^{-18})$ by \eqref{rbound1}.  Consequently, with this same probability, we note that $\| \xtilde_i \| \gtrsim  \theta_i$, so that
\begin{align*}
   \bigg\| \frac{1}{\|\xtilde_{i\cdot} \|} \bigg( \mathbf{I} - \frac{ \xtilde_{i\cdot} \xtilde\t_{i\cdot}}{\|\xtilde_{i\cdot}\|^2} \bigg)\big(  \mathbf{R} \big)_{i\cdot} \bigg\| &\lesssim \frac{1}{\theta_i} \bigg\| \bigg( \mathbf{I} - \frac{ \xtilde_i \xtilde\t_i}{\|\xtilde_i\|^2} \bigg)\big( \mathbf{R} \big)_{i\cdot} \bigg\| \\
   &\lesssim  \alpha_{\mathbf{R}},
\end{align*}
since the term $\mathbf{I} - \frac{ \xtilde_{i\cdot} \xtilde\t_{i\cdot}}{\|\xtilde_{i\cdot}\|^2}$ is a projection matrix. We therefore need only bound the term $e_i\t \mathbf{\tilde R}_Y$ which satisfies
\begin{align*}
    \| e_i\t \mathbf{\tilde R}_Y \| &\lesssim r^2 \max_{|\alpha|=2} \sup_{\|c - \xtilde_{i\cdot}\|} \| \mathbf{D}^{\alpha}(c) \|.
\end{align*}
We now note that the mixed partials of the mapping $x\mapsto \frac{x}{\|x\|}$ are given by
\begin{align*}
    \frac{\partial^2}{\partial x_i \partial x_j} \frac{x_k}{\|x\|} &=  \frac{3x_ix_jx_k}{\|x\|^5} - \frac{\delta_{ik} x_j + \delta_{ij} x_k + \delta_{jk} x_i}{\|x\|^3}.
\end{align*}
We evaluate this in a neighborhood of $\xtilde_{i\cdot}$ of radius at most $r$ where $r$ satisfies the inequality in \eqref{eq:residual-taylor-bound}.
It is straightforward to observe that since $r \lesssim \| \xtilde_{i\cdot} \|$, we have 
\begin{align*}
    \max_{|\alpha|=2} \sup_{\|c - \xtilde_i\|\leq r} \| \mathbf{D}^{\alpha}(c) \| &\lesssim \frac{1}{\|\xtilde_i\|^2}.
\end{align*}
Therefore,
\begin{align*}
    \| e_i\t \mathbf{\tilde R}_Y \| &\lesssim \frac{r^2}{\|\xtilde_i\|^2} \\
    &\lesssim \frac{ \theta_i^2 }{\|\xtilde_{i\cdot}\|^2}
    \Bigg\{  \bigg( \frac{\theta_{\max}}{\theta_i} \bigg)^{1/2} \bigg[ \frac{K \sqrt{\theta_{\max} \|\theta\|_1 \log(n)}}{\lambda_{\min}^{1/2}\|\theta\|^2} \bigg] +  \alpha_{\mathbf{R}} \Bigg\}^2\\
    &\lesssim  \Bigg\{  \bigg( \frac{\theta_{\max}}{\theta_i} \bigg)^{1/2} \bigg[ \frac{K \sqrt{\theta_{\max} \|\theta\|_1 \log(n)}}{\lambda_{\min}^{1/2}\|\theta\|^2} \bigg] +  \alpha_{\mathbf{R}} \Bigg\}^2 \\
    &\lesssim \bigg( \frac{\theta_{\max}}{\theta_i} \bigg)  \frac{K^2 \theta_{\max} \|\theta\|_1 \log(n)}{\lambda_{\min}\|\theta\|^4}  +  \alpha_{\mathbf{R}}, 
\end{align*}
which holds as long as $C$ in \cref{ass:networklevel} is larger than t he universal constants above, and hence both terms will be smaller than one. 
Therefore, we obtain that 
\begin{align*}
    \| e_i\t\mathcal{R}_{\mathrm{Stage \ I}} \| &\lesssim \bigg( \frac{\theta_{\max}}{\theta_{\min}} \bigg)  \frac{K^2 \theta_{\max} \|\theta\|_1 \log(n)}{\lambda_{\min}\|\theta\|^4}  +  \alpha_{\mathbf{R}} \\
    &\asymp \bigg( \frac{\theta_{\max}}{\theta_{\min}} \bigg)  \frac{K^2 \theta_{\max} \|\theta\|_1 \log(n)}{\lambda_{\min}\|\theta\|^4}  +  {\sqrt{K}}\frac{K^2 \theta_{\max} \|\theta\|_1}{\|\theta\|^4 \lambda_{\min}^{2}} + \bigg( \frac{\theta_{\max}}{\theta_{\min}} \bigg)^{1/2} \bigg( \frac{K^{9/2} \theta_{\max} \|\theta\|_1 \log(n)}{\|\theta\|^4 \lambda_{\min}^{3/2}}\bigg) \\
    &\lesssim  \frac{K^2 \theta_{\max} \|\theta\|_1}{\lambda_{\min} \|\theta\|^4} \bigg( \log(n) \frac{\theta_{\max}}{\theta_{\min}} + \frac{{\sqrt{K}}}{\lambda_{\min}} +\bigg(\frac{\theta_{\max}}{\theta_{\min}} \bigg)^{1/2}\frac{K^{5/2}\log(n)}{ \lambda_{\min}^{1/2}} \bigg) 
\end{align*}
which holds with probability at least $1- O(n^{-18})$.  This is the advertised bound, which completes the proof. \end{proof}

\subsection{Proofs of Lemmas~\ref{lem:rowwiseerror1} and \ref{lem:rowwiseerror2}} \label{sec:rowwiseproof}
To prove these lemmas we require leave-one-out sequences, similar to \citet{abbe_entrywise_2020}.  First we state the following lemma concerning the leave-one-out sequences.  The proof is deferred to \cref{sec:loolems}.

\begin{restatable}[Good properties of Leave-one-out sequences]{lemma}{loosequences} \label{lem:goodloo}
Let $\A\mi$ denote the matrix $\A\m$ with its $i$'th row and column replaced with $\bP \m$.  Let $\uhat^{(-i)}$ denote the leading $K$ eigenvectors of $\A\mi$.  Suppose that $\lambda \gtrsim \sqrt{\theta_{\max} \|\theta\|_1 \log(n)}$ and $\min_{i} \theta_i \|\theta\|_1 \gtrsim\log(n)$.  Then the following hold with probability at least $1 - O(n^{-20})$:
\begin{align*}
     |\lambda_{K}( \A\m ) - \lambda_{K+1}(\A\mi)| &\gtrsim \lambda\m; \\ 
     \| e_i\t ( \A\m - \bP \m) \uhat^{(-i)} \| &\lesssim \sqrt{\theta_i \|\theta\|_1 \log(n)} \| \uhat\|_{2,\infty}; \\
     \|  \uhat \uhat\t - \uhat^{(-i)} (\uhat^{(-i)})\t  \| &\lesssim \frac{\sqrt{\theta_i \|\theta\|_1 \log(n)}}{\lambda} \| \uhat\|_{2,\infty}.
\end{align*}
\end{restatable}
\noindent
We now prove \cref{lem:rowwiseerror1}.  The statement is repeated for convenience.

\rowwiseerror*

\begin{proof}[Proof of \cref{lem:rowwiseerror1}]
First, let $\uhat^{(-i)}$ denote the eigenvectors of $\A\m$ with the $i$'th row and column replaced with the corresponding row and column of $\bP\m$.  Observe that
\begin{align*}
    \| e_i\t (\mathbf{A} - \E \mathbf{A}) \uhat \| &=  \| e_i\t (\mathbf{A} - \E \mathbf{A}) \uhat \uhat\t \| \\
    &\leq \| e_i\t ( \A\m - \bP \m) \uhat^{(-i)} (\uhat^{(-i)})\t \| + \| e_i\t (\A \m - \bP \m ) \big( \uhat \uhat\t - \uhat^{(-i)} (\uhat^{(-i)})\t \big) \| \\
    &\leq  \| e_i\t ( \A\m - \bP \m) \uhat^{(-i)} \| + \| e_i\t (\A \m - \bP \m ) \| \|  \uhat \uhat\t - \uhat^{(-i)} (\uhat^{(-i)})\t  \| \\
    &\leq \| e_i\t ( \A\m - \bP \m) \uhat^{(-i)} \| + \| (\A \m - \bP \m ) \| \|  \uhat \uhat\t - \uhat^{(-i)} (\uhat^{(-i)})\t  \| \\
    &\lesssim \sqrt{\theta_i \| \theta\|_1 \log(n)} \| \uhat\|_{2,\infty} + \sqrt{\theta_{\max} \| \theta\|_1 } \frac{\sqrt{\theta_i \|\theta\|_1 \log(n)}}{\lambda} \| \uhat\|_{2,\infty},
\end{align*}
where the final inequality holds with probability at least $1 - O(n^{-20})$ by \cref{lem:goodloo} and \cref{lem:spectralnormconcentration}.  Consequently, since $\lambda \gtrsim \sqrt{\theta_{\max} \|\theta\|_1}$, we obtain that
\begin{align*}
     \| e_i\t (\mathbf{A} - \E \mathbf{A}) \uhat \| &\lesssim \sqrt{\theta_i \| \theta\|_1 \log(n)} \| \uhat\|_{2,\infty}
\end{align*}
with probability at least $1 - O(n^{-20})$ which completes the proof.
\end{proof}

We now restate \cref{lem:rowwiseerror2} for convenience.

\rowwiseerrortwo*

\begin{proof}
First we will argue that 
\begin{align}
     \| e_i\t (\mathbf{A} - \E \mathbf{A}) (\uhat \uhat \t \U - \U) \| &\lesssim \sqrt{\theta_i  \|\theta\|_1 \log(n)} \frac{\sqrt{\theta_i \|\theta\|_1 \log(n)}}{\lambda} \|\uhat\|_{2,\infty} \nonumber \\
     &\quad + \sqrt{\theta_i \|\theta\|_1 \log(n)}\| \uhat \uhat\t \U - \U \|_{2,\infty} \label{claim1}
\end{align}
with probability at least $1 - O(n^{-20})$.  Provided this is true, by \cref{lem:uhatucloseness}, we have that 
\begin{align*}
    \|\uhat\|_{2,\infty} &\lesssim \| \U \|_{2,\infty}; \\
    \| \uhat \uhat\t \U - \U \|_{2,\infty} &\lesssim \frac{\sqrt{\theta_{\max} \|\theta\|_1 \log(n)}}{\lambda} \| \U \|_{2,\infty},
\end{align*}
with probability at least $1 - O(n^{-19})$.  Plugging these in yields
\begin{align*}
     \| e_i\t (\mathbf{A} - \E \mathbf{A}) (\uhat \uhat \t \U - \U) \| &\lesssim \sqrt{\theta_i  \|\theta\|_1 \log(n)} \frac{\sqrt{\theta_i \|\theta\|_1 \log(n)}}{\lambda} \|\U\|_{2,\infty} \\
     &\quad + \sqrt{\theta_i \|\theta\|_1 \log(n)}\|\frac{\sqrt{\theta_{\max} \|\theta\|_1 \log(n)}}{\lambda} \| \U \|_{2,\infty} \\
     &\lesssim  \sqrt{\theta_i \|\theta\|_1 \log(n)}\|\frac{\sqrt{\theta_{\max} \|\theta\|_1 \log(n)}}{\lambda} \| \U \|_{2,\infty},
\end{align*}
which is the desired bound.  Therefore, it remains to prove the claim \eqref{claim1}.

Proceeding similarly  to the proof of \cref{lem:rowwiseerror1},
\begin{align*}
     \| e_i\t (\mathbf{A} -& \E \mathbf{A}) (\uhat\uhat\t \U - \U) \| \\
     &\leq \| e_i\t \big( \A - \bP \big) (\uhat \uhat\t \U - \uhat^{(-i)} (\uhat^{(-i)})\t \U) \| + \| e_i\t \big( \A - \bP \big) \big( \uhat^{(-i)} (\uhat^{(-i)})\t \U - \U \big) \| \\
     &\leq \| e_i\t \big(\A - \bP \big) \| \| \uhat \uhat\t  - \uhat^{(-i)} (\uhat^{(-i)})\t \| + \|  e_i\t \big( \A - \bP \big) \big( \uhat^{(-i)} (\uhat^{(-i)})\t \U - \U \big) \|.
\end{align*}
First, we note that the matrix $ \big( \uhat^{(-i)} (\uhat^{(-i)})\t \U - \U \big)$ is independent from the $i$'th row of $\A - \bP$.  The matrix Bernstein inequality (Corollary 3.3 of \cite{chen_spectral_2021}) shows that
\begin{align*}
     \|  e_i\t \big( \A - \bP \big) \big( \uhat^{(-i)} (\uhat^{(-i)})\t \U - \U \big) \| &\leq \sqrt{42 v \log(n)} + \frac{42}{3} w \log(n)
\end{align*}
with probability at least $1 - 2 n^{-20}$, where we have defined
\begin{align*}
    v :&= \max\bigg\{ \bigg\| \sum_{j=1}^{n} \E \big[ ( \A_{ij} - \bP_{ij} ) \big( \uhat^{(-i)} (\uhat^{(-i)})\t \U - \U \big)_{j\cdot} \big( \uhat^{(-i)} (\uhat^{(-i)})\t \U - \U \big)_{j\cdot}\t (\A_{ij} - \bP_{ij} ) \bigg\|, \\
    &\qquad \bigg\| \sum_{j=1}^{n} \E \bigg[ \big( \uhat^{(-i)} (\uhat^{(-i)})\t \U - \U \big)_{j\cdot}\t (\A_{ij} - \bP_{ij})^2 \big( \uhat^{(-i)} (\uhat^{(-i)})\t \U - \U \big)_{j\cdot} \bigg] \bigg\|\bigg\}; \\
    w :&= \max_{1\leq j \leq n} \| (\A_{ij} - \bP_{ij} ) \big( \uhat^{(-i)} (\uhat^{(-i)})\t \U - \U \big)_{j\cdot} \| \\
    &\leq \|\uhat^{(-i)} (\uhat^{(-i)})\t \U - \U \|_{2,\infty}.
\end{align*}
For the term $v$, we recognize that $\A_{ij} - \bP_{ij}$ is a scalar, yielding
\begin{align*}
    v &\leq \theta_i \| \theta\|_1 \|\uhat^{(-i)} (\uhat^{(-i)})\t \U - \U \|_{2,\infty}^2
\end{align*}
(for details on this calculation, see the proof of \cref{lem:goodloo}).  Consequently,
\begin{align*}
     \|  e_i\t \big( \A - \bP \big) \big( \uhat^{(-i)} (\uhat^{(-i)})\t \U - \U \big) \| &\leq \sqrt{42 v \log(n)} + \frac{42}{3} w \log(n) \\
     &\lesssim \sqrt{\theta_i \| \theta\|_1 \log(n)} \| \uhat^{(-i)} (\uhat^{(-i)})\t \U - \U \|_{2,\infty},
\end{align*}
as long as $\min_i \theta_i \|\theta\|_1 \gtrsim \log(n)$.  Moreover, a straightforward Bernstein inequality argument shows that $\| e_i\t \big(\A - \bP \big) \| \lesssim \sqrt{\theta_i \|\theta\|_1 \log(n)}$ with probability at least $1 - O(n^{-20})$.  Consequently, by \cref{lem:goodloo} and \cref{lem:spectralnormconcentration}, with probability at least $1 - O(n^{-20})$ it holds that
\begin{align*}
    \| e_i\t (\mathbf{A} -& \E \mathbf{A}) (\uhat\uhat\t \U - \U) \| \\
    &\leq \| e_i\t \big(\A - \bP \big) \| \| \uhat \uhat\t  - \uhat^{(-i)} (\uhat^{(-i)})\t \| \\
    &\quad + \|  e_i\t \big( \A - \bP \big) \big( \uhat^{(-i)} (\uhat^{(-i)})\t \U - \U \big) \| \\
    &\lesssim  \sqrt{\theta_i  \|\theta\|_1 \log(n)} \frac{\sqrt{\theta_i \|\theta\|_1 \log(n)}}{\lambda} \|\uhat\|_{2,\infty} \\
    &\quad + \sqrt{\theta_i \|\theta\|_1 \log(n)}\| \uhat^{(-i)} (\uhat^{(-i)})\t \U - \U \|_{2,\infty} \\
    &\lesssim  \sqrt{\theta_i  \|\theta\|_1 \log(n)} \frac{\sqrt{\theta_i \|\theta\|_1 \log(n)}}{\lambda} \|\uhat\|_{2,\infty} \\
    &\quad + \sqrt{\theta_i \|\theta\|_1 \log(n)} \bigg( \|\uhat^{(-i)} (\uhat^{(-i)})\t \U - \uhat \uhat\t \U \|_{2,\infty} + \| \uhat \uhat\t \U - \U \|_{2,\infty} \bigg) \\
    &\lesssim \sqrt{\theta_i  \|\theta\|_1 \log(n)} \frac{\sqrt{\theta_i \|\theta\|_1 \log(n)}}{\lambda} \|\uhat\|_{2,\infty} + \sqrt{\theta_i \|\theta\|_1 \log(n)}\| \uhat \uhat\t \U - \U \|_{2,\infty}.
\end{align*} 
\end{proof}

\subsubsection{Proof of Lemma~\ref{lem:goodloo}} \label{sec:loolems}

We restate \cref{lem:goodloo} for convenience.

\loosequences*

\begin{proof}[Proof of \cref{lem:goodloo}]
First, by \cref{lem:spectralnormconcentration}, it holds that
\begin{align*}
    \| \A \m - \bP \m \| &\lesssim \sqrt{\theta_{\max} \| \theta\|_1} \\
    &\leq \lambda/\sqrt{\log(n)}.
\end{align*}
Therefore, Weyl's inequality shows that
\begin{align*}
    |\lambda_K (\A\m)| &\geq |\lambda_K| -  \| \A\m - \bP \m  \| \\
    &\geq \lambda - \lambda/\sqrt{\log(n)} \\
    &\geq \lambda/2 \gtrsim \lambda,
\end{align*}
and that $|\lambda_{K+1}(\A\m) | \leq  \| \A \m - \bP \m \| \leq |\lambda_K|/\sqrt{\log(n)}$.  Therefore, $ |\lambda_K (\A\m)| -  |\lambda_{K+1}(\A\m)| \gtrsim \lambda$.  Furthermore,
\begin{align*}
    \| e_i\t \big( \A\m - \bP \m \big) \| &\leq \| \A\m - \bP \m  \|.
\end{align*}
Observe that $\A\m = \A\mi + e_i e_i\t\big( \A\m - \bP \m \big) + \big( \A\m - \bP \m \big) e_i e_i\t - e_i e_i\t \big( \A\m - \bP \m \big) e_i e_i\t$. Consequently, by Weyl's inequality,
\begin{align*}
    |\lambda_K( \A\mi) - \lambda_{K+1}(\A\m) | &\geq |\lambda_K(\A \m) |- |\lambda_{K+1}(\A\m)| \\
    &-  \bigg\| e_i e_i\t\big( \A\m - \bP \m \big) + \big( \A\m - \bP \m \big) e_i e_i\t - e_i e_i\t \big( \A\m - \bP \m \big) e_i e_i\t \bigg\| \\
    &\gtrsim |\lambda_K|.
\end{align*}
This proves the first assertion.  As a byproduct, we are free to apply the Davis-Kahan Theorem to $\uhat$ and $\uhat^{(-i)}$ to observe that
\begin{align*}
    \| \uhat \uhat\t - \uhat^{(-i)} (\uhat^{(-i)})\t \| &\lesssim \frac{ \| e_i\t \big( \A - \bP\m \big) \uhat^{(-i)} \| + \| \big( \A - \bP \big) e_i e_i\t \uhat^{(-i)} \|}{\lambda} \\
    &\lesssim \frac{ \| e_i\t \big( \A - \bP\m \big) \uhat^{(-i)} \|}{\lambda} +\frac{ \| e_i\t \big( \A - \bP \big) \| \| e_i\t \uhat^{(-i)} \|}{\lambda}.
\end{align*}
Consequently, we need only bound the numerators above; however, a bound on the first term will also prove the second assertion of this lemma. Note that
\begin{align*}
    e_i\t \big( \A - \bP \big) \uhat^{(-i)} &= \sum_{j=1}^{n} \big( \A_{ij} - \bP_{ij} \big) \uhat^{(-i)}_{j\cdot}.
\end{align*}
Since $\uhat^{(-i)}$ is independent from the $i$'th row of $\A_{ij}$, this is a sum of $n$ independent random matrices condition on $\uhat^{(-i)}$.  Therefore, the matrix Bernstein inequality (Corollary 3.3 of \cite{chen_spectral_2021}) reveals that
\begin{align*}
    \| e_i\t \big( \A - \bP \big) \uhat^{(-i)} \| &\leq \sqrt{42  v \log(n)} + \frac{42}{3} w \log(n)
\end{align*}
with probability at least $1 - 2 n^{-20}$.  Here we note that
\begin{align*}
    v :&= \max\bigg\{ \bigg\| \sum_{j=1}^n \E \big[ ( \A_{ij} - \bP_{ij}) \uhat^{(-i)}_{j\cdot} (\uhat^{(-i)}_{j\cdot})\t ( \A_{ij} - \bP_{ij}) \big] \bigg\|, \bigg\| \sum_{j=1}^n \E \big[(\uhat^{(-i)}_{j\cdot})\t  ( \A_{ij} - \bP_{ij})^2 \uhat^{(-i)}_{j\cdot}   \big] \bigg\|\bigg\}; \\
    w :&= \max_{1\leq j \leq n} \| (\A_{ij} - \bP_{ij}) \uhat^{(-i)}_{j\cdot} \| \\
&\leq \| \uhat^{(-i)} \|_{2,\infty},
\end{align*}
where the expectation in the first term is conditional on $\uhat^{(-i)}$.
Observing that $\A_{ij} - \bP_{ij}$ is a scalar reveals that
\begin{align*}
    v &\leq \max\bigg\{ \bigg\| \sum_{j=1}^n  \uhat^{(-i)}_{j\cdot} (\uhat^{(-i)}_{j\cdot})\t \E  ( \A_{ij} - \bP_{ij})^2  \bigg\|, \bigg\| \sum_{j=1}^n (\uhat^{(-i)}_{j\cdot})\t   \uhat^{(-i)}_{j\cdot}  \E ( \A_{ij} - \bP_{ij})^2  \bigg\|\bigg\} \\
    &\leq \sum_{j=1}^{n} \| \uhat^{(-i)}_{j\cdot} \|^2 \theta_i \theta_j \\
    &\leq \| \uhat^{(-i)} \|_{2,\infty}^2 \theta_i \| \theta \|_1.
\end{align*}
Therefore, it holds that
\begin{align*}
     \| e_i\t \big( \A - \bP \big) \uhat^{(-i)} \| &\leq \sqrt{42  v \log(n)} + \frac{42}{3} L \log(n) \\
     &\leq \sqrt{42 \theta_i \|\theta\|_1 \log(n)} \| \uhat^{(-i)} \|_{2,\infty} + \frac{42}{3} \log(n) \| \uhat^{(-i)} \|_{2,\infty} \\
     &\lesssim \sqrt{\theta_i \| \theta\|_1 \log(n)} \| \uhat^{(-i)} \|_{2,\infty},
\end{align*}
which holds as long as $\min_i \theta_i \| \theta \|_1 \gtrsim \log(n)$.  Moreover, we have that $\| e_i\t (\A - \bP ) \| \lesssim \sqrt{\theta_i \| \theta \|_1 \log(n)}$ by a direct application of matrix Bernstein again.  Consequently, applying these bounds yields that
\begin{align*}
     \| \uhat \uhat\t - \uhat^{(-i)} (\uhat^{(-i)})\t \| &\lesssim \frac{ \| e_i\t \big( \A - \bP\m \big) \uhat^{(-i)} \|}{\lambda} +\frac{ \| e_i\t \big( \A - \bP \big) \| \| e_i\t \uhat^{(-i)} \|}{\lambda} \\
     &\lesssim \frac{ \sqrt{\theta_i \| \theta\|_1 \log(n)} }{\lambda}  \| \uhat^{(-i)} \|_{2,\infty}.
\end{align*}
As a byproduct, we also have that
\begin{align*}
    \| \uhat^{(-i)} \|_{2,\infty} &= \| \uhat^{(-i)} (\uhat^{(-i)})\t \|_{2,\infty} \\
    &\leq \| \uhat^{(-i)} (\uhat^{(-i)})\t - \uhat \uhat\t \|_{2,\infty}  + \| \uhat \uhat\t \|_{2,\infty} \\
    &\leq \frac{1}{2} \| \uhat^{(-i)} \|_{2,\infty} + \| \uhat \|_{2,\infty},
\end{align*}
which holds as long as $\lambda \gtrsim \sqrt{\theta_{\max} \| \theta\|_1 \log(n)}$.  Consequently, by rearranging, we have that $\|\uhat^{(-i)} \|_{2,\infty} \lesssim \| \uhat\|_{2,\infty}$ which yields the inequality
\begin{align*}
    \| \uhat \uhat\t - \uhat^{(-i)} (\uhat^{(-i)})\t \| &\lesssim\frac{ \sqrt{\theta_i \| \theta\|_1 \log(n)} }{\lambda}  \| \uhat \|_{2,\infty},
\end{align*}
which holds with probability at least $1 - O(n^{-20})$.  Moreover, with this same probability, we have that
\begin{align*}
    \|e_i\t (\A - \bP) \uhat^{(-i)} \| &\lesssim \sqrt{\theta_i \|\theta\|_1 \log(n)} \| \uhat\|_{2,\infty}.
\end{align*}
This completes the proof.
\end{proof}

\section{Proof of Second Stage \texorpdfstring{$\sin\bTheta$}{sin(Theta)} Bound (Theorem~\ref{thm:step2sintheta})}
First we will restate \cref{thm:step2sintheta}.

\sintheta*

In what follows we give a high-level overview of the proof. Define the matrix $\mathcal{ Y} := [\ytilde\one, \cdots, \ytilde\M] \in \mathbb{R}^{n \times LK}$, and let $\mathcal{\hat Y}$ be defined similarly.  Since we consider the singular vectors of $\mathcal{ Y}$ and $\mathcal{\hat Y}$, we will examine the \emph{eigenvectors} of their associated $n \times n$ Gram matrices, or the matrices $\ycal \ycal\t$ and $\yhatcal \yhatcal\t$ respectively.  Therefore, we will view $\yhatcal \yhatcal\t$ as a perturbation of matrix $\ycal \ycal\t$.  We expand via
\begin{align*}
    \ycal \ycal\t -\yhatcal \yhatcal\t &=  \mathcal{L}(\mathcal{E})\mathcal{ Y}\t + \mathcal{ Y} \mathcal{L}(\mathcal{E})\t + \mathcal{R}_{\mathrm{all}},
\end{align*}
where we define
\begin{align*}
    \mathcal{R}_{\mathrm{all}} :&= \sum_{l} \mathcal{L}(\mathbf{E}\m)\mathcal{L}(\mathbf{E}\m)\t + \mathcal{L}(\mathbf{E}\m) (\mathcal{R}_{\mathrm{Stage \ I}}\m)\t+ \mathcal{R}_{\mathrm{Stage \ I}}\m\mathcal{L}(\mathbf{E}\m)\t + \mathcal{R}_{\mathrm{Stage \ I}}\m (\mathcal{R}_{\mathrm{Stage \ I}}\m)\t \\
        &\quad + \sum_l \mathcal{R}_{\mathrm{Stage \ I}}\m (\ytilde\m)\t + \ytilde\m (\mathcal{R}_{\mathrm{Stage \ I}}\m)\t,
\end{align*}
and
\begin{align*}
\mathcal{L}(\mathcal{E}) :&= \big[ \mathcal{L}(\mathbf{E}\one), \cdots , \mathcal{L}(\mathbf{E}\M) \big],
\end{align*}
where we have defined $\mathbf{E}\l$ as the mean-zero random matrix $\mathbf{E}\m := \A\m - \bP\m$.  Hence,
\begin{align*}
    \mathcal{L}(\mathcal{E})\mathcal{ Y}\t + \mathcal{ Y} \mathcal{L}(\mathcal{E})\t &= \sum_l \mathcal{L}(\mathbf{E}\m) (\ytilde\m)\t + \ytilde\m \mathcal{L}(\mathbf{E}\m)\t .
\end{align*}
By virtue of the tight characterization for each $\mathcal{R}\m_Y$ in \cref{thm:firststep}, we can see that $\ycal\ycal\t$ is \emph{nearly} a linear perturbation of $\yhatcal\yhatcal\t$.  The proof of \cref{thm:step2sintheta} makes this rigorous.

\subsection{Preliminary Lemmas: Spectral Norm Concentration Bounds} \label{sec:step2prelim}
Throughout this section we use the notation $\mathbf{E}\l := \A\m - \bP\m$.  The following lemma bounds several terms involving $\mathcal{L}(\mathbf{E}\m)$ in spectral norm.
\begin{lemma}[Linear Term Spectral Norm Concentration] \label{lem:LEmspectral}
It holds that
\begin{align*}
    \| \mathcal{L}(\mathbf{E}\m) \| &\lesssim \frac{K \sqrt{n\log(n)}}{(\lambda_{\min}\m)^{1/2} \|\theta\m\|} \bigg( \frac{\theta_{\max}\l}{\theta_{\min}\l} \bigg)^{1/2};\\
    \| \sum_{l} \mathcal{L}(\mathbf{E}\m) (\ytilde\m)\t \| &\lesssim  K n \sqrt{L\log(n)} \Bigg[ \frac{1}{L} \sum_{l} \bigg( \frac{\theta_{\max}\m}{\theta_{\min}\m}\bigg) \frac{1}{\lambda_{\min}\m \|\theta\m\|^2} \Bigg]^{1/2}.
\end{align*}
with probability at least $1 - O(n^{-15}).$
\end{lemma}

\begin{proof}[Proof of \cref{lem:LEmspectral}]
We recall that
\begin{align*}
    \mathcal{L}(\mathbf{E}\m)_{i\cdot} &= \mathbf{J}(\xtilde_{i\cdot}) \bigg( (\A\m - \bP \m) \U\m |\Lambda\m|^{-1/2} \ipq\m \bigg)_{i\cdot}.
\end{align*}
Therefore, we can write this matrix via
\begin{align*}
    \mathcal{L}(\mathbf{E}\m) &= \bigg( \sum_{i,j} \mathbf{E}\m_{ij} e_i e_j\t \bigg) \bigg( \U\m |\Lambda\m|^{-1/2} \ipq \m \mathbf{J}(\xtilde_{i\cdot} ) \bigg) \\
    &= \sum_{i\leq j} \mathbf{E}\m_{ij} e_i e_j\t \bigg( \U\m |\Lambda\m|^{-1/2} \ipq \m \mathbf{J}(\xtilde_{i\cdot} ) \bigg) + \sum_{j < i} \mathbf{E}\m_{ij} e_j e_i\t \bigg( \U\m |\Lambda\m|^{-1/2} \ipq \m \mathbf{J}(\xtilde_{j\cdot} )\bigg),
\end{align*}
both of which are a sum of independent random matrices.  Without loss of generality we bound the first term; the second is similar.  We will apply the matrix Bernstein inequality (\citet{chen_spectral_2021}, Corollary 3.3).  We need to bound:
\begin{align*}
    v :&= \max \bigg\{ \bigg\| \sum_{i\leq j} \E (\mathbf{E}\m_{ij})^2  \bigg( \U\m |\Lambda\m|^{-1/2} \ipq \m \mathbf{J}(\xtilde_{i\cdot} ) \bigg)\t e_i e_j\t e_j e_i\t \bigg( \U\m |\Lambda\m|^{-1/2} \ipq \m \mathbf{J}(\xtilde_{i\cdot} ) \bigg)\t \bigg\|,\\
    &\qquad \bigg\| \sum_{i\leq j} \E (\mathbf{E}\m_{ij})^2 e_i e_j\t  \bigg( \U\m |\Lambda\m|^{-1/2} \ipq \m \mathbf{J}(\xtilde_{i\cdot} ) \bigg) \bigg( \U\m |\Lambda\m|^{-1/2} \ipq \m \mathbf{J}(\xtilde_{i\cdot} ) \bigg)\t e_j e_i\t \bigg\|\bigg\}; \\
    w:&= \max_{i,j} \| \mathbf{E}\m_{ij} e_i e_j\t \bigg( \U\m |\Lambda\m|^{-1/2} \ipq \m \mathbf{J}(\xtilde_{i\cdot} ) \bigg) \|.
\end{align*}
For $v$, we note that 
\begin{align*}
    \bigg\| \sum_{i\leq j} \E & (\mathbf{E}\m_{ij})^2  \bigg( \U\m |\Lambda\m|^{-1/2} \ipq \m \mathbf{J}(\xtilde_{i\cdot} ) \bigg)\t e_i e_j\t e_j e_i\t \bigg( \U\m |\Lambda\m|^{-1/2} \ipq \m \mathbf{J}(\xtilde_{i\cdot} ) \bigg)\t \bigg\| \\
    &\leq \sum_{i\leq j} \E (\mathbf{E}\m_{ij})^2 \bigg\|  \bigg( \U\m |\Lambda\m|^{-1/2} \ipq \m \mathbf{J}(\xtilde_{i\cdot} ) \bigg)\t e_i e_j\t e_j e_i\t \bigg( \U\m |\Lambda\m|^{-1/2} \ipq \m \mathbf{J}(\xtilde_{i\cdot} ) \bigg)\t \bigg\| \\
    &\leq \sum_{i\leq j} \theta_i\l \theta_j\l \| e_j\t \U\m \|^2 \| |\Lambda\m|^{-1/2} \|^2 \| \mathbf{J}(\mathbf{X}_{i\cdot}) \|^2 \\
    &\leq \frac{K }{\lambda_{\min}\m\|\theta\m\|^2}\sum_{i\leq j} \theta_i\l \theta_j\l \| e_j\t \U\m \|^2\| \mathbf{J}(\mathbf{X}_{i\cdot}) \|^2 \\
    &\leq  \frac{K }{\lambda_{\min}\m\|\theta\m\|^2}\sum_{i\leq j} \theta_i\l \theta_j\l \frac{(\theta_j\l)^2 K}{\|\theta\m\|^2} \frac{1}{(\theta_i\l)^2} \\
    &\leq \frac{K^2}{\lambda_{\min}\m \|\theta\m\|^4} \sum_{i\leq j} \frac{\theta_j\l}{\theta_i\l} (\theta_j\l)^2 \\
    &\leq \frac{K^2n }{\lambda_{\min}\m \|\theta\m\|^4} \bigg( \frac{\theta_{\max}\l}{\theta_{\min}\l}\bigg) \sum_{j}  \theta_j^2 \\ 
    &\leq \frac{K^2n }{\lambda_{\min}\m \|\theta\m\|^4} \bigg( \frac{\theta_{\max}\l}{\theta_{\min}\l}\bigg) \|\theta\m\|^2\\ 
    &\leq \frac{K^2n }{\lambda_{\min}\m\|\theta\m\|^2} \bigg( \frac{\theta_{\max}\l}{\theta_{\min}\l}\bigg).
\end{align*}
The other term satisfies the same upper bound.  In addition,
\begin{align*}
    w &= \max_{i,j} \| \mathbf{E}\m_{ij} e_i e_j\t \bigg( \U\m |\Lambda\m|^{-1/2} \ipq \m \mathbf{J}(\xtilde_{i\cdot} ) \bigg) \| \\
    &\leq  \| \U\m \|_{2,\infty} \| |\Lambda\m|^{-1/2}\| \max_{i} \|  \mathbf{J}(\xtilde_{i\cdot} ) \| \\
    &\leq \frac{K }{\|\theta\|^2 \lambda_{\min}\m}  \bigg( \frac{\theta_{\max}\l}{(\theta_{\min}\l)^{1/2}} \bigg).
\end{align*}
Therefore, by the Matrix Bernstein inequality, with probability at least $1 - O(n^{-20})$ it holds that
\begin{align*}
    \| \mathcal{L}(\mathbf{E}\m) \| &\lesssim \sqrt{v \log(n)} + w \log(n) \\
    &\lesssim \frac{K \sqrt{n\log(n)} }{(\lambda_{\min}\m)^{1/2} \|\theta\m\|} \bigg( \frac{\theta_{\max}\l}{\theta_{\min}\l}\bigg)^{1/2} + \frac{K \log(n) }{\|\theta\l\|^2 (\lambda_{\min}\m)^{1/2}}  \bigg( \frac{\theta_{\max}\l}{\theta_{\min}\l} \bigg) \\
    &\leq \frac{K \sqrt{\log(n)}}{(\lambda_{\min}\m)^{1/2} \|\theta\m\|} \bigg( \frac{\theta_{\max}\l}{\theta_{\min}\l} \bigg)^{1/2} \max\bigg\{ \sqrt{n}, \bigg(\frac{\theta_{\max}\l}{\theta_{\min}\l} \bigg)^{1/2} \frac{\sqrt{\log(n)}}{\|\theta\m\|} \bigg\}.
\end{align*}
Finally, we note that by \cref{ass:networklevel}, it holds that  $\frac{\theta_{\max}}{\theta_{\min}} \lesssim \sqrt{n}$, which implies that $\sqrt{n}$ is the maximum of the term above.  Therefore,
\begin{align*}
     \| \mathcal{L}(\mathbf{E}\m) \| &\lesssim \frac{K \sqrt{n\log(n)}}{(\lambda_{\min}\m)^{1/2} \|\theta\m\|} \bigg( \frac{\theta_{\max}\l}{\theta_{\min}\l} \bigg)^{1/2},
\end{align*}
which completes the proof of the first statement.  

For the next statement, we proceed similarly, only now streamlining the analysis.  Representing the sum similarly, we have that
\begin{align*}
    \sum_{l} \mathcal{L}(\mathbf{E}\m) (\ytilde\m)\t &= \sum_{l} \bigg( \sum_{i\leq j} \mathbf{E}\m_{ij} e_i e_j\t \bigg( \U\m |\Lambda\m|^{-1/2} \ipq \m \mathbf{J}(\xtilde_{i\cdot} )  \bigg)(\ytilde\m)\t \\
    &\qquad +  \sum_l \bigg( \sum_{j < i} \mathbf{E}\m_{ij} e_j e_i\t \bigg( \U\m |\Lambda\m|^{-1/2} \ipq \m \mathbf{J}(\xtilde_{j\cdot} ) \bigg) (\ytilde\m)\t.
\end{align*}
We focus again on the first term.  Since it holds that $\|\ytilde\m\| \leq \|\ytilde\m \|_F = \sqrt{n}$, we have that
\begin{align*}
    v &\leq \sum_{l} \frac{K^2n^2 }{\lambda_{\min}\m \|\theta\m\|^2} \bigg( \frac{\theta_{\max}\m}{\theta_{\min}\m}\bigg) \\
    &= K^2 n^2 \sum_{l} \bigg( \frac{\theta_{\max}\m}{\theta_{\min}\m}\bigg) \frac{1}{\lambda_{\min}\m \|\theta\m\|^2}
\end{align*}
and
\begin{align*}
     w &= \max_{i,j,m} \| \mathbf{E}\m_{ij} e_i e_j\t \bigg( \U\m |\Lambda\m|^{-1/2} \ipq \m \mathbf{J}(\xtilde_{i\cdot} ) \bigg) (\ytilde\m)\t \| \\
    &\leq  \max_l \frac{K \sqrt{n} }{\|\theta\m\|^2 (\lambda_{\min}\m)^{1/2}}  \bigg( \frac{\theta_{\max}\m}{\theta_{\min}\m} \bigg).
\end{align*}
Therefore, with probability at least $1 - O(n^{-15})$,
\begin{align*}
   \| \sum_{l} \mathcal{L}(\mathbf{E}\m) (\ytilde\m)\t \| &\lesssim  K n \sqrt{\log(n)} \Bigg[\sum_{l} \bigg( \frac{\theta_{\max}\m}{\theta_{\min}\m}\bigg) \frac{1}{\lambda_{\min}\m \|\theta\m\|^2} \Bigg]^{1/2} \\
   &\qquad + K \sqrt{n} \log(n) \max_l \frac{1 }{\|\theta\m\|^2 (\lambda_{\min}\m)^{1/2}}  \bigg( \frac{\theta_{\max}\m}{\theta_{\min}\m} \bigg).
\end{align*}
Finally, we note that as long as $\frac{\theta_{\max}}{\theta_{\min}} \lesssim \sqrt{n}$, the first term dominates.  Therefore,
\begin{align*}
     \| \sum_{l} \mathcal{L}(\mathbf{E}\m) (\ytilde\m)\t \| &\lesssim  K n \sqrt{L\log(n)} \Bigg[ \frac{1}{L} \sum_{l} \bigg( \frac{\theta_{\max}\m}{\theta_{\min}\m}\bigg) \frac{1}{\lambda_{\min}\m \|\theta\m\|^2} \Bigg]^{1/2}.
\end{align*}
\end{proof}

Next, we bound residual term $\mathcal{R}_{\mathrm{all}}$ in spectral norm.
\begin{lemma}[Residual Term Spectral Concentration]\label{lem:step2rallspectralbound}
The residual term $\mathcal{R}_{\mathrm{all}}$ satisfies
\begin{align*}
    \| \mathcal{R}_{\mathrm{all}} \| &\lesssim L K^2 n \log(n) \| \snr\inv \|_{\infty}^2 + K L n \sqrt{\log(n)} \alpha_{\max} \| \snr\inv \|_{\infty} + n L \alpha_{\max}
\end{align*}
with probability at least $1 -O(n^{-15})$.  
\end{lemma}

\begin{proof}[Proof of \cref{lem:step2rallspectralbound}]
Recall that
\begin{align*}
       \mathcal{R}_{\mathrm{all}} :&= \sum_{l} \mathcal{L}(\mathbf{E}\m)\mathcal{L}(\mathbf{E}\m)\t + \mathcal{L}(\mathbf{E}\m) (\mathcal{R}_{\mathrm{Stage \ I}}\m)\t+ \mathcal{R}_{\mathrm{Stage \ I}}\m\mathcal{L}(\mathbf{E}\m)\t + \mathcal{R}_{\mathrm{Stage \ I}}\m (\mathcal{R}_{\mathrm{Stage \ I}}\m)\t \\
        &\quad + \sum_l \mathcal{R}_{\mathrm{Stage \ I}}\m (\ytilde\m)\t + \ytilde\m (\mathcal{R}_{\mathrm{Stage \ I}}\m)\t \\
        :&= (I) + (II) + (III) + (IV),
\end{align*}
where
\begin{align*}
    (I) :&= \sum_{l} \mathcal{L}(\mathbf{E}\m)\mathcal{L}(\mathbf{E}\m)\t; \\
    (II) :&=\sum_{l} \mathcal{L}(\mathbf{E}\m) (\mathcal{R}_{\mathrm{Stage \ I}}\m)\t+ (\mathcal{R}_{\mathrm{Stage \ I}}\m)\mathcal{L}(\mathbf{E}\m)\t ; \\
    (III) :&= \sum_l \mathcal{R}_{\mathrm{Stage \ I}}\m (\mathcal{R}_{\mathrm{Stage \ I}}\m)\t; \\
    (IV) :&= \sum_l \mathcal{R}_{\mathrm{Stage \ I}}\m (\ytilde\m)\t + \ytilde\m (\mathcal{R}_{\mathrm{Stage \ I}}\m)\t.
\end{align*}
We bound each term separately.
\\ \ \\ \noindent
\textbf{The Term $(I)$:} We note that by \cref{lem:LEmspectral} we have the bound
\begin{align*}
    \| \mathcal{L}(\mathbf{E}\m) \| &\lesssim \frac{K \sqrt{n\log (n)}}{(\lambda_{\min }\l)^{1/2}\left\|\theta\l\right\|}\left(\frac{\theta_{\max }\m}{\theta_{\min }\m}\right)^{1 / 2}.
\end{align*}
Therefore, 
\begin{align*}
    \bigg\| \sum_l \mathcal{L}(\mathbf{E}\m) \mathcal{L}(\mathbf{E}\m) \bigg\| &\lesssim L \max_l \bigg(\frac{K \sqrt{n\log (n)}}{(\lambda_{\min}\l)^{1/2}\left\|\theta^{(l)}\right\|}\left(\frac{\theta_{\max }\m}{\theta_{\min }\m}\right)^{1 / 2} \bigg)^2 \\
    &= L K^2 n\log(n) \max_l \bigg( \frac{\theta_{\max}\m}{\theta_{\min}\m} \bigg) \frac{1}{\lambda_{\min}\m \|\theta\m\|^2} \\
    &\asymp L K^2 n\log(n) \| \snr\inv \|_{\infty}^2.
    \numberthis \label{lelebound}
\end{align*}
\textbf{The term $(II):$} without loss of generality we consider the first term. By \cref{lem:step2resspectral}, it holds that
\begin{align*}
    \| \mathcal{R}_{\mathrm{Stage \ I}}\m \| &\lesssim \sqrt{n} \alpha\m,
\end{align*}
where $\alpha\m$ is the residual bound from \cref{thm:firststep}. 
Therefore, 
\begin{align*}
    \sum_l \| \mathcal{L}(\mathbf{E}\m) \mathcal{R}_{\mathrm{Stage \ I}}\m \| &\lesssim L\sqrt{n} \max_l \alpha\l \max_l \| \mathcal{L}(\mathbf{E}\m) \| \\
    &\asymp KL n \sqrt{\log(n)} \alpha_{\max} \| \snr\inv \|_{\infty} , \numberthis \label{lerybound}
\end{align*}
where we set $\alpha_{\max} := \max_l \alpha\m$.   \\ \ \\ \noindent \textbf{The Term $(III)$:} By a similar argument,
\begin{align*}
    (III) &\lesssim  n L \max_l (\alpha\m)^2 \\
    &\lesssim nL \alpha_{\max}^2. \numberthis \label{RyRybound}
\end{align*}
\textbf{The term $(IV)$:}
Finally, it holds that
\begin{align*}
    \sum_l \| \mathcal{R}_{\mathrm{Stage \ I}}\m \| \| \mathbf{Y}\m \| &\lesssim L \sqrt{n} \alpha_{\max} \max_l \| \mathbf{Y}\m \| \\
    &\lesssim L n \alpha_{\max}. \numberthis \label{RyYbound}
\end{align*}
\textbf{Putting it all together:}
Combining \eqref{lelebound}, \eqref{lerybound}, \eqref{RyRybound}, and \eqref{RyYbound}, we have that
\begin{align*}
    \| \mathcal{R}_{\mathrm{all}} \| &\lesssim L K^2 n \log(n) \| \snr\inv \|_{\infty}^2 + K L n \sqrt{\log(n)} \alpha_{\max} \| \snr\inv \|_{\infty} + nL \alpha_{\max}^2 + n L \alpha_{\max} \\
    &\asymp L K^2 n \log(n) \| \snr\inv \|_{\infty}^2 + K L n \sqrt{\log(n)} \alpha_{\max} \| \snr\inv \|_{\infty} + n L \alpha_{\max},
    \end{align*}
    since $\alpha_{\max} < 1$ by \cref{ass:networklevel} (as shown in the proof of \cref{thm:firststep}). 
\end{proof}

\subsection{Proof of Theorem~\ref{thm:step2sintheta}}

\begin{proof}[Proof of \cref{thm:step2sintheta}]
First, by \cref{lem:LEmspectral}, we have the bound
\begin{align*}
    \| \sum_l \mathcal{L}(\mathbf{E}\m) (\ytilde\m)\t \| &\lesssim K n \sqrt{L\log(n)} \Bigg[ \frac{1}{L} \sum_{l} \bigg( \frac{\theta_{\max}\m}{\theta_{\min}\m}\bigg) \frac{1}{\lambda_{\min}\m \|\theta\m\|^2} \Bigg]^{1/2}.
\end{align*}
Recall we define
\begin{align*}
  \big( \frac{1}{L} \| \snr\inv \|_2^2\big) :&= \frac{1}{L} \sum_l \bigg( \frac{\theta_{\max}\m}{\theta_{\min}\m} \bigg) \frac{1}{\lambda_{\min}\m \|\theta\m\|^2} .
\end{align*}
Then the bound can be concisely written as
\begin{align*}
    \| \sum_l \mathcal{L}(\mathbf{E}\m) (\ytilde\m)\t \| &\lesssim K n \sqrt{L \log(n) } \big( \frac{1}{L} \| \snr\inv \|_2^2\big)^{1/2}.
\end{align*}
In addition, by \cref{lem:step2rallspectralbound}, we have that
\begin{align*}
    \mathcal{R}_{\mathrm{all}} &\lesssim L K^2 n \log(n) \| \snr\inv \|_{\infty}^2 + K L n \sqrt{\log(n)} \alpha_{\max} \| \snr\inv \|_{\infty} + n L \alpha_{\max}.
\end{align*}
Therefore, it holds that
\begin{align*}
    \| \ycal\ycal\t - \yhatcal\yhatcal\t \| &\lesssim K n \sqrt{L \log(n) } \big( \frac{1}{L} \| \snr\inv \|_2^2\big)^{1/2} + L K^2 n \log(n) \| \snr\inv \|_{\infty}^2 \\
    &\quad + KL n \sqrt{\log(n)} \alpha_{\max} \| \snr\inv \|_{\infty} + n L \alpha_{\max}.
\end{align*}
Recall that $\lambda^2_{Y} \gtrsim \frac{n}{K} L \bar \lambda$ by \cref{lem:step2popprop}. Therefore, as long as
\begin{align*}
    n L \bar \lambda &\gtrsim K^2 n \sqrt{L \log(n) } \big( \frac{1}{L} \| \snr\inv \|_2^2\big)^{1/2} + LK^3 n \log(n) \| \snr\inv \|_{\infty}^2 \\
    &\quad + K^2 Ln \sqrt{\log(n)} \alpha_{\max} \| \snr\inv \|_{\infty} + n L K \alpha_{\max} \numberthis \label{toverify}
\end{align*}
it holds that
\begin{align*}
    \| \sin\bTheta(\uhat,\U) \| &\lesssim   K^2 \sqrt{\log(n)} \frac{\big( \frac{1}{L} \| \snr\inv \|_2^2\big)^{1/2}}{\sqrt{L}\bar\lambda} + K^3 \log(n) \frac{\| \snr\inv \|_{\infty}^2}{\bar \lambda}  \\
    &\quad + K^2 \sqrt{ \log(n)} \frac{\alpha_{\max} \| \snr\inv \|_{\infty}}{\bar \lambda} + \frac{K \alpha_{\max}}{\bar \lambda}. \numberthis \label{sinthetabound2}
\end{align*}
Since the events listed above hold together with probability at least $1 - O(Ln^{-15})$, we see that the whole event holds with probability at least $1 - O(n^{-10})$ by the assumption that $L \lesssim n^{5}$.

We now verify \eqref{toverify}.  It is sufficient to check that the $\sin\bTheta$ bound in \eqref{sinthetabound2} is less than one (which is equivalent to checking \eqref{toverify}). In fact, we will show that each term is less than (in order) $\frac{1}{K}$, which is the second statement of the result. 

\cref{ass:networklevel} requires that
\begin{align*}
   C \bigg( \frac{\theta_{\max}\l}{\theta_{\min}\l} \bigg) \frac{K^{8} \theta_{\max}\l \|\theta\l\|_1 \log(n)}{\|\theta\l\|^4 (\lambda_{\min}\l)^2} \leq \bar \lambda.
\end{align*}
This immediately implies that $\frac{\alpha_{\max}}{\bar \lambda} \lesssim \frac{1}{K}$ from the definition of $\alpha_{\max}$.
By plugging in the definition of $\snr_l\inv$, we see that we require
\begin{align*}
   C \frac{K^{8} \theta_{\max}\l \|\theta\l\|_1 \log(n)}{\|\theta\l\|^2\snr_l^2} &\leq \bar \lambda \lambda_{\min}\l.
\end{align*}
Therefore the final three terms being are less than $\frac{1}{K}$ since $\frac{\theta_{\max}\l \|\theta\l\|_1}{\|\theta\l\|^2}$ is always larger than one.  For the remaining term, we observe that by averaging the above equation over $L$, we require that
\begin{align*}
C\frac{K^{8} \log(n)}{L}     \sum_{l}   \frac{\theta_{\max}\l \|\theta\l\|_1 }{\|\theta\l\|^2\snr_l^2} &\leq \bar \lambda^2. \numberthis \label{condcond}
\end{align*}
By squaring the first term, we see that we need the first term to satisfy 
\begin{align*}
    \frac{K^2 \log(n) }{L^2 \bar \lambda} \|\snr\inv\|_2^2 &\lesssim \frac{1}{K^2}.
\end{align*}
This is weaker than the condition \eqref{condcond}.
 The proof is now complete.
\end{proof}

\section{Proof of Second Stage Asymptotic Expansion (Theorem~\ref{thm:step2asympexp})}
\label{sec:step2expansion}
First we will restate \cref{thm:step2asympexp}.

\asymptoticexpansion*

 To prove \cref{thm:step2asympexp} we first state and prove several $\|\cdot\|_{2,\infty}$ concentration results for the residual terms that arise in the asymptotic expansion, and we prove \cref{thm:step2asympexp} in \cref{sec:steptwoasympexpproof}.

\subsection{Preliminary Lemmas: \texorpdfstring{$\ell_{2,\infty}$}{l\_\{2,Infinity\}} Residual Concentration Bounds} \label{sec:step2twoinftyprelim}

The following lemma bounds each of these residual terms in $\|\cdot\|_{2,\infty}.$
\begin{lemma}[Second Stage Residual Bounds] \label{lem:step2residualtwoinfty}
The following bounds hold with probability at least $1 - O(n^{-10})$:
\begin{align*}
    \| \U \U\t \mathcal{L}(\mathcal{E})\mathcal{ Y}\t \uhat \hat \Sigma^{-2} \|_{2,\infty} &\lesssim \frac{K^3 \sqrt{\log(n)}}{n \sqrt{L} \bar \lambda} \big( \frac{1}{L} \| \snr\inv \|_2^2\big)^{1/2}; \\
    \| (\I - \U \U\t ) \mathcal{R}_{\mathrm{all}} \uhat \hat \Sigma^{-2} \|_{2,\infty} &\lesssim  \frac{K^{3/2} \alpha_{\max}}{\sqrt{n} \bar \lambda} + \frac{K^{7/2} \log(n) \| \snr\inv \|_{\infty}^2}{\sqrt{n}\bar \lambda} + \frac{K^{5/2} \sqrt{\log(n)} \alpha_{\max} \| \snr\inv \|_{\infty}}{\sqrt{n \bar \lambda}}. 
\end{align*}
\end{lemma}

\begin{proof}[Proof of \cref{lem:step2residualtwoinfty}]
At the outset, we note that Weyl's inequality and the condition in \cref{thm:step2sintheta} implies that $\|\hat \Sigma^{-2}\| \lesssim K(\bar\lambda n L)\inv$ with high probability. 

We analyze each term separately.  First, we observe that
\begin{align*}
    \| \U \U\t \mathcal{L}(\mathcal{E})\mathcal{ Y}\t \uhat \hat \Sigma^{-2} \|_{2,\infty} &\lesssim \| \U \|_{2,\infty} \| \| \U\t \mathcal{L}(\mathcal{E}) \ycal\t \| \| \hat \Sigma^{-2} \| \\
    &\lesssim \frac{K^{3/2}}{\sqrt{n}} \frac{1}{nL \bar \lambda} \| \U\t \mathcal{L}(\mathcal{E}) \ycal\t \|.
\end{align*}
We now establish a concentration inequality for the term $\U\t \mathcal{L}(\mathcal{E}) \ycal\t$.  The result is similar to the proof of \cref{lem:LEmspectral}, so we postpone it to the end.  For now, we simply state that with probability at least $1 - O(n^{-20})$,
\begin{align*}
    \| \U\t \mathcal{L}(\mathcal{E}) \ycal\t \| &\lesssim K^{3/2} \sqrt{n L\log(n)} \big( \frac{1}{L} \| \snr\inv \|_2^2\big)^{1/2} + K^{3/2} \log(n)  \max_l \bigg( \frac{\theta_{\max}}{\theta_{\min}} \bigg) \frac{1}{\|\theta\m\|^2(\lambda_{\min}\m)^{1/2}} \numberthis\label{toprove}\\
    &\lesssim K^{3/2} \sqrt{n L\log(n)} \big( \frac{1}{L} \| \snr\inv \|_2^2\big)^{1/2},
\end{align*}
as long as $\max_l \frac{\theta_{\max}\m}{\theta_{\min}\m} \lesssim \sqrt{n/\log(n)}$, which holds under \cref{ass:networklevel}.  Putting it together, we obtain
\begin{align*}
     \| \U \U\t \mathcal{L}(\mathcal{E})\mathcal{ Y}\t \uhat \hat \Sigma^{-2} \|_{2,\infty} &\lesssim\frac{K^{3/2}}{\sqrt{n}} \frac{1}{nL \bar \lambda}  K^{3/2} \sqrt{n L\log(n)} \big( \frac{1}{L} \| \snr\inv \|_2^2\big)^{1/2} \\
     &\asymp \frac{K^3 \sqrt{\log(n)}}{n \sqrt{L} \bar \lambda} \big( \frac{1}{L} \| \snr\inv \|_2^2\big)^{1/2}.
\end{align*}
For the next term, we note that
\begin{align*}
    \|(\mathbf{I}- \U \U\t) \mathcal{R}_{\mathrm{all}} \uhat \hat \Sigma^{-2} \|_{2,\infty} &\leq \| \mathcal{R}_{\mathrm{all}} \|_{2,\infty} \| \hat \Sigma^{-2} \| + \|\U \|_{2,\infty} \| \mathcal{R}_{\mathrm{all}} \| \|\hat \Sigma^{-2} \| \\
    &\lesssim K\frac{\| \mathcal{R}_{\mathrm{all}} \|_{2,\infty}}{nL\bar \lambda} + \frac{K^{3/2}}{n^{3/2}L \bar \lambda} \|\mathcal{R}_{\mathrm{all}} \|.  \numberthis \label{ralltwoinfty}
\end{align*}
By \cref{lem:lineartermtwoinfty}, \cref{lem:LEmspectral}, and \cref{lem:step2resspectral}, we have the bounds
\begin{align*}
    \| \mathcal{L}(\mathbf{E}\m) \|_{2,\infty} &\lesssim \bigg( \frac{\theta_{\max}\l}{\theta_{\min}\l} \bigg)^{1/2} \frac{K\sqrt{\log(n)}}{(\lambda_{\min}\l)^{1/2} \|\theta\l\|}; \\
    \| \mathcal{L}(\mathbf{E}\m) \| &\leq \frac{K \sqrt{n\log (n)}}{(\lambda_{\min} \l)^{1/2}\left\|\theta^{(l)}\right\|}\left(\frac{\theta_{\max }\m}{\theta_{\min }\m}\right)^{1 / 2}; \\
    \| \mathcal{R}_{\mathrm{all}} \| &\lesssim L K^2 n \log(n) \| \snr\inv \|_{\infty}^2 + K L n \sqrt{\log(n)} \| \snr\inv \|_{\infty} \alpha_{\max} + n L \alpha_{\max}; \\
     \| \mathcal{R}_{\mathrm{Stage \ I}}\m \|_{2,\infty} &\lesssim \alpha_{\max}; \\
     \| \mathcal{R}_{\mathrm{Stage \ I}}\m \| &\lesssim \sqrt{n} \alpha_{\max}
\end{align*} 
Therefore, we obtain
\begin{align*}
    \| \mathcal{R}_{\mathrm{all}} \|_{2,\infty} &\lesssim L \max_l \| \mathcal{L}(\mathbf{E}\m) \|_{2,\infty} \| \mathcal{L}(\mathbf{E}\m) \| + L \max_l\| \mathcal{L}(\mathbf{E}\m) \|_{2,\infty} \| \mathcal{R}_{Y}\m \| \\
     &\qquad +L \max_l \| \mathcal{R}_{\mathrm{Stage \ I}}\m \|_{2,\infty} \| \mathcal{L}(\mathbf{E}\m) \| +L \max_l \| \mathcal{R}_{\mathrm{Stage \ I}}\m \|_{2,\infty} \|  \mathcal{R}_{\mathrm{Stage \ I}}\m\| \\
     &\qquad + L \max_l \|\mathcal{R}_{\mathrm{Stage \ I}}\m \|_{2,\infty} \| \mathbf{Y}\m \| +L \max_l \| \ytilde\m \|_{2,\infty} \| \mathcal{R}_{\mathrm{Stage \ I}}\m \| \\
     &\lesssim  L \bigg( \frac{\theta_{\max}\l}{\theta_{\min}\l} \bigg)^{1/2} \frac{K\sqrt{\log(n)}}{(\lambda_{\min}\l)^{1/2} \|\theta\l\|}  \frac{K \sqrt{n\log (n)}}{(\lambda_{\min}\l)^{1/2}\left\|\theta^{(m)}\right\|}\left(\frac{\theta_{\max }\m}{\theta_{\min }\m}\right)^{1 / 2}
       \\
     &\qquad + L \max_l \bigg( \frac{\theta_{\max}\l}{\theta_{\min}\l}\bigg)^{1/2} \frac{K \sqrt{\log(n)}}{(\lambda_{\min}\l)^{1/2} \|\theta\l\|^2} \sqrt{n} \alpha_{\max} \\
     &\qquad+ L  \alpha_{\max} \max_l \bigg( \frac{\theta_{\max}\m}{\theta_{\min}\m} \bigg)^{1/2} \frac{K \sqrt{n\log(n)}}{(\lambda_{\min}\m)^{1/2} \|\theta\m\|}  +L \sqrt{n} \alpha_{\max}^2 + \sqrt{n} L\alpha_{\max} \\
     &\asymp L \sqrt{n} \bigg( \frac{\theta_{\max}\l}{\theta_{\min}\l} \bigg) \frac{K^2 \log(n)}{\lambda_{\min}\l \|\theta\l\|^2}  + L \sqrt{n} \alpha_{\max},
     \end{align*}
     where we have used the assumption that      \begin{align*}
\bigg( \frac{\theta_{\max}\l}{\theta_{\min}\l} \bigg)^{1/2} \frac{K \sqrt{\log(n)}}{(\lambda_{\min}\l)^{1/2} \|\theta\|} \lesssim 1. \numberthis \label{toverify2}
     \end{align*}
 We will verify this momentarily. 
Plugging this into \eqref{ralltwoinfty}, we obtain that
\begin{align*}
     \|(\mathbf{I}- \U \U\t) \mathcal{R}_{\mathrm{all}} \uhat \hat \Sigma^{-2} \|_{2,\infty} &\lesssim K\frac{\| \mathcal{R}_{\mathrm{all}} \|_{2,\infty}}{nL\bar \lambda} + \frac{K^{3/2}}{n^{3/2} L \bar \lambda} \|\mathcal{R}_{\mathrm{all}} \| \\
     &\lesssim \frac{K}{nL \bar \lambda} \bigg( L \sqrt{n}  \max_l\bigg( \frac{\theta_{\max}\l}{\theta_{\min}\l} \bigg) \frac{K^2 \log(n)}{\lambda_{\min}\l \|\theta\l\|^2}  + L \sqrt{n} \alpha_{\max} \bigg) \\
     &\quad +  \frac{K^{3/2}}{n^{3/2}L \bar \lambda} \bigg\{  L K^2 n \log(n) \| \snr\inv \|_{\infty}^2 + K L n \sqrt{\log(n)} \| \snr\inv \|_{\infty} \alpha_{\max}\\
     & \quad\quad+ n L \alpha_{\max} \bigg\} \\
     &\asymp 
     \frac{K^{3/2} \alpha_{\max}}{\sqrt{n} \bar \lambda} + \frac{K^{7/2} \log(n) \| \snr\inv \|_{\infty}^2}{\sqrt{n}\bar \lambda} + \frac{K^{5/2} \sqrt{\log(n)} \alpha_{\max} \| \snr\inv \|_{\infty}}{\sqrt{n \bar \lambda}},
\end{align*}
since
\begin{align*}
    \| \snr\inv \|_{\infty} &= \max_l \bigg( \frac{\theta_{\max}\l}{\theta_{\min}\l} \bigg)^{1/2} \frac{1}{(\lambda_{\min}\l)^{1/2} \|\theta\l\|}.
\end{align*}
which holds with probability at least $1 - O(n^{-15})$. 

We now verify \eqref{toverify2}.     By  \cref{ass:networklevel}, the definition of $\snr$, and the fact that $\frac{\theta_{\max}\l \|\theta\l\|_1}{\|\theta\l\|^2 } \geq 1$, it holds that $\bar \lambda \geq K^{5} \log(n) \| \snr\inv \|_{\infty}^2$, which in particular implies that $K^2 \log(n) \| \snr\|_{\infty}^2 \leq 1$ since $\bar \lambda \leq 1$. This verifies \eqref{toverify2}.

Therefore, we will have completed the proof provided we can establish the bound \eqref{toprove}.  Observe that 
\begin{align*}
    \U\t \mathcal{L}(\mathcal{E}) \ycal\t &= \sum_{l} \U\t \bigg[ \sum_{i\leq j} \mathbf{E}_{ij}\m e_i e_j\t \bigg]\big( \U\m |\Lambda\m|^{-1/2} \ipq \m \mathbf{J}(\xtilde_{i\cdot}) \big) (\ytilde\m)\t \\
    &\qquad + \sum_{l} \U\t \bigg[ \sum_{j < i} \mathbf{E}_{ij}\m e_j e_i\t \bigg]\big( \U\m |\Lambda\m|^{-1/2} \ipq \m \mathbf{J}(\xtilde_{j\cdot}) \big) (\ytilde\m)\t \\
    &= \sum_{l} \sum_{i\leq j} \mathbf{E}_{ij}\m \U\t e_i e_j\t \big( \U\m |\Lambda\m|^{-1/2} \ipq \m \mathbf{J}(\xtilde_{i\cdot}) \big) (\ytilde\m)\t \\
    &\qquad + \sum_{l} \sum_{j < i} \mathbf{E}_{ij}\m \U\t e_j e_i\t \big( \U\m |\Lambda\m|^{-1/2} \ipq \m \mathbf{J}(\xtilde_{j\cdot}) \big) (\ytilde\m)\t,
\end{align*}
both of which are a sum of independent random matrices.  We bound the first term now; the second is similar.  We will apply Matrix Bernstein  (Corollary 3.3 of \cite{chen_spectral_2021}).  To wit, we need to bound
\begin{align*}
    v :&= \sum_{l} \sum_{i\leq j} \E \big( \mathbf{E}\m_{ij}\big)^2 \| \U\t e_i e_j\t \big( \U\m |\Lambda\m|^{-1/2} \ipq \m \mathbf{J}(\xtilde_{i\cdot}) \big) (\ytilde\m)\t \|^2; \\
    w :&= \max_{m,i,j} \| \U\t e_i e_j\t \big( \U\m |\Lambda\m|^{-1/2} \ipq \m \mathbf{J}(\xtilde_{i\cdot}) \big) (\ytilde\m)\t \|.
\end{align*}
We observe that
\begin{align*}
    v &\leq \sum_{l} \sum_{i\leq j} \theta_i\m \theta_j\m \| \U\t e_i e_j\t \big( \U\m |\Lambda\m|^{-1/2} \ipq \m \mathbf{J}(\xtilde_{i\cdot}) \big) (\ytilde\m)\t \|^2\\
    &\leq \sum_l \sum_{i\leq j} \theta_i\m \theta_j\m \| \U \|_{2,\infty}^2 \| e_j\t \big( \U\m |\Lambda\m|^{-1/2} \ipq \m \mathbf{J}(\xtilde_{i\cdot}) \big) (\ytilde\m)\t \|^2 \\
    &\lesssim \frac{K}{n} \sum_l \sum_{i\leq j} \theta_i\m \theta_j\m
 \| e_j\t \U\m \|^2 \| |\Lambda\m|^{-1/2} \|^2 \| \mathbf{J}(\xtilde_{i\cdot}) \|^2 \| \ytilde\m \|^2 \\
 &\lesssim K^3 \sum_l \sum_{i\leq j} \theta_i\m \theta_j\m\frac{(\theta_{j}\l)^2}{\|\theta\m\|^2} \frac{1}{\lambda_{\min}\m \|\theta\m\|^2} \frac{1}{(\theta_i\m)^2 } \\
 &\lesssim K^3 \sum_l \frac{1}{\lambda_{\min}\m \|\theta\m\|^4} \sum_{i,j} \frac{\theta_j\m}{\theta_i\m} (\theta_j\m)^2 \\
 &\lesssim K^3 n \sum_l \frac{1}{\lambda_{\min}\m \|\theta\m\|^2} \bigg( \frac{\theta_{\max}\m}{\theta_{\min}\m} \bigg) \\
 &= K^3 n L \big( \frac{1}{L} \| \snr\inv \|_2^2\big),
  \end{align*}
  where we have implicitly used \cref{lem:popprop}.  In addition, via similar arguments,
  \begin{align*}
      w &\lesssim K^{3/2} \max_l \bigg( \frac{\theta_{\max}}{\theta_{\min}} \bigg) \frac{1}{\|\theta\m\|^2(\lambda_{\min}\m)^{1/2}}. 
  \end{align*}
  Therefore, the result is completed by applying Matrix Bernstein.  This completes the proof.
\end{proof}

The following result bounds several additional ``approximate commutation'' terms, analogous to \cref{lem:approximatecommutation} for Stage 1.  
\begin{lemma}[Second Stage Approximate Commutation] \label{lem:step2approximatecommutation}
The following bounds hold with probability at least $1 - O(n^{-10}):$
\begin{align*}
    \| \U\t \uhat -  \wstar  \| &\lesssim \bigg( K^2 \sqrt{\log(n)} \frac{\big( \frac{1}{L} \| \snr\inv \|_2^2\big)^{1/2}}{\sqrt{L}\bar\lambda} + K^3 \log(n) \frac{\| \snr\inv \|_{\infty}^2}{\bar \lambda} \\
    &\qquad \qquad + K^2 \sqrt{ \log(n)} \frac{\alpha_{\max} \| \snr\inv \|_{\infty}}{\bar \lambda} + \frac{K\alpha_{\max}}{\bar \lambda} \bigg)^2; \\
    \|  \Sigma^{-2} \U\t \uhat - \U\t \uhat \hat \Sigma^{-2} \| &\lesssim\frac{K^3 \sqrt{\log(n)} \big( \frac{1}{L} \| \snr\inv \|_2^2\big)^{1/2}}{n L^{3/2} \bar \lambda^2} + \frac{K^4 \log(n) \| \snr\inv \|_{\infty}^2}{n L \bar \lambda^2} \\
    &\qquad + \frac{K^3 \sqrt{\log(n)} \alpha_{\max} \| \snr\inv \|_{\infty}}{n L \bar \lambda^2} + \frac{K^2\alpha_{\max}}{n L \bar \lambda^2},; \\
    \| \sum_{l}\mathcal{L}(\mathbf{E}\m) ( \ytilde\m)\t \|_{2,\infty} &\lesssim K \sqrt{Ln\log(n)} \big( \frac{1}{L} \| \snr\inv \|_2^2\big)^{1/2}; \\
\end{align*}
\end{lemma}

\begin{proof}
For the first bound, we observe that
\begin{align*}
    \| \U\t \uhat - \wstar  \| &\lesssim  \| \sin\bTheta(\uhat,\U) \|^2 \\
    &\lesssim  \bigg( K^2 \sqrt{\log(n)} \frac{\big( \frac{1}{L} \| \snr\inv \|_2^2\big)^{1/2}}{\sqrt{L}\bar\lambda} + K^3 \log(n) \frac{\| \snr\inv \|_{\infty}^2}{\bar \lambda} \\
    &\qquad + K^2 \sqrt{ \log(n)} \frac{\alpha_{\max} \| \snr\inv \|_{\infty}}{\bar \lambda} + \frac{K\alpha_{\max}}{\bar \lambda} \bigg)^2,
\end{align*}
where the final inequality holds by \cref{thm:step2sintheta}, with probability at least $1 - O(n^{-10})$.  

For the second bound, we observe that
\begin{align*}
    \| \Sigma^{-2} &\U\t \uhat - \U\t \uhat \hat \Sigma^{-2} \| \\
    &= \| \Sigma^{-2} \big( \U\t \uhat \hat \Sigma^2 - \Sigma^2 \U\t \uhat \big) \hat \Sigma^{-2} \| \\
    &\lesssim \frac{K^2}{n^2L^2 \bar \lambda^2} \| \U\t \uhat \hat \Sigma^2 - \Sigma^2 \U\t \uhat \| \\
    &\lesssim \frac{K^2}{n^2 L^2 \bar \lambda^2} \| \U\t (\yhatcal\yhatcal\t -\ycal \ycal\t )\uhat  \| \\
    &\lesssim \frac{K^2}{n^2 L^2 \bar \lambda^2} \bigg\{ \| \mathcal{L}(\mathcal{E}) \ycal\t \| +  \|\mathcal{R}_{\mathrm{all}}\| \bigg\} \\
    &\lesssim \frac{K^2}{n^2 L^2 \bar\lambda^2} \bigg\{ K n \sqrt{L\log(n)} \big( \frac{1}{L} \| \snr\inv \|_2^2\big)^{1/2}  + L K^{2} n \log (n) \| \snr\inv \|_{\infty}^{2}\\
    &\qquad \qquad +K L n \sqrt{\log (n)} \alpha_{\max } \| \snr\inv \|_{\infty}+n L \alpha_{\max } \bigg\} \\
    &\asymp  \frac{K^3 \sqrt{\log(n)} \big( \frac{1}{L} \| \snr\inv \|_2^2\big)^{1/2}}{n L^{3/2} \bar \lambda^2} + \frac{K^4 \log(n) \| \snr\inv \|_{\infty}^2}{n L \bar \lambda^2} \\
    &\qquad + \frac{K^3 \sqrt{\log(n)} \alpha_{\max} \| \snr\inv \|_{\infty}}{n L \bar \lambda^2} + \frac{K^2\alpha_{\max}}{n L \bar \lambda^2},
\end{align*}
which holds with probability at least $1 - O(n^{-10})$ by \cref{lem:LEmspectral} and \cref{lem:step2rallspectralbound}.  

For the third term, we note that we can write the $i$'th row of the matrix in question via
\begin{align*}
\bigg(\sum_{l} \mathcal{L}(\mathbf{E}\m) (\ytilde\m)\t \bigg)_{i\cdot} &= 
    \sum_{l} \sum_{j} \mathbf{E}\m_{ij} \bigg( \U\m |\Lambda\m|^{-1/2} \ipq\m \mathbf{J}(\xtilde_{i\cdot}) (\ytilde\m)\t  \bigg)_{j\cdot},
\end{align*}
which is a sum of independent random matrices.  To wit, we bound via the Matrix Bernstein inequality (Corollary 3.3 of \cite{chen_spectral_2021}).  The proof is similar to \cref{lem:LEmspectral} (amongst others), so we omit the detailed proof for brevity.  Matrix Bernstein then implies that with probability at least $1 - O(n^{-11})$ that
\begin{align*}
 \bigg\|    \sum_{l} &\sum_{j} \mathbf{E}\m_{ij} \bigg( \U\m |\Lambda\m|^{-1/2} \ipq\m \mathbf{J}(\xtilde_{i\cdot}) (\ytilde\m)\t \bigg)_{j\cdot} \bigg\| \\
 &\lesssim K \sqrt{\log(n)} \max_l \| (\ytilde\m)\t  \| \bigg( \sum_{l} \bigg( \frac{\theta_{\max}\m}{\theta_{\min}\m} \bigg) \frac{1}{\lambda_{\min}\m \|\theta\m\|^2} \bigg)^{1/2} \\
 &\lesssim K \sqrt{L\log(n)} \big( \frac{1}{L} \| \snr\inv \|_2^2\big)^{1/2}  \max_l \| (\ytilde\m)\t \| \\
 &\lesssim 
 K \sqrt{Ln\log(n)} \big( \frac{1}{L} \| \snr\inv \|_2^2\big)^{1/2}.
\end{align*}
Taking a union bound over all $n$ rows completes the proof of this bound.  
\end{proof}

\subsection{Proof of Theorem~\ref{thm:step2asympexp}} \label{sec:steptwoasympexpproof}

\begin{proof}[Proof of \cref{thm:step2asympexp}]
First, recall we have the expansion 
\begin{align*}
      \mathcal{\hat Y} \mathcal{\hat Y}\t -  \mathcal{ Y}\mathcal{ Y}\t &= \sum_{l=1}^L \yhat\m (\yhat\m)\t - (\ytilde\m) (\ytilde\m)\t \\
     :&=   \mathcal{L}(\mathcal{E})\mathcal{ Y}\t + \mathcal{ Y} \mathcal{L}(\mathcal{E})\t + \mathcal{R}_{\mathrm{all}},
\end{align*}
where recall we define
\begin{align*}
    \mathcal{R}_{\mathrm{all}} :&= \sum_{l} \mathcal{L}(\mathbf{E}\m)\mathcal{L}(\mathbf{E}\m)\t + \mathcal{L}(\mathbf{E}\m) (\mathcal{R}^{(l)})\t+ \mathcal{R}^{(l)}\mathcal{L}(\mathbf{E}\m)\t + \mathcal{R}^{(l)} (\mathcal{R}^{(l)})\t \\
        &\quad + \sum_l \mathcal{R}^{(l)} (\ytilde\m)\t + \ytilde\m (\mathcal{R}^{(l)})\t,
\end{align*}
and
\begin{align*}
\mathcal{L}(\mathcal{E}) :&= \big[ \mathcal{L}(\mathbf{E}\one), \cdots , \mathcal{L}(\mathbf{E}\M) \big],
\end{align*}
and hence that
\begin{align*}
    \mathcal{L}(\mathcal{E})\mathcal{ Y}\t + \mathcal{ Y} \mathcal{L}(\mathcal{E})\t &= \sum_l \mathcal{L}(\mathbf{E}\m) (\ytilde\m)\t + \ytilde\m \mathcal{L}(\mathbf{E}\m)\t .
\end{align*}
We now study how well $\uhat$ approximates $\U$ in an entrywise sense.  We start with the expansion:
\begin{align*}
    \uhat - \U \wstar &= (\I - \U \U\t) \big( \mathcal{ Y}\mathcal{ Y}\t - \mathcal{\hat Y}\mathcal{\hat Y}\t \big) \uhat \hat \Sigma^{-2} + \U (\U\t \uhat - \wstar) \\
    &=  (\I - \U \U\t) \big( \mathcal{L}(\mathcal{E})\mathcal{ Y}\t + \mathcal{ Y} \mathcal{L}(\mathcal{E})\t + \mathcal{R}_{\mathrm{all}} \big) \uhat \hat \Sigma^{-2}+ \U (\U\t \uhat - \wstar) \\
    &= \mathcal{L}(\mathcal{E})\mathcal{ Y}\t \uhat \hat \Sigma^{-2}  - \U \U\t \mathcal{L}(\mathcal{E})\mathcal{ Y}\t \uhat \hat \Sigma^{-2} \\
    &\quad + (\I - \U \U\t ) \mathcal{ Y} \mathcal{L}(\mathcal{E})\t\uhat \hat \Sigma^{-2} + (\I - \U \U\t ) \mathcal{R}_{\mathrm{all}}  \uhat \hat \Sigma^{-2}+ \U (\U\t \uhat - \wstar) \\
    &= \mathcal{L}(\mathcal{E})\mathcal{ Y}\t \uhat \hat \Sigma^{-2}  - \U \U\t \mathcal{L}(\mathcal{E})\mathcal{ Y}\t \uhat \hat \Sigma^{-2} \\
    &\quad +(\I - \U \U\t ) \mathcal{R}_{\mathrm{all}}  \uhat \hat \Sigma^{-2}+ \U (\U\t \uhat - \wstar) , \numberthis \label{initialexpansion}
\end{align*}
where we have observed that the term
\begin{align*}
    (\I - \U \U\t) \mathcal{ Y} \mathcal{L}(\mathcal{E})\t \uhat \hat \Sigma^{-2} \equiv 0,
\end{align*}
since $\mathcal{ Y}$ has left singular vectors $\U$.  We now expand the first-order term out further.  Observe that
\begin{align*}
    \sum_{l} \mathcal{L}(\mathbf{E}\m) ( \ytilde\m)\t \uhat\hat\Sigma^{-2} &= \sum_{l}\mathcal{L}(\mathbf{E}\m) ( \ytilde\m)\t \U \Sigma^{-2} \wstar + \sum_{l}\mathcal{L}(\mathbf{E}\m) ( \ytilde\m)\t \U \Sigma^{-2} \big( \wstar - \U\t \uhat \big) \\
    &\quad + \sum_{l}\mathcal{L}(\mathbf{E}\m) ( \ytilde\m)\t \U\big( \Sigma^{-2} \U\t \uhat - \U\t \uhat \hat \Sigma^{-2} \big) \\
    &\quad + \sum_{l}\mathcal{L}(\mathbf{E}\m) ( \ytilde\m)\t \big( \uhat - \U \U\t \uhat \big) \hat \Sigma^{-2}. \numberthis \label{firstorderterm}
\end{align*}
Plugging \eqref{firstorderterm} into \eqref{initialexpansion} yields the full expansion
\begin{align*}
     \uhat - \U \wstar &=\sum_{l}\mathcal{L}(\mathbf{E}\m) ( \ytilde\m)\t \U \Sigma^{-2} \wstar + \sum_{l}\mathcal{L}(\mathbf{E}\m) ( \ytilde\m)\t \U \Sigma^{-2} \big( \wstar - \U\t \uhat \big) \\
    &\quad + \sum_{l}\mathcal{L}(\mathbf{E}\m) ( \ytilde\m)\t \U\big( \Sigma^{-2} \U\t \uhat - \U\t \uhat \hat \Sigma^{-2} \big) \\
    &\quad + \sum_{l}\mathcal{L}(\mathbf{E}\m) ( \ytilde\m)\t \big( \uhat - \U \U\t \uhat \big) \hat \Sigma^{-2}\\
    &\qquad - \U \U\t \mathcal{L}(\mathcal{E})\mathcal{ Y}\t \uhat \hat \Sigma^{-2} \\
    &\quad +(\I - \U \U\t ) \mathcal{R}_{\mathrm{all}}  \uhat \hat \Sigma^{-2}+ \U (\U\t \uhat - \wstar).
\end{align*}
Multiplying through by $\wstar\t$ yields
\begin{align*}
    \uhat \wstar\t - \U &= \sum_{l}\mathcal{L}(\mathbf{E}\m) ( \ytilde\m)\t \U \Sigma^{-2} + \sum_{l}\mathcal{L}(\mathbf{E}\m) ( \ytilde\m)\t \U \Sigma^{-2} \big( \wstar - \U\t \uhat \big) \wstar\t \\
    &\quad + \sum_{l}\mathcal{L}(\mathbf{E}\m) ( \ytilde\m)\t \U\big( \Sigma^{-2} \U\t \uhat - \U\t \uhat \hat \Sigma^{-2} \big) \wstar\t \\
    &\quad + \sum_{l}\mathcal{L}(\mathbf{E}\m) ( \ytilde\m)\t \big( \uhat - \U \U\t \uhat \big) \hat \Sigma^{-2} \wstar\t \\
    &\qquad - \U \U\t \mathcal{L}(\mathcal{E})\mathcal{ Y}\t \uhat \hat \Sigma^{-2} \wstar\t \\
    &\quad +(\I - \U \U\t ) \mathcal{R}_{\mathrm{all}}  \uhat \hat \Sigma^{-2}+ \U (\U\t \uhat - \wstar) \wstar\t \\
    :&= \sum_{l}\mathcal{L}(\mathbf{E}\m) ( \ytilde\m)\t \U \Sigma^{-2} + \mathbf{R}_1 + \mathbf{R}_2 + \mathbf{R}_3 + \mathbf{R}_4 + \mathbf{R}_5 + \mathbf{R}_6,
\end{align*}
where
\begin{align*}
    \mathbf{R}_1 :&= \sum_{l}\mathcal{L}(\mathbf{E}\m) ( \ytilde\m)\t \U \Sigma^{-2} \big( \wstar - \U\t \uhat \big) \wstar\t; \\
    \mathbf{R}_2 :&= \sum_{l}\mathcal{L}(\mathbf{E}\m) ( \ytilde\m)\t \U\big( \Sigma^{-2} \U\t \uhat - \U\t \uhat \hat \Sigma^{-2} \big) \wstar\t; \\
    \mathbf{R}_3 :&= \sum_{l}\mathcal{L}(\mathbf{E}\m) ( \ytilde\m)\t \big( \uhat - \U \U\t \uhat \big) \hat \Sigma^{-2} \wstar\t; \\
    \mathbf{R}_4:&= -\U \U\t \mathcal{L}(\mathcal{E})\mathcal{ Y}\t \uhat \hat \Sigma^{-2} \wstar\t ; \\   \mathbf{R}_5 :&= (\I - \U \U\t ) \mathcal{R}_{\mathrm{all}}  \uhat \hat \Sigma^{-2}; \\
      \mathbf{R}_6 :&=\U (\U\t \uhat - \wstar) \wstar\t.
\end{align*}
By \cref{lem:step2residualtwoinfty},  we have the bounds
\begin{align*}
    \| \mathbf{R}_4 \|_{2,\infty} &\lesssim  \frac{K^3 \sqrt{\log(n)}}{n\sqrt{L} \bar \lambda} \big( \frac{1}{L} \| \snr\inv \|_2^2\big)^{1/2}; \\
    \| \mathbf{R}_5 \|_{2,\infty} &\lesssim \frac{K^{3/2} \alpha_{\max}}{\sqrt{n} \bar \lambda} + \frac{K^{7/2} \log(n) \| \snr\inv \|_{\infty}^2}{\sqrt{n}\bar \lambda} + \frac{K^{5/2} \sqrt{\log(n)} \alpha_{\max} \| \snr\inv \|_{\infty}}{\sqrt{n \bar \lambda}}. 
\end{align*}
In addition, by properties of the $\ell_{2,\infty}$ norm and \cref{lem:step2approximatecommutation}, it holds that
\begin{align*}
    \| \mathbf{R}_6 \|_{2,\infty} &\leq \| \U \|_{2,\infty} \| \U\t \uhat - \wstar \| \\
    &\lesssim \sqrt{\frac{K}{n}} \bigg( K^2 \sqrt{\log(n)} \frac{\big( \frac{1}{L} \| \snr\inv \|_2^2\big)^{1/2}}{\sqrt{L}\bar\lambda} + K^3 \log(n) \frac{\| \snr\inv \|_{\infty}^2}{\bar \lambda} \\
    &\qquad \qquad + K^2 \sqrt{ \log(n)} \frac{\alpha_{\max} \| \snr\inv \|_{\infty}}{\bar \lambda} + \frac{K\alpha_{\max}}{\bar \lambda} \bigg)^2 \\
    &\lesssim \frac{K^{9/2} \log(n) \big( \frac{1}{L} \| \snr\inv \|_2^2\big)}{\sqrt{n} L \bar \lambda^2} + \frac{K^{7/2} \log(n) \| \snr\inv\|_{\infty}^2}{\sqrt{n} \bar \lambda} \\
    &\quad + \frac{K^{5/2} \sqrt{\log(n)} \alpha_{\max}\| \snr\inv \|_{\infty}}{\sqrt{n} \bar \lambda} + \frac{K^{3/2} \alpha_{\max}^2}{\sqrt{n} \bar \lambda^2}, 
\end{align*}
where we have used the fact that each of the terms inside of the parentheses on the bound for $\mathbf{R}_6$ is less than one, which was verified in the proof of \cref{thm:step2sintheta} (note that these terms in parentheses are simply the $\sin\bTheta$ upper bound).

Combining these, we obtain that with probability at least $1 - O(n^{-10})$,
\begin{align*}
    \| \mathbf{R}_4 \|_{2,\infty} + \|\mathbf{R}_5\|_{2,\infty} + \|\mathbf{R}_6\|_{2,\infty} &\lesssim \frac{K^3 \sqrt{\log(n)}}{n\sqrt{L} \bar \lambda} \big( \frac{1}{L} \| \snr\inv \|_2^2\big)^{1/2} + \frac{K^{7/2} \log(n) \big( \frac{1}{L} \| \snr\inv \|_2^2\big)}{\sqrt{n} L \bar \lambda^2} \\
    &\quad + \frac{K^{7/2} \log(n) \| \snr\inv \|_{\infty}^2}{\sqrt{n} \bar \lambda}  + \frac{K^{5/2} \sqrt{\log(n)} \alpha_{\max}\| \snr\inv \|_{\infty}}{\sqrt{n}\bar \lambda} + \frac{ \alpha_{\max}}{\sqrt{n} \bar \lambda},
\end{align*}
where we have used the fact that $\frac{K^{3/2} \alpha_{\max}}{\bar \lambda} \lesssim 1$.

For the terms $\mathbf{R}_1$ through $\mathbf{R}_3$, we observe that
\begin{align*}
    \| \mathbf{R}_1 \|_{2,\infty} &\lesssim \frac{K}{nL\bar \lambda}\| \sum_l \mathcal{L}(\mathbf{E}\m) (\ytilde\m)\t \|_{2,\infty} \| \wstar - \U\t \uhat \|; \\
    \| \mathbf{R}_2 \|_{2,\infty} &\lesssim \| \sum_l \mathcal{L}(\mathbf{E}\m) (\ytilde\m)\t \|_{2,\infty} \| \Sigma^{-2} \U\t \uhat - \U\t \uhat \hat \Sigma^{-2} \|; \\
    \| \mathbf{R}_3 \|_{2,\infty} &\lesssim \frac{K}{nL\bar \lambda} \| \sum_{l} \mathcal{L}(\mathbf{E}\m) (\ytilde\m)\t \|_{2,\infty} \| \sin\bTheta(\uhat,\U) \| .
    \end{align*}
\cref{lem:step2approximatecommutation} shows that
 with probability at least $1 - O(n^{-10})$ that
\begin{align*}
 \| \U\t \uhat -  \wstar  \| &\lesssim \bigg( K^2 \sqrt{\log(n)} \frac{\big( \frac{1}{L} \| \snr\inv \|_2^2\big)^{1/2}}{\sqrt{L}\bar\lambda} + K^3 \log(n) \frac{\| \snr\inv \|_{\infty}^2}{\bar \lambda} \\
    &\qquad \qquad + K^2 \sqrt{ \log(n)} \frac{\alpha_{\max} \| \snr\inv \|_{\infty}}{\bar \lambda} + \frac{K\alpha_{\max}}{\bar \lambda} \bigg)^2; \\
    \|  \Sigma^{-2} \U\t \uhat - \U\t \uhat \hat \Sigma^{-2} \| &\lesssim\frac{K^3 \sqrt{\log(n)} \big( \frac{1}{L} \| \snr\inv \|_2^2\big)^{1/2}}{n L^{3/2} \bar \lambda^2} + \frac{K^4 \log(n) \| \snr\inv \|_{\infty}^2}{n L \bar \lambda^2} \\
    &\qquad + \frac{K^3 \sqrt{\log(n)} \alpha_{\max} \| \snr\inv \|_{\infty}}{n L \bar \lambda^2} + \frac{K^2\alpha_{\max}}{n L \bar \lambda^2}; \\
    \| \sum_{l}\mathcal{L}(\mathbf{E}\m) ( \ytilde\m)\t \|_{2,\infty} &\lesssim K \sqrt{Ln\log(n)} \big( \frac{1}{L} \| \snr\inv \|_2^2\big)^{1/2}. 
%
\end{align*}
In addition, by \cref{thm:step2sintheta}, we have that
\begin{align*}
    \|\sin\bTheta(\uhat,\U)\| &\lesssim  K^2 \sqrt{\log(n)} \frac{\big( \frac{1}{L} \| \snr\inv \|_2^2\big)^{1/2}}{\sqrt{L}\bar\lambda} + K^3 \log(n) \frac{\| \snr\inv \|_{\infty}^2}{\bar \lambda} \\
    &\quad + K^2 \sqrt{ \log(n)} \frac{\alpha_{\max} \| \snr\inv \|_{\infty}}{\bar \lambda} + \frac{K\alpha_{\max}}{\bar \lambda}.
\end{align*}
Plugging these bounds in yields that
\begin{align*}
     \| \mathbf{R}_1 \|_{2,\infty} &\lesssim \frac{K}{nL\bar \lambda} K \sqrt{Ln\log(n)} ( \frac{1}{L} \|\snr\inv\|_2^2 )^{1/2} \\
     &\quad \times \bigg( K^2 \sqrt{\log(n)} \frac{\big( \frac{1}{L} \| \snr\inv \|_2^2\big)^{1/2}}{\sqrt{L}\bar\lambda} + K^3 \log(n) \frac{\| \snr\inv \|_{\infty}^2}{\bar \lambda} \\
    &\qquad \qquad + K^2 \sqrt{ \log(n)} \frac{\alpha_{\max} \| \snr\inv \|_{\infty}}{\bar \lambda} + \frac{K\alpha_{\max}}{\bar \lambda} \bigg)^2 \\
     &\asymp \frac{K^4 \log(n) \big( \frac{1}{L} \| \snr\inv \|_2^2\big)}{L \sqrt{n} \bar \lambda^2} + \frac{K^5 \log^{3/2}(n) \big( \frac{1}{L} \| \snr\inv \|_2^2\big)^{1/2} \| \snr\inv \|_{\infty}^2}{\sqrt{nL} \bar \lambda^2} \\
     &\quad + \frac{K^4 \log(n) \alpha_{\max}\| \snr\inv \|_{\infty} \big( \frac{1}{L} \| \snr\inv \|_2^2\big)^{1/2}}{\sqrt{nL} \bar \lambda^2} + \frac{K^3 \sqrt{\log(n)}\big( \frac{1}{L} \| \snr\inv \|_2^2\big)^{1/2} \alpha_{\max}}{\sqrt{nL} \bar \lambda^2}; \\
     \| \mathbf{R}_2\|_{2,\infty} &\lesssim K \sqrt{Ln\log(n)} \big( \frac{1}{L} \| \snr\inv \|_2^2\big)^{1/2}  \\
     &\quad \times \bigg( \frac{K^3 \sqrt{\log(n)} \big( \frac{1}{L} \| \snr\inv \|_2^2\big)^{1/2}}{n L^{3/2} \bar \lambda^2} + \frac{K^4 \log(n) \| \snr\inv \|_{\infty}^2}{n L \bar \lambda^2} \\
    &\qquad + \frac{K^3 \sqrt{\log(n)} \alpha_{\max} \| \snr\inv \|_{\infty}}{n L \bar \lambda^2} + \frac{K^2\alpha_{\max}}{n L \bar \lambda^2} \bigg) \\
     &\asymp \frac{K^4 \log(n) \big( \frac{1}{L} \| \snr\inv \|_2^2\big)}{\sqrt{n} L \bar \lambda^2} + \frac{K^5 \log^{3/2}(n) \| \snr\inv \|_{\infty}^2 \big( \frac{1}{L} \| \snr\inv \|_2^2\big)^{1/2}}{\sqrt{nL} \bar\lambda^2} \\
     &\quad + \frac{K^4 \log(n) \alpha_{\max}\| \snr\inv \|_{\infty} \big( \frac{1}{L} \| \snr\inv \|_2^2\big)^{1/2}}{\sqrt{nL}\bar\lambda^2} + \frac{K^3 \sqrt{\log(n)} \alpha_{\max} \big( \frac{1}{L} \| \snr\inv \|_2^2\big)^{1/2}}{\sqrt{nL} \bar \lambda^2}; \\
     \| \mathbf{R}_3 \|_{2,\infty} &\lesssim \frac{K}{nL \bar \lambda}  K \sqrt{Ln\log(n)} \big( \frac{1}{L} \| \snr\inv \|_2^2\big)^{1/2} \\
     &\times \bigg( K^2 \sqrt{\log(n)} \frac{\big( \frac{1}{L} \| \snr\inv \|_2^2\big)^{1/2}}{\sqrt{L}\bar\lambda} + K^3 \log(n) \frac{\| \snr\inv \|_{\infty}^2}{\bar \lambda} \\
    &\quad + K^2 \sqrt{ \log(n)} \frac{\alpha_{\max} \| \snr\inv \|_{\infty}}{\bar \lambda} + \frac{K\alpha_{\max}}{\bar \lambda} \bigg) \\
     &\asymp \frac{K^4 \log(n) \big( \frac{1}{L} \| \snr\inv \|_2^2\big)}{L \sqrt{n} \bar \lambda^2} + \frac{K^5 \log^{3/2}(n) \big( \frac{1}{L} \| \snr\inv \|_2^2\big)^{1/2} \| \snr\inv \|_{\infty}^2}{\sqrt{nL} \bar \lambda^2} \\
     &\quad + \frac{K^4 \log(n) \alpha_{\max}\| \snr\inv \|_{\infty} \big( \frac{1}{L} \| \snr\inv \|_2^2\big)^{1/2}}{\sqrt{nL} \bar \lambda^2} + \frac{K^3 \sqrt{\log(n)}\big( \frac{1}{L} \| \snr\inv \|_2^2\big)^{1/2} \alpha_{\max}}{\sqrt{nL} \bar \lambda^2}.
\end{align*}
We note that we have used the fact that 
\begin{align*}
     K^2 \sqrt{\log(n)} \frac{\big( \frac{1}{L} \|\snr\inv\|_2^2 \big)^{1/2}}{\sqrt{L} \bar \lambda} + K^3 \log(n) \frac{\|\snr\inv\|^2_{\infty}}{\bar \lambda} + K^2 \sqrt{\log(n)} \frac{\alpha_{\max} \|\snr\inv \|_{\infty}}{\bar \lambda} + \frac{K \alpha_{\max}}{\bar \lambda} \lesssim 1,
\end{align*}
as was verified in the proof of \cref{thm:step2sintheta} (observe that this term matches the $\sin\bTheta$ upper bound, and hence is less than one by assumption). Consequently, since each term is the same, we obtain
\begin{align*}
    \| \mathbf{R}_1 \|_{2,\infty} + &\|\mathbf{R}_2 \|_{2,\infty} + \|\mathbf{R}_3\|_{2,\infty}\\
    &\lesssim \frac{K^4 \log(n) \big( \frac{1}{L} \| \snr\inv \|_2^2\big)}{L \sqrt{n} \bar \lambda^2} + \frac{K^5 \log^{3/2}(n) \big( \frac{1}{L} \| \snr\inv \|_2^2\big)^{1/2} \| \snr\inv \|_{\infty}^2}{\sqrt{nL} \bar \lambda^2} \\
     &\quad + \frac{K^4 \log(n) \alpha_{\max}\| \snr\inv \|_{\infty} \big( \frac{1}{L} \| \snr\inv \|_2^2\big)^{1/2}}{\sqrt{nL} \bar \lambda^2} + \frac{K^3 \sqrt{\log(n)}\big( \frac{1}{L} \| \snr\inv \|_2^2\big)^{1/2} \alpha_{\max}}{\sqrt{nL} \bar \lambda^2} \\
    &\asymp \frac{K^4 \log(n) \big( \frac{1}{L} \| \snr\inv \|_2^2\big)}{L \sqrt{n} \bar \lambda^2} + \frac{K^5 \log^{3/2}(n) \big( \frac{1}{L} \| \snr\inv \|_2^2\big)^{1/2} \| \snr\inv \|_{\infty}^2}{\sqrt{nL} \bar \lambda^2} \\
    &\quad + \frac{K^3 \sqrt{ \log(n)} \alpha_{\max} \big( \frac{1}{L} \| \snr\inv \|_2^2\big)^{1/2}}{\sqrt{nL} \bar \lambda^2},
\end{align*}
    where we have used the assumption that that $K \sqrt{\log(n)}\| \snr\inv \|_{\infty} \lesssim 1$, which follows immediately the fact that $\frac{\theta_{\max}\l\|\theta\l\|_1}{\|\theta\l\|^2 \lambda_{\min}\l}\geq 1$ and from \cref{ass:networklevel}, which requires that $K^{8} \log(n) \| \snr\inv\|_{\infty}^2 \frac{\theta_{\max}\l \|\theta\l\|_1}{\|\theta\l\|^2 \lambda_{\min}\l }\lesssim \bar \lambda$.  Therefore, we have shown that
\begin{align*}
    \uhat \wstar\t - \U &= \sum_{l}\mathcal{L}(\mathbf{E}\m) ( \ytilde\m)\t \U \Sigma^{-2}  + \mathcal{R}_{\mathrm{Stage \ II}},
    \end{align*}
    with
    \begin{align*}
    \|\mathcal{R}_{\mathrm{Stage \ II}}\|_{2,\infty} &\lesssim  \frac{K^3 \sqrt{\log(n)}}{n\sqrt{L} \bar \lambda} \big( \frac{1}{L} \| \snr\inv \|_2^2\big)^{1/2} + \frac{K^{7/2} \log(n) \big( \frac{1}{L} \| \snr\inv \|_2^2\big)}{\sqrt{n} L \bar \lambda^2} \\
    &\quad + \frac{K^{7/2} \log(n) \| \snr\inv \|_{\infty}^2}{\sqrt{n} \bar \lambda}  + \frac{K^{5/2} \sqrt{\log(n)} \alpha_{\max}\| \snr\inv \|_{\infty}}{\sqrt{n}\bar \lambda} + \frac{ \alpha_{\max}}{\sqrt{n} \bar \lambda} \\
    &\quad + \frac{K^4 \log(n) \big( \frac{1}{L} \| \snr\inv \|_2^2\big)}{L \sqrt{n} \bar \lambda^2} + \frac{K^5 \log^{3/2}(n) \big( \frac{1}{L} \| \snr\inv \|_2^2\big)^{1/2} \| \snr\inv \|_{\infty}^2}{\sqrt{nL} \bar \lambda^2} \\
    &\quad + \frac{K^3 \sqrt{ \log(n)} \alpha_{\max} \big( \frac{1}{L} \| \snr\inv \|_2^2\big)^{1/2}}{\sqrt{nL} \bar \lambda^2} \\
&\asymp \frac{K^3 \sqrt{\log(n)}}{n L \bar \lambda} \| \snr\inv\|_2 + \frac{K^4 \log(n)}{L^2 \sqrt{n} \bar \lambda^2} \| \snr\inv\|_2^2 + \frac{K^{7/2} \log(n)}{\sqrt{n} \bar \lambda} \| \snr\inv\|_{\infty}^2 \\
&\quad + \frac{K^{5/2 } \sqrt{\log(n)}}{\sqrt{n} \bar \lambda} \alpha_{\max} \|\snr\inv\|_{\infty} + \frac{\alpha_{\max}}{\sqrt{n} \bar \lambda} \\
&\quad + \frac{K^5 \log^{3/2}(n)}{L \sqrt{n} \bar \lambda^2} \|\snr\inv\|_2 \|\snr\inv\|_{\infty}^2 + \frac{K^3 \sqrt{\log(n)}}{L \sqrt{n} \bar \lambda^2} \alpha_{\max} \|\snr\inv\|_2^2 \\
&= \frac{K^3 \sqrt{\log(n)}}{nL \bar \lambda} \|\snr\inv\|_2 + \frac{K^4 \log(n)}{L^2 \sqrt{n} \bar \lambda^2} \|\snr\inv\|_2^2 \\
&\quad + \frac{K^{7/2} \log(n)}{\sqrt{n} \bar \lambda} \| \snr\inv\|_{\infty}^2 \bigg( 1 + \frac{K^{3/2} \sqrt{\log(n)}}{L \bar \lambda} \|\snr\inv\|_2 \bigg) \\
&\quad + \frac{\alpha_{\max}}{\sqrt{n} \bar \lambda} \bigg( 1 + K^{5/2} \sqrt{\log(n)} \|\snr\inv\|_{\infty} + \frac{K^3 \sqrt{\log(n)}}{L \bar \lambda} \|\snr\inv\|_2^2 \bigg) \\
&\lesssim \frac{K^3 \sqrt{\log(n)}}{nL \bar \lambda} \|\snr\inv\|_2 + \frac{K^4 \log(n)}{L^2 \sqrt{n} \bar \lambda^2} \|\snr\inv\|_2^2 \\
&\quad + \frac{K^{7/2} \log(n)}{\sqrt{n} \bar \lambda} \| \snr\inv\|_{\infty}^2 + \frac{\alpha_{\max}}{\sqrt{n} \bar \lambda},
\end{align*}
where the final inequality holds as long as
\begin{align}
    \frac{K^{3/2} \sqrt{\log(n)}}{L \bar \lambda} \|\snr\inv\|_2 &\lesssim 1; \label{verifybb}\\
    K^{5/2} \|\snr\inv\|_{\infty} &\lesssim 1 \label{verifybb1} \\
    \frac{K^3 \sqrt{\log(n)}}{L \bar \lambda} \|\snr\inv\|_2^2 &\lesssim 1. \label{verifybb2}
\end{align}
We will verify these bounds now.  First, \cref{ass:networklevel} implies that
\begin{align*}
    \frac{K^8 \log(n) \theta_{\max}\l \|\theta\l\|_1}{\|\theta\l\|_2^2} (\snr_l\inv)^2 \lesssim \bar \lambda \lambda_{\min}\l,
\end{align*}
as long as $C$ in the assumption is sufficiently large. 
Observe that this immediately implies equation \eqref{verifybb1} since $\frac{\theta_{\max}\l \|\theta\l\|_1}{\|\theta\l\|_2^2} \geq 1$ and $\lambda_{\min}\l \in (0,1)$ by assumption.  For the other two terms, by averaging this condition over $l$, we see that \cref{ass:networklevel} implies
\begin{align*}
     \frac{K^8 \log(n)}{L} \|\snr\inv\|_2^2 \lesssim \bar \lambda^2. \numberthis \label{hamlet}
\end{align*}
This implies \eqref{verifybb} and \eqref{verifybb2}.
Hence, we have shown so far that 
\begin{align*}
  \|\mathcal{R}_{\mathrm{Stage \ II}}\|_{2,\infty} &\lesssim   \frac{K^3 \sqrt{\log(n)}}{nL \bar \lambda} \|\snr\inv\|_2 + \frac{K^4 \log(n)}{L^2 \sqrt{n} \bar \lambda^2} \|\snr\inv\|_2^2 \\
&\quad + \frac{K^{7/2} \log(n)}{\sqrt{n} \bar \lambda} \| \snr\inv\|_{\infty}^2 + \frac{\alpha_{\max}}{\sqrt{n} \bar \lambda}.
\end{align*}
This holds cumulatively with probability at least $1 - O(n^{-10})$.  We now verify that the sum of these terms is  less than $\frac{1}{16\sqrt{n_{\max}}}$.  Since $n_{\max} \leq  n$, it suffices to show that this upper bound is at most $\frac{1}{16\sqrt{n}}$.  By pulling out a factor of  $1/\sqrt{n}$  it suffices to show that 
\begin{align*}
    \frac{K^3 \sqrt{\log(n)}}{\sqrt{n} L \bar \lambda} \|\snr\inv\|_2 + \frac{K^4 \log(n)}{L^2 \bar \lambda^2} \|\snr\inv\|_2^2 + \frac{K^{7/2} \log(n)}{ \bar \lambda} \| \snr\inv\|_{\infty}^2 + \frac{\alpha_{\max}}{ \bar \lambda} \lesssim 1.  
\end{align*}
By similar manipulations as in verifying the bounds \eqref{verifybb}, \eqref{verifybb1}, and \eqref{verifybb2}, it is straightforward to check the condition above holds, except for the condition $\frac{\alpha_{\max}}{\bar \lambda} \lesssim 1$.  Plugging in the definition for $\alpha_{\max}$, we see that we require
\begin{align*}
    \frac{1}{\bar \lambda} \max_{l} \frac{K^2 \theta_{\max}\l \|\theta\l\|_1}{\lambda_{\min}\l \|\theta\l\|^4} \bigg( \log(n) \frac{\theta_{\max}\l}{\theta_{\min}\l} + \frac{{\sqrt{K}}}{\lambda_{\min}\l} + \bigg( \frac{\theta_{\max}\l}{\theta_{\min}\l}\bigg)^{1/2} \frac{K^{5/2} \log(n)}{(\lambda_{\min}\l)^{1/2}} \bigg) \lesssim 1.
\end{align*}
This is covered by \cref{ass:networklevel}.  Therefore, this completes the proof.
\end{proof}

\section{Proof of Extension to Different Network Setting (Theorem \ref{thm:newtheorem})}
Let $\tilde{\ycal}$ be the same as the matrix $\ycal$, except constructed using the ``true'' membership matrix $\bZ$.  First we study the spectral structure of $\tilde{\ycal}$ and $\ycal$, yielding a $\sin\Theta$ bound between their respective eigenvectors.  We then use this result together with a deterministic $\ell_{2,\infty}$ bound from \citet{cape_two--infinity_2019} bound the $\ell_{2,\infty}$ difference.  Combining these results we provide a modified proof of \cref{thm:clusteringerror}. 
\\ \ \\
\noindent \textbf{Step 1: Spectral Structure of $\ycal\ycal\t$}.  First, we note that \cref{lem:step2popprop} applies to $\tilde{\ycal}\tilde{\ycal}\t$, and hence it holds that
\begin{align*}
    \lambda^2_Y \gtrsim \frac{n}{K} L \bar \lambda.
\end{align*}
Furthermore, since the first phase population analysis continues to hold without modification for $\ycal\ycal\t$, it holds that $\ycal\ycal\t = \sum_l \bZ^{(l)} \bM^{(l)} (\bZ^{(l)})\t$ for some (positive definite) matrices $\bM^{(l)}$.  The proof of this same lemma reveals further that $\max_l \|\bM^{(l)} \| \lesssim 1$.  We also have that $\tilde{\ycal} \tilde{\ycal} = \bZ \big( \sum_l \bM^{(l)} \big) \bZ\t$.  Therefore, it holds that
\begin{align*}
    \| \tilde{\ycal}\tilde{\ycal}\t - \ycal \ycal\t \| &= \bigg\| \sum_l \bZ^{(l)} \bM^{(l)} (\bZ^{(l)})\t - \bZ \big( \sum_l \bM^{(l)} \big) \bZ\t \bigg\| \\
    &\leq \bigg\| \sum_l \big[ \bZ^{(l)} - \bZ ] \bM^{(l)} (\bZ^{(l)})\t \bigg\| + \bigg\| \sum_l \big[ \bZ^{(l)} - \bZ ] \bM^{(l)} \bZ\t \bigg\| \\
    &\lesssim \sqrt{n} \max_{i} \bigg\| e_i\t \sum_l \big[ \bZ^{(l)} - \bZ\big] \bM^{(l)} (\bZ^{(l)})\t \bigg\| + n \max_i \bigg\| e_i\t \sum_l \big[ \bZ^{(l)} - \bZ\big] \bM^{(l)} \bigg\|\\
     &\lesssim n \max_{i}\sum_l \bigg\| e_i\t  \big[ \bZ^{(l)} - \bZ\big]  \bigg\| + n \max_i\sum_l \bigg\| e_i\t  \big[ \bZ^{(l)} - \bZ\big] \bigg\|\\
    &\lesssim n L \delta. 
\end{align*}
Therefore, under the assumption  $\delta  \ll \frac{\bar \lambda}{K}$, letting $\tilde \lambda^2_{Y}$ denote the $K$-th eigenvalue of $\tilde{\ycal}\tilde{\ycal}\t$, Weyl's inequality implies that $\tilde \lambda^2_{Y} \gtrsim \frac{n}{K} L \bar \lambda$.  Therefore, the Davis-Kahan Theorem implies that
\begin{align*}
    \| \sin\Theta(\bU, \tilde{\bU}) \| &\lesssim K\frac{\delta}{\bar \lambda}. \numberthis \label{shakespeare}
\end{align*}
Note that under our assumptions the quantity on the right hand side above is $o(1)$.  
\\ \ \\
\noindent
\textbf{Step 2: Bounding The $\ell_{2,\infty}$ difference}.  We apply Theorem 3.7 of \citet{cape_two--infinity_2019} to reveal that
\begin{align*}
    \| \tilde{\bU} - \bU \mathbf{W}_U \|_{2,\infty} &\lesssim \underbrace{\frac{ \|(\mathbf{I} - \tilde{\U} \tilde{\U}\t) ( \tilde{\ycal}\tilde{\ycal}\t - \ycal \ycal\t ) \bU \bU\t \|_{2,\infty}}{\lambda^2_Y}}_{\alpha_1} \\
    &\quad + \underbrace{\frac{ \|(\mathbf{I} - \tilde{\U} \tilde{\U}\t) ( \tilde{\ycal}\tilde{\ycal}\t - \ycal \ycal\t ) (\mathbf{I} - \tilde{\U} \tilde{\U}\t )\|_{2,\infty}}{\lambda^2_Y} \| \sin\Theta(\U,\tilde{\bU})\|}_{\alpha_2} \\
    &\quad + \underbrace{\| \sin\Theta(\U,\tilde{\bU})\|^2 \| \bU \|_{2,\infty}}_{\alpha_3}.
\end{align*}
We will bound each term above separately.  However, before doing so we observe that
\begin{align*}
     \| \tilde{\ycal}\tilde{\ycal}\t - \ycal \ycal\t \|_{2,\infty} &\leq \max_{1\leq i \leq n}  \bigg\|e_i\t    \sum_l \big[ \bZ^{(l)} - \bZ ] \bM^{(l)} (\bZ^{(l)})\t \bigg\| + \bigg\| e_i\t \bZ \sum_l   \bM^{(l)}\big[ \bZ^{(l)} - \bZ ]  \t \bigg\|. 
   \end{align*}
   For a fixed $i$ it holds that
   \begin{align*}
         \sum_l \| e_i\t& \big( \bZ^{(l)} - \bZ \big) \| \| \bM^{(l)} \| \| \bZ^{(l)} \| + \| \bZ \|_{\infty,\infty}  \max_i \bigg\|  e_i\t \sum_l \bM^{(l)}[\bZ^{(l)} - \bZ]\t  \bigg\|     \\
         &\lesssim L\sqrt{n}  \delta +       \max_i \bigg\|  e_i\t \sum_l \bM^{(l)}[\bZ^{(l)} - \bZ]\t  \bigg\| \\   
         &\lesssim L \sqrt{n} \delta + \max_l \| \bM^{(l)} \|_{\max} \max_i \sum_l \bigg\| e_i\t \big[ \bZ^{(l)} - \bZ] \bigg\| \\
         &\lesssim L \sqrt{n} \delta. \numberthis \label{debussy}
\end{align*}
We now bound $\alpha_i$ in turn.
\begin{itemize}
    \item \textbf{The term $\alpha_1$}:  By \cref{debussy},
    \begin{align*}
         \|(\mathbf{I} - \tilde{\U} \tilde{\U}\t) ( \tilde{\ycal}\tilde{\ycal}\t - \ycal \ycal\t ) \tilde{\U} \tilde{\U}\t \|_{2,\infty} &\leq \|(\mathbf{I} - \tilde{\U} \tilde{\U}\t) \|_{\infty,\infty} \|\tilde{\ycal}\tilde{\ycal}\t - \ycal \ycal\t  \|_{2,\infty} \\
         &\leq \sqrt{K} L \sqrt{n} \delta,          
    \end{align*}
    where $\|\cdot\|_{\infty,\infty}$ denotes the $\ell_{\infty}$ operator norm on matrices, and the bound $\|(\mathbf{I} - \tilde{\U} \tilde{\U}\t) \|_{\infty,\infty} \lesssim \sqrt{K}$ comes from the fact that $\|\tilde{\U}\|_{2,\infty} \lesssim \sqrt{\frac{K}{n}}$.  As a consequence, recalling that $\lambda^2_Y \gtrsim \frac{n}{K} L \bar \lambda$, we have that
    \begin{align*}
        \alpha_1 &\lesssim \frac{\sqrt{K} L \sqrt{n} \delta}{\frac{n}{K} L \bar \lambda} \asymp \frac{K^{3/2}\delta}{\sqrt{n} \bar \lambda}. \numberthis \label{alpha1bound}
    \end{align*}
    \item \textbf{The term $\alpha_2$}:  By a similar argument as above,
    \begin{align*}
        \alpha_2 &\lesssim \| \sin\Theta(\U,\tilde{\bU})\|  \frac{K^{3/2}\delta}{\sqrt{n} \bar \lambda} \lesssim K^{5/2} \frac{\delta^2}{\bar \lambda^2 \sqrt{n}}, \numberthis \label{alpha2bound}
    \end{align*}
    where the final inequality follows from \cref{shakespeare}.
    \item \textbf{The term $\alpha_3$}:  By \cref{shakespeare}, we have that
    \begin{align*}
        \alpha_3 &\lesssim K \frac{\delta}{\bar \lambda} \| \tilde{\U} \|_{2,\infty} \lesssim K^{3/2} \frac{\delta}{\bar \lambda \sqrt{n}}. \numberthis\label{alpha3bound}
    \end{align*}
\end{itemize}
Therefore, combining \cref{alpha1bound,alpha2bound,alpha3bound}, we have that
\begin{align*}
    \| \bU - \tilde{\bU} \mathbf{W}_U \|_{2,\infty} &\lesssim \frac{K^{3/2}\delta}{\sqrt{n} \bar \lambda} + K^{5/2} \frac{\delta^2}{\bar \lambda^2 \sqrt{n}} + K^{3/2} \frac{\delta}{\bar \lambda \sqrt{n}} \\
    &\lesssim \frac{K^{3/2}\delta}{\sqrt{n} \bar \lambda},
\end{align*}
where the final inequality follows from the assumption $\delta \ll \frac{\bar \lambda}{K}$.  Note that this result also implies that
\begin{align*}
    \|\bU \|_{2,\infty} &\lesssim \sqrt{\frac{K}{n}}. \numberthis \label{Uisincoherent}
\end{align*}
\\ \ \\
\textbf{Step 3: Modifying the proof of \cref{thm:clusteringerror}}:
First, the proof of \cref{thm:firststep} goes through without modification since the proof only relies on the fact that each probability matrix is rank $K$, which continues to hold.  The proof of \cref{thm:step2sintheta} also continues to hold without modification as it relies only on eigengap assumptions, which hold from the argument in Step 1 of this proof. As a result, the proof of \cref{thm:step2asympexp} nearly holds, except one now has the asymptotic expansion
\begin{align*}
    \uhat  - \bU\mathbf{W}_* = \sum_l \mathcal{L}(\bA^{(l)} - \bP^{(l)}) ( \mathbf{Y}^{(l)})\t \bU \Sigma^{-2}\mathbf{W}_* \t + ( \mathbf{I} - \bU \bU\t ) \mathcal{Y} \mathcal{L}(\mathcal{E})\t \uhat \hat \Sigma^{-2} + \mathcal{R}_{\mathrm{Stage \ II}},
\end{align*}
where $\mathcal{R}_{\mathrm{Stage \ II}}$ continues to have the same upper bound as before.  In contrast to the case where all the $\bZ$'s are the same, the additional quantity $ ( \mathbf{I} - \bU \bU\t ) \mathcal{Y} \mathcal{L}(\mathcal{E})\t \uhat \hat \Sigma^{-2}$ does not vanish as $\mathcal{Y}$ is no longer exactly rank $K$, but instead only approximately so. The following lemma controls this additional term.

\begin{lemma} \label{lem:lastlemma}
    Under the conditions of \cref{thm:newtheorem}, it holds that
    \begin{align*}
        \| (\mathbf{I} - \U \U\t) \mathcal{Y} \mathcal{L}(\mathcal{E})\t \uhat \hat \Sigma^{-2} \|_{2,\infty} \leq c_0\sqrt{\frac{K}{n}},
        \end{align*}
        where $c_0$ is some sufficiently small constant.  
\end{lemma}
Explicitly, these results in tandem imply that
\begin{align*}
    \uhat - \mathbf{\tilde U} \mathbf{W}_U \mathbf{W}_* &= \sum_l \mathcal{L}(\bA^{(l)} - \bP^{(l)}) ( \mathbf{Y}^{(l)})\t \bU \Sigma^{-2}\mathbf{W}_* \t  + O\bigg( \frac{K^{3/2}\delta}{\sqrt{n}\bar\lambda}\bigg) + c \sqrt{\frac{K}{n}} \\
    &= \sum_l \mathcal{L}(\bA^{(l)} - \bP^{(l)}) ( \mathbf{Y}^{(l)})\t \bU \Sigma^{-2}\mathbf{W}_* \t  +  c \sqrt{\frac{K}{n}},
\end{align*}
where the constant $c$ is sufficiently small, where we have used the assumption $\delta \ll \frac{\bar \lambda}{K}$.  We may therefore modify the proof of \cref{thm:clusteringerror}.  First, we still have the bound $\|\sin\Theta(\uhat,\tilde{\U})\| \leq \frac{\beta}{8 K \sqrt{C_{\eps}}}$ on the event $\mathcal{E}_{\sin\Theta}$.  Replacing the appearances of $\U$ with $\tilde{\bU}$ we see that step two of the argument remains valid without further modification, which relies only on the $\sin\Theta$ bound.  Step three also remains the same by appealing to the fact that $\| \U \|_{2,\infty} \lesssim \sqrt{\frac{K}{n}}$ by \eqref{Uisincoherent} (which is a deterministic bound).  Therefore, with these modifications, the conclusion of \cref{thm:clusteringerror} continues to hold under the assumptions of \cref{thm:newtheorem}, which completes the proof.

\subsection{Proof of Lemma \ref{lem:lastlemma}}
Before proving this lemma, we introduce the following lemma establishing concentration inequalities for two terms that appear in the analysis.
\begin{lemma} \label{lem:secondtolastlemma}
With probability at least $1 - O(n^{-10})$ it holds that
    \begin{align*}
        \| (\mathbf{I} - \U \U\t) \mathcal{Y} \mathcal{L}(\mathcal{E})\t \U \|_{2,\infty} &\lesssim   \frac{K^{2}}{\sqrt{n}} \max_l \frac{\|\theta\l\|_1}{\|\theta\l\|^2 (\lambda_{\min}\l)^{1/2}} \sqrt{L\log(n)}; \\
          \| (\mathbf{I} - \U \U\t) \mathcal{Y} \mathcal{L}(\mathcal{E})\t \|_{2,\infty} &\lesssim K n \sqrt{L \log(n)} \bigg( \frac{1}{L} \|\snr\inv \|_2^2 \big)^{1/2}.
    \end{align*}
\end{lemma}

\begin{proof}[Proof of \cref{lem:secondtolastlemma}]
The proof is similar to the lemmas in \cref{sec:secondstageasympexp}, where we write everything as a sum of independent random variables and apply Bernstein's inequality.  First, fix a row $i$ of the matrix $(\mathbf{I} - \U \U\t) \mathcal{Y} \mathcal{L}(\mathcal{E})\t \U$.  Observe that
\begin{align*}
    \bigg\| e_i\t (\mathbf{I} - \U \U\t) \mathcal{Y} \mathcal{L}(\mathcal{E})\t \U \bigg\| &\leq \sqrt{K} \max_{1\leq \nu \leq K} \bigg| e_i\t (\mathbf{I} - \U \U\t) \mathcal{Y} \mathcal{L}(\mathcal{E})\t \U e_{\nu} \bigg|. \numberthis \label{kfactor}
\end{align*}
Let $\bY^{(l)}_{\perp}$ denote the $l$'th matrix $(\mathbf{I} - \U\U\t)\bY^{(l)}$.  Then the $i,\nu$ entry above can be written via
\begin{align*}
   e_{i}\t  \sum_l &(\bY^{(l)}_{\perp}) \mathcal{L}\big( \bA\l - \bP\l \big)\t \bU e_{\nu}\\
   &=     \sum_l \sum_{k=1}^{K} \sum_{c = 1}^{n} (\bY^{(l)}_{\perp})_{i k} \mathcal{L}\big( \bA\l - \bP\l \big)_{ck} \bU_{c \nu} \\
   &= \sum_l \sum_{k=1}^{K} \sum_{c = 1}^{n} \sum_{f=1}^{n} (\bY^{(l)}_{\perp})_{i k} \bigg( \bA\l - \bP\l \bigg)_{cf} \bigg(  \bU^{(l)} |\Lambda^{(l)}|^{-1/2} \ipq^{(l)} \mathbf{J}( \mathbf{X}_{f \cdot} \big)  \bigg)_{fk} \bU_{c \nu} \\
   &= \sum_l\sum_{c = 1}^{n} \sum_{f=1}^{n}  \bigg( \bA\l - \bP\l \bigg)_{cf} \Bigg[ \sum_{k=1}^{K}  (\bY^{(l)}_{\perp})_{i k}  \bigg(  \bU^{(l)} |\Lambda^{(l)}|^{-1/2} \ipq^{(l)} \mathbf{J}( \mathbf{X}_{f \cdot} \big)  \bigg)_{fk} \bU_{c \nu} \Bigg]\\
   &= \sum_{l} \sum_{c=1}^{n} \sum_{f\leq c,f =1}^{n}  \bigg( \bA\l - \bP\l \bigg)_{cf} \Bigg[ \sum_{k=1}^{K}  (\bY^{(l)}_{\perp})_{i k}  \bigg(  \bU^{(l)} |\Lambda^{(l)}|^{-1/2} \ipq^{(l)} \mathbf{J}( \mathbf{X}_{f \cdot} \big)  \bigg)_{fk} \bU_{c \nu} \Bigg] \\
   &\quad + \sum_{l} \sum_{c=1}^{n} \sum_{f=c+1}^{n}  \bigg( \bA\l - \bP\l \bigg)_{cf} \Bigg[ \sum_{k=1}^{K}  (\bY^{(l)}_{\perp})_{i k}  \bigg(  \bU^{(l)} |\Lambda^{(l)}|^{-1/2} \ipq^{(l)} \mathbf{J}( \mathbf{X}_{f \cdot} \big)  \bigg)_{fk} \bU_{c \nu} \Bigg].
\end{align*}
Observe that each term above is a sum over the independent random variables $\{(\bA\l - \bP\l)_{cf}\}_{l=1}^{L}$.  We bound the first term as the second term is similar.  Bernstein's inequality shows us that for fixed indices $i$ and $\nu$ it holds that with probability at least $1 - O(n^{-20})$,
\begin{align*}
    \sum_{l} \sum_{c=1}^{n} \sum_{f\leq c,f =1}^{n}  \bigg( \bA\l - \bP\l \bigg)_{cf} \Bigg[ \sum_{k=1}^{K}  (\bY^{(l)}_{\perp})_{i k}  \bigg(  \bU^{(l)} |\Lambda^{(l)}|^{-1/2} \ipq^{(l)} \mathbf{J}( \mathbf{X}_{f \cdot} \big)  \bigg)_{fk} \bU_{c \nu} \Bigg] &\lesssim \sqrt{v\log(n)} + w \log(n),
\end{align*}
where $v$ is the sum of the variances and $w$ is a bound on the maximum value for all $c,f$ and $l$.  We therefore bound directly, observing that by \cref{lem:popprop},
\begin{align*}
    v &\leq \sum_l \sum_{c=1}^{n} \sum_{f\leq c f=1}^{n} \theta\l_{c} \theta\l_f \Bigg[ \sum_{k=1}^{K}  (\bY^{(l)}_{\perp})_{i k}  \bigg(  \bU^{(l)} |\Lambda^{(l)}|^{-1/2} \ipq^{(l)} \mathbf{J}( \mathbf{X}_{f \cdot} \big)  \bigg)_{fk} \bU_{c \nu} \Bigg]^2 \\
    &\leq  \sum_l \sum_{c=1}^{n} \sum_{f\leq c f=1}^{n} \theta\l_{c} \theta\l_f \bU_{c \nu}^2  \Bigg[ \sum_{k=1}^{K}  (\bY^{(l)}_{\perp})_{i k}  \bigg(  \bU^{(l)} |\Lambda^{(l)}|^{-1/2} \ipq^{(l)} \mathbf{J}( \mathbf{X}_{f \cdot} \big)  \bigg)_{fk}  \Bigg]^2 \\
    &\leq \sum_l \sum_{c=1}^{n} \sum_{f\leq c f=1}^{n} \theta\l_{c} \theta\l_f \bU_{c \nu}^2  \|  (\bY^{(l)}_{\perp})_{i\cdot} \|^2  \Bigg\|    \bigg(  \bU^{(l)} |\Lambda^{(l)}|^{-1/2} \ipq^{(l)} \mathbf{J}( \mathbf{X}_{f \cdot} \big)  \bigg)_{f\cdot}  \Bigg\|^2 \\
    &\lesssim \sum_l \sum_{c=1}^{n} \sum_{f=1}^{n}\theta\l_{c} \theta\l_f \bU_{c \nu}^2 \|  (\bY^{(l)}_{\perp})_{i\cdot} \|^2   \frac{K (\theta_f\l)^2}{\|\theta\l\|^2} \frac{K}{\|\theta\l\|^2 \lambda_{\min}\l } \frac{1}{(\theta_f\l)^2} \\
    &\lesssim K^2 \|  (\bY^{(l)}_{\perp})_{i\cdot} \|^2 \sum_l \sum_{c=1}^{n} \sum_{f=1}^{n}\theta\l_{c} \theta\l_f \bU_{c \nu}^2    \frac{ 1}{\|\theta\l\|^4\lambda_{\min}\l } \\
    &\lesssim \frac{K^3}{n} \|  (\bY^{(l)}_{\perp})_{i\cdot} \|^2 L \max_l  \frac{ 1}{\|\theta\l\|^4\lambda_{\min}\l }\sum_{c=1}^{n} \sum_{f=1}^{n}\theta\l_{c} \theta\l_f     \\
    &\lesssim \frac{K^3}{n} \|  (\bY^{(l)}_{\perp})_{i\cdot} \|^2 L \max_l  \frac{ \| \theta\l \|_1^2}{\|\theta\l\|^4\lambda_{\min}\l } \\
    &\lesssim \frac{K^3}{n}  L \max_l  \frac{ \| \theta\l \|_1^2}{\|\theta\l\|^4\lambda_{\min}\l } \numberthis \label{variancebound1}
\end{align*}
where we have used the inequality \eqref{Uisincoherent} to bound $ \U_{c\nu}^2\lesssim \frac{K}{n},$ and the final inequality follows from the fact that $\|(\bY^{(l)}_{\perp})_{i\cdot} \| \leq 1$.  The maximum over $l,c,$ and $f$ is bounded by
\begin{align*}
    \max_{l,c,f} \Bigg| \sum_{k=1}^{K}  (\bY^{(l)}_{\perp})_{i k}  \bigg(  \bU^{(l)} |\Lambda^{(l)}|^{-1/2} \ipq^{(l)} \mathbf{J}( \mathbf{X}_{f \cdot} \big)  \bigg)_{fk} \bU_{c \nu} \Bigg| &\lesssim   \max_{l,f}  \| \bU \|_{2,\infty} \|(\bY^{(l)}_{\perp})_{i\cdot} \| \bigg\|  \bigg(  \bU^{(l)} |\Lambda^{(l)}|^{-1/2} \ipq^{(l)} \mathbf{J}( \mathbf{X}_{f \cdot} \big)  \bigg)_{f\cdot} \bigg\|  \\
    &\lesssim \sqrt{\frac{K}{n}} \frac{K}{\|\theta\l\|^2 (\lambda_{\min}\l)^{1/2}} \\
    &\asymp \frac{K^{3/2}}{\sqrt{n}} \max_l \frac{1}{\|\theta\l\|^2(\lambda_{\min}\l)^{1/2}}\numberthis \label{secondguy}
\end{align*} 
Therefore, combining \eqref{variancebound1} and \eqref{secondguy}, it holds that with probability at least $1 - O(n^{-20})$,
\begin{align*}
 \Bigg|    \sum_{l} &\sum_{c=1}^{n} \sum_{f\leq c,f =1}^{n}  \bigg( \bA\l - \bP\l \bigg)_{cf} \Bigg[ \sum_{k=1}^{K}  (\bY^{(l)}_{\perp})_{i k}  \bigg(  \bU^{(l)} |\Lambda^{(l)}|^{-1/2} \ipq^{(l)} \mathbf{J}( \mathbf{X}_{f \cdot} \big)  \bigg)_{fk} \bU_{c \nu} \Bigg] \Bigg|  \\
 &\lesssim    \sqrt{\frac{K^3}{n}  L \max_l  \frac{ \| \theta\l \|_1^2}{\|\theta\l\|^4\lambda_{\min}\l }\log(n)} + \frac{K^{3/2}}{\sqrt{n}} \max_l \frac{1}{\|\theta\l\|^2(\lambda_{\min}\l)^{1/2}} \log(n) \\
 &\lesssim \frac{K^{3/2}}{\sqrt{n}} \max_l \frac{\|\theta\l\|_1}{\|\theta\l\|^2 (\lambda_{\min}\l)^{1/2}} \sqrt{L\log(n)},
\end{align*}
where the final inequality follows from the assumption that $\|\theta\l\|_1 \gtrsim c \log(n)$.  Therefore, combining this bound with the bound \eqref{kfactor}, we obtain the desired result.  

We now bound the second quantity, though it is significantly easier due to the previous arguments.  By \cref{lem:approximatecommutation}, with probability at least $1 - O(n^{-10})$,
\begin{align*}
     \bigg\| e_i\t (\mathbf{I} - \U \U\t) \mathcal{Y} \mathcal{L}(\mathcal{E})\t \bigg\| &\leq \sqrt{n} \max_{1 \leq \nu \leq n } \bigg |e_i\t (\mathbf{I} - \U \U\t) \mathcal{Y} \mathcal{L}(\mathcal{E})\t e_{\nu} \bigg| \\
     &\leq \sqrt{n} \|  \mathcal{L}(\mathcal{E}) \mathcal{Y}\t (\mathbf{I} - \U \U\t ) \|_{2,\infty}  \\
     &\lesssim K n \sqrt{L \log(n)} \bigg( \frac{1}{L} \|\snr\inv \|_2^2 \big)^{1/2}.
\end{align*}
This completes the proof.
\end{proof}

\begin{proof}[Proof of \cref{lem:lastlemma}]
First, we decompose via:
\begin{align*}
     (\mathbf{I} - \U \U\t) \mathcal{Y} \mathcal{L}(\mathcal{E})\t \uhat \hat \Sigma^{-2} &=  \underbrace{(\mathbf{I} - \U \U\t) \mathcal{Y} \mathcal{L}(\mathcal{E})\t \bU \Sigma^{-2} \wstar}_{\alpha_1} \\
     &\quad + \underbrace{(\mathbf{I} - \U \U\t) \mathcal{Y} \mathcal{L}(\mathcal{E})\t \bU \Sigma^{-2} ( \wstar - \U\t \uhat )}_{\alpha_2} \\
     &\quad + \underbrace{(\mathbf{I} - \U \U\t) \mathcal{Y} \mathcal{L}(\mathcal{E})\t \bU (\Sigma^{-2} \U\t\uhat - \U\t \uhat \hat \Sigma^{-2})}_{\alpha_3} \\
     &\quad + \underbrace{(\mathbf{I} - \U \U\t) \mathcal{Y} \mathcal{L}(\mathcal{E})\t (\uhat - \U \U\t \uhat) \hat \Sigma^{-2}}_{\alpha_4}.
\end{align*}
We will bound each sequentially.  Before moving on we note that by \cref{lem:step2approximatecommutation} it holds that
\begin{align*}
    \| \U\t \uhat -  \wstar  \| &\lesssim \bigg( K^2 \sqrt{\log(n)} \frac{\big( \frac{1}{L} \| \snr\inv \|_2^2\big)^{1/2}}{\sqrt{L}\bar\lambda} + K^3 \log(n) \frac{\| \snr\inv \|_{\infty}^2}{\bar \lambda} \\
    &\qquad \qquad + K^2 \sqrt{ \log(n)} \frac{\alpha_{\max} \| \snr\inv \|_{\infty}}{\bar \lambda} + \frac{K\alpha_{\max}}{\bar \lambda} \bigg)^2 \\
    :&= \mathcal{A}^2 \numberthis \label{romeo} \\
    \|  \Sigma^{-2} \U\t \uhat - \U\t \uhat \hat \Sigma^{-2} \| &\lesssim\frac{K^3 \sqrt{\log(n)} \big( \frac{1}{L} \| \snr\inv \|_2^2\big)^{1/2}}{n L^{3/2} \bar \lambda^2} + \frac{K^4 \log(n) \| \snr\inv \|_{\infty}^2}{n L \bar \lambda^2} \\
    &\qquad + \frac{K^3 \sqrt{\log(n)} \alpha_{\max} \| \snr\inv \|_{\infty}}{n L \bar \lambda^2} + \frac{K^2\alpha_{\max}}{n L \bar \lambda^2}; \\ 
    :&= \mathcal{B} \numberthis \label{juliet} \\
    \| \sin\Theta(\uhat,\U) \| &\lesssim  K^2 \sqrt{\log(n)} \frac{\big( \frac{1}{L} \| \snr\inv \|_2^2\big)^{1/2}}{\sqrt{L}\bar\lambda} + K^3 \log(n) \frac{\| \snr\inv \|_{\infty}^2}{\bar \lambda}  \\
    &\quad + K^2 \sqrt{ \log(n)} \frac{\alpha_{\max} \| \snr\inv \|_{\infty}}{\bar \lambda} + \frac{K \alpha_{\max}}{\bar \lambda} \\
    &= \mathcal{A} \numberthis \label{macbeth}.
\end{align*}
We now bound $\alpha_i$.  
\begin{itemize}
    \item \textbf{The term $\alpha_1$}: 
    By \cref{lem:secondtolastlemma}, with probability at least $1 - O(n^{-10})$,
    \begin{align*}
    \| \alpha_1 \|_{2,\infty} &= \|(\mathbf{I} - \U \U\t) \mathcal{Y} \mathcal{L}(\mathcal{E})\t \bU \Sigma^{-2} \wstar  \|_{2,\infty} \\
    &\leq \|(\mathbf{I} - \U \U\t) \mathcal{Y} \mathcal{L}(\mathcal{E})\t \bU  \|_{2,\infty} \| \Sigma^{-2} \| \\
    &\lesssim \frac{K^{2}}{\sqrt{n}} \max_l \frac{\|\theta\l\|_1}{\|\theta\l\|^2 (\lambda_{\min}\l)^{1/2}} \sqrt{L\log(n)} \| \Sigma^{-2} \| \\
    &\lesssim \frac{K^{2}}{\sqrt{n}} \max_l \frac{\|\theta\l\|_1}{\|\theta\l\|^2 (\lambda_{\min}\l)^{1/2}} \sqrt{L\log(n)} \frac{K}{nL \bar \lambda} \\
    &\asymp \frac{K^3 \sqrt{\log(n)}}{n^{3/2} \sqrt{L} \bar \lambda} \max_{l} \frac{\|\theta\l\|}{\|\theta\l\|^2 (\lambda_{\min}\l)^{1/2}}. \numberthis \label{radiohead}
       \end{align*}
    \item \textbf{The term $\alpha_2$}:  Again by \cref{lem:secondtolastlemma}, 
    \begin{align*}
          \|\alpha_2\|_{2,\infty} &= \| (\mathbf{I} - \U \U\t) \mathcal{Y} \mathcal{L}(\mathcal{E})\t \bU \Sigma^{-2} ( \wstar - \U\t \uhat )\|_{2,\infty} \\
          &\leq  \| (\mathbf{I} - \U \U\t) \mathcal{Y} \mathcal{L}(\mathcal{E})\t \bU  \|_{2,\infty} \|\Sigma^{-2}\| \| ( \wstar - \U\t \uhat )\| \\
          &\lesssim \| (\mathbf{I} - \U \U\t) \mathcal{Y} \mathcal{L}(\mathcal{E})\t \bU  \|_{2,\infty} \|\Sigma^{-2}\| \\
          &\lesssim \frac{K^3\sqrt{\log(n)}}{n^{3/2} \sqrt{L} \bar \lambda} \max_{l} \frac{\|\theta\l\|}{\|\theta\l\|^2 (\lambda_{\min}\l)^{1/2}}. \numberthis \label{fakeplastictrees} 
       \end{align*}
       which satisfies the same bound as $\alpha_1$.  Here we have noted that $\|\wstar - \U\t \uhat \| \lesssim \mathcal{A}^2 \lesssim \frac{1}{K^2}$ by \eqref{macbeth} and the second conclusion of  \cref{thm:step2sintheta}.
    \item \textbf{The term $\alpha_3$}: By \eqref{lem:secondtolastlemma} and \eqref{juliet}, with high probability we have that
    \begin{align*}
         \|\alpha_3\|_{2,\infty} &= \| (\mathbf{I} - \U \U\t) \mathcal{Y} \mathcal{L}(\mathcal{E})\t \bU (\Sigma^{-2} \U\t\uhat - \U\t \uhat \hat \Sigma^{-2}) \|_{2,\infty} \\
         &\leq \| (\mathbf{I} - \U \U\t) \mathcal{Y} \mathcal{L}(\mathcal{E})\t \bU  \|_{2,\infty} \| (\Sigma^{-2} \U\t\uhat - \U\t \uhat \hat \Sigma^{-2}) \| \\
         &\lesssim \frac{K^{2}}{\sqrt{n}} \max_l \frac{\|\theta\l\|_1}{\|\theta\l\|^2 (\lambda_{\min}\l)^{1/2}} \sqrt{L\log(n)} \| (\Sigma^{-2} \U\t\uhat - \U\t \uhat \hat \Sigma^{-2}) \| \\
         &\lesssim \frac{K^{2} \sqrt{L\log(n)}}{\sqrt{n}} \max_l \frac{\|\theta\l\|_1}{\|\theta\l\|^2 (\lambda_{\min}\l)^{1/2}} \mathcal{B} . \numberthis \label{idioteque}
        \end{align*}
        \item \textbf{The term $\alpha_4$}: We note that
        \begin{align*}
         \| \alpha_4 \|_{2,\infty} &= \| (\mathbf{I} - \U \U\t) \mathcal{Y} \mathcal{L}(\mathcal{E})\t (\uhat - \U \U\t \uhat) \hat \Sigma^{-2} \| \\
         &\lesssim \| (\mathbf{I} - \U \U\t) \mathcal{Y} \mathcal{L}(\mathcal{E})\t \|_{2,\infty} \| \sin\Theta(\uhat,\U) \|  \hat \Sigma^{-2} \|\\
         &\lesssim K n \sqrt{L \log(n)} \bigg( \frac{1}{L} \|\snr\inv \|_2^2 \bigg)^{1/2}\| \sin\Theta(\uhat,\U) \|  \hat \Sigma^{-2} \| \\ 
         &\lesssim K n \sqrt{L \log(n)} \bigg( \frac{1}{L} \|\snr\inv \|_2^2 \bigg)^{1/2} \mathcal{A} \frac{K}{n L \bar \lambda } \\
         &\asymp \mathcal{A}\frac{K^2 \sqrt{\log(n)}}{\sqrt{L} \bar \lambda} \bigg( \frac{1}{L} \|\snr\inv \|_2^2 \bigg)^{1/2} . \numberthis \label{15step}
        \end{align*}
\end{itemize}
Combining \eqref{radiohead}, \eqref{fakeplastictrees}, \eqref{idioteque}, and \eqref{15step}, we have that with probability at least $1 - O(n^{-10})$, 
    \begin{align*}
         \|  (\mathbf{I} - \U \U\t) \mathcal{Y} \mathcal{L}(\mathcal{E})\t \uhat \hat \Sigma^{-2} \|_{2,\infty} &\lesssim \underbrace{\frac{K^3\sqrt{\log(n)}}{n^{3/2} \sqrt{L} \bar \lambda} \max_{l} \frac{\|\theta\l\|_1}{\|\theta\l\|^2 (\lambda_{\min}\l)^{1/2}}}_{\beta_1} \\
           &\quad +\underbrace{\frac{K^{2} \sqrt{L\log(n)}}{\sqrt{n}} \max_l \frac{\|\theta\l\|_1}{\|\theta\l\|^2 (\lambda_{\min}\l)^{1/2}} \mathcal{B}}_{\beta_2} \\
           &\quad + \underbrace{\frac{K^2 \sqrt{\log(n)}}{\sqrt{L} \bar \lambda} \mathcal{A} \bigg( \frac{1}{L} \|\snr\inv \|_2^2 \bigg)^{1/2}}_{\beta_3}.
    \end{align*}
    The proof is complete if we can argue that each of the three quantities $\beta_i$ above are smaller than $c_0 \sqrt{\frac{K}{n}}$ for some sufficiently small constant $c_0$.  
    \begin{itemize}
        \item \textbf{The quantity $\beta_1$:} First, note that \cref{ass:networklevel} implies that
        \begin{align*}
            \frac{\|\theta\l\|_1}{\|\theta\l\|^2 (\lambda_{\min}\l)^{1/2} \bar \lambda} &\leq c_0 \frac{\|\theta\l\|^2}{\log(n) K^8},
        \end{align*}
        provided the constant $C$ in the assumption is sufficiently large.  
        Therefore,
        \begin{align*}
            \beta_1 &= \frac{K^3\sqrt{\log(n)}}{n^{3/2} \sqrt{L} \bar \lambda} \max_{l} \frac{\|\theta\l\|_1}{\|\theta\l\|^2 (\lambda_{\min}\l)^{1/2}}  \\
            &\leq \frac{K^3\sqrt{\log(n)}}{n^{3/2} \sqrt{L} } c_0 \frac{\max_l\|\theta\l\|^2}{\log(n) K^8} \\
            &\leq c_0 \frac{n \max_l (\theta\l_{\max})^2  }{n^{3/2} \sqrt{L} \sqrt{\log(n)} K^5} \\
            &\leq c_0  \frac{1}{\sqrt{nL\log(n)} K^5} \\
            &\leq c_0 \sqrt{\frac{K}{n}}.  
        \end{align*}
                \item \textbf{The quantity $\beta_2$:} By \cref{ass:networklevel},
                \begin{align*}
                    \beta_2 &\leq c_0 \frac{K^2 \sqrt{L\log(n)}}{\sqrt{n}} \mathcal{B} \frac{\|\theta\l\|^2 \bar \lambda}{\log(n) K^8} \\
                    &\lesssim c_0 \frac{ \sqrt{Ln }\bar \lambda}{\sqrt{\log(n)} K^6}  \bigg( \frac{K^3 \sqrt{\log(n)} \big( \frac{1}{L} \| \snr\inv \|_2^2\big)^{1/2}}{n L^{3/2} \bar \lambda^2} + \frac{K^4 \log(n) \| \snr\inv \|_{\infty}^2}{n L \bar \lambda^2} \\
    &\qquad \qquad \qquad  + \frac{K^3 \sqrt{\log(n)} \alpha_{\max} \| \snr\inv \|_{\infty}}{n L \bar \lambda^2} + \frac{K^2\alpha_{\max}}{n L \bar \lambda^2} \bigg) \\
    &\asymp c_0 \frac{( \frac{1}{L} \| \snr\inv \|_2^2)^{1/2}}{K^3 L \sqrt{n} \bar \lambda} + c_0 \frac{\sqrt{\log(n)} \| \snr\inv \|_{\infty}^2}{K^2\sqrt{nL} \bar \lambda} + \frac{\alpha_{\max} \| \snr\inv \|_{\infty}}{K^3\sqrt{nL} \bar \lambda} + \frac{\alpha_{\max}}{\sqrt{\log(n)} \sqrt{nL} \bar \lambda K^4} \\ \\
    &\leq c_0 \frac{1}{K^3L \sqrt{n}} + c_0 \frac{1}{K^2 \sqrt{nL}} + \frac{1}{K^3 \sqrt{nL}} + \frac{1}{K^4\sqrt{nL}} \\
    &\leq c_0 \sqrt{\frac{K}{n}},
                \end{align*}
                which follows from algebra similar to the proof of \cref{thm:step2asympexp} (e.g.,\eqref{verifybb}) which shows that
        \begin{align*}
            \max\bigg\{ \frac{1}{L \bar \lambda} \| \snr\inv \|_2^2, \frac{ \sqrt{\log(n)} \| \snr\inv\|_{\infty}^2}{\bar \lambda}, \frac{\alpha_{\max}}{\bar \lambda} \bigg\} \lesssim 1.
        \end{align*}
        \item \textbf{The quantity $\beta_3$:}  By \eqref{hamlet}, it holds that
        \begin{align*}
            \frac{K^8 \log(n)}{L} \| \snr\inv \|_2^2 \leq c_0\bar \lambda^2, 
        \end{align*}
        and hence,
        \begin{align*}
            \beta_3 &= \frac{K^2 \sqrt{\log(n)}}{\sqrt{L} \bar \lambda} \mathcal{A} \bigg( \frac{1}{L} \| \snr\inv \|_2^2 \bigg)^{1/2} \\
            &\lesssim c_0 \frac{1}{K^2 \sqrt{L}}\mathcal{A} \\
            &\lesssim c_0 \frac{1}{K^3 \sqrt{L}},
        \end{align*}
        which follows from the fact that $\mathcal{A} \lesssim \frac{1}{K}$ as discussed previously.  From our assumption that $L \geq n/ K^7$, we see that $\frac{1}{K^3 \sqrt{L}} \leq c_0 \sqrt{\frac{K}{n}}$.  
    \end{itemize}
    Therefore, combining these arguments and reassigning constants if necessary, the proof is complete.   
\end{proof}

\section{Proof of Minimax Lower Bound (Theorem \ref{thm:minimax})}
\begin{proof}
    Our proof mimics that of \citet{gao_community_2018,han_exact_2021}.  For technical convenience we assume that $n$ is divisible by $KL$ and that the communities are equal-sized.
 
    We now proceed in steps.
\begin{itemize}
    \item \textbf{Step 1: Reduction to fundamental testing problem.} First, let $z \in [K]^{n}$ be such that 
    \begin{align*}
        c \frac{n}{K} = n_1 = n_2 \leq n_3 \leq \cdots \leq n_{K},
    \end{align*}
    with $n_k = |\mathcal{C}_k|$.  In addition, let $\mathbf{B}^{(l)}$ be the matrix $\mathbf{I} +\frac{ \lambda_{\min} - 1}{K}  \mathbf{1} \mathbf{1}\t$, which ensures that $\lambda_{\min}\l = \lambda_{\min}$.  We claim that there exists $\theta\l$ such that $\frac{\|\theta\l_{\mathcal{C}_k}\|^2}{K}$ is the same for all $l \in [L]$ and $k \in [K]$.  In this manner we have that $\mathbf{P}^{(1)},\cdots \mathbf{P}^{(L)} \in \mathcal{P}(\lambda_{\min},K,n,\theta,L)$.  We will verify the existence of $\theta\l$ at the end of the proof.

    For each $a \in [K]$, let $T_a$ be a subset of $\mathcal{C}_a$ with cardinality $|T_a| = \lceil n_a - \frac{c n}{4 K^2} \rceil$.  Let $T = \cup_{a \in [K]} T_a$ and define
    $$\mathcal{Z}_T:=\{z': c\frac{n}{K}\leq |\{j\in[n]:z_j'=a\}|\leq C\frac{n}{K}\text{ for all }a\in[K], z_j'=z_j\text{ for all }j\in T\}.$$
    
    If $\tilde z \neq z$ and $\tilde z \in \mathcal{Z}_T$, we have
    \begin{align*}
        \frac{1}{n} \sum_{j=1}^{n} \mathbb{I}\{ \tilde z_j \neq z_j \} \leq \frac{1}{n} |T^c| \leq \frac{K}{n} \frac{c n}{4 K^2} = \frac{c}{4 K}.
    \end{align*}
    Similarly, if $\pi$ is any non-identity permutation on $[K]$ it holds that 
    \begin{align*}
         \frac{1}{n} \sum_{j=1}^{n} \mathbb{I}\{ \pi(\tilde z_j) \neq z_j \} &\geq \frac{1}{n} \min_a |T_a| \geq \frac{1}{n} \bigg( \frac{c n}{K} - \frac{c n}{4 K^2} \bigg) \geq \frac{3 c}{4 K}.
    \end{align*}
    Therefore for any $\tilde z \in \mathcal{Z}_T$, it holds that the identity permutation is the optimal permutation.  

    Therefore, following \citet{gao_community_2018}, Theorem 2, it holds that
    \begin{align*}
        \inf_{\hat z} \sup_{z} \mathbb{E} \ell(\hat z, z) &\geq \frac{c}{6 K |T^c|} \sum_{j \in T^c} \bigg[ \frac{1}{2 K^2} \inf_{\hat z_j}\big( \p_{H_0}( \hat z_j = 1) + \p_{H_1}(\hat z_j = 0) \big) \bigg] \\
        &\geq C' \frac{1}{K^3 |T^c|} \sum_{j \in T^c} \bigg[  \inf_{\hat z_j}\big( \p_{H_0}( \hat z_j = 1) + \p_{H_1}(\hat z_j = 0) \big) \bigg]
    \end{align*}
    where $H_0$ and $H_1$ are the distributions given $z_j$ is  in community 1 or 2 respectively.      
    \item \textbf{Step 2: Lower bounding the node-wise Type I and Type II error}.  We next note that the quantity  $\inf_{\hat z_j}\big( \p_{H_0}( \hat z_j = 1) + \p_{H_1}(\hat z_j = 0) \big)$ is the sum of the type I and type II error for the simple-simple hypothesis test.  By standard testing results, it holds that
    \begin{align*}
        \p\big( \hat z(i) = 2| z(i) = 1 \big) + \p\big[ \hat z(i) =1 | z(i) = 2\big] &\geq  \frac{1}{2} - \frac{1}{2}\sqrt{\frac{1}{2} d_{KL}(P_0|| P_1)},
    \end{align*}
    where $d_{KL}$ is the K-L divergence between $H_0$ and $H_1$. We claim that with our particular choice of $\theta\l$ we can demonstrate that
    \begin{align*}
        d_{KL}(P_0||P_1) &\leq c_0 \lambda_{\min}^2 \theta_{\max}^2 \frac{n}{L} < 1.
    \end{align*}
    If this is the case,  by combining our arguments and noting that the errors are all the same across all nodes, we complete the proof.  Therefore, it remains to verify the existence of $\theta\l$ satisfying the requisite inequalities.    
    \item \textbf{Details of the construction.} We now explain the construction of $\theta\l$, which we will assign differently according to each node and network.  Recall that we assume that each community is of exact size $\frac{n}{K}$.  Within each community, divide the nodes into subsets of equal size $\frac{n}{KL}$, so that each community $k$ is partitioned into $L$ different subsets $\mathcal{S}_1^k, \cdots \mathcal{S}^k_L$.  We will assign each degree correction within each subset to be either $\theta_{\max}$ or $\theta_{\min} = \frac{\theta_{\max}}{L}$ depending on which network it belongs to such that the degree correction parameter for the nodes in $\mathcal{S}_l^k$ are equal to $\theta_{\max}$ within network $l$ and otherwise are equal to $\theta_{\min}$.  
%
In this manner each vertex has $\theta_{\max}$ as its degree correction parameter exactly once across all the networks, and $\theta_{\min} = \frac{\theta_{\max}}{L}$ as its degree correction parameter $L-1$ times.  Note that 
    \begin{align*}
        \|\theta\|_1 &= K\bigg(\frac{n}{KL} \theta_{\max} + \frac{n(L-1)}{K} \theta_{\min}\bigg) \\
        &= \frac{n}{L} \theta_{\max} + \frac{n(L-1)}{L} \frac{\theta_{\max}}{L} \\
        &= \theta_{\max} \frac{n}{L} \bigg( 2 - \frac{1}{L} \bigg); \\
        \|\theta\|_2^2 &= K\bigg(\frac{n}{KL} \theta_{\max}^2 + \frac{n(L-1)}{KL^2} \theta_{\min}^2\bigg) \\
        &= K \bigg(\frac{n}{KL} \theta_{\max}^2 + \frac{n(L-1)}{KL^3} \theta_{\max}^2\bigg) \\
        &= \theta_{\max}^2 \frac{n}{L} \bigg( 1 + \frac{L-1}{L^2} \bigg).
    \end{align*}
    Consequently, our assumption implies that
    \begin{align*}
    c_1 \leq \frac{\theta_{\max} \|\theta\|_1}{\|\theta\|^4 \lambda_{\min}^2} &= \frac{\theta_{\max}^2 \frac{n}{L}}{\theta_{\max}^4 \frac{n^2}{L^2} \lambda_{\min}^2} \frac{2 - \frac{1}{L}}{(1 + \frac{L-1}{L^2})^2} = \frac{1}{\theta_{\max}^2 \frac{n}{L} \lambda_{\min}^2 } \big( 2 + o(1) \big).
    \end{align*}
    This implies that
    \begin{align*}
        \theta_{\max}^2 \frac{n}{L} \lambda_{\min}^2 \leq \frac{2 + o(1)}{c_1}.
    \end{align*}
    As long as $c_1 \geq 10$, we have that $\theta_{\max}^2 \frac{n}{L} \lambda_{\min}^2 \leq \frac{1}{4}$.  Next, the KL-divergence between two Bernoulli distributions with parameters $a$ and $b$ is governed by
    \begin{align*}
        d_{KL}(a,b) &= a \log(a/b) + (1- a) \log\frac{1-a}{1-b} \leq \frac{(a-b)^2}{b(1-b)}.
    \end{align*}
By the product property of independent tests, we have that as long as $\lambda_{\min} \leq \frac{1}{2}$,
      \begin{align*}
        d_{KL}(P_0||P_1) &\leq \sum_{l=1}^{L} \sum_{i=1}^n \frac{(\theta_i\l \theta_1\l  - \theta_i\l \theta_1\l (1 - \lambda_{\min}))^2}{\theta_i\l \theta_1\l  (1 - \lambda_{\min}) (1 - \theta_i\l \theta_1\l(1 - \lambda_{\min}))}        
        \\
        &\leq \sum_{l=1}^{L} \sum_{i=1}^n \frac{(\theta_i\l \theta_1\l)^2  \lambda_{\min}^2}{\theta_i\l \theta_1\l  (1 - \lambda_{\min}) (1 - \theta_i\l \theta_1\l (1 - \lambda_{\min}))}   \\
        &\leq 2 \sum_{l=1}^{L} \sum_{i=1}^n \theta_i^{(l)} \theta_1^{(l)} \frac{(\lambda_{\min})^2}{1 - \lambda_{\min}} \\
        &\leq  4 \| \theta\|_1 \lambda_{\min}^2 \sum_{l=1}^{L} \theta_{1}^{(l)} \\
        &\leq 4  \| \theta\|_1 \lambda_{\min}^2 \bigg( \theta_{\max} + (L - 1) \theta_{\min} \bigg) \\
        &=4  \theta_{\max}^2 \frac{n}{L} \big(2 - \frac{1}{L} \big)  \lambda_{\min}^2 \bigg(  1 + \frac{(L - 1)}{L} \bigg) \\
        &\leq 16 \theta_{\max}^2 \lambda_{\min}^2 \frac{n}{L} \\
        &\leq \frac{1}{2}. 
    \end{align*}
This verifies our main condition.  

    We also need to check that our condition on the sparsity holds.  We have that
    \begin{align*}
        \frac{\sqrt{\log(n+L)}}{\theta_{\max}^2 \frac{n}{L} \lambda_{\min}^2} \ll \sqrt{L} \iff \theta_{\max} \lambda_{\min} \gg \frac{L^{1/4}}{\sqrt{n}}\log^{1/4}(n+L).
    \end{align*} 
    Therefore,    since $\lambda_{\min} \asymp 1$ we may take $\theta_{\max} \asymp \frac{L^{1/4 +\eps}}{\sqrt{n}} \log^{1/4}(n+L)$. 
\end{itemize}
 This completes the proof.
 \end{proof}
\section{Further simulations and theory for single network spherical clustering} 
\label{sec:singlenetwork}

In the main paper, we have compared our results to the best-known expected misclustering error for spectral clustering without refinement for degree-corrected stochastic blockmodels; i.e., the result in \citet{jin_improvements_2022}. 
However, DC-MASE uses the spherical normalization, and the result in \citet{jin_improvements_2022} uses the SCORE normalization.  While \citet{jin_improvements_2022} demonstrate that the SCORE procedure exhibits an exponential misclustering rate, to the best of our knowledge there is no similarly strong error rate for vanilla spectral clustering with the spherical normalization, though there are  polynomial upper bounds \citep{lei_consistency_2015,qin_regularized_2013}, as well as some perfect clustering results \citep{lyzinski2014perfect,su_strong_2020}.  
Conveniently, as a byproduct of our analysis we characterize the rows of $\yhat\l$,
and we are able to apply the same proof strategy for \cref{thm:clusteringerror} to analyze the result of running $K$-means on these rows.
The following theorem demonstrates an exponential error rate for single network clustering. For simplicity, we suppress the dependence of the parameters on the index $l$.
\begin{theorem}[Single Network Misclustering Rate: Spherical Normalization] \label{thm:singlenetwork}
 Assume that \cref{ass:communityass} and \cref{ass:networklevel} hold (with $\bar \lambda =\lambda_{\min}$).  
Then the output of $(1 + \eps)$ $K$-means on the rows of $\yhat$ satisfies
\begin{align*}
    \E \ell(\hat z, z) &\leq \frac{2K}{n} \sum_{i=1}^{n} \exp\bigg( - c \theta_i \min\bigg\{ \frac{\|\theta\|^4 \lambda_{\min}^2}{K^3\|\theta\|_3^3 }, \frac{\|\theta\|^2 \lambda_{\min}}{K^{3/2} \theta_{\max}} \bigg\} \bigg) + O(n^{-10}).
\end{align*}
\end{theorem}
This rate exactly matches the rate obtained in \citet{jin_improvements_2022}, but the assumptions are somewhat different, which we now describe. First, the signal-strength assumption in \citet{jin_improvements_2022} requires that
\begin{align*}
    \frac{K^{8} \theta_{\max} \|\theta\|_1 \log(n)}{\|\theta\|^4 \lambda_{\min}^{2}} \bigg( \frac{\theta_{\max}}{\theta_{\min}}  \bigg)^2  \lesssim 1.
\end{align*}
In contrast, we require that
\begin{align*}
    \frac{K^{8} \theta_{\max} \|\theta\|_1 \log(n)}{\|\theta\|^4 \lambda_{\min}^{3}}  \frac{\theta_{\max}}{\theta_{\min}}  &\lesssim 1,
\end{align*}
which is weaker whenever $
\frac{1}{\lambda_{\min}} \lesssim \bigg( \frac{\theta_{\max}}{\theta_{\min}} \bigg). 
$
This regime corresponds to high degree heterogeneity relative to the community separation. For example, if the network is sparse (e.g. $\|\theta\| \asymp \sqrt{\log(n)}$), then it must be that $\lambda_{\min} \asymp 1$ (or else the assumption fails).

One possible reason that the spherical normalization requires weaker conditions on the degree heterogeneity is that the spherical normalization is more ``robust'' to severe degree heterogeneity as it uses all the eigenvectors simultaneously to normalize, whereas the SCORE-based approach only uses a single eigenvector.    In essence, the standard deviation of the leading eigenvector exhibits additional dependence on $\frac{\theta_{\max}}{\theta_{\min}}$, but no dependence on $\lambda_{\min}$, whereas the standard deviation of the spherical normalization does not depend as strongly on  $\theta_{\max}/\theta_{\min}$, but has additional dependence on $\lambda_{\min}$.

Furthermore, \citet{jin_improvements_2022} impose an additional (perhaps artificial) assumption on the leading eigenvalue and eigenvector of the matrix $K\| \theta\|^{-2} (\bZ\t \bTheta^2 \bZ)^{1/2} \bB (\bZ\t \bTheta^2 \bZ)^{1/2}$; namely that the leading eigenvalue is well separated from the remaining eigenvalues and the leading eigenvector has entries of similar order. From a technical perspective, such a condition is required so that the SCORE procedure (which uses the entries of the leading eigenvector of the adjacency matrix for its normalization) does not ``explode.''  However, such a condition may be unintuitive.  Consider, for example, the case that all $\theta_i$'s are equal and that $\mathbf{B}$ is of the form
\begin{align*}
    \mathbf{B} &= \begin{pmatrix}
        1 & .5 & 0 \\
        .5 & 1 & b  \\
        0 & b & 1 
    \end{pmatrix}.
\end{align*}
Assuming that all the communities are of equal size, the leading eigenvector is of the form 
\begin{align*}
    \mathbf{v} &= \big( \frac{1}{2b}, \frac{\sqrt{1 + 4b^2}}{2b}, 1 \big)\t
\end{align*}
with corresponding eigenvalue $\frac{1}{2} \big( \sqrt{4b^2 +1} + 2\big)$.  As $b\to 0$, the assumption in \citet{jin_improvements_2022} is violated as the leading two entries of $\mathbf{v}$ diverge.  While such a setting is perhaps slightly contrived (as the third community is significantly easier to separate), the assumption imposed on the leading eigenvector in \citet{jin_improvements_2022} fails to accommodate this scenario.  In contrast, the spherical normalization is able to handle such scenarios.  In essence, the reason for this difference is that the SCORE normalization requires estimating the leading eigenvector with high fidelity, which depends on the gap between the leading eigenvalue and the bottom $K-1$ eigenvalues.  While the Perron-Frobenius Theorem shows that there is  necessarily \emph{some} separation, the additional assumption imposed by \citet{jin_improvements_2022} ensures that that this separation is sufficiently strong.  The spherical normalization does not require such a separation as it uses all eigenvectors simultaneously.

In summary, we see that the assumptions imposed by the spherical normalization require a) slightly weaker assumptions on the degree heterogeneity, b) slightly stronger assumptions on the smallest eigenvalue $\lambda_{\min}$, and c) no additional assumptions ensuring that the leading eigenvector is well-separated.

\subsection{Simulations for single network clustering}
In this section, we study the role of spherical normalization versus other normalization procedures in the DCSBM via simulated data.  

Given an adjacency matrix $\bA$, write its  eigendecomposition as
$$\bA = \widehat{\bU}\widehat{\bLambda}\widehat{\bU} + \widehat{\bU}_\perp\widehat{\bLambda}_\perp\widehat{\bU}_\perp,$$ where $\widehat{\bLambda}\in\real^{K\times K}$ is a diagonal matrix containing the $K$ leading eigenvalues of $\bA$ (in magnitude) and $\widehat{\bU}\in\real^{n\times K}$ is a matrix containing the corresponding $K$ leading eigenvectors. The methods we consider in the simulations are (1) \textit{spherical} spectral clustering using the scaled matrix of eigenvectors $\widehat{\bU}|\widehat{\bLambda}|^{1/2}$, (2) \textit{spherical unscaled} spectral clustering, where we consider the matrix of eigenvectors $\widehat{\bU}$, (3) the \textit{SCORE} normalization as proposed by \cite{jin_fast_2015}, (4) the SCORE+ method of \cite{jin_improvements_2022}, and (5) the \textit{unthresholded SCORE} method of \cite{jin_fast_2015}. We note that \cite{jin_fast_2015} proposes a thresholding approach to remove the low-degree vertices, and provides a theoretical analysis of this step. In principle, implementing a thresholding approach like this can potentially improve methods that use the spherical normalization as well, but we decided to include the un-thresholded SCORE in the simulations to observe the effect of the embedding methodology directly without further removal of low-degree nodes. All variations of the SCORE methods are computed using the \texttt{ScorePlus} R package \citep{jin2022package}.

All networks in the simulations have $n=300$ vertices and $K=3$ communities. Unless explicitly indicated, we consider the following parameters:
$$\bB = \left(\begin{array}{ccc}
     1 & \gamma & \gamma  \\
     \gamma & 1 & \gamma \\
     \gamma & \gamma & 1  
\end{array}\right), \quad\quad \bZ = \mathbf{1}_{n/K} \otimes \mathbf{I}_{K} , \quad\quad \quad \theta_1, \ldots,\theta_n \overset{\text{i.i.d.}}{\sim}\text{Uniform}(0.1, 1),$$
$$\bTheta = \text{diag}(\theta_1, \ldots,\theta_n).  $$
The entries of the symmetric adjacency matrix are sampled independently with probabilities given by $\bP = \alpha \bTheta\bZ\bB\bZ^\top \bTheta$, where $\alpha$ is a constant adjusted to make the average expected degree of $\bP$ equal to 15.
The particular simulation scenarios considered are as follows:
\begin{itemize}
    \item \textit{Between-community connectivity}: the off-diagonal values of $\bB$ are varied from $0$ to $0.7$.
    \item \emph{Community imbalance:} we change the value of $\bB_{11}$ to increase the connectivity of the first community.
    \item \textit{Community sizes:} vertex memberships are assigned independently at random with probabilities $(1/3 + \epsilon, 1/3 - \frac{\epsilon}{3}, 1/3 - \frac{2\epsilon}{3})$, with $\epsilon \in [$
    \item \emph{Degree distribution power}: the degree-correction parameters are simulated as $\theta_1, \ldots, \theta_n\overset{\text{i.i.d.}}{\sim} [U(0.1,1)]^p$, i.e., uniformly distributed random variables raised to the power of $p\geq 0$.
\end{itemize}

The simulation results are shown in Figure~\ref{fig:single-network-comparisons}. The results confirm the theoretical analysis, showing that the spherical normalization is more robust to degree heterogeneity than the SCORE normalization, as demonstrated in the performance with respect to changes in the degree distribution power, but this later one has better performance in terms of community magnitudes and community sizes.

\begin{figure}
    \centering
    \includegraphics[width=\textwidth]{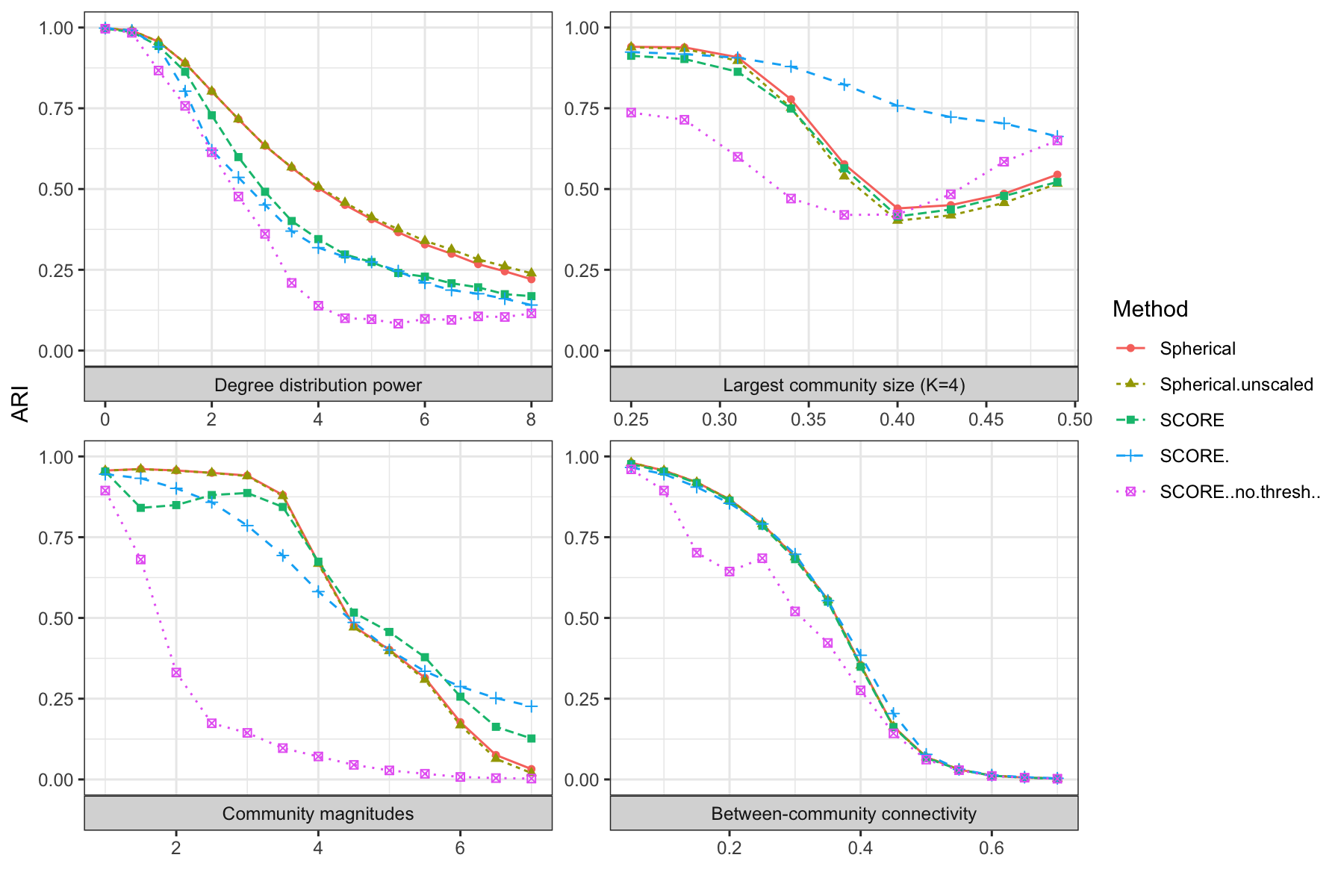}
    \caption{Adjusted Rand index (ARI) of different normalization strategies for single network spectral clustering. Values close to one indicate agreement with the true communities}
    \label{fig:single-network-comparisons}
\end{figure}

\subsection{Scaled vs. unscaled spherical spectral clustering}
This paper mainly considers scaled eigenvectors for spectral clustering in single and multilayer networks. Namely, step 1(a) of Algorithm~\ref{alg:dcmase} uses the matrix of scaled eigenvectors given by $\widehat{\bX}^{(l)} = \widehat{\bU}^{(l)}|\widehat{\Lambda}^{(l)}|^{1/2}$. \cref{thm:clusteringerror} provides an upper bound for the misclustering error rate of a version of single-layer spherical spectral clustering using this scaled matrix
before normalizing its rows and clustering via $K$-means. 
Alternatively, one might consider the unscaled eigenvector matrix $\widehat{\bU}^{(l)}$ followed by row-normalization and clustering \citep{lei_consistency_2015,qin_regularized_2013}. While both approaches can estimate the clusters consistently, we argue that the scaled matrix can alleviate the effect of different community sizes in the performance of spectral clustering. This property can be explained by the interpretation of the embeddings as latent positions of a generalized random dot product graph, which are invariant to community sizes.

To explain the intuition behind the use of scaled eigenvectors, let $\bP = \bTheta \bZ\bB\bZ^\top \bTheta$ with $\bZ\in\{0,1\}^{n\times K}$ be the probability matrix of a DCSBM with $K$ communities, and connectivity matrix $\bB\in\real^{K\times K}$ with $\bB_{kk}=1, k\in[K]$. It can be checked that the matrix of $K$ leading eigenvectors of $\bP$, denoted by $\bU\in\real^{n\times K}$, has the form
$$\bU = \bTheta \bZ \mathbf{T}^{-1/2}\bQ$$ 
for some orthogonal matrix $\mathbf{Q}\in\real^{K\times K}$ and with $\mathbf{T} = \bZ^\top\bTheta^2\bZ$ a $K\times K$ diagonal matrix. In this case, for a given row of $\bU$, say $i\in[n]$, if $\bZ_{ik}=1$, then $\|\bU_{i,\cdot}\| = \theta_i\mathbf{T}^{-1/2}_{kk}= \theta_i/\|\bTheta\bZ e_k\|$, where $e_k\in\real^{K}$ is the standard unit vector. Thus, the norms of the rows of $\bU$ depend on the community size, $n_k$, as well as the magnitude of the degree correction parameters for vertices in the community. On the other hand, the $i$-th row of $\bX = \bU|\bLambda|^{1/2}$ has norm given by $\|\bX_{i,\cdot}\| = \theta_i$, which can be verified by observing that $\bX_{i,\cdot} \bX_{i,\cdot}^\top = \bP_{ii} = \theta_i^2$. Thus, we argue that the row-normalization in spherical spectral clustering is more prone to affect the clustering error in the unscaled eigenvector case when the communities have different sizes, as the variance of this normalization has a different order.

Figure ~\ref{fig:Unorm-vs-Xnorm} shows an illustration of the effect of different community sizes in the embeddings obtained by the normalized rows of $\widehat{\bU}$ (left panel)
and $\widehat{\bX}$ (right panel), demonstrating that the variance of the point clouds of vectors corresponding to the larger communities can be much larger than the 
ones in smaller communities, which might yield poor clustering performance. More specifically, a single graph is generated from the DCSBM model with $n=800$, $K=4$ and connectivity matrix $\bB = (1-\gamma)\bI_K + \gamma\mathbf{1}_K\mathbf{1}_K^\top$, for $\gamma =0.4$. The degree correction parameters are generated at random as $\theta_1, \ldots, \theta_n\overset{\text{i.i.d.}}{\sim}U(0.8, 0.1)$ and the community memberships are assigned at random with probabilities  $\left(\frac{5}{11}, \frac{5}{11}, \frac{1}{22}, \frac{1}{22}\right)$. In the left panel (row-normalized eigenvectors of $\bA$), the spread of the larger clusters (corresponding to communities 1 and 2) dominates, which results in these communities being partitioned into half by the $K$-means algorithm. The right panel (row-normalized scaled eigenvectors of $\bA$) shows that scaling alleviates the effect of different community sizes, as all point clouds show a similar spread. In this particular simulation, the misclustering error in the unscaled eigenvectors is 0.205, whereas the scaled eigenvectors recover the communities perfectly. We repeated this experiment for different values of $\gamma$ in $[0.1, 0.8]$ and averaged the results of 10 simulations (see Table~\ref{table:scaled-vs-unscaled}, observing a superior performance in clustering using the scaled eigenvectors for a wide range of values of $\gamma$.

\begin{figure}
    \centering
    \includegraphics[width=0.48\textwidth]{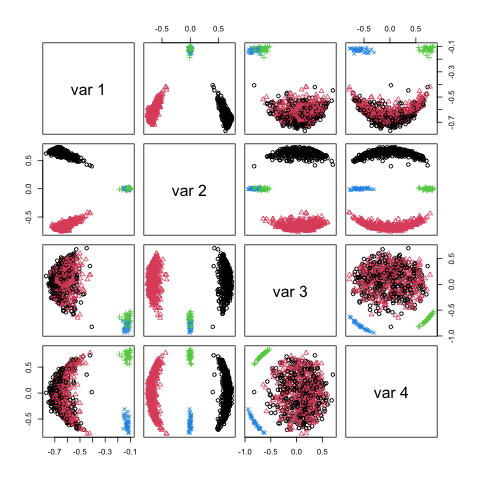}
    \includegraphics[width=0.48\textwidth]{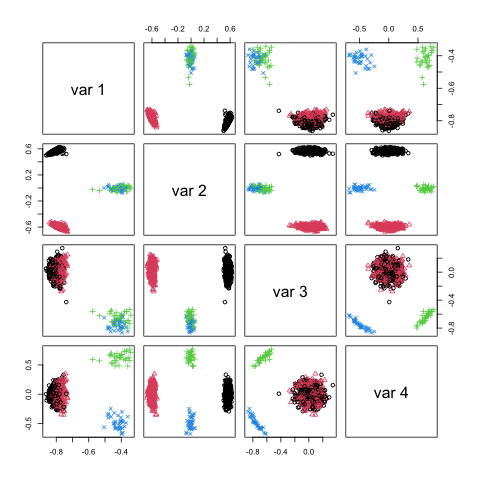}
    \caption{Row-normalized embeddings obtained from the eigenvectors of the adjacency matrix (left) and the scaled eigenvectors (right). The network is generated from a DCSBM with four communities with the fraction of vertices on each of them given by $\left(\frac{5}{11}, \frac{5}{11}, \frac{1}{22}, \frac{1}{22}\right)$.  }
    \label{fig:Unorm-vs-Xnorm}
\end{figure}
\begin{table}[ht]
\centering

\begin{tabular}{r|rrrrrrrr}
  \hline
$\gamma$ & 0.10 & 0.20 & 0.30 & 0.40 & 0.50 & 0.60 & 0.70 & 0.80 \\ 
\hline
  Unscaled eigenvectors & 0.00 & 0.00 & 0.12 & 0.20 & 0.23 & 0.27 & 0.39 & 0.49 \\ 
  Scaled eigenvectors & 0.00 & 0.00 & 0.02 & 0.13 & 0.22 & 0.25 & 0.34 & 0.49 \\ 
   \hline
\end{tabular}
\caption{{Misclustering error rate of spherical spectral clustering using the unscaled vs. scaled eigenvectors as a function of $\gamma$, the connectivity between communities. Graphs are generated from a DCSBM with $K=4$ unbalanced communities and $n=800$.\label{table:scaled-vs-unscaled}}}
\end{table}

\subsection{Proof of Theorem~\ref{thm:singlenetwork}}
\begin{proof}[Proof of \cref{thm:singlenetwork}]
The proof of this result is similar to the proof of the main result.  First we demonstrate the initial error implies that each community contains at least $\frac{3}{4}$  
of its true members, whereupon we study the empirical centroids and show that they are closer to their true cluster centroid than they are to each other.  Finally, instead of applying \cref{thm:step2asympexp} to obtain the exponential error rate we apply \cref{thm:firststep}.  As this result only involves a single network, we suppress the dependency on $l$ for ease of notation.
\\ \ \\ 
\noindent
\textbf{Step 1: Initial Hamming Error}\\
Observe that $\ytilde = \mathbf{Z} \mathbf{M_Y},$ where it straightforward to check that
\begin{align*}
    \sqrt{\lambda_{\min}} \leq \| \big(\mathbf{M_Y}\big)_{r\cdot} - \big(\mathbf{M_Y}\big)_{s\cdot} \| \leq 2.  
\end{align*}
The upper bound is immediate; as for the lower bound, we may apply the same argument as in the proof of \cref{lem:step2popprop}.  Let the matrix $\mathbf{\hat X} := \mathbf{\hat Z}\mathbf{\hat M}_{\mathbf{Y}}$, where $\mathbf{\hat Z}$ and $\mathbf{\hat M}_{\mathbf{Y}}$ are the outputs of $(1 + \eps)$ $K$-means on the rows of $\yhat$, and let $S_r := \{ i \in \mathcal{C}(r): \| \wstar \mathbf{\hat X}_{i\cdot} - \ytilde_{i\cdot} \| \geq \delta_r/2\}$, where $\delta_r = \sqrt{\lambda_{\min}}$.  By Lemma 5.3 of \cite{lei_consistency_2015} and a similar argument as in the proof of \cref{thm:clusteringerror}, it holds that
\begin{align*}
     \inf_{\mathcal{P}} \sum_{i=1}^{n} \mathbb{I}\{ \hat z(i) \neq \mathcal{P}( z(i)) \}  &\leq \frac{C_{\eps}}{\lambda_{\min}} \| \yhat \wstar\t - \ytilde \|_F^2 \\
    &\leq \frac{C_{\eps} n}{\lambda_{\min}}  \| \yhat \wstar\t - \ytilde \|_{2,\infty}^2.
\end{align*}
By \cref{cor:step1twoinfty}, with probability at least $1 - O(n^{-15})$ it holds that
\begin{align*}
  \| \yhat \wstar\t - \ytilde \|_{2,\infty} &\lesssim   \bigg( \frac{\theta_{\max}}{\theta_{\min}} \bigg)^{1/2} \frac{K\sqrt{\log(n)}}{\|\theta\|\lambda_{\min}^{1/2} } \\
     &\quad +\frac{K^2 \theta_{\max}\m \|\theta\m\|_1}{\lambda_{\min}\m \|\theta\m\|^4} \bigg( \log(n) \frac{\theta_{\max}\m}{\theta_{\min}\m} + \frac{1}{\lambda_{\min}\m} +\bigg(\frac{\theta_{\max}\m}{\theta_{\min}\m} \bigg)^{1/2}\frac{K^{5/2}\log(n)}{ (\lambda_{\min}\m)^{1/2}} \bigg) \\
     &\leq \frac{\beta}{8 \sqrt{C_{\eps} K}}\lambda_{\min},
\end{align*}
where $\beta \in (0,1]$  
is such that $n_{\min} \geq \beta n_{\max}$, and where the final bound holds under the conditions of \cref{thm:singlenetwork}. 
Let this event be denoted $\mathcal{E}$.  By squaring the above bound we arrive at
\begin{align*}
 \inf_{\mathcal{P}}\sum_{i=1}^{n} \mathbb{I}\{ \hat z(i) \neq \mathcal{P}( z(i)) \}  &\leq 
   n \frac{\beta^2}{64K } \lambda_{\min} \\
   &\leq \frac{\beta}{64}\lambda_{\min} n_{\min}.
\end{align*}
Therefore, each cluster is associated to a true cluster, denoted as $\mathcal{\hat C}(r)$, where $|\mathcal{\hat C}(r)| \geq (1 - \beta \lambda_{\min}  /64)n_{\min}$ and $|\mathcal{\hat C}(r) \setminus \mathcal{C}(r) | \leq \frac{\beta \lambda_{\min}}{64}  n_{\min}$. Note that since $\beta\in(0,1)$ and $\lambda_{\min} \in (0,1)$, then $\beta \lambda_{\min}/64 < 1$  
this is a well-defined fraction. \\ \ \\ 
\noindent
\textbf{Step 2: Properties of Empirical Centroids}\\
Recall that we denote $(\mathbf{\hat M}_{\mathbf{Y}})_{r\cdot}$ and $(\mathbf{M}_{\mathbf{Y}})_{r\cdot}$ as the cluster centroids for $\mathcal{\hat C}(r)$ and $\mathcal{C}(r)$ respectively. 
Then by a similar argument as in the proof of \cref{thm:clusteringerror}, we have that
\begin{align*}
    \| \wstar (\mathbf{\hat M}_{\mathbf{Y}})_{r\cdot} - (\mathbf{M}_{\mathbf{Y}})_{r\cdot} \| &\leq \frac{1}{|\mathcal{\hat C}(r)|^{1/2}} \| \yhat \wstar\t - \ytilde \|_F + 2 \frac{|\mathcal{\hat C}(r) \setminus \mathcal{C}(r) |}{|\mathcal{\hat C}(r)|}  \\
    &\leq \frac{1}{\sqrt{n_{\min}(1 - \beta \lambda_{\min}/64)}} \sqrt{n} \|\yhat \wstar\t - \ytilde \|_{2,\infty} + 2\frac{\beta n_{\min} \lambda_{\min}}{64 (1 - \beta \lambda_{\min}/64)n_{\min}}  \\
    &\leq \frac{1}{\sqrt{n_{\min}\beta (1 - \beta \lambda_{\min}/64)}} \sqrt{n} \frac{\beta}{8 \sqrt{C_{\eps} K}}\lambda_{\min} + \frac{\lambda_{\min}}{32 (1 - \beta\lambda_{\min}/64)} \\
    &\leq \frac{1}{\sqrt{n_{\min}\beta (1 - \beta \lambda_{\min}/64)}}\sqrt{\frac{K n_{\min}}{\beta}} \frac{\beta}{8 \sqrt{C_{\eps} K}}\lambda_{\min} + \frac{\lambda_{\min}}{32 (1 - \beta\lambda_{\min}/64)} \\
    &\leq \frac{1}{8} \sqrt{\lambda_{\min}},
\end{align*}
since $\lambda_{\min} \in (0,1)$ by assumption.  The above bound holds on the event $\mathcal{E}$. 
\\ \ \\ \noindent \textbf{Step 3: Applying The Asymptotic Expansion} \\
Arguing similarly as in the proof of \cref{thm:clusteringerror}, it holds that
\begin{align*}
    \E \ell(\hat z, z) &\leq \frac{1}{n} \sum_{i=1}^{n} \p\big( \mathbf{Z}_{i\cdot} \neq \mathbf{\hat Z}_{i\cdot}, \mathcal{E}\big) + O(n^{-15}).
\end{align*}
On the event $\mathcal{E}$, it holds that  
\begin{align*}
    \| \mathcal{R}_{\mathrm{Stage \ I}} \|_{2,\infty} &\leq \frac{1}{8} \sqrt{\lambda_{\min}},
\end{align*}
and hence by repeating the arguments in the proof of \cref{thm:clusteringerror},
\begin{align*}
    \p( \mathbf{Z}_{i\cdot} \neq \mathbf{\hat Z}_{i\cdot}, \mathcal{E}) &\leq \p\big( \| (\yhat \wstar\t )_{i\cdot} - \ytilde_{i\cdot} \| \geq \frac{1}{4} \sqrt{\lambda_{\min}}, \mathcal{E} \big) \\
    &\leq K \max_k \p\bigg\{ \big| e_i\t \mathcal{L}(\mathbf{E} ) \U\l |\Lambda\l|^{-1/2} \ipq e_k \big| \geq \frac{1}{4} \sqrt{\lambda_{\min}/K}\bigg\}.
\end{align*}
Here $\mathcal{L}(\mathbf{E})$ is the linear term from \cref{thm:firststep} with $\mathbf{E} = \bA - \bP$.  We now apply Bernstein's inequality.  The variance $v$ is upper bounded by
\begin{align*}
    v &\leq \sum_{j} \theta_i \theta_j \| e_j\t \U\l \|^2 \| | \Lambda\l|^{-1/2} \|^2 \| \mathbf{J}(\mathbf{X}_{i\cdot}) \|^2 \\
    &\lesssim \sum_{j} \theta_i \theta_j \frac{\theta_j^2 K}{\|\theta\|^2} \frac{K}{\|\theta\|^2 \lambda_{\min}} \frac{1}{\theta_i^2} \\
    &\lesssim \frac{K^2 \|\theta\|_3^3}{\|\theta\|^4 \lambda_{\min} \theta_i}.
\end{align*}
Similarly,
\begin{align*}
    \max_{j} \| e_j\t \U\l \| \| |\Lambda\l|^{-1/2} \| \| \mathbf{J}(\mathbf{X}_{i\cdot}) \| &\lesssim \theta_j \frac{K}{\|\theta\|^2 \lambda_{\min}^{1/2}} \frac{1}{\theta_i} \\
    &\lesssim   \frac{K \theta_{\max}}{\|\theta\|^2 \lambda_{\min}^{1/2} \theta_i}.
\end{align*}
By Bernstein's inequality,
\begin{align*}
    \p\bigg\{ \big| e_i\t \mathcal{L}(\mathbf{E} ) \U\l |\Lambda\l|^{-1/2} \ipq e_k \big| \geq \frac{1}{8} \sqrt{\lambda_{\min}/K}\bigg\} 
    &\leq 2 \exp\bigg( - \frac{\frac{1}{128} \frac{\lambda_{\min}}{K}}{C\frac{K^2 \|\theta\|_3^3}{\|\theta\|^4 \lambda_{\min} \theta_i}+ \frac{\lambda_{\min}^{1/2}}{24 \sqrt{K}} C  \frac{K \theta_{\max}}{\|\theta\|^2 \lambda_{\min}^{1/2} \theta_i}} \bigg) \\
    &\leq 2 \exp\bigg( - c \theta_i \min\bigg\{ \frac{ \|\theta\|^4 \lambda_{\min}^2}{K^3 \|\theta\|_3^3}, \frac{\|\theta\|^2 \lambda_{\min}}{K^{3/2}\theta_{\max}} \bigg\} \bigg).
\end{align*}
Assembling everything together completes the proof.
\end{proof}


\section{Additional simulation experiments for multilayer networks}
\label{sec:additionalsims}

{We evaluate the performance of different multilayer community detection methods in terms of the sparsity of the networks. We use the same simulation settings described in \cref{sec:sims}, but here we fix the number of layers as $L = 20$ and $n=150$, and we change the value of the edge density, which is controlled by $\alpha^{(l)}$. In particular, we change the value of this parameter in order to obtain a specific expected edge density, defined as $\frac{1}{n(n-1)}\sum_{i,j}\bP^{(l)}$. This value of the edge density is changed in the range $[2/150, 24/150]$. }

{The results of this simulation are shown in Figure~\ref{fig:simulation-sparsity}. For a given edge density (x-axis), a point representing the average misclustering error of 100 simulation results using a given method is plotted. The results show that DC-MASE is always able to estimate the communities correctly if the networks are sufficiently dense. For the scenarios considered and the range of edge density values, there is no other method that is able to always perform perfect clustering. The only other method that always improves its performance with more density is graph-tool. This might be expected from the fact that this method uses the correct likelihood for the model, but DC-MASE substantially outperforms the method in the last column. Besides, DC-MASE is computationally more scalable than graph-tool.}
\begin{figure}
    \centering
    \includegraphics[width=\textwidth]{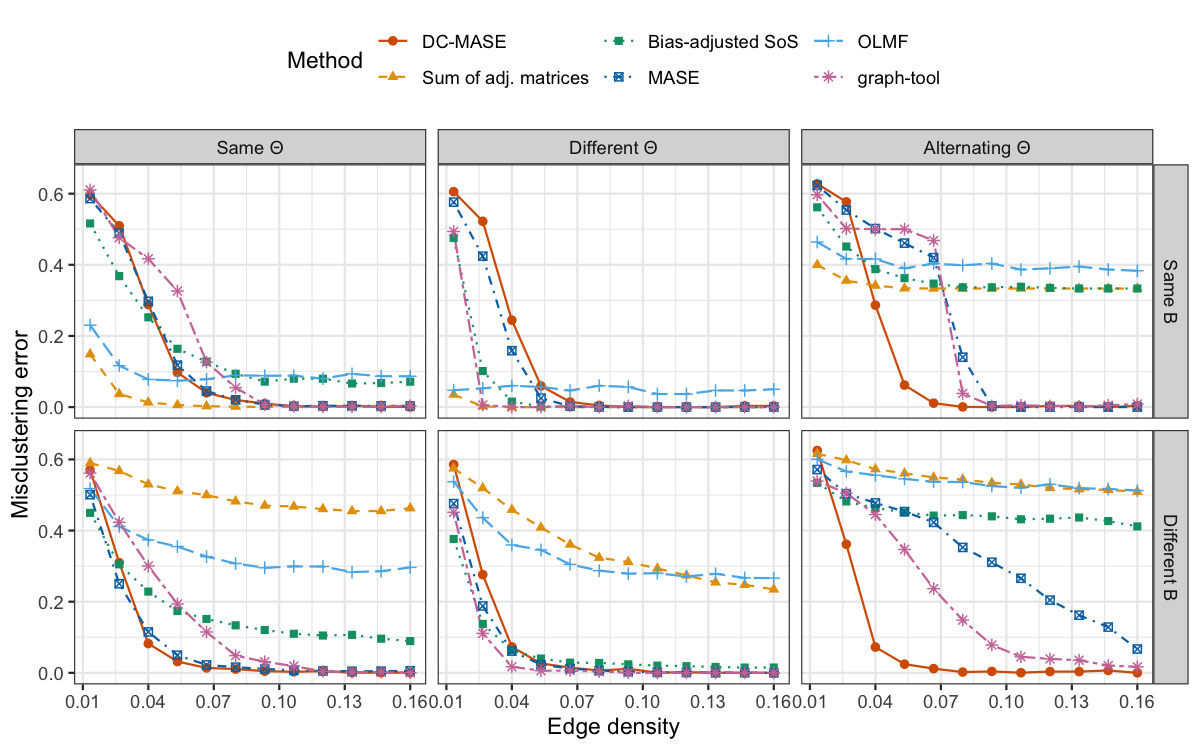}
    \caption{{Misclustering error rate of different community detection methods as a function of the edge density of the networks. The number of layers is fixed as $L=20$, and the simulation scenarios contemplate different types of heterogeneity in the parameters, as described in Section~\ref{sec:additionalsims}.}}
    \label{fig:simulation-sparsity}
\end{figure}

\section{Additional data results} \label{sec:additionaldatanalysis}

\subsection{Out-of-sample performance and robustness to choice of $K$}
{We compared the performance of DC-MASE with the other spectral clustering algorithms considered in \cref{sec:sims}. In the absence of ground truth communities, we measure the performance in terms of out-of-sample mean squared error (MSE) for a given graph $l$ and some number of communities $K$, defined as
$$\text{MSE}(K, l) = \frac{1}{n^2}\|\bA\m - \widehat{\bP}^{(l)}_{\widehat{\bZ}^{(-l,K)}}\|_F^2.$$
Here, $\widehat{\bZ}^{(-l,K)}$ indicates the estimated community memberships obtained from a particular method fitted on the set of graphs indexed by $[L]\setminus\{l\}$ with $K$ communities. Given $\widehat{\bZ}$, the value of the expected adjacency matrix is estimated as $\widehat{\bP}^{(l)}_{\widehat{\bZ}} = \widehat{\bTheta}\m_{\widehat{\bZ}} {\widehat{\bZ}} \widehat{\bB}\m_{\widehat{\bZ}} {\widehat{\bZ}}^\top \widehat{\bTheta}\m_{\widehat{\bZ}}$, where $\widehat{\bTheta}\m_{\widehat{\bZ}}$ and $\widehat{\bB}\m_{\widehat{\bZ}}$ are the plug-in estimates defined via Eq.~\eqref{eq:plug-in-parameters} using the communities defined by $\widehat{\bZ}$. As the expected value of the average MSE is minimized by the  expected adjacency matrices calculated with the correct communities, small values of this quantity are a proxy for the quality of the community estimates.}

After calculating the MSE for all the graphs in the data and for different values of $K$, we performed a paired comparison via the MSE difference between the results for a given method and DC-MASE for each value of $K$ and $l$. Figure~\ref{fig:airport-paired-dff} shows boxplots of these differences across all values of $l\in[L]$ and as a function of the number of communities. Notably, the MSE differences are positive for almost all  graphs in the data and all values of $K$, indicating that the communities obtained by DC-MASE generally have smaller generalization error than the ones obtained by the other spectral methods considered.

\begin{figure}
    \centering
    \includegraphics[width=0.5\textwidth,keepaspectratio]{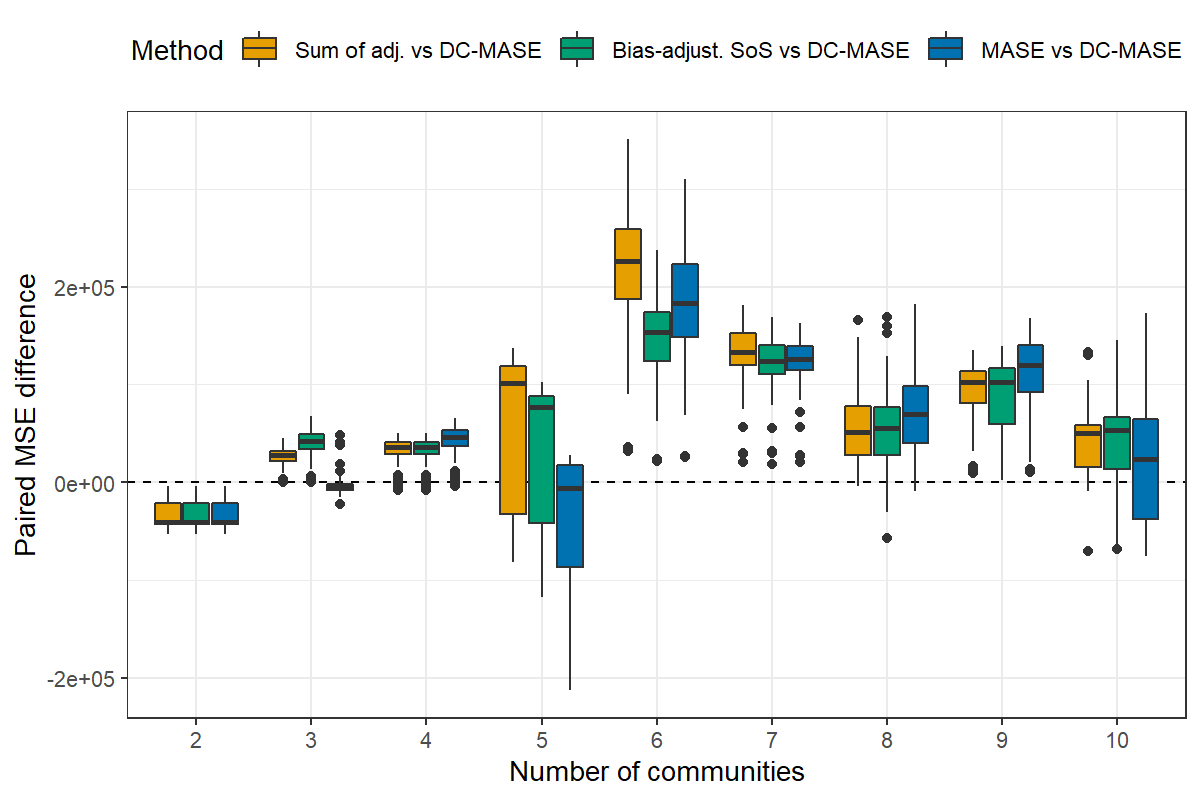}
    \caption{Paired out-of-sample mean squared error (MSE) difference for the Frobenius error of the estimated expected adjacency matrices obtained by each method and DC-MASE. Positive values indicate that the MSE of the respective method is larger than the MSE of DC-MASE.}
    \label{fig:airport-paired-dff}
\end{figure}

\subsection{On the common community membership assumption}

To validate a multilayer DCSBM with common community memberships across time in the airport network data, we compare the community memberships in each month. For this goal, we fit community memberships for each layer (month) independently using our method on a single network (assuming $K=4$, as estimated before). We then compare the community memberships recovered by each layer with the overall community memberships using our joint spectral clustering algorithm with all the layers together. Figure~\ref{fig:before-after} shows the percentage of nodes with different memberships on each clustering result (monthly vs overall) for each month. As can be noticed, most of the months before the start of the pandemic (month 50, corresponding to February 2020) were in close agreement with the communities recovered by the joint clustering. During this period, the percentage of difference usually ranges between 5\% to 10\%. However, when the pandemic started, this percentage rose up to 33\%, and went down again by the end of the period of study. This suggests that a model with constant communities across the layers is at least reasonable before the pandemic. After the pandemic started, the number of flights generally decreased, making the networks sparser and, hence, increasing the uncertainty in the membership estimates. Nevertheless, Theorem \ref{thm:newtheorem} suggests that even if there are changes in the community memberships, our method is still able to recover a common clustering structure for the majority of the networks.

\begin{figure}
    \centering
    \includegraphics[width=0.5\textwidth]{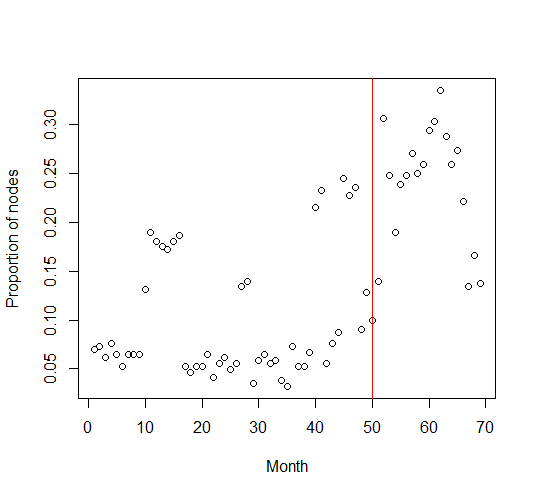}
    \caption{{Proportion of nodes with different community memberships on a model fit in a single month compared with the community memberships obtained by our algorithm for joint spectral clustering with all months. Overall, most months before February 2020 (red vertical line) show agreement in terms of community memberships.}}
    \label{fig:before-after}
\end{figure}



\bibliography{dc_mase2,DC_MASE}

\end{document}